\documentclass{book}
\usepackage[finnish]{babel}
\usepackage[T1]{fontenc}
\usepackage[utf8]{inputenc}
\usepackage{index}
\usepackage{epsfig}
\usepackage{amsmath,amsthm,amssymb}

\makeindex

\title{Johdatus funktionaalianalyysiin}
\author{Tommi Höynälänmaa}

\newcommand{\defterm}[1]{\textbf{#1}}
\newcommand{\setsep}{\;\vert\;}
\newcommand{\bigsetsep}{\;\Big\vert\;}
\newcommand{\naturalnumbers}{\mathbb{N}}
\newcommand{\integernumbers}{\mathbb{Z}}
\newcommand{\rationalnumbers}{\mathbb{Q}}
\newcommand{\realnumbers}{\mathbb{R}}
\newcommand{\complexnumbers}{\mathbb{C}}
\newcommand{\positiveintegers}{\mathbb{Z}_+}
\newcommand{\extrealnumbers}{\realnumbers_*}
\newcommand{\powerset}[1]{2^{#1}}
\newcommand{\union}{\cup}
\newcommand{\intersection}{\cap}
\newcommand{\onto}{\twoheadrightarrow}
\newcommand{\setimage}[2]{#1[#2]}
\newcommand{\abs}[1]{\vert{#1}\vert}
\newcommand{\szabs}[1]{\left\vert{#1}\right\vert}
\newcommand{\norm}[1]{\Vert{#1}\Vert}
\newcommand{\sznorm}[1]{\left\Vert{#1}\right\Vert}
\newcommand{\bignorm}[1]{\big\Vert{#1}\big\Vert}
\newcommand{\norminspace}[2]{\norm{{#1}\;\vert\;{#2}}}
\newcommand{\bignorminspace}[2]{\sznorm{{#1}\;\Big\vert\;{#2}}}
\newcommand{\norminfty}[1]{\norm{#1}_\infty}
\newcommand{\rn}{\realnumbers^n}
\newcommand{\inverseimage}[2]{{#1}^{-1}[#2]}
\newcommand{\eqclass}[2]{E(#1;#2)}

\newcommand{\eqc}[1]{\overline{#1}}
\newcommand{\seq}[1]{\mathbf{#1}}
\newcommand{\rnvec}[1]{\mathbf{#1}}
\newcommand{\cnvec}[1]{\mathbf{#1}}
\newcommand{\seqelem}[2]{{#1}[#2]}
\newcommand{\vecelem}[2]{{#1}[#2]}
\newcommand{\vect}[1]{\mathbf{#1}}
\newcommand{\imagunit}{\mathsf{i}}
\newcommand{\seqalt}[2]{\left({#1}_{#2}\right)_{{#2}=0}^\infty}
\newcommand{\vectorsubspace}{\subset_{\textrm{v.s.}}}
\newcommand{\unitvec}[2]{\mathbf{e}^{[#1]}_{#2}}
\newcommand{\ip}[2]{\langle{#1},{#2}\rangle}
\newcommand{\szip}[2]{\left\langle{#1},{#2}\right\rangle}
\newcommand{\distrappl}[2]{\langle{#1},{#2}\rangle}
\newcommand{\szdistrappl}[2]{\left\langle{#1},{#2}\right\rangle}
\newcommand{\columnvector}[2]{\left(\begin{array}{c}{#1}\\ \vdots\\{#2}\end{array}\right)}
\newcommand{\spaceafter}{\;\;\;\;}
\newcommand{\openball}[3]{B({#1};{#2},{#3})}
\newcommand{\szopenball}[3]{B\left({#1};{#2},{#3}\right)}
\newcommand{\closedball}[3]{\bar{B}({#1};{#2},{#3})}
\newcommand{\Lp}[2]{L^{#1}({#2})}
\newcommand{\Cb}[1]{C_b({#1})}
\newcommand{\Co}[1]{C_0({#1})}
\newcommand{\Cu}[1]{C_{\textrm{u}}({#1})}
\newcommand{\Ccom}[1]{C_{\textrm{com}}({#1})}
\newcommand{\equaltvs}{=_{\textrm{tvs}}}
\newcommand{\mca}{\mathcal{A}}
\newcommand{\mcr}{\mathcal{R}}
\newcommand{\sigmaalgebra}{\(\sigma\)-algebra}
\newcommand{\besovspace}[4]{B^{#1}_{#2,#3}({#4})}
\newcommand{\besovaltnorm}[6]{\norm{{#6}\;\vert\;\besovspace{#1}{#2}{#3}{#4};{#5}}}
\newcommand{\tlspace}[4]{F^{#1}_{#2,#3}({#4})}
\newcommand{\holderspace}[2]{C^{#1}({#2})}
\newcommand{\zygmundspace}[2]{\mathcal{Z}^{#1}({#2})}
\newcommand{\sobolevspace}[3]{W^{#1}_{#2}({#3})}
\newcommand{\slobspace}[3]{W^{#1}_{#2}({#3})}
\newcommand{\besselpotspace}[3]{H^{#1}_{#2}({#3})}
\newcommand{\localhardyspace}[2]{h_{#1}({#2})}
\newcommand{\bmospace}[1]{\mathrm{bmo}({#1})}
\newcommand{\ft}{\mathbf{F}}
\newcommand{\contemb}{\subset_{\textrm{c}}}
\newcommand{\diff}[1]{\Delta_{#1}}
\newcommand{\higherdiff}[2]{\Delta^{#1}_{#2}}
\newcommand{\genmodcont}[4]{\omega^{#1}_{#2}({#3},{#4})}
\newcommand{\modcont}[2]{\omega({#1},{#2})}
\newcommand{\intervaloc}[2]{]{#1},{#2}]}
\newcommand{\intervalco}[2]{[{#1},{#2}[}
\newcommand{\intervaloo}[2]{]{#1},{#2}[}
\newcommand{\intervalcc}[2]{[{#1},{#2}]}
\newcommand{\testfunctionspace}[1]{\mathcal{D}({#1})}
\newcommand{\tfs}[1]{\testfunctionspace{#1}}
\newcommand{\distrspace}[1]{\mathcal{D}'({#1})}

\newcommand{\tempdistrspace}[1]{S'({#1})}

\newcommand{\quotientspace}[2]{{#1}\;/{#2}}
\newcommand{\vquotientspace}[2]{{#1}/{#2}}
\newcommand{\directsum}{\oplus}
\newcommand{\zeroset}{\{0\}}
\newcommand{\herm}{\dagger}
\newcommand{\mideriv}[1]{\partial^{#1}}

\newcommand{\szceil}[1]{\left\lceil{}{#1}\right\rceil}
\newcommand{\cp}[2]{\left({#1}\times{#2}\right)}
\newcommand{\rnb}{\rnvec{b}}
\newcommand{\rnd}{\rnvec{d}}
\newcommand{\rnh}{\rnvec{h}}
\newcommand{\rnx}{\rnvec{x}}
\newcommand{\rny}{\rnvec{y}}
\newcommand{\rnz}{\rnvec{z}}
\newcommand{\rnxo}{\rnvec{x}_0}
\newcommand{\id}[1]{\mathrm{id}_{#1}}
\newcommand{\pd}[1]{\frac{\partial}{\partial{}#1}}
\newcommand{\isomemb}{\subset_1}
\newcommand{\sqordfunction}{{\sigma_{\mathrm{sq}}}}
\newcommand{\sqordfirstfunction}{{\sigma_{\mathrm{sq1}}}}
\newcommand{\sqordsecondfunction}{{\sigma_{\mathrm{sq2}}}}
\newcommand{\sqord}[1]{{\sqordfunction({#1})}}
\newcommand{\sqordfirst}[1]{{\sqordfirstfunction({#1})}}
\newcommand{\sqordsecond}[1]{{\sqordsecondfunction({#1})}}
\newcommand{\mytopspace}{T}
\newcommand{\mytop}{\tau}
\newcommand{\firsttopspace}{T}
\newcommand{\firsttop}{\tau}
\newcommand{\secondtopspace}{U}
\newcommand{\secondtop}{\upsilon}

\DeclareMathOperator{\bilop}{bil}
\newcommand{\bil}[2]{\bilop({#1}\times{#2})}

\DeclareMathOperator{\ReAlt}{Re}
\DeclareMathOperator{\ImAlt}{Im}
\DeclareMathOperator{\spanop}{span}
\DeclareMathOperator{\clos}{clos}
\DeclareMathOperator{\interior}{int}
\DeclareMathOperator{\exterior}{ext}
\DeclareMathOperator{\esssup}{ess\;sup}
\DeclareMathOperator{\essinf}{ess\;inf}
\DeclareMathOperator{\bor}{bor}
\DeclareMathOperator{\stsupp}{supp_{set}}
\DeclareMathOperator{\supp}{supp}
\DeclareMathOperator{\image}{im}
\DeclareMathOperator{\dist}{dist}

\newenvironment{examples}{\vspace{0.5cm}\noindent\textbf{Esimerkkejä:}\begin{enumerate}}{\end{enumerate}}
\newenvironment{example}{\item}{}

\newenvironment{exercises}{\vspace{0.5cm}\noindent\textbf{Tehtäviä:}\begin{enumerate}}{\end{enumerate}}
\newenvironment{exercise}[1]{\item[#1]}{}

\newenvironment{solution}[1]{\vspace{0.5cm}\noindent\textbf{Ratkaisu tehtävään {#1}:}}{}

\newtheorem{theorem}{Lause}[section]
\newtheorem{corollary}[theorem]{Seurauslause}

\theoremstyle{remark}
\newtheorem{remark}[theorem]{Huomautus}

\theoremstyle{definition}
\newtheorem{definition}[theorem]{Määritelmä}

\begin{document}

\frontmatter

\maketitle

\begin{center}
  \textbf{Abstract}
\end{center}

This book concentrates on functional analysis. The text is written so
that it can be followed on the basis of high school mathematics. The
book introduces the set theoretical foundations of mathematics, the
basic theories of linear algebra and topology, and the theory of
topological vector spaces. Distributions, measure theory, and some
common spaces of complex-valued functions of real variables are
handled, too.

\vspace{1cm}

\noindent
-----

\begin{center}
  \textbf{Tiivistelmä}
\end{center}

Tämä kirja keskittyy funktionaalianalyysiin. Teksti on kirjoitettu
niin, että sitä pystyy seuraamaan lukiomatematiikan
pohjalta. Kirjassa käydään läpi matematiikan joukko-opilliset
perusteet, lineaarialgebran ja topologian perusteoriat ja topologisten
vektoriavaruuksien teoria. Lisäksi kirjassa käsitellään distribuutiot,
mittateoria, ja joitakin yleisiä reaalimuuttujien kompleksiarvoisia
funktioavaruuksia.

\tableofcontents

\listoffigures

\chapter{Esipuhe}

Funktionaalianalyysi voidaan määritellä topologisia vektoriavaruuksia
tutkivaksi matematiikan haaraksi. Tässä mielessä funktionaalianalyysi
on topologian ja lineaarialgebran synteesi.

Luku \ref{ch:johdanto} esittelee tarvittavia peruskäsitteitä, kuten
joukko, relaatio ja funktio. Luku \ref{ch:ryhmat-ja-kunnat} käsittelee
myöhemmin tarvittavia algebrallisia struktuureja: ryhmiä ja
kuntia. Luvussa \ref{ch:lukujoukot} konstruoidaan lukujoukot:
luonnolliset luvut, kokonaisluvut, rationaaliluvut, reaaliluvut ja
kompleksiluvut. Luku \ref{ch:vektoriavaruudet} esittelee
vektoriavaruudet ja luku \ref{ch:topologiset-avaruudet} jonot, verkot,
filtterit, topologiset avaruudet ja topologisen avaruuden
erikoistapauksen metrisen avaruuden. Luku
\ref{ch:topologiset-vektoriavaruudet} käsittelee topologisia
vektoriavaruuksia, jotka ovat karkeasti vektoriavaruuksia, joissa on
määritelty topologia. Tässä luvussa määritellään myös lokaalikonveksit
avaruudet, sisätuloavaruudet ja normiavaruudet, jotka ovat
topologisten vektoriavaruuksien erikoistapauksia.
Luvussa \ref{ch:distribuutiot} määritellään
distribuutiot ja temperoidut distribuutiot.
Luku
\ref{ch:mittateoriaa} käsittelee mittateoriaa ja siinä määritellään
Lebesguen integraali.
Luvussa
\ref{ch:yleisia-funktioavaruuksia} määritellään usein tarvittavia
funktioavaruuksia.

Yleisesti kukin luku nojaa kaikkiin sitä edeltäviin
lukuihin. Poikkeuksena mittateoriaa käsittelevä luku
\ref{ch:mittateoriaa} ei tarvitse vektoriavaruuksien, topologisten
avaruuksien eikä topologisten vektoriavaruuksien teoriaa.
%%Kuvassa \ref{fig:riippuvuudet}
Kuvassa 1
on esitetty
kaavio lukujen välisistä riippuvuuksista: \(A \rightarrow B\)
tarkoittaa, että luku \(B\) riippuu luvusta \(A\).

\begin{figure}[htb]
\label{fig:riippuvuudet}
\caption{Lukujen väliset riippuvuudet.}
\includegraphics{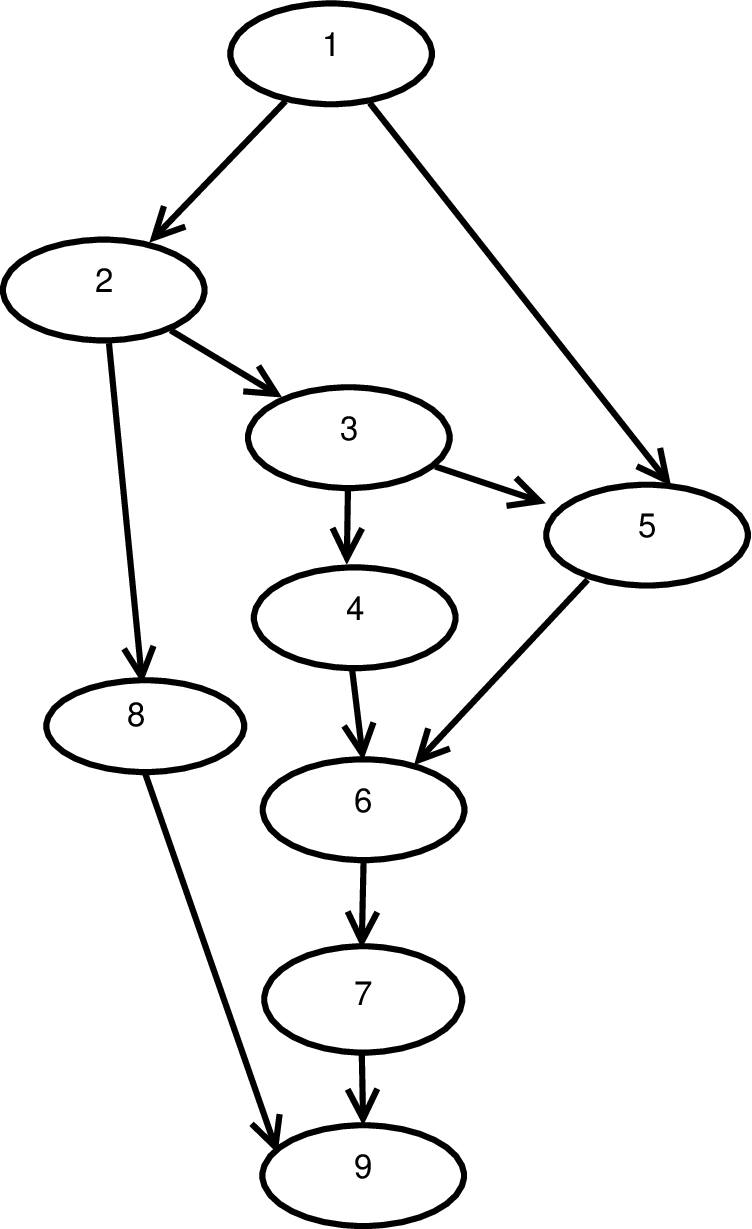}
\end{figure}

\mainmatter

\chapter{Johdanto}
\label{ch:johdanto}

\section{Matematiikan perusteista}

Matemaattinen teoria koostuu joukosta aksioomia eli perusoletuksia,
joista kaikki teorian lauseet johdetaan. Tosin Gödelin lauseen mukaan
jokaisessa teoriassa, joka sisältää luonnolliset luvut, on olemassa
tosia lauseita, joita ei voida johtaa teorian aksioomista.
Luvut, joukot ja muut matemaattiset käsitteet ovat ongelmallisia
materialistisen todellisuuskäsityksen kannalta, koska ne eivät ole
ainetta. Ilmeisesti todellisten olioiden joukon (ei
matemaattisessa mielessä) on oltava laajempi kuin aineellisten
olioiden joukko.

\index{algoritmi}
\defterm{Algoritmi} on lista toimintaohjeista jonkin ongelman
ratkaisemiseksi. On olemassa myös algoritmeja, jotka eivät pääty
äärellisessä ajassa, ja ongelmia, joihin ei ole algoritmia ollenkaan.
\index{Turingin kone}
\defterm{Turingin kone} on kuvitteellinen laite joka koostuu
lukupäästä ja sen läpi kulkevasta nauhasta. Kone lukee vuorollaan
merkin nauhasta ja päättää sen perusteella, miten se siirtää nauhaa
seuraavaksi. Algoritmisesti ratkeavien ongelmien luokka on sama, kuin
ne ongelmat, jotka Turingin kone pystyy ratkaisemaan. Tämä luokka on
myös sama kuin ne ongelmat, jotka tietokone pystyy ratkaisemaan.  Jos
algoritmin nopeus on syötteen määrän polynomiaalinen funktio, niin
algoritmi kuuluu luokkaan P. Muussa tapauksessa algoritmi kuuluu
luokkaan NP. On avoin kysymys, ovatko P ja NP samat.

\section{Joukot}

\index{joukko}
Matematiikan peruskäsite on \defterm{joukko}. Joukko on olio, jolle
jokaisesta toisesta oliosta voidaan sanoa, kuuluuko se joukkoon vai
\index{alkio}
ei. Joukkoon kuuluvia olioita sanotaan sen \defterm{alkioiksi}.
Aikoinaan Bertrand Russell keksi määritellä joukon, joka sisältää
kaikki ne joukot, jotka eivät ole itsensä alkioita, siis
\begin{math}
P := \{ x \setsep x \not\in x \}
\end{math}.
Merkintä \(A := B\) tarkoittaa, että muuttujan \(A\) arvoksi
määritellään \(B\).
Nyt voidaan kysyä, onko \(P\) itsensä alkio. Jos se on, se ei ole, ja
\index{Russellin paradoksi}
jos se ei ole, se on. Tämä on \defterm{Russellin
  paradoksi}. Paradoksin johdosta Frege katsoi koko luomansa
matematiikan aksiomaattisen järjestelmän romahtavan, mutta näin synkkä
ei tilanne ollut. Paradoksi voitiin kiertää rajoittamalla sitä, mitä
joukon jäsenenä voi olla. Olio, joka ei voi olla minkään joukon alkio,
\index{luokka}
on nimeltään \defterm{luokka}. Yleisesti käytetty joukko-opin
\index{Zermelo-Fraenkelin aksioomajärjestelmä}
Zermelo-Fraenkelin aksioomajärjestelmä ei sisällä luokkia, mutta von
\index{von Neumann-Bernays-Gödelin aksioomajärjestelmä}
Neumann-Bernays-Gödelin aksioomajärjestelmä sisältää ne.

Joukko voidaan määritellä esim. luettelemalla sen alkiot: \(A := \{ 1,
3, 5, 6 \}\) tai antamalla jokin predikaatti (ehto), joka kertoo,
kuuluuko kukin olio joukkoon vai ei:
\(B := \{ x \in \realnumbers \setsep x > 5 \}\).
Joukko \(B\) sisältää kaikki reaaliluvut, jotka ovat suurempia kuin
\index{tyhjä joukko}
5. \defterm{Tyhjää joukkoa} merkitään \(\emptyset\).
Jos \(A\) ja \(B\) ovat joukkoja, niin sanomme, että \(A\) on
\index{osajoukko}
joukon \(B\) \defterm{osajoukko}, jos jokainen joukon \(A\) alkio
kuuluu joukkoon \(B\). Tätä merkitään \(A \subset B\). Jos \(A \subset
B\) mutta \(A \not= B\), niin sanomme, että \(A\) on joukon
\index{aito osajoukko}
\(B\) \defterm{aito osajoukko}.
\index{potenssijoukko}
Joukon \(A\) \defterm{potenssijoukko} \(\powerset{A}\) on joukon
\(A\) kaikkien osajoukkojen joukko. Esim. jos \(A := \{ 1, 2, 3 \}\),
niin
\[\powerset{A} = \{ \emptyset, \{1\}, \{2\}, \{3\}, \{1, 2\},
\{1, 3\}, \{2, 3\}, \{1, 2, 3\}\}.\]

Kun alkio \(x\) kuuluu joukkoon \(A\), niin merkitään \(x \in A\). Jos
\(x\) ei kuulu joukkoon \(A\), niin merkitään \(x \not\in
\index{unioni}
A\). Joukkojen \(A\) ja \(B\) \defterm{unioni} \(A \union B\) koostuu molempien
joukkojen kaikista alkioista. On huomattava, ettei mikään alkio voi
esiintyä unionissa tai muussakaan joukossa kahta kertaa. Esim. jos \(A
:= \{1, 2, 3\}\) ja \(B := \{ 3, 5 \}\), niin \(A \union B = \{ 1, 2,
\index{leikkaus}
3, 5 \}\).  Joukkojen \(A\) ja \(B\) \defterm{leikkaus} \(A \intersection B\)
koostuu alkioista, jotka kuuluvat sekä joukkoon \(A\) että joukkoon
\(B\). Edellisen esimerkin joukoilla \(A \intersection B = \{ 3 \}\).
\index{joukkojen erotus}
\index{erotus}
Joukkojen \(A\) ja \(B\) \defterm{erotus} määritellään
\begin{displaymath}
  A \setminus B := \{ x \in A \setsep x \not\in B \}
\end{displaymath}
\index{komplementti}
ja sitä sanotaan myös joukon \(B\) \defterm{komplementiksi} joukossa
\(A\).

\index{karteesinen tulo}
Joukkojen \(A\) ja \(B\) \defterm{karteesinen tulo} määritellään
\begin{displaymath}
  A \times B := \{ (x,y) \setsep x \in A, y \in B \} .
\end{displaymath}
Esim. jos \(A := \{1,2,3\}\) ja \(B := \{4,5\}\), niin
\begin{displaymath}
  A \times B = \{
  (1, 4), (1, 5), (2, 4), (2, 5), (3, 4), (3, 5)\} .
\end{displaymath}
Huomaa, että parit (x,y) ovat järjestettyjä: esim. \((1,2)
\not= (2,1)\).
Kun \(A\) on joukko ja \(n\) on positiivinen kokonaisluku, niin
määritellään karteesinen tulo
%% Tark. termi
\begin{displaymath}
  A^n := \{ (a_1,\ldots,a_n) \setsep a_1,\ldots,a_n \in A \} .
\end{displaymath}
Myös vektorit \((a_1,\ldots,a_n)\) ovat järjestettyjä.
Kun \(A\) on joukko, niin määritellään
\begin{displaymath}
  \index{$\#$}
  \# A := \text{joukon } A \text{ alkioiden lukumäärä.}
\end{displaymath}
\(\# A\) voi olla ääretön.

\begin{definition}
\index{erilliset joukot}
Olkoon \(I\) (mahdollisesti ääretön) epätyhjä joukko. Olkoot \(A_i\),
\(i \in I\), joukkoja. Sanomme, että joukot \(A_i\) ovat erillisiä,
jos \(A_i \intersection A_j = \emptyset\) kaikilla \(i,j \in I\), \(i
\not= j\).
\end{definition}

\index{valinta-aksiooma}

\defterm{Valinta-aksiooma} on aksiooma, jonka mukaan minkä tahansa joukon \(A\)
potenssijoukosta \(\powerset{A}\) voidaan muodostaa uusi joukko
poimimalla kustakin joukon \(\powerset{A}\) alkiosta yksi alkio ja
muodostamalla niistä joukko. Edellisen esimerkin joukossa tällainen
joukko on esim. \(\{1,2\}\). 
Äärettömien joukkojen tapauksessa tilanne ei ole näin selvä, eikä
valinta-aksiooma kuulukaan kaikkiin joukko-opin aksioomajärjestelmiin.
\index{Banach-Tarskin paradoksi}
Valinta-aksioomasta seuraa \defterm{Banach-Tarskin paradoksi}
\cite{wpbanachtarski}: Jos on annettu täytetty pallo kolmiulotteisessa
avaruudessa, niin tämä pallo voidaan jakaa äärelliseen määrään
erillisiä joukkoja, joista voidaan sitten muodostaa kaksi alkuperäisen
pallon identtistä kopiota.

\section{Logiikkaa}

Logiikka operoi väitteillä, jotka ovat joko tosia tai
epätosia. Yleensä 1 tarkoittaa totta väitettä ja 0 epätotta. Loogiset
operaatiot ovat:
\begin{itemize}
  \item
    \index{looginen ei-operaatio}
    \(\neg P\) : ei \(P\)
  \item
    \index{looginen ja-operaatio}
    \(P \land Q\) : \(P\) ja \(Q\)
  \item
    \index{looginen tai-operaatio}
    \(P \lor Q\) : \(P\) tai \(Q\)
\end{itemize}

JA- ja TAI-operaatioiden totuustaulu on taulukossa \ref{tab:ja-tai}.
Huomaa, että \(P \lor Q\) on tosi myös silloin, kun sekä \(P\) että
\(Q\) ovat tosia. EI-operaation totuustaulu on taulukossa \ref{tab:ei}.

\begin{table}
\label{tab:ja-tai}
\caption{Loogisten JA- ja TAI-operaatioiden totuustaulu}
\begin{displaymath}
  \begin{array}{rrrr}
    P & Q & P \land Q & P \lor Q \\
    \hline
    0 & 0 & 0 & 0 \\
    0 & 1 & 0 & 1 \\
    1 & 0 & 0 & 1 \\
    1 & 1 & 1 & 1
  \end{array}
\end{displaymath}
\end{table}

\begin{table}
\label{tab:ei}
\caption{Loogisen EI-operaation totuustaulu}
\begin{displaymath}
  \begin{array}{rr}
    P & \neg P \\
    \hline
    0 & 1 \\
    1 & 0
  \end{array}
\end{displaymath}
\end{table}

Jokin väite \(P\) voi sisältää muuttujia, esim. ``\(x > 3\)''. Tällöin
väitettä voidaan merkitä \(P(x)\).

Merkintä \(\forall x : P(x)\) tarkoittaa, että väite \(P(x)\) on tosi
kaikille olioille \(x\). Merkintä \(\exists x : P(x)\) tarkoittaa,
että on olemassa ainakin yksi olio \(x\), jolle \(P(x)\) on tosi.
Merkintä \(\exists! x : P(x)\) tarkoittaa, että on olemassa tasan yksi
olio \(x\), jolle \(P(x)\) on tosi. Symbolia \(\forall\) sanotaan
\index{universaalikvanttori}
\defterm{universaalikvanttoriksi} ja symbolia \(\exists\)
\index{eksistenssikvanttori}
\defterm{eksistenssikvanttoriksi}.

Jos \(A\) on joukko, niin \(\forall x \in A : P(x)\) tarkoittaa, että
väite \(P(x)\) on tosi kaikille \(x \in A\) ja \(\exists x \in A :
P(x)\), että on olemassa ainakin yksi \(x \in A\), jolle \(P(x)\) on
tosi. Jos \(A = \emptyset\), niin \(\forall x \in A : P(x)\) on tosi
ja \(\exists x \in A : P(x)\) on epätosi riippumatta väitteestä \(P(x)\).
Merkintä \(\exists! x \in A : P(x)\) tarkoittaa, että on olemassa tasan yksi
olio \(x \in A\), jolle \(P(x)\) on tosi.

Jos \(P\) ja \(Q\) ovat väitteitä, niin \(P \iff Q\) tarkoittaa, että
\(P\) ja \(Q\) ovat yhtäpitäviä eli ekvivalentteja. Tämä tarkoittaa,
että \(P\) on tosi, jos ja vain jos \(Q\) on tosi. Merkintä \(P
\implies Q\) tarkoittaa, että väitteestä \(P\) seuraa väite \(Q\).
\index{ekvivalenssi}
Symboli ``\(\iff\)'' on nimeltään \defterm{ekvivalenssi} ja
\index{implikaatio}
ja symboli ``\(\implies\)'' \defterm{implikaatio}.

\section{Relaatiot ja funktiot}

Epätyhjien joukkojen \(A\) ja \(B\) karteesisen tulon \(A \times B\)
\index{relaatio}
osajoukkoa \(R\) kutsutaan joukkojen \(A\) ja \(B\)
väliseksi \defterm{relaatioksi}. Alkiot \(x \in A\) ja \(y \in B\)
toteuttavat relaation \(R\), jos ja vain jos \((x,y) \in R\). Jos
\(x\) ja \(y\) toteuttavat relaation \(R\), niin merkitään \(xRy\).
Jos \(A=B\) niin kutsumme relaatiota \(R\) relaatioksi joukossa \(A\).

Olkoot \(A\) ja \(B\) epätyhjiä joukkoja ja \(f\) niiden välinen
relaatio. Jos kullekin \(x \in A\) on olemassa tasan yksi \(y \in B\),
\index{funktio}
jolle \(xRy\), niin sanotaan, että \(f\) on \defterm{funktio} ja
merkitään \(f(x) = y\). Tällöin funktiota voidaan merkitä \(f: A \to
B\). Jos \(f: A \to B\) on funktio, niin sitä voidaan merkitä
\begin{displaymath}
  x \in A \mapsto f(x) ,
\end{displaymath}
esim. \(x \in \realnumbers \mapsto 5x\).
\index{määrittelyjoukko}
Joukkoa \(A\) sanotaan funktion \(f\) \defterm{määrittelyjoukoksi} ja
\index{maalijoukko}
joukkoa \(B\) funktion \(f\) \defterm{maalijoukoksi}.
Kun \(C \subset A\), niin määritellään
\begin{displaymath}
  \setimage{f}{C} := \{ f(x) \setsep x \in C \} .
\end{displaymath}
Joukkoa \(\setimage{f}{A}\) sanotaan funktion \(f\)
\index{kuvajoukko}
\defterm{kuvajoukoksi}.  Sanomme, että funktio \(f\) on
\index{injektio}
\defterm{injektio}, jos se kuvaa eri alkiot aina eri alkioille eli \(x,
y \in A \;\;\textrm{ja}\;\; x \not= y \implies f(x) \not=
\index{surjektio}
f(y)\). Funktio \(f : A \to B\) on \defterm{surjektio}, jos jokaiselle
\(b \in B\) on olemassa \(a \in A\) siten, että \(f(a) =
b\). Surjektiota voidaan merkitä \(f : A \onto B\). Jos \(f\) on sekä
injektio että surjektio, niin sanomme että \(f\) on
\index{bijektio}
\defterm{bijektio}. Havainnollisesti bijektio tarkoittaa, että
joukkojen \(A\) ja \(B\) alkiot kuvautuvat toisilleen yksi yhteen.
Jos \(f: A \to B\) on surjektio, niin sanomme, että \(f\) on funktio
joukolta \(A\) \textit{joukolle} \(B\). Jos \(f: A \to B\) ei
välttämättä ole
surjektio, niin sanomme, että \(f\) on funktio joukolta \(A\)
\textit{joukkoon} \(B\).  Jos \(A\) ja \(B\) ovat epätyhjiä joukkoja,
niin määritellään
\begin{displaymath}
  A^B := \{ \textrm{kaikki funktiot joukolta}\; B \;\textrm{joukkoon} A \} .
\end{displaymath}
\index{identtinen funktio}
Kun \(X\) on epätyhjä joukko, niin määritellään \defterm{identtinen
  funktio} \(\id{X} : X \onto X\) asettamalla
\begin{displaymath}
  \id{X}(x) := x
\end{displaymath}
kaikille \(x \in X\).
Jos \(I\) ja \(A\) ovat epätyhjiä joukkoja ja \(f\) funktio joukolta
\(I\) joukkoon \(A\), niin voimme merkitä tätä funktiota
\((x_\alpha)_{\alpha \in I}\), missä \(x_\alpha = f(\alpha)\), jolloin
sanomme sitä \defterm{perheeksi}.

\begin{definition}
  Olkoot \(A\) ja \(B\) joukkoja ja \(f : A \to B\) funktio. Olkoon
  \index{alkukuva}
  \(C \subset B\). Joukon \(C\) \defterm{alkukuva} määritellään
  \(\inverseimage{f}{C} := \{ x \in A \setsep f(x) \in C \}\).
\end{definition}

\begin{definition}
  Olkoot \(A\) ja \(B\) joukkoja ja \(f : A \to B\) funktio.
  Jos jokaiselle \(y \in B\) on olemassa tasan yksi \(x \in A\) jolle
  \(f(x) = y\), niin määritellään funktion \(f\)
  \index{käänteisfunktio}
  \defterm{käänteisfunktio} \(f^{-1} : B \to A\) asettamalla \(x :=
  f^{-1}(y)\).
\end{definition}

\begin{remark}
  Alkukuva \(\inverseimage{f}{C}\) on aina määritelty, vaikka
  käänteisfunktio \(f^{-1}\) ei olisikaan.
\end{remark}

\begin{definition}
  Olkoon \(A\) joukko ja \(\leq\) relaatio joukossa \(A\)
  \index{osittainen järjestys}
  \index{järjestys}
  Sanomme, että \(\leq\) on \defterm{osittainen järjestys}, jos
  seuraavat aksioomat ovat voimassa:
  \begin{itemize}
    \item[(OJ1)] antisymmetrisyys: \(x \leq y \land y \leq x \implies x = y\)
    \item[(OJ2)] transitiivisuus: \(x \leq y \land y \leq z \implies x \leq z\)
    \item[(OJ3)] refleksiivisyys: \(x \leq x\)
  \end{itemize}
  kaikille \(x, y, z \in A\).
  Jos joukossa \(A\) on määritelty osittainen järjestys, niin sanomme,
  \index{osittain järjestetty joukko}
  että joukko \(A\) on \defterm{osittain järjestetty}.
\end{definition}

\begin{definition}
  Olkoon \(A\) joukko ja \(\leq\) relaatio joukossa \(A\)
  \index{totaalinen järjestys}
  \index{järjestys}
  Sanomme, että \(\leq\) on \defterm{totaalinen järjestys}, jos
  seuraavat aksioomat ovat voimassa:
  \begin{itemize}
    \item[(J1)] antisymmetrisyys: \(x \leq y \land y \leq x \implies x = y\)
    \item[(J2)] transitiivisuus: \(x \leq y \land y \leq z \implies x \leq z\)
    \item[(J3)] totaalisuus: \(x \leq y \lor y \leq x\)
  \end{itemize}
  kaikille \(x, y, z \in A\).
  Jos joukossa \(A\) on määritelty totaalinen järjestys, niin sanomme,
  \index{totaalisesti järjestetty joukko}
  että joukko \(A\) on \defterm{totaalisesti järjestetty}.
  Jos lisäksi jokaisella joukon \(A\) epätyhjällä osajoukolla on
  pienin (relaation \(\leq\) suhteen) alkio, niin sanomme, että \(A\) on
  \index{hyvin järjestetty joukko}
  \defterm{hyvin järjestetty}.
\end{definition}

\begin{definition}
  Olkoon \(A\) totaalisesti järjestetty joukko ja \(B \subset
  A\). Olkoon \(r \in A\).  Jos \(x \leq r\) kaikilla \(x \in B\),
  niin sanomme, että \(r\) on \index{yläraja} joukon \(B\)
  \defterm{yläraja}. Jos joukolla \(B\) on yläraja, niin
  \index{ylhäältä rajoitettu joukko} sanomme, että \(B\) on
  \defterm{ylhäältä rajoitettu}.  Jos \(r \leq x\) kaikilla \(x \in
  B\), niin sanomme, että \(r\) on \index{alaraja} joukon \(B\)
  \defterm{alaraja}. Jos joukolla \(B\) on alaraja, niin
  \index{alhaalta rajoitettu joukko} sanomme, että \(B\) on
  \defterm{alhaalta rajoitettu}.
\end{definition}

\begin{definition}
  Olkoon \(A\) totaalisesti järjestetty joukko ja \(B \subset A\).
  \index{$\sup$} Jos \(B\) on ylhäältä rajoitettu, määritellään \(\sup
  B\) joukon \(B\) pienimmäksi ylärajaksi, eli jokaiselle joukon \(B\)
  ylärajalle \(r \in A\) on \(\sup B \leq r\). Jos \(B\) ei ole
  ylhäältä rajoitettu, määritellään \(\sup B := \infty\).
  \index{$\inf$} Jos \(B\) on alhaalta rajoitettu, määritellään \(\inf
  B\) joukon \(B\) suurimmaksi alarajaksi, eli jokaiselle joukon \(B\)
  alarajalle \(q \in A\) on \(q \leq \inf B\). Jos \(B\) ei ole
  alhaalta rajoitettu, määritellään \(\inf B := -\infty\).
  Määritellään lisäksi \(\sup \emptyset := \infty\) ja \(\inf
  \emptyset := -\infty\).  Määritellään edelleen
  \begin{displaymath}
    \sup_{x \in A} f(x) := \sup \{ f(x) \setsep x \in A \}
  \end{displaymath}
  ja
  \begin{displaymath}
    \inf_{x \in A} f(x) := \inf \{ f(x) \setsep x \in A \} .
  \end{displaymath}
\end{definition}

\begin{definition}
  Olkoon \(A\) joukko ja \(\sim\) joukkojen \(A\) ja \(A\) välinen
  \index{ekvivalenssirelaatio}
  relaatio. Sanomme, että \(\sim\) on \defterm{ekvivalenssirelaatio}, jos
  seuraavat aksioomat ovat voimassa:
  \begin{itemize}
    \item[(E1)] refleksiivisyys: \(x \sim x\)
    \item[(E2)] symmetrisyys: \(x \sim y \iff y \sim x\)
    \item[(E3)] transitiivisuus: \(x \sim y \land y \sim
      z \implies x \sim z\)
  \end{itemize}
  kaikille \(x, y, z \in A\).
\end{definition}

\begin{definition}
  Olkoon \(A\) joukko, \(\sim\) ekvivalenssirelaation joukossa
  \(A\) ja \(x \in A\). Sanomme joukkoa \(\{y \in A \setsep x \sim
  \index{ekvivalenssiluokka}
  y\}\) alkion \(x\) määräämäksi \defterm{ekvivalenssiluokaksi} ja
  merkitsemme sitä \(\eqclass{\sim}{x}\).
\end{definition}

\begin{definition}
  Olkoon \(A\) joukko ja \(\sim\) ekvivalenssirelaation joukossa
  \(A\).
  \index{jäännösluokka-avaruus}
  Määritellään \defterm{jäännösluokka-avaruus}
  \begin{displaymath}
    \quotientspace{A}{\sim} := \{ \eqclass{\sim}{x} \setsep x \in A \}.
  \end{displaymath}
\end{definition}

\begin{theorem}
  \label{th:ekv-avaruus}
  %% L13052015-1
  Olkoon \(A\) joukko ja \(\sim\) ekvivalenssirelaatio joukossa
  \(A\).
  Nyt
  \begin{itemize}
    \item[(1)]
      \[\bigcup_{J \in \quotientspace{A}{\sim}} J = A\]
    \item[(2)]
      Joukon \(\quotientspace{A}{\sim}\) alkiot ovat erillisiä.
  \end{itemize}
\end{theorem}

\begin{proof}
  \mbox{ }
  \begin{itemize}
  \item[(1)]
    Jokainen \(x \in A\) kuuluu johonkin \(L \in \quotientspace{A}{\sim}\) ja
    ekvivalenssirelaation määritelmän nojalla \(L' \subset A\)
    kaikille \(L' \in \quotientspace{A}{\sim}\), joten
    \begin{displaymath}
      \bigcup_{J \in \quotientspace{A}{\sim}} J = A .
    \end{displaymath}
  \item[(2)]
    Nyt \(\quotientspace{A}{\sim} \; = \{ \eqclass{\sim}{x} \setsep x \in A\}\) ja
    \(\eqclass{\sim}{x} \not= \emptyset\) kaikille \(x \in A\).

    Oletetaan, että joukossa \(\quotientspace{A}{\sim}\) on tasan yksi alkio
    \(L_1\). Nyt \(\eqclass{\sim}{x} = L_1\) kaikille \(x \in A\),
    joten \(L_1 = A\). Täten (2) on tosi.

    Oletetaan seuraavassa, että joukossa \(\quotientspace{A}{\sim}\)
    on enemmän kuin yksi alkio. Olkoot \(M, N
    \in \quotientspace{A}{\sim}\) ja \(M \not= N\).  Nyt \(M =
    \eqclass{\sim}{x}\) ja \(N = \eqclass{\sim}{y}\) joillekin \(x,y
    \in A\), \(x \not= y\).

    Oletetaan, että \(M \intersection N \not= \emptyset\)
    (vastaoletus). Olkoon \(z \in M \intersection N\). Jos \(m \in
    M\), niin \(z \sim m\). Koska \(z \in N\), niin \(m \in N\). Siis
    \(M \subset N\). Jos \(n \in N\), niin \(z \sim n\). Koska \(z \in
    M\), niin \(n \in M\). Siis \(N \subset M\). Täten \(M = N\), mikä
    on ristiriidassa vastaoletuksen kanssa. Siis on oltava \(M
    \intersection N = \emptyset\).
  \end{itemize}
\end{proof}

\begin{definition}
  \cite{wpdirectedset}
  \index{suunnattu joukko}
  Olkoon \(A\) epätyhjä joukko ja \(\geq\) relaatio joukossa
  \(A\). Sanomme, että \((A, \geq)\) on \defterm{suunnattu joukko},
  jos
  \begin{itemize}
    \item[(1)] \(a \geq a\) kaikille \(a \in A\).
    \item[(2)] \(a \geq b \land b \geq c \implies a \geq c\)
      kaikille \(a, b, c \in A\).
    \item[(3)] Jos \(a, b \in A\), niin on olemassa \(c \in A\) siten,
      että \(c \geq a\) ja \(c \geq b\).
  \end{itemize}
\end{definition}

\begin{remark}
  Jos järjestysrelaatio on asiayhteydestä selvä, niin voimme puhua
  myös suunnatusta joukosta \(A\). Sovellamme samanlaista käytäntöä
  myös muihin matemaattisiin struktuureihin, kuten vektoriavaruuksiin
  ja topologisiin avaruuksiin.
\end{remark}

\section{Kardinaali- ja ordinaaliluvut}

Ks. \cite{wpcardinal,wpordinal}.

%% MIETI tyhjä joukko
\begin{definition}
  \index{mahtavuus}
  Olkoot \(A\) ja \(B\) joukkoja. Sanomme, että joukot \(A\) ja \(B\)
  ovat \defterm{yhtä mahtavia}, jos on olemassa bijektio joukolta
  \(A\) joukolle \(B\).
\end{definition}

\begin{definition}
  Olkoon \(A\) epätyhjä joukko. Sanomme, että \(A\) on
  \index{ääretön joukko}
  \defterm{ääretön}, jos on olemassa bijektio joukolta \(A\) jollekin
  sen aidolle osajoukolle. Muussa tapauksessa sanomme, että joukko on
  \index{äärellinen joukko}
  \defterm{äärellinen}.
  Sanomme myös tyhjää joukkoa äärelliseksi.
\end{definition}

Äärelliset joukot ovat yhtä mahtavia, jos ja vain jos niissä on yhtä
monta alkiota. Jos joukon \(A\) mahtavuus on sama kuin luonnollisten
lukujen joukon \(\naturalnumbers\), niin sanomme, että \(A\) on
\index{numeroituva joukko}
\defterm{numeroituva}.
\index{korkeintaan numeroituva joukko}
Sanomme, että \(A\) on \defterm{korkeintaan numeroituva}, jos \(A\) on
äärellinen tai numeroituva. Jos \(A\) on ääretön mutta ei numeroituva,
\index{ylinumeroituva joukko}
niin sanomme, että \(A\) on \defterm{ylinumeroituva}.

\index{kardinaliteetti}
Kardinaaliluvut kuvaavat joukkojen mahtavuuksia. Äärellisen joukon
kardinaliteetti on sama kuin sen alkioiden lukumäärä. Kaikkien
äärettömien joukkojen kardinaliteetit eivät kuitenkaan ole samat. Näin
on esim. luonnollisille luvuille ja reaaliluvuille.
Luonnollisten lukujen joukon \(\naturalnumbers\) mahtavuutta merkitään
\(\aleph_0\) ja reaalilukujen joukon mahtavuus on \(2^{\aleph_0}\),
ks. luku \ref{ch:lukujoukot}.

\hyphenation{järjestys-isomorfisia}
Jos joukossa \(S\) on määritelty osittainen järjestys \(<\) ja
joukossa \(S'\) osittainen järjestys \(<'\), niin sanomme, että joukot
\index{järjestysisomorfisuus}
\((S, <)\) ja \((S', <')\) ovat \defterm{järjestysisomorfisia}, jos on
olemassa bijektio \(f : S \onto S'\), joka säilyttää järjestyksen, eli
\(f(a) <' f(b)\), jos ja vain jos \(a < b\).

\hyphenation{kardinaali-luvut}

Jokainen hyvin järjestetty joukko \((S, <)\) on järjestysisomorfinen
joukon kanssa, joka koostuu kaikista ordinaaliluvuista vähemmän kuin
eräs tietty ordinaaliluku.  Sanomme tätä ordinaalilukua joukon \((S,
<)\) järjestystyypiksi.
Jokainen ordinaaliluku \(\alpha\) on järjestystyyppi joukolle, joka
koostuu kaikista ordinaaliluvuista vähemmän kuin \(\alpha\).
Äärelliset ordinaaliluvut ovat samat kuin luonnolliset luvut. Pienin
ääretön ordinaaliluku on \(\omega\), joka samaistetaan
kardinaalilukuun \(\aleph_0\). Kuitenkin transfiniittisessa
tapauksessa ordinaaliluvut erottelevat joukkoja tarkemmin kuin
kardinaaliluvut.

%%\hyphenation{ekvivalenssi-luokkina}
\hyphenation{ordinaali-luvulla}

Ordinaaliluvut määriteltiin alunperin hyvin järjestettyjen joukkojen
ekvivalenssiluokkina. Tästä määrittelystä täytyy kuitenkin luopua
\index{Zermelo-Fraenkelin aksioomajärjestelmä}
joukko-opin Zermelo-Fraenkel ja siihen liittyvissä
aksioomajärjestelmissä, koska nämä ekvivalenssiluokat ovat liian
laajoja ollakseen joukkoja. Sen sijaan, että määrittelisimme
ordinaaliluvun hyvin järjestettyjen joukkojen ekvivalenssiluokkana,
määrittelemme sen yhtenä tiettynä hyvin määriteltynä joukkona, joka
edustaa tätä luokkaa. Siten ordinaaliluku on eräs hyvin määritelty
joukko, ja jokainen hyvin järjestetty joukko on järjestysisomorfinen
tasan yhden ordinaaliluvun kanssa.

John von Neumannin esittämä standardi määritelmä on:
Jokainen ordinaaliluku on sitä pienempien ordinaalilukujen muodostama
hyvin järjestetty joukko. Muodollisemmin tämä kuuluu:
Joukko \(S\) on ordinaaliluku, jos ja vain jos \(S\) on aidosti hyvin
järjestetty joukon jäsenyysrelaation suhteen ja jokainen joukon \(S\)
alkio on myös joukon \(S\) osajoukko.
Luonnolliset luvut ovat ordinaalilukuja tämän määritelmän
mukaan. Esim. 2 on luvun \(4 = \{0, 1, 2, 3\}\) alkio ja \(2 = \{0,
1\}\) on joukon \(\{0, 1, 2, 3\}\) osajoukko.

Jokaisella nollasta poikkeavalla ordinaaliluvulla on pienin alkio
\(0\). Ordinaaliluvulla joko on tai ei ole suurinta alkiota. Jos
ordinaaliluvulla \(\beta\) on maksimi \(\alpha\), niin ordinaalilukua
\index{seuraaja}
\index{seuraajaordinaali}
\(\beta\) kutsutaan seuraajaordinaaliksi (ordinaaliluvun \(\alpha\)
seuraajaksi) ja sitä merkitään \(\alpha + 1\). Von Neumannin
määritelmän mukaan ordinaaliluvun \(\alpha\) seuraaja on \(\alpha \cup
\{\alpha\}\). Jos ordinaaliluku \(\alpha\) ei ole minkään toisen
ordinaaliluvun seuraaja, niin ordinaalilukua \(\alpha\) kutsutaan
\index{rajaordinaali}
\defterm{rajaordinaaliksi}.

Jos oletetaan, että valinta-aksiooma on voimassa, niin joukon \(X\)
kardinaliteetti on pienin ordinaaliluku \(\alpha\) siten, että on
olemassa bijektio joukolta \(X\) joukolle \(\alpha\). Tätä määrittelyä
kutsutaan von Neumannin kardinaalilukujen määrittelyksi.

\hyphenation{kompleksi-lukujen}

\begin{examples}
  \begin{example}
    Joukko \(\{1, 2, 3\}\) on äärellinen.
  \end{example}
  \begin{example}
    Luonnollisten lukujen joukko \(\naturalnumbers\) on
    numeroituvuuden määritelmän nojalla numeroituva.
  \end{example}
  \begin{example}
    Kokonaislukujen joukko \(\integernumbers\) ja rationaalilukujen
    joukko \(\rationalnumbers\) ovat numeroituvia.
  \end{example}
  \begin{example}
    Reaalilukujen joukko \(\realnumbers\), irrationaalilukujen joukko
    \(\realnumbers \setminus \rationalnumbers\) ja kompleksilukujen
    joukko \(\complexnumbers\) ovat ylinumeroituvia.
  \end{example}
\end{examples}

\begin{definition}
  Olkoot \(A\) ja \(B\) joukkoja. Sanomme, että joukon \(A\) mahtavuus
  on \defterm{aidosti suurempi} kuin joukon \(B\), jos on olemassa
  bijektio joukolta \(B\) jollekin joukon \(A\) aidolle osajoukolle.
  Sanomme, että joukon \(A\) mahtavuus on \defterm{aidosti pienempi}
  kuin joukon \(B\) mahtavuus, jos on olemassa bijektio joukolta \(A\)
  jollekin joukon \(B\) aidolle osajoukolle.
\end{definition}

\index{kontinuumihypoteesi}

\defterm{Kontinuumihypoteesi} väittää, ettei ole olemassa joukkoa,
jonka mahtavuus on aidosti suurempi kuin joukon \(\naturalnumbers\) ja
aidosti pienempi kuin joukon \(\powerset{\naturalnumbers}\).

\begin{exercises}
  \begin{exercise}{1.1}
    %% L14052015-2
    Määritellään relaatio \(\sim\) joukossa \(\realnumbers\)
    asettamalla
    \begin{displaymath}
      x \sim y \iff x^2 = y^2 , \spaceafter x, y \in \realnumbers .
    \end{displaymath}
    \begin{enumerate}
      \item
        Osoita, että \(\sim\) on ekvivalenssirelaatio.
      \item
        Konstruoi jäännösluokka-avaruus
        \(\quotientspace{\realnumbers}{\sim}\).
      \item
        Osoita, että \(\sim\) ei ole funktio.
    \end{enumerate}
  \end{exercise}
\end{exercises}

\chapter{Ryhmät ja kunnat}
\label{ch:ryhmat-ja-kunnat}

\section{Ryhmät}

\begin{definition}
  Olkoon \(G\) epätyhjä joukko ja \(\circ : G \times G \to G\)
  \index{ryhmä}
  funktio. Sanomme, että \((G, \circ)\) on \defterm{ryhmä}, jos
  seuraavat aksioomat ovat voimassa:
  \begin{itemize}
    \item[(G1)]
      \(a \circ (b \circ c) = (a \circ b) \circ c\)
      kaikille \(a,b,c \in G\).
    \item[(G2)] On olemassa alkio \(e \in G\) siten, että \(e \circ a
      = a \circ e = a\) kaikille \(a \in G\).
    \item[(G3)]
      Jokaisella \(a \in G\) on olemassa alkio \(b \in G\) siten, että
      \(a \circ b = b \circ a = e\). Määritellään \(a^{-1} := b\).
  \end{itemize}
  Jos lisäksi on voimassa aksiooma
  \begin{itemize}
    \item[(G4)] \(a \circ b = b \circ a\) kaikille \(a, b \in G\).
  \end{itemize}
  \index{kommutatiivinen ryhmä}
  niin sanomme, että \((G, \circ)\) on \defterm{kommutatiivinen ryhmä}
  \index{Abelin ryhmä}
  eli \defterm{Abelin ryhmä}.
\end{definition}

Alkio \(e \in G\) on nimeltään \index{neutraalialkio}
\defterm{neutraalialkio}. Alkiota \(a^{-1}\) kutsutaan alkion \(a\)
\index{käänteisalkio}
\defterm{käänteisalkioksi}.  Toinen tapa sanoa,
että \((G, \circ)\) on ryhmä on, että \(G\) on ryhmä operaation
\(\circ\) suhteen. Alkioiden \(a\) ja \(b\) välistä ryhmäoperaatiota
voidaan merkitä myös \(ab\). Tämä on yleistä silloin, kun
ryhmäoperaationa on kertolasku.  Kommutatiivisen ryhmän
ryhmäoperaatiota merkitään usein ``+'' ja käänteisalkiota \(-a\).

\begin{figure}[htb]
\label{fig:vesimolekyyli}
\caption{Symmetriaryhmän $C_{2v}$ operaatiot vesimolekyylille.}
\includegraphics[scale=1.0]{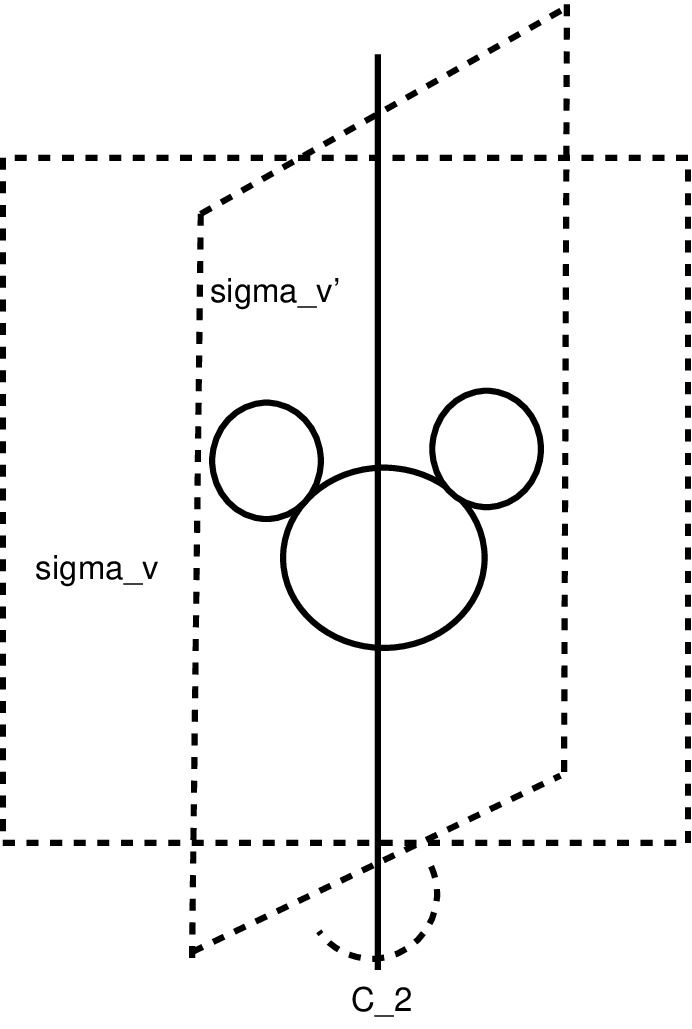}
\end{figure}

\begin{examples}
  \begin{example}
    \((\realnumbers, +)\) on ryhmä.
  \end{example}
  \begin{example}
    \((\realnumbers \setminus \{0\}, \cdot)\) on ryhmä
  \end{example}
  \begin{example}
    Vastaavasti \(\complexnumbers\) on ryhmä yhteenlaskun suhteen, ja
    \(\complexnumbers \setminus \{0\}\) on ryhmä kertolaskun suhteen.
  \end{example}
  \begin{example}
    \index{permutaatio}
    Bijektiota \(p : \{1,\ldots,n\} \onto \{1,\ldots,n\}\) sanotaan
    asteen \(n\) \defterm{permutaatioksi}.
    Permutaatiot tarkoittavat äärellisen joukon
    uudelleenjärjestelyjä. Esim. jos luvut 1, 2 ja 3 kuvataan
    luvuille 2, 3 ja 1 niin vastaavaa permutaatio on
    \begin{displaymath}
      \left(
      \begin{array}{ccc}
        1 & 2 & 3 \\
        2 & 3 & 1
      \end{array}
      \right)
    \end{displaymath}
    asteen \(n\) permutaatioiden neutraalialkio on
    \begin{displaymath}
      \left(
      \begin{array}{ccc}
        1 & \ldots & n \\
        1 & \ldots & n
      \end{array}
      \right)
    \end{displaymath}
    Permutaatioiden tulo \(a \circ b\) määritellään tavallisena
    kuvausten yhdistämisenä (ensin sovelletaan \(b\) ja sitten \(a\)).
    Esim.
    \begin{displaymath}
      \left(
      \begin{array}{ccc}
        1 & 2 & 3 \\
        2 & 3 & 1
      \end{array}
      \right)
      \circ
      \left(
      \begin{array}{ccc}
        1 & 2 & 3 \\
        2 & 1 & 3
      \end{array}
      \right)
      =
      \left(
      \begin{array}{ccc} 
        1 & 2 & 3 \\
        1 & 3 & 2
      \end{array}
      \right)
    \end{displaymath}
    \index{permutaatioryhmä}
    Edellä annetuilla määritelmillä kaikki \(n \in \naturalnumbers\), \(n >
    0\), alkion permutaatiot muodostavat ryhmän, jota merkitsemme \(P_n\).
  \end{example}
  \begin{example}
    \index{vesimolekyyli}
    Vesimolekyylin symmetriaryhmä on \(C_{2v}\) \cite{af2005}, ks. kuva
    \ref{fig:vesimolekyyli}.
    Ryhmän \(C_{2v}\) operaatiot ovat:
    \begin{enumerate}
      \item \(e\): ei tehdä mitään
      \item \(C_2\): \(180^\circ\) asteen kierto pääakselin ympäri
      \item \(\sigma_v\): peilaus sen tason suhteen, jossa molekyyli on
      \item \(\sigma_v'\): peilaus molekyylin läpi menevän tason suhteen
    \end{enumerate}
  \end{example}
\end{examples}

\begin{definition}
  Olkoon \((G, \circ)\) ryhmä. Jos \(A \subset G\) ja \(A\) on myös
  itse ryhmä, niin sanomme, että \(A\) on ryhmän \(G\)
  \index{aliryhmä}
  \defterm{aliryhmä}.
\end{definition}

\begin{theorem}
  \label{th:aliryhma}
  Olkoon \((G, \circ)\) ryhmä ja \(A \subset G\), \(A \not=
  \emptyset\). \(A\) on ryhmän \(G\) aliryhmä, jos ja vain jos \(a
  \circ b \in A\) kaikilla \(a, b \in A\) ja \(a^{-1} \in A\) kaikilla
  \(a \in A\).
\end{theorem}

Todistus harjoitustehtävänä.

\begin{examples}
  \begin{example}
    Kun \(A\) on ryhmä, niin \(\{e\}\) ja koko ryhmä \(A\) ovat ryhmän
    \(A\) aliryhmiä.
  \end{example}
  \begin{example}
    Joukko
    \begin{displaymath}
      B := \left\{
      \left(
      \begin{array}{ccc}
        1 & 2 & 3 \\
        1 & 2 & 3
      \end{array}
      \right) ,
      \left(
      \begin{array}{ccc}
        1 & 2 & 3 \\
        2 & 1 & 3
      \end{array}
      \right)
      \right\}
    \end{displaymath}
    on permutaatioryhmän \(P_3\) aliryhmä.
  \end{example}
\end{examples}

\begin{definition}
  Olkoon \((G,\circ)\) ryhmä, \(H\) ryhmän \(G\) aliryhmä ja \(a \in
  G\).
  \index{sivuluokka}
  \index{vasen sivuluokka}
  Määritellään aliryhmän \(H\) \defterm{vasen sivuluokka} asettamalla
  \begin{displaymath}
    a \circ H := \{ a \circ x \setsep x \in H \}
  \end{displaymath}
  \index{oikea sivuluokka}
  ja \defterm{oikea sivuluokka} asettamalla
  \begin{displaymath}
    H \circ a := \{ x \circ a \setsep x \in H \}.
  \end{displaymath}
\end{definition}

Kommutatiivisessa ryhmässä vasemmat ja oikeat sivuluokat ovat samat.

%% Sivuluokat jakavat ryhmän G erillisiksi joukoiksi

\begin{definition}
  Olkoot \((G,\circ)\) ja \((H,\star)\) ryhmiä ja \(f: G \to H\)
  \index{homomorfismi}
  funktio. Sanomme, että \(f\) on \defterm{homomorfismi} jos
  \begin{math}
    f(x \circ y) = f(x) \star f(y)
  \end{math}
  kaikille \(x, y \in G\).
  Jos lisäksi \(f\) on bijektio, niin sanomme, että \(f\) on
  \index{isomorfismi}
  \defterm{isomorfismi}.
\end{definition}

\section{Kunnat}

\begin{definition}
  \cite{wpfield}
  Olkoon \(K\) epätyhjä joukko ja \(+ : K \times K \to K\), \(\cdot :
  K \times K \to K\) funktioita. Sanomme, että \((K, +, \cdot)\) on
  \index{kunta}
  \defterm{kunta}, jos seuraavat aksioomat ovat voimassa:
  \begin{itemize}
    \item[(K1)] \(a + ( b + c ) = ( a + b ) + c\).
    \item[(K2)] \(a \cdot ( b \cdot c ) = ( a \cdot b ) \cdot c\).
    \item[(K3)] \(a+b=b+a\).
    \item[(K4)] \(a \cdot b = b \cdot a\).
    \item[(K5)] On olemassa \(0_K \in K\) siten, että \(0_K+a=a\) kaikille
      \(a \in K\).
    \item[(K6)] On olemassa \(1_K \in K\) siten, että \(1_K \cdot a = a\)
      kaikille \(a \in K\).
    \item[(K7)] \(0_K \not= 1_K\).
    \item[(K8)] Jokaiselle \(a \in K\) on olemassa \(-a \in K\) siten,
      että \(a+(-a)=0_K\).
    \item[(K9)] Kun \(a \in K \setminus \{0_K\}\), niin on olemassa
      \(a^{-1} \in K\) siten, että \(a \cdot a^{-1} = 1_K\).
    \item[(K10)] \(a \cdot (b+c) = (a \cdot b) + (a \cdot c)\).
  \end{itemize}
\end{definition}

\index{ykkösalkio}
Alkiota \(1_K\) sanotaan kunnan \(K\) \defterm{ykkösalkioksi} ja
\index{nolla-alkio}
alkiota \(0_K\) kunnan \(K\) \defterm{nolla-alkioksi}.

\begin{examples}
  \begin{example}
    Joukko \(\{0,1\} \subset \naturalnumbers\) on kunta. Tämä on
    pienin mahdollinen kunta.
  \end{example}
  \begin{example}
    Rationaalilukujen joukko \(\rationalnumbers\), reaalilukujen
    joukko \(\realnumbers\) ja kompleksilukujen joukko
    \(\complexnumbers\) ovat kuntia.
  \end{example}
\end{examples}

\begin{definition}
  Olkoon \(K\) kunta. Sanomme, että funktio \(\abs{\cdot} : K \to
  \realnumbers_0\), missä \(\realnumbers_0 := \{ r \in \realnumbers
  \index{itseisarvo}
  \setsep r \geq 0 \}\), on \defterm{itseisarvo}, jos
  \begin{itemize}
    \item[(i)] \(\abs{\lambda} = 0 \iff \lambda = 0\) ;
    \item[(ii)] \(\abs{\lambda + \mu} \leq \abs{\lambda} + \abs{\mu}\).
  \end{itemize}
\end{definition}

\begin{exercises}
  \begin{exercise}{2.1}
    Osoita, että ryhmässä on ainoastaan yksi neutraalialkio.
  \end{exercise}
  \begin{exercise}{2.2}
    Osoita, että kullakin ryhmän alkiolla on ainoastaan yksi
    käänteisalkio.
  \end{exercise}
  \begin{exercise}{2.3}
    %% L18052015-1
    Kun \(n \in \positiveintegers\), osoita että permutaatioryhmä
    \(P_n\) on ryhmä.
  \end{exercise}
  \begin{exercise}{2.4}
    Todista lause \ref{th:aliryhma}.
  \end{exercise}
\end{exercises}

\chapter{Lukujoukot}
\label{ch:lukujoukot}

\section{Luonnolliset luvut $\naturalnumbers$}

\index{luonnolliset luvut}
\index{$\naturalnumbers$}

Määritellään \(S(X) := X \union \{X\}\), missä \(X\) on joukko.

Luonnolliset luvut määritellään seuraavasti \cite{wpnaturalnumber}:
\begin{itemize}
\item
  \(0 := \emptyset \in \naturalnumbers\)
\item
  Jos \(k \in \naturalnumbers\), niin \(S(k) \in \naturalnumbers\).
\end{itemize}

Saadaan seuraavanlainen jono:

\(0 = \emptyset\)

\(1 = S(0) = S(\emptyset) = \{\emptyset\} = \{0\}\)

\(2 = S(1) = \{\emptyset, \{ \emptyset \} \} = \{0, 1\}\)

\(3 = S(2) = \{ \emptyset, \{ \emptyset \}, \{ \emptyset,
\{\emptyset\}\}\} = \{ 0, 1, 2 \}\)

\(\vdots\)

\vspace{0.5cm}

%%Huom! Kuitenkaan lukua \(0\) ei yleensä pidetä samana kuin tyhjää joukkoa
%%\(\emptyset\).

\index{seuraaja}
Operaatiota \(S(k)\) kutsutaan luvun \(k\) \defterm{seuraajaksi}.
Yhteenlasku määritellään rekursiivisesti:
\begin{itemize}
  \item \(a + 0 = a\) ;
  \item \(a + S(b) = S(a+b)\) ,
\end{itemize}
missä \(a, b \in \naturalnumbers\).
Kertolasku määritellään rekursiivisesti:
\begin{itemize}
  \item \(a \cdot 0 = 0\) ;
  \item \(a \cdot S(b) = (a \cdot b) + a\) ,
\end{itemize}
missä \(a, b \in \naturalnumbers\).
Kertolaskua \(a \cdot b\) voidaan merkitä myös \(ab\).
\index{totaalinen järjestys}
Luonnollisten lukujen joukossa määritellään totaalinen järjestys
asettamalla
\begin{displaymath}
  a \leq b \iff \exists c \in \naturalnumbers : a + c = b .
\end{displaymath}

Luonnolliset
luvut noudattavat Peanon aksioomajärjestelmää.

%% Peanon aksioomat tähän ?

\section{Kokonaisluvut $\integernumbers$}

Määritellään ensin ekvivalenssirelaatio luonnollisten lukujen
muodostamille (järjestetyille) pareille:
\begin{displaymath}
  (a,b) \sim (c,d) \iff a+d = c+b
\end{displaymath}
kaikille \((a,b), (c,d) \in \naturalnumbers \times \naturalnumbers\).
\index{kokonaisluvut}
\index{$\integernumbers$}
Kokonaislukujen joukko määritellään jäännösluokka-avaruutena
\begin{displaymath}
  \integernumbers := \cp{\naturalnumbers}{\naturalnumbers} / \sim .
\end{displaymath}
Kun \(a, b \in \naturalnumbers\), niin parin \((a, b)\) määräämää
ekvivalenssiluokkaa merkitään \(\eqc{(a, b)}\).
Määritellään kokonaislukujen yhteenlasku
\begin{displaymath}
  \eqc{(a, b)} + \eqc{(c, d)} := \eqc{(a + c, b + d)}
\end{displaymath}
kaikille \(a,b,c,d \in \naturalnumbers\),
\defterm{vastaluku}
\begin{displaymath}
  -\eqc{(a, b)} := \eqc{(b, a)}
\end{displaymath}
kaikille \(a, b \in \integernumbers\),
vähennyslasku
\begin{displaymath}
  \eqc{(a, b)} - \eqc{(c, d)} := \eqc{(a, b)} + (-\eqc{(c, d)})
\end{displaymath}
kaikille \(a,b,c,d \in \naturalnumbers\)
ja kertolasku
\begin{displaymath}
  \eqc{(a,b)} \cdot \eqc{(c,d)} := \eqc{(ac+bd,ad+bc)}
\end{displaymath}
kaikille \(a,b,c,d \in \naturalnumbers\).
Kokonaislukujen joukossa määritellään totaalinen järjestys asettamalla
\begin{displaymath}
  \eqc{(a,b)} \leq \eqc{(c,d)} \iff
  \eqc{(c+b,a+d)} = \eqc{(n, 0)} \; \text{jollekin}\; n \in
  \naturalnumbers .
\end{displaymath}
Luonnolliset luvut \(a \in \naturalnumbers\) samaistetaan
kokonaislukuihin \(\eqc{(a,0)}\).

\index{$\positiveintegers$}

Määritellään \(\positiveintegers := \{ x \in \integernumbers \setsep x
> 0\}\).

\section{Rationaaliluvut $\rationalnumbers$}

Määritellään ekvivalenssirelaatio \(\sim\) joukossa \(\integernumbers
\times \positiveintegers\) asettamalla
\begin{displaymath}
  (a,b) \sim (c,d) \iff ad = bc
\end{displaymath}
kaikille \(a,c \in \integernumbers\), \(b,d \in \positiveintegers\).
\index{rationaaliluvut}
\index{$\rationalnumbers$}
Määritellään rationaalilukujen joukko
\begin{displaymath}
  \rationalnumbers := \cp{\integernumbers}{\positiveintegers} / \sim .
\end{displaymath}
Kun \(a \in \integernumbers\) ja \(b \in \positiveintegers\), niin
parin \((a, b)\)
määräämää ekvivalenssiluokkaa merkitään \(\eqc{(a, b)}\).
Kunnan \(\rationalnumbers\) nolla-alkio on \(0_\rationalnumbers :=
\eqc{(0, 1)}\) ja
ykkösalkio \(1_\rationalnumbers := \eqc{(1, 1)}\).
Oletetaan seuraavassa, että \(\eqc{(a,b)}, \eqc{(c,d)} \in
\rationalnumbers\).
Määritellään rationaalilukujen yhteenlasku
\begin{displaymath}
  \eqc{(a,b)} + \eqc{(c,d)} := \eqc{(ad + bc, bd)} ,
\end{displaymath}
vastaluku
\begin{displaymath}
  -\eqc{(a,b)} := \eqc{(-a,b)} ,
\end{displaymath}
vähennyslasku
\begin{displaymath}
  \eqc{(a,b)} - \eqc{(c,d)} = \eqc{(a,b)} + (-\eqc{(c,d)})
\end{displaymath}
ja kertolasku
\begin{displaymath}
  \eqc{(a,b)} \cdot \eqc{(c,d)} := \eqc{(ac,bd)}
\end{displaymath}
Kun \(\eqc{(a,b)} \in \rationalnumbers\) ja \(\eqc{(a,b)} \not=
0_\rationalnumbers\), niin määritellään rationaaliluvun
\index{käänteisluku}
\(\eqc{(a,b)}\) \defterm{käänteisluku}
\begin{displaymath}
  \eqc{(a,b)}^{-1} := \eqc{(b,a)}
\end{displaymath}
ja jakolasku
\begin{displaymath}
  \eqc{(c,d)} / \eqc{(a,b)} :=
  \eqc{(c,d)} \cdot \eqc{(a,b)}^{-1} .
\end{displaymath}
Määritellään rationaalilukujen joukossa totaalinen järjestys
asettamalla
\begin{displaymath}
  \eqc{(a, b)} \leq \eqc{(c, d)} \iff ad-bc \leq 0
\end{displaymath}
Rationaalilukujen joukko \(\rationalnumbers\) on kunta.
Kokonaisluvut \(a \in \integernumbers\) samaistetaan
rationaalilukuihin \(\eqc{(a,1)}\).
Määritellään
\index{$\rationalnumbers_+$}
\begin{math}
  \rationalnumbers_+ := \{ q \in \rationalnumbers \setsep q > 0 \}
\end{math}.

\section{Reaaliluvut $\realnumbers$}

Reaaliluvut voidaan konstruoida rationaaliluvuista joko Cauchyn
jonoina tai Dedekindin leikkauksina. Seuraamme ensimmäistä tapaa.
Jonot on esitelty luvussa \ref{sec:j-v-f}.

\begin{definition}
  Olkoon \(\seq{a} := (a_k)_{k=0}^\infty \subset \rationalnumbers\). Sanomme,
  \index{rationaalinen Cauchyn jono}
  \index{Cauchyn jono}
  että \(\seq{a}\) on \defterm{rationaalinen Cauchyn jono}, jos
  \begin{displaymath}
    \forall \varepsilon \in \rationalnumbers_+ : \exists N \in
    \naturalnumbers : \forall n, m \in \naturalnumbers :
    n > N \land m > N \implies \abs{a_n - a_m} < \varepsilon .
  \end{displaymath}
\end{definition}

Määritellään
\begin{displaymath}
  A := \{ \seq{a} \setsep \seq{a} := (a_k)_{k=0}^\infty \subset
  \rationalnumbers \;\text{on rationaalinen Cauchyn
    jono} \}.
\end{displaymath}
Määritellään relaatio \(\sim\) joukossa \(A\) asettamalla
\begin{displaymath}
  \seq{a} \sim \seq{b} \iff \lim_{k \to \infty}
  \abs{\seqelem{\seq{a}}{k} - \seqelem{\seq{b}}{k}} = 0 . 
\end{displaymath}
Määritellään reaalilukujen joukko
\index{reaaliluvut}
\index{$\realnumbers$}
\begin{displaymath}
  \realnumbers := A / \sim .
\end{displaymath}
Kun \(\seq{a} \in A\), niin alkion \(\seq{a}\) määräämää
ekvivalenssiluokkaa merkitään \(\eqc{\seq{a}}\).

Määritellään
\begin{displaymath}
  0_\realnumbers := \eqc{(0)_{k=0}^\infty}
\end{displaymath}
ja
\begin{displaymath}
  1_\realnumbers := \eqc{(1)_{k=0}^\infty} .
\end{displaymath}
Määritellään reaaliluvun vastaluku
\begin{displaymath}
  -\eqc{\seq{a}} := \eqc{-\seq{a}}
\end{displaymath}
missä \(\eqc{\seq{a}} \in \realnumbers\)
ja yhteen- ja vähennyslasku
\begin{eqnarray*}
  \eqc{\seq{a}} + \eqc{\seq{b}} & := &
  \eqc{\seq{a} + \seq{b}} \\
  \eqc{\seq{a}} - \eqc{\seq{b}} & := &
  \eqc{\seq{a} - \seq{b}}
\end{eqnarray*}
missä \(\eqc{\seq{a}}, \eqc{\seq{b}} \in \realnumbers\).
Määritellään reaalilukujen kertolasku
\begin{displaymath}
  \eqc{\seq{a}} \cdot \eqc{\seq{b}}
  :=
  \left( \seqelem{\seq{a}}{k} \cdot \seqelem{\seq{b}}{k} \right)_{k=0}^\infty .
\end{displaymath}
Määritellään nollasta poikkeavan reaaliluvun käänteisluku seuraavasti:
Olkoon \(\eqc{\seq{a}} \in \realnumbers, \eqc{\seq{a}} \not= 0_\realnumbers\).
Nyt \(\lim_{k \to \infty} \seqelem{\seq{a}}{k} \not= 0\), joten on
olemassa luku \(k_0 \in \naturalnumbers\) siten, että
\(\seqelem{\seq{a}}{k} \not= 0\) kaikille \(k \in \naturalnumbers\),
\(k \geq k_0\). Asetetaan
\begin{displaymath}
  \left(\eqc{\seq{a}}\right)^{-1} :=
  \eqc{\left(\frac{1}{\seqelem{\seq{a}}{k_0+k}}\right)_{k=0}^\infty} .
\end{displaymath}
Määritellään reaalilukujen jakolasku
\begin{displaymath}
  \eqc{\seq{a}} / \eqc{\seq{b}} :=
  \eqc{\seq{a}} \cdot \left(\eqc{\seq{b}}\right)^{-1}
\end{displaymath}
missä \(\eqc{\seq{a}}, \eqc{\seq{b}} \in \realnumbers, \eqc{\seq{b}}
\not= 0_\realnumbers\).
Määritellään totaalinen järjestys reaalilukujen joukossa asettamalla
\(\eqc{\seq{a}} \leq \eqc{\seq{b}}\) jos ja vain jos
on olemassa \(k_0 \in \naturalnumbers\) siten, että
\(\seqelem{\seq{a}}{k} \leq \seqelem{\seq{b}}{k}\) kaikilla \(k \in
\naturalnumbers\), \(k \geq k_0\).

\begin{figure}[htb]
\label{fig:sup}
\caption{Täydellisyysaksiooma: jokaisella ylhäältä rajoitetulla
  joukolla $A \subset \realnumbers$ on pienin yläraja $\sup A \in
  \realnumbers$.}
\includegraphics{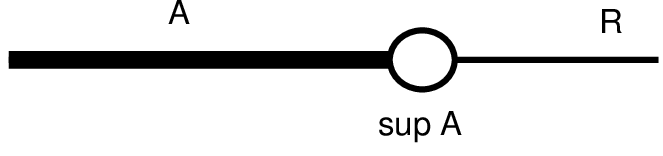}
\end{figure}

Reaalilukujen joukko \(\realnumbers\) on kunta ja totaalisesti järjestetty
joukko. Lisäksi reaalilukujen joukossa on voimassa
\index{täydellisyysaksiooma}
\defterm{täydellisyysaksiooma}:
Jokaisella ylhäältä rajoitetulla joukolla \(J \subset \realnumbers\)
on pienin yläraja. Täydellisyysaksioomaa on havainnollistettu kuvassa
\ref{fig:sup}.
Aksiooman todistus on seuraavassa \cite{wpconstrrealnumbers}:
Olkoon \(S\) jokin joukon \(\realnumbers\) epätyhjä osajoukko.
Jos kaikki joukon \(S\) alkiot ovat ei-negatiivisia, niin olkoon \(U\)
jokin joukon \(S\) rationaalinen yläraja
ja valitaan luku \(L \in \rationalnumbers\) siten, että \(L < s\)
jollekin \(s \in S\).
Määritellään jonot \((u_n)\) ja \((l_n)\)
seuraavasti:
\begin{itemize}
  \item Asetetaan \(u_0 := U\) ja \(l_0 := L\).
  \item Määritellään \(m_n := (u_n + l_n)/2\).
  \item Jos \(m_n\) on joukon \(S\) yläraja, asetetaan \(u_{n+1} :=
    m_n\) ja \(l_{n+1} := l_n\). Muuten asetetaan \(l_{n+1} := m_n\)
    ja \(u_{n+1} :=u_n\).
\end{itemize}
Näin määritellään kaksi rationaalista Cauchyn jonoa, joten saamme
reaaliluvut \(u := \eqc{(u_n)}\) ja \(l := \eqc{(l_n)}\). Induktiolla
voidaan osoittaa, että \(u_n\) on joukon \(S\) yläraja kaikille \(n\)
ja \(l_n\) ei ole joukon \(S\) yläraja millekään \(n\). Meillä on
\begin{displaymath}
  \lim_{n \to \infty} (u_n-l_n) = 0 ,
\end{displaymath}
joten \(l = u\). Oletetaan, että \(b < u = l\) olisi pienempi yläraja
joukolle \(S\) kuin \(u\). Koska \((l_n)\) on kasvava, niin \(b <
l_k \) jollekin \(k\). Mutta \(l_n\) ei ole joukon \(S\) yläraja ja
siten ei myöskään \(b\) ole yläraja. Täten \(u\) on joukon \(S\)
pienin yläraja.

Huomaa, että täydellisyysaksiooma ei ole voimassa rationaalilukujen
joukossa: esim. joukolla
\begin{math}
  \{ x \in \rationalnumbers \setsep x < \sqrt{2} \}
\end{math}
ei ole pienintä ylärajaa rationaalilukujen joukossa.

\index{$\realnumbers_+$}
\index{$\realnumbers_0$}
Määritellään \(\realnumbers_+ := \{ r \in \realnumbers \setsep r > 0
\}\) ja \(\realnumbers_0 := \{ r \in \realnumbers \setsep r \geq 0
\index{laajennettu reaalilukusuora}
\index{$\extrealnumbers$}
\}\). Määritellään laajennettu reaalilukusuora \(\extrealnumbers :=
\realnumbers \union \{ -\infty, +\infty \}\).
%% laajennetun reaalilukusuoran laskutoimitukset

\begin{definition}
  Kun \(a, b \in \extrealnumbers\) ja \(a \leq b\), niin määritellään
  reaalilukuvälit
  \begin{eqnarray*}
    \intervalcc{a}{b} & := & \{ x \in \extrealnumbers \setsep a \leq x \leq b \} \\  
    \intervalco{a}{b} & := & \{ x \in \extrealnumbers \setsep a \leq x < b \} \\  
    \intervaloc{a}{b} & := & \{ x \in \extrealnumbers \setsep a < x \leq b \} \\  
    \intervaloo{a}{b} & := & \{ x \in \extrealnumbers \setsep a < x < b \}
  \end{eqnarray*}
  Väliä \([a,b]\) kutsutaan \defterm{suljetuksi väliksi} ja väliä
  \(]a,b[\) \defterm{avoimeksi väliksi}.
\end{definition}

\section{Kompleksiluvut $\complexnumbers$}

Reaalilukujen joukossa ei voi ratkaista yhtälöä
\begin{displaymath}
  x^2 = -1 .
\end{displaymath}
Tästä nousee esiin tarve lukualueen laajentamiselle reaaliluvuista, ja
näin syntyvät kompleksiluvut.

\index{kompleksiluvut}
\index{$\complexnumbers$}
Määritellään kompleksiluvut reaalilukupareina:
\begin{displaymath}
  \complexnumbers := \{ (x,y) \setsep x, y \in \realnumbers \} .
\end{displaymath}
Kun \(z := (x,y) \in \complexnumbers\), niin lukua \(x\) sanotaan
\index{reaaliosa}
kompleksiluvun \(z\) \defterm{reaaliosaksi} ja lukua \(y\)
\index{imaginaariosa}
kompleksiluvun \(z\) \defterm{imaginaariosaksi}. Määritellään
\begin{math}
  \ReAlt z := x
\end{math}
ja
\begin{math}
  \ImAlt z := y
\end{math}.
Määritellään kompleksiluvuille
\begin{itemize}
  \item yhteenlasku: \((a,b)+(c,d) := (a+c,b+d)\)
  \item vastaluku: \(-(a,b) := (-a,-b)\)
  \item vähennyslasku: \((a,b)-(c,d) := (a-b,c-d)\)
  \item kertolasku: \((a,b) \cdot (c,d) := (ac-bd,ad+bc)\)
\end{itemize}
kaikille \((a,b), (c,d) \in \complexnumbers\).
Kun \(((a,b), (c,d) \in \complexnumbers\) ja \((c,d) \not= (0,0)\),
niin määritellään kompleksiluvun käänteisluku
\begin{displaymath}
  (c,d)^{-1} := \frac{(c,-d)}{c^2+d^2}
\end{displaymath}
ja kompleksilukujen jakolasku
\begin{displaymath}
  \frac{(a,b)}{(c,d)} := \frac{(ac+bd,bc-ad)}{c^2+d^2} .
\end{displaymath}

Kompleksilukujen joukko \(\complexnumbers\) on kunta, jonka ykkösalkio
on \((1,0)\) ja nolla-alkio \((0,0)\). Kompleksilukujen joukko ei
kuitenkaan ole totaalisesti järjestetty joukko.

\index{imaginaariyksikkö}
Määritellään \defterm{imaginaariyksikkö} \(\imagunit := (0,1)\).
Nyt jokainen kompleksiluku \((a,b)\) voidaan kirjoittaa muodossa
\(a + \imagunit b\), missä \(a,b \in \realnumbers\).
%% Tod. HT
Kun \(z := a + \imagunit b \in \complexnumbers\), niin määritellään
\index{liittoluku}
\index{konjugaatti}
kompleksiluvun \(z\) \defterm{liittoluku} eli \defterm{konjugaatti}
\(z^* := a - \imagunit b\).

Tarkastellaan kompleksilukuja \((a,0)\) ja \((b,0)\), \(a,b \in
\realnumbers\).
Nyt
\begin{eqnarray*}
  (a,0)+(b,0) & = & (a+b,0) \\
  -(a,0) & = & (-a,0) \\
  (a,0) \cdot (b,0) & = & (ab,0) \\
  (a,0)^{-1} & = & \left( \frac{1}{a}, 0 \right), \;\;\;\; \textrm{kun}
  \; a \not= 0
\end{eqnarray*}
Siis kompleksiluvut \((a,0)\), \(a \in \realnumbers\), voidaan
samaistaa reaalilukuihin \(a\). Täten kompleksilukujen kunta on
reaalilukujen kunnan kuntalaajennus. Voidaan osoittaa, että
kompleksilukujen kunta on laajin kunta, joka sisältää reaaliluvut.

\begin{definition}
  \label{def:joukko-opillinen-kantaja}
  Olkoon \(X\) epätyhjä joukko ja \(f : X \to
  \complexnumbers\). Määrittelemme funktion \(f\)
  \index{joukko-opillinen kantaja}
  \index{kantaja}
  \defterm{joukko-opillisen kantajan} asettamalla
  \begin{displaymath}
    \stsupp f := \{ x \in X \setsep f(x) \not= 0 \}.
  \end{displaymath}
\end{definition}

Vertaa topologisen kantajan määritelmään
\ref{def:topologinen-kantaja}.

\begin{exercises}
  \begin{exercise}{3.1}
    %% L16052015-1
    Osoita, että
    \begin{displaymath}
      \left(\frac{1}{k+1}\right)_{k=0}^\infty
    \end{displaymath}
    on rationaalinen Cauchyn jono.
  \end{exercise}
  \begin{exercise}{3.2}
    Osoita, että jokainen kompleksiluku \((a, b)\), missä \(a, b \in
    \realnumbers\), voidaan esittää muodossa \(a + \imagunit b\).
  \end{exercise}
\end{exercises}

%%\input Jonot-ja-verkot.tex
%% yhdistetty lukuun Topologiset avaruudet

\chapter{Vektoriavaruudet}
\label{ch:vektoriavaruudet}

\section{Matriisit}

\begin{definition}
  Määritellään
  \index{Kroneckerin delta}
  \defterm{Kroneckerin delta}
  \begin{displaymath}
    \delta_{i,j} :=
    \left\{
    \begin{array}{ll}
      1; \;\;\;\; i = j \\
      0; \;\;\;\; i \not= j
    \end{array}
    \right.
  \end{displaymath}
  kaikille \(i,j \in \naturalnumbers\).
\end{definition}

\begin{definition}
  Kun \(n \in \positiveintegers\), niin määritellään
  \begin{displaymath}
    Z(n) := \{ k \in \positiveintegers \setsep k \leq n \} .
  \end{displaymath}
\end{definition}

\begin{definition}
  Olkoot \(n, m \in \positiveintegers\).
  Kutsumme funktiota \(A: Z(n) \times Z(m) \to \complexnumbers\)
  \index{matriisi}
  \defterm{matriisiksi}. Matriisia voidaan merkitä
  \((a_{i,j})_{i=1,j=1}^{n,m}\), missä \(a_{i,j} = A(i, j)\) kaikille
  \(i \in Z(n)\) ja \(j \in Z(m)\). Funktiota \(A\) kutsutaan myös
  \(n\)-riviseksi ja \(m\)-sarakkeiseksi matriisiksi.
  \index{vaakavektori}
  Jos \(n=1\), niin sanomme, että \(A\) on \defterm{vaakavektori}.
  \index{pystyvektori}
  Jos \(m=1\), niin sanomme, että \(A\) on \defterm{pystyvektori}.
  \index{skalaari}
  Jos \(n=m=1\), niin sanomme, että \(A\) on \defterm{skalaari}.
  \index{neliömatriisi}
  Jos \(n=m\), niin sanomme, että \(A\) on \defterm{neliömatriisi}.
\end{definition}

Pystyvektoria voidaan merkitä myös \((a_k)_{k=1}^n\). Kompleksiset (ja
vastaavasti reaaliset) \(n\) alkion pystyvektorit voidaan samaistaa
karteesisen tulon \(\complexnumbers^n\) (ja vastaavasti
\(\realnumbers^n\)) kanssa. Vaakavektorit voidaan samaistaa
karteesisiin tuloihin samalla lailla.

\begin{definition}
  Olkoon \(n \in \positiveintegers\). Määritellään
  \index{yksikkövektori}
  \defterm{yksikkövektorit}
  \begin{displaymath}
    \unitvec{n}{k} := \left( \delta_{j,k} \right)_{j=1}^n ,
  \end{displaymath}
  missä \(k \in Z(n)\).
\end{definition}

\begin{definition}
  Olkoon \(n \in \positiveintegers\). Määritellään 
  \index{yksikkömatriisi}
  \defterm{yksikkömatriisi}
  \begin{displaymath}
    I_n := \left( \delta_{i,j} \right)_{i=1,j=1}^{n,n} .
  \end{displaymath}
\end{definition}

\begin{definition}
  Olkoot \(n, m \in \positiveintegers\). Määritellään
  \index{nollamatriisi}
  \defterm{nollamatriisi}
  \begin{displaymath}
    0_{n \times m} := \left( 0 \right)_{i=1,j=1}^{n,m} .
  \end{displaymath}
\end{definition}

Matriisi voidaan konstruoida kirjoittamalla sen alkiot sulkeiden sisälle,
esim.
\begin{displaymath}
  A := (a_{i,j})_{i=1,j=1}^{3,4} :=
  \left(
  \begin{array}{rrrr}
    1 & 7 & 6.5 & 0 \\
    3 & 9 & -\imagunit & 0 \\
    5 & -5 & 1 & 2+4\imagunit
  \end{array}
  \right)
\end{displaymath}
Ensimmäinen indeksi \(i\) ilmoittaa alkion rivin ja toinen indeksi
\(j\) sen sarakkeen. Esim. \(a_{1,2} = 7\) ja \(a_{3,2} = -5\).
Esimerkki vaakavektorista on
\begin{displaymath}
  \left(
  \begin{array}{rrrr}
    1 & 5 & 8 & -10
  \end{array}
  \right)
\end{displaymath}
ja pystyvektorista
\begin{displaymath}
  \left(
  \begin{array}{r}
    5 \\
    -1 \\
    0 \\
    8.5 \\
    4
  \end{array}
  \right) .
\end{displaymath}

Matriisien \((a_{i,j})_{i=1,j=1}^{n,m}\) ja
\index{matriisien summa}
\((b_{i,j})_{i=1,j=1}^{n,m}\) \defterm{summa} määritellään
\begin{displaymath}
  (a_{i,j})_{i=1,j=1}^{n,m} + (b_{i,j})_{i=1,j=1}^{n,m} := (a_{i,j} +
  b_{i,j})_{i=1,j=1}^{n,m}
\end{displaymath}
\index{matriisien erotus}
ja \defterm{erotus}
\begin{displaymath}
  (a_{i,j})_{i=1,j=1}^{n,m} - (b_{i,j})_{i=1,j=1}^{n,m} := (a_{i,j} -
  b_{i,j})_{i=1,j=1}^{n,m} .
\end{displaymath}
\index{vastamatriisi}
Matriisin \((a_{i,j})_{i=1,j=1}^{n,m}\) \defterm{vastamatriisi}
määritellään
\begin{displaymath}
  -(a_{i,j})_{i=1,j=1}^{n,m} := (-a_{i,j})_{i=1,j=1}^{n,m} .
\end{displaymath}
Matriisien \(A := (a_{i,j})_{i=1,j=1}^{n,m}\) ja
\index{matriisitulo}
\(B := (b_{i,j})_{i=1,j=1}^{m,p}\) \defterm{matriisitulo} määritellään
\begin{displaymath}
  AB := \left( \sum_{k=1}^m a_{i,k} b_{k,j} \right)_{i=1,j=1}^{n,p} .
\end{displaymath}

\begin{definition}
  Olkoon \(n \in \positiveintegers\) ja \(A \in
  \complexnumbers^{Z(n) \times Z(n)}\). Jos on olemassa neliömatriisi
  \(B \in \complexnumbers^{Z(n) \times Z(n)}\) siten, että
  \(AB = I_n\) niin sanomme, että \(B\) on matriisin \(A\)
  \index{käänteismatriisi}
  \defterm{käänteismatriisi}, ja merkitsemme \(A^{-1} := B\).
\end{definition}

\begin{remark}
  Käänteismatriisia ei välttämättä ole olemassa.
\end{remark}

\begin{definition}
  Olkoon \(n \in \positiveintegers\) ja \(A \in
  \complexnumbers^{Z(n) \times Z(n)}\). Jos matriisin \(A\)
  käänteismatriisi on olemassa, niin sanomme, että \(A\) on
  \index{ei-singulaarinen matriisi}
  \defterm{ei-singulaarinen}. Jos käänteismatriisia ei ole olemassa,
  \index{singulaarinen matriisi}
  niin sanomme, että \(A\) on \defterm{singulaarinen}.
\end{definition}

\begin{definition}
  \label{def:matr-tr-hk}
  Olkoot \(n, m \in \positiveintegers\) ja \(A \in
  \complexnumbers^{Z(n) \times Z(m)}\). Matriisin \(A\)
  \index{transpoosi}
  \index{matriisin transpoosi}
  \defterm{transpoosi} \(A^T \in \complexnumbers^{Z(m) \times Z(n)}\)
  määritellään asettamalla \(A^T(i,j) := A(j,i)\) kaikille \(i \in
  \index{Hermiten konjugaatti}
  \index{matriisin Hermiten konjugaatti}
  Z(m)\) ja \(j \in Z(n)\). Matriisin \(A\) \defterm{Hermiten
    konjugaatti} \(A^\herm \in \complexnumbers^{Z(m) \times Z(n)}\)
  määritellään asettamalla \(A^\herm (i,j) := A(j,i)^*\) kaikille \(i \in
  Z(m)\) ja \(j \in Z(n)\).
\end{definition}

Vertaa vastaavaan määritelmään lineaarisille kuvauksille,
\ref{def:transpoosi-ja-h-konjugaatti}.

\begin{definition}
  \label{def:matriisilajeja}
  Olkoon \(n \in \positiveintegers\) ja \(A \in
  \complexnumbers^{Z(n) \times Z(n)}\). Sanomme, että
  \begin{itemize}
    \index{symmetrinen matriisi}
  \item \(A\) on \defterm{symmetrinen}, jos \(A = A^T\).
    \index{hermiittinen matriisi}
  \item \(A\) on \defterm{hermiittinen}, jos \(A = A^\herm\).
    \index{ortogonaalinen matriisi}
  \item \(A\) on \defterm{ortogonaalinen}, jos \(A^{-1} = A^T\).
    \index{unitaarinen matriisi}
  \item \(A\) on \defterm{unitaarinen}, jos \(A^{-1} = A^\herm\).
  \end{itemize}
\end{definition}

Vertaa vastaavaan määritelmään lineaarisille funktioille,
\ref{def:lin-kuv-lajeja}.

%% Lause: jokainen hermiittinen matriisi on unitaarisesti
%% diagonalisoituva

\index{ominaisarvoyhtälö matriiseille}

\begin{definition}
  \label{def:matriisien-oay}
  Olkoon \(n \in \positiveintegers\) ja \(A \in
  \complexnumbers^{Z(n) \times Z(n)}\). 
  Jos
  \begin{displaymath}
    A \vect{v} = \lambda \vect{v}
  \end{displaymath}
  missä \(\vect{v} \in \complexnumbers^{Z(n) \times Z(1)}\) on
  pystyvektori, \(\vect{v} \not= 0_{n \times 1}\) ja \(\lambda \in
  \complexnumbers\), niin sanomme, että
  \index{ominaisarvo}
  \(\lambda\) on matriisin \(A\) \defterm{ominaisarvo} ja \(\vect{v}\)
  \index{ominaisvektori}
  sitä vastaava \defterm{ominaisvektori}.
\end{definition}

Vertaa määritelmään \ref{def:oay}.

\begin{theorem}
  Olkoon \(n \in \positiveintegers\) ja \(A \in
  \complexnumbers^{Z(n) \times Z(n)}\). Tällöin matriisi \(A\) on
  singulaarinen, jos ja vain jos sillä on ominaisarvo 0.
\end{theorem}

%% TOD.

\section{Vektoriavaruudet}

\begin{definition}
  Olkoon \(V\) epätyhjä joukko, \(K\) kunta, ja \(+ : V \times V \to
  V\) ja \(\cdot : K \times V \to V\) funktioita. Tuloa \(c \cdot v\)
  merkitään myös \(cv\).
  Sanomme, että \((V,K,+,\cdot)\)
  \index{vektoriavaruus}
  \index{kerroinkunta}
  on \defterm{vektoriavaruus} ja \(K\) sen \defterm{kerroinkunta},
  jos seuraavat aksioomat ovat voimassa:
  \begin{itemize}
    \item[(V1)] \(\forall v,w,z \in V : v + ( w + z ) = ( v + w ) +
      z\)
    \item[(V2)] \(\forall v, w \in V : v+w=w+v\)
    \item[(V3)] \(\exists 0_V \in V : \forall v \in V : v + 0_V = 0_V
      + v = v\)
    \item[(V4)] \(\forall x \in V : \exists y \in V : x + y = 0_V\)
    \item[(V5)] \(\forall a \in K, v, w \in V : a ( v + w ) = av + aw\)
    \item[(V6)] \(\forall a, b \in K, x \in V: (a + b) x = ax + bx\)
    \item[(V7)] \(\forall a, b \in K, x \in V: a(bx) = (ab)x\)
    \item[(V8)] \(\forall v \in V: 1_K v = v\)
  \end{itemize}
\end{definition}

\index{vektori}
Vektoriavaruuden alkioita kutsutaan \defterm{vektoreiksi}.
\index{lineaarinen avaruus}
Vektoriavaruudesta käytetään myös nimeä \defterm{lineaarinen
  avaruus}. Nollavektoria \(0_V\) merkitään yleensä \(0\).
\index{vektorien summa}
\index{summa}
Operaatiota \(+\) kutsutaan vektorien \defterm{summaksi} ja
\index{skalaarin ja vektorin tulo}
\index{tulo}
operaatiota \(\cdot\) skalaarin ja vektorin \defterm{tuloksi}.
\index{vastavektori}
Aksiooman (V4) mukaista alkiota \(y\) kutsutaan vektorin \(x\)
\defterm{vastavektoriksi} ja sitä merkitään \(-x\).

\begin{examples}
  \begin{example}
    Reaalilukujen joukko \(\realnumbers\) ja kompleksilukujen joukko
    \(\complexnumbers\) varustettuna normaalilla summalla ja tulolla
    ovat vektoriavaruuksia. Ne ovat itse itsensä kerroinkuntia.
  \end{example}
  \begin{example}
    Joukko \(\{0\}\) on vektoriavaruus.
  \end{example}
  \begin{example}
    Kun \(n,m \in \positiveintegers\), niin kaikki \(n\)-riviset ja
    \(m\)-sarakkeiset kompleksiset (ja vastaavasti reaaliset) matriisit
    muodostavat vektoriavaruuden kerroinkunnalla \(\complexnumbers\)
    (vastaavasti \(\realnumbers\)).
  \end{example}
  \begin{example}
    Karteesinen tulo \(K^n\), missä \(n \in \positiveintegers\) ja \(K
    = \realnumbers\) tai \(K = \complexnumbers\) on
    vektoriavaruus. Summa määritellään alkioittain:
    \begin{displaymath}
      \left(
      \begin{array}{c}
        a_1 \\
        \vdots \\
        a_n
      \end{array}
      \right)
      +
      \left(
      \begin{array}{c}
        b_1 \\
        \vdots \\
        b_n
      \end{array}
      \right)
      =
      \left(
      \begin{array}{c}
        a_1 + b_1 \\
        \vdots \\
        a_n + b_n
      \end{array}
      \right) ,
      \;\;\;\;
      \left(
      \begin{array}{c}
        a_1 \\
        \vdots \\
        a_n
      \end{array}
      \right) ,
      \left(
      \begin{array}{c}
        b_1 \\
        \vdots \\
        b_n
      \end{array}
      \right)
      \in K^n
    \end{displaymath}
    ja skalaarin ja vektorin tulo:
    \begin{displaymath}
      c
      \left(
      \begin{array}{c}
        a_1 \\
        \vdots \\
        a_n
      \end{array}
      \right)
      =
      \left(
      \begin{array}{c}
        ca_1 \\
        \vdots \\
        ca_n
      \end{array}
      \right) ,
      \;\;\;\;
      c \in K,
      \left(
      \begin{array}{c}
        a_1 \\
        \vdots \\
        a_n
      \end{array}
      \right)
      \in K^n .
    \end{displaymath}
  \end{example}
  \begin{example}
    Olkoon \(V := \realnumbers^\realnumbers\) kaikkien
    reaalifunktioiden joukko. Määritellään funktioiden yhteenlasku
    asettamalla
    \begin{displaymath}
      (f+g)(x) := f(x) + g(x)
    \end{displaymath}
    kaikille \(f,g \in V, x \in \realnumbers\)
    ja skalaarin ja funktion tulo asettamalla
    \begin{displaymath}
      (c \cdot f)(x) := (cf)(x) := cf(x)
    \end{displaymath}
    kaikille \(f \in V, c, x \in \realnumbers\).
    Nyt \((V,\realnumbers,+,\cdot)\) on vektoriavaruus.
  \end{example}
\end{examples}

\begin{theorem}
  \label{th:vekt-av}
  Olkoon \(V\) vektoriavaruus kerroinkunnalla \(K\) ja \(v \in
  V\). Nyt
  \begin{itemize}
  \item[(i)] \(0_K v = 0_V\) ;
  \item[(ii)] \((-1)v = -v\).
  \end{itemize}
\end{theorem}

Todistus harjoitustehtävänä.

\begin{definition}
  Olkoon \(V\) vektoriavaruus kerroinkunnalla \(K\) ja \(L \subset V\)
  (mahdollisesti ääretön) joukko avaruuden \(V\) vektoreita. Sanomme,
  \index{lineaarinen riippumattomuus}
  että \(L\) on \defterm{lineaarisesti riippumaton}, jos
  \begin{displaymath}
    \sum_{k=1}^n a_kx_k = 0 \implies \forall k \in Z(n) : a_k = 0
  \end{displaymath}
  kaikille joukon \(L\) äärellisille osajoukoille
  \(\{x_1,\ldots,x_n\}\) ja kaikille skalaareille \(a_k \in K\), \(k
  \in Z(n)\).
  Jos \(L\) ei ole lineaarisesti riippumaton, niin sanomme, että se on
  \index{lineaarinen riippuvuus}
  \defterm{lineaarisesti riippuva}.
\end{definition}

\begin{definition}
  Olkoon \(V\) vektoriavaruus kerroinkunnalla \(K\), \(n \in
  \positiveintegers\), \(x_1,\ldots,x_n \in V\) ja \(a_1,\ldots,a_n
  \in K\). Lauseketta
  \begin{displaymath}
    \sum_{k=1}^n a_k x_k
  \end{displaymath} 
  \index{lineaarikombinaatio}
  sanotaan \defterm{lineaarikombinaatioksi}.
\end{definition}

\begin{definition}
  Olkoon \(V\) vektoriavaruus kerroinkunnalla \(K\) ja \(A \subset
  V\). Sanomme, että \(A\) on vektoriavaruuden \(V\)
  \index{vektorialiavaruus}
  \defterm{vektorialiavaruus}, jos \(A\) on itse vektoriavaruus.
  Tällöin merkitään \(A \vectorsubspace V\).
\end{definition}

\begin{theorem}
  \label{th:vektorialiavaruus}
  %% L13052015-4
  Olkoon \(V\) vektoriavaruus kerroinkunnalla \(K\) ja \(A \subset
  V\), \(A \not= \emptyset\). Joukko \(A\) on vektoriavaruuden \(V\)
  vektorialiavaruus jos ja vain jos seuraavat ehdot ovat voimassa:
  \begin{itemize}
    \item[(A1)] \(\forall x,y \in A: x + y \in A\)
    \item[(A2)] \(\forall x \in A, a \in K: ax \in A\)
  \end{itemize}
\end{theorem}

%% TARKISTA
\begin{proof}
  Jos \(A \vectorsubspace V\), niin \(A\) on vektoriavaruus, joten se
  on suljettu summan ja skalaarilla kertomisen suhteen.

  Oletetaan että lauseen ehdot (A1) ja (A2) ovat voimassa. Siis \(A\)
  on suljettu summan ja skalaarilla kertomisen
  suhteen. Vektoriavaruuden aksioomat (V1), (V2), (V5), (V6), (V7) ja (V8)
  ovat tosia, koska ne ovat tosia koko avaruudessa \(V\) mukaanlukien
  joukon \(A\).  Olkoon \(a \in A\) mielivaltainen. Nyt \(0_V = a - a
  \in A\). Täten aksiooma (V3) on voimassa. Aksiooma (V4) on voimassa
  ehdon (A2) lauseen \ref{th:vekt-av} kohdan (ii) nojalla. Siis \(A\)
  on vektoriavaruus ja \(A \vectorsubspace V\).
\end{proof}

\begin{examples}
  \begin{example}
    Olkoon \(V = \realnumbers^2\) (kaikkien tasovektoreiden joukko),
    kerroinkunta \(\realnumbers\), ja \(\rnvec{v} \in V\). Asetetaan
    \begin{displaymath}
      A := \{ a \rnvec{v} \setsep a \in \realnumbers \} .
    \end{displaymath}
    Nyt \(A\) on vektoriavaruuden \(V\) vektorialiavaruus.
  \end{example}
\end{examples}

\begin{definition}
  Kun \(V\) on vektoriavaruus kerroinkunnalla \(K\) ja \(A\)
  avaruuden \(V\) vektorialiavaruus,
  \index{jäännösluokka}
  niin määritellään \defterm{jäännösluokka} 
  \begin{displaymath}
    v + A := \{ v + a \setsep a \in A \}, \spaceafter v \in V .
  \end{displaymath}
  Määritellään
  %% vai sivuluokka-avaruus?
  \index{jäännösluokka-avaruus}
  \defterm{jäännösluokka-avaruus} asettamalla
  \begin{displaymath}
    \vquotientspace{V}{A} := \{ v + A \setsep v \in V\}
  \end{displaymath}
  ja laskutoimitukset
  \begin{eqnarray}
    \label{eq:quotient-sum}
    (u + A) + (v + A) & = & (u + v) + A , \\
    \label{eq:quotient-product}
    \alpha (u + A) & = & \alpha u + A ,
  \end{eqnarray}
  missä \(u, v \in V\) ja \(\alpha \in K\).
  Olkoon \(0_{\vquotientspace{V}{A}} := 0 + A = A\)
 \end{definition}

\begin{theorem}
  %% L13052015-3
  Olkoon \(V\) on vektoriavaruus kerroinkunnalla \(K\) ja \(A\)
  avaruuden \(V\) vektorialiavaruus. Nyt
  \begin{itemize}
    \item[(1)]
      \begin{displaymath}
        \bigcup_{J \in \vquotientspace{V}{A}} J = V .
      \end{displaymath}
    \item[(2)]
      Joukon \(\vquotientspace{V}{A}\) alkiot ovat erillisiä.
    \item[(3)]
      Laskutoimitukset \eqref{eq:quotient-sum} ja
      \eqref{eq:quotient-product} ovat hyvin määriteltyjä.
    \item[(4)] Joukko \(\vquotientspace{V}{A}\) varustettuna
      laskutoimituksilla \eqref{eq:quotient-sum} ja
      \eqref{eq:quotient-product} on vektoriavaruus, jonka
      kerroinkunta on \(K\).
  \end{itemize}
\end{theorem}

\begin{proof}
  Vertaa lauseen \ref{th:ekv-avaruus} todistukseen.
  \item[(1)]
    Kun \(x \in V\), niin \(x + A \in \vquotientspace{V}{A}\) ja \(x +
    A \subset V\), joten
    \begin{displaymath}
      \bigcup_{J \in \vquotientspace{V}{A}} J = V
    \end{displaymath}
  \item[(2)]
    Jos \(\quotientspace{V}{A} = \{ A \}\) (eli \(A = V\)), niin (2)
    on tosi. Oletetaan, että joukossa \(\quotientspace{V}{A}\) on
    enemmän kuin yksi alkio. Olkoon \(M := x + A\), \(x \in V\), \(N
    := y + A\), \(y \in V\) ja \(M \not= N\).
    
    Oletamme nyt, että olisi \(M \intersection N \not= \emptyset\)
    (vastaoletus). Olkoon \(z \in M \intersection N\). Nyt \(z = x +
    a_1 = y + a_2\) joillekin \(a_1, a_2 \in A\).  Jos \(m \in M\),
    niin \(m = x + a_3\) eräälle \(a_3 \in A\).  Olkoon \(a_4 := a_2 -
    a_1 + a_3 \in A\).  Nyt \(y + a_4 = y + a_2 - a_1 + a_3 = x + a_3
    = m\), joten \(m = y + a_4 \in y + A = N\).  Täten \(M \subset
    N\). Vastaavasti osoitetaan, että \(N \subset M\), joten \(M =
    N\). Tämä on ristiriita, joten vastaoletus on väärä, ja \(M
    \intersection N = \emptyset\).
  \item[(3)]
    \mbox{ }
    \begin{itemize}
      \item[(3.1)]
        Osoitetaan, että \((u+v)+A = (u'+v')+A\), kun \(u + A = u' +
        A\) ja \(v + A = v' + A\), \(u,v,u',v' \in V\).
        Oletetaan, että \(p \in (u+v)+A\). Nyt \(p = u+v+a\) eräälle
        \(a \in A\). Edelleen \(u+a=u'+a'\) ja \(v+a=v'+a''\), missä
        \(a', a'' \in A\). Nyt
        \begin{eqnarray*}
          p & = & u + v + a = u' + a' - a + v' + a'' - a = u' + v' + a' +
          a'' - 2a \\
          & \in & (u'+v')+A .
        \end{eqnarray*}
        Siis \((u+v)+A \subset (u'+v')+A\). Vastaavasti osoitetaan,
        että \((u'+v')+A \subset (u+v)+A\). Täten
        \begin{displaymath}
          (u+v)+A = (u'+v')+A .
        \end{displaymath}
      \item[(3.2)]
        Osoitetaan, että \(\alpha u + A = \alpha u' + A\), kun
        \(u+A=u'+A\), \(u,u' \in V\).
        Olkoon \(p = \alpha u + b\) eräälle \(b \in A\). Nyt
        \(u+b=u'+b'\) eräälle \(b' \in A\). Edelleen \(u = u' + b' -
        b\). Siten
        \begin{displaymath}
          \alpha u + a = \alpha u' + \alpha b' - \alpha b + b = \alpha
          u' + c \in \alpha u' + A ,
        \end{displaymath}
        missä \(c := \alpha b' - \alpha b + b \in A\).
        Siten \(\alpha u + A \subset \alpha u' + A\). Vastaavasti
        osoitetaan, että \(\alpha u' + A \subset \alpha u + A\).
        Täten \(\alpha u + A = \alpha u' + A\).
    \end{itemize}
  \item[(4)]
    Osoitetaan, että \(\vquotientspace{V}{A}\) on vektoriavaruus. Kun
    \(x+A,y+A \in \vquotientspace{V}{A}\), niin
    \begin{displaymath}
      (x+y)+A = \{ x+y+a \setsep a \in A \} \in \vquotientspace{V}{A} .
    \end{displaymath}
    Kun \(\alpha \in K\) ja \(x \in \vquotientspace{V}{A}\), niin
    \begin{displaymath}
      \alpha x + A =  \{ \alpha x + a \setsep a \in A \}
      \in \vquotientspace{V}{A} .
    \end{displaymath}
    \begin{itemize}
      \item[(V1):] \((a+A)+((b+A)+(c+A)) = (a+b+c)+A = ((a+A) +
        (b+A)) + (c+A)\),
      \item[(V2):] \((a+A)+(b+A) = (a+b)+A = (b+a) + A = (b+A)+(a+A)\),
      \item[(V3):] \((a+A) + 0_{\vquotientspace{V}{A}} = (a+A) + (0_K +
        A) = (a+0_K)+A = a+A\),
      \item[(V4):] \((a+A) + ((-a)+A) = (a+(-a)) + A = 0 + A = A =
        0_{\vquotientspace{V}{A}}\)
      \item[(V5):] \(c((a+A)+(b+A)) = (ca+cb)+A = c(a+A)+c(b+A)\),
      \item[(V6):] \((c + d)(a + A) = (ca + da) + A\).
      \item[(V7):] \(c(d(a+A)) = c(da + A) = (cd)(a+A)  \).
      \item[(V8):] \(1(a + A) = a + A\).
    \end{itemize}
    missä \(a,b \in V\) ja \(c,d \in K\).
    Täten \(\vquotientspace{V}{A}\) on vektoriavaruus.
\end{proof}

\begin{figure}[htb]
\label{fig:jaannosluokat}
\caption{Suora $A := \{ (x,ax) \setsep x \in \realnumbers \}$ on
  vektoriavaruuden $\realnumbers^2$ aliavaruus. Sen jäännösluokat ovat
  suorat $b + A = \{(x,ax+b) \setsep x \in \realnumbers \}$.}
\includegraphics{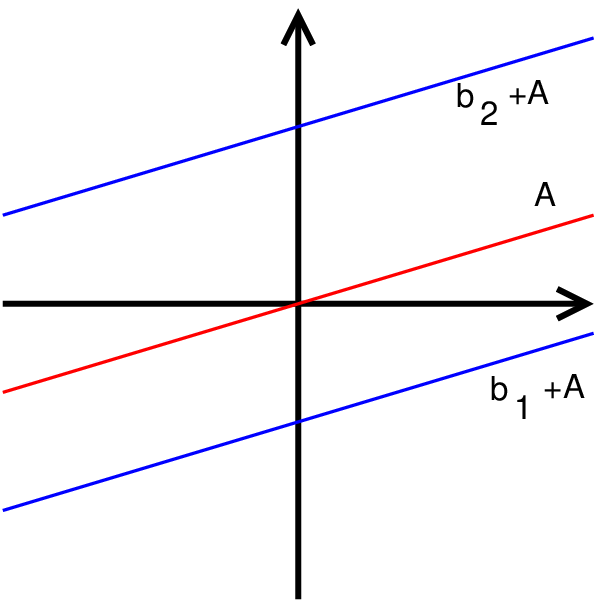}
\end{figure}

Jäännösluokkia on havainnollistettu kuvassa \ref{fig:jaannosluokat}.

%% \begin{definition}
%%   Olkoon \(V\) vektoriavaruus ja \(A\) avaruuden \(V\)
%%   vektorialiavaruus.
%%   \index{jäännösluokkafunktio}
%%   Määritellään \defterm{jäännösluokkafunktio} \(\pi : V \to V / A\)
%%   asettamalla \(\pi(x) = x + A\) kaikille \(x \in V\).
%% \end{definition}

\begin{definition}
  Olkoon \(V\) vektoriavaruus kerroinkunnalla \(K\) ja \(A \subset
  V\), \(A \not= \emptyset\). Olkoon
  \begin{displaymath}
    S := \{ B \setsep B \vectorsubspace V \land A \subset B \} .
  \end{displaymath}
  Määritellään
  \begin{displaymath}
    \spanop A := \bigcap_{J \in S} J .
  \end{displaymath}
\end{definition}

\hyphenation{vektori-avaruudeksi}

%% HT: os. että \span A on vektoriavaruus
\index{joukon virittämä vektoriavaruus}
Nyt \(\spanop A\) on vektoriavaruus, jota kutsutaan joukon \(A\)
\defterm{virittämäksi vektoriavaruudeksi}.

\begin{theorem}
  Olkoon \(V\) vektoriavaruus kerroinkunnalla \(K\) ja \(A \subset
  V\), \(A \not= \emptyset\). Nyt
  \begin{equation}
    \label{eq:span-formula}
    \spanop A = \left\{ \sum_{k=1}^n a_k x_k \bigsetsep n \in
    \positiveintegers, x_k \in V, a_k \in K, k \in Z(n) \right\} .
  \end{equation}
\end{theorem}

\begin{remark}
  Lineaarikombinaatio yhtälössä \eqref{eq:span-formula} on aina
  äärellinen vaikka joukko \(A\) olisi ääretön.
\end{remark}

\begin{definition}
  Olkoon \(V\) vektoriavaruus ja \(L \subset V\) lineaarisesti
  riippumaton. Jos \(\spanop L = V\), niin sanomme, että \(L\) on
  \index{kanta}
  \index{Hamelin kanta}
  vektoriavaruuden \(V\) \defterm{(Hamelin) kanta}.
\end{definition}

Voidaan osoittaa, että jokaisella vektoriavaruudella on kanta ja
kaikkien tietyn vektoriavaruuden kantojen mahtavuus on sama.
Jos vektoriavaruudella on äärellinen kanta, niin sanomme, että
\index{äärellisulotteinen vektoriavaruus}
vektoriavaruus on \defterm{äärellisulotteinen}. Muussa tapauksessa
\index{ääretönulotteinen vektoriavaruus}
sanomme, että vektoriavaruus on \defterm{ääretönulotteinen}.

\index{lineaarinen funktio}

\begin{definition}
  Olkoot \(V\) ja \(W\) vektoriavaruuksia kerroinkunnalla \(K\) ja
  olkoon \(f : V \to W\) funktio. Sanomme, että \(f\) on
  \defterm{lineaarinen}, jos
  \begin{itemize}
    \item[(1)] \(f(x+y) = f(x) + f(y)\)
    \item[(2)] \(f(ax) = af(x)\)
  \end{itemize}
  kaikille \(x, y \in V\) ja \(a \in K\).
\end{definition}

\index{matriisi lineaarisena funktiona}

\begin{remark}
  \(n\) rivin ja \(m\) sarakkeen matriisi kunnassa \(K\) voidaan tulkita
  lineaariseksi funktioksi vektoriavaruudelta \(K^m\)
  vektoriavaruuteen \(K^n\) asettamalla
  \begin{displaymath}
    A(v) := Av , \;\;\;\; v \in K^n ,
  \end{displaymath}
  missä yhtälön oikea puoli on matriisitulo ja \(v\) tulkitaan \(n\)
  alkion pystyvektoriksi.
\end{remark}

\begin{definition}
  Olkoon \(V\) vektoriavaruus kerroinkunnalla \(K\).
  \index{duaali}
  \index{algebrallinen duaali}
  Määritellään vektoriavaruuden \(V\) \defterm{algebrallinen duaali}
  asettamalla
  \begin{displaymath}
    V^\sharp := \{ f : V \to K \setsep f \;\textrm{on lineaarinen}\} .
  \end{displaymath}
\end{definition}

\begin{definition}
  Olkoot \(V\) ja \(W\) vektoriavaruuksia kerroinkunnalla
  \(K\). Olkoon \(f : V \onto W\) lineaarinen surjektio.
%% ???
  \index{homomorfismi}
  Sanomme, että \(f\) on \defterm{homomorfismi} avaruudelta \(V\)
  avaruudelle \(W\) ja avaruudet \(V\) ja \(W\) ovat
  \index{homomorfiset vektoriavaruudet}
  \defterm{homomorfiset}.
  Jos lisäksi \(f\) on bijektio, niin sanomme, että \(f\) on
  \index{isomorfismi}
  \index{algebrallinen isomorfismi}
  \index{isomorfiset vektoriavaruudet}
  \defterm{(algebrallinen) isomorfismi} ja avaruudet \(V\) ja \(W\)
  ovat \defterm{isomorfisia}.
\end{definition}

Isomorfisten vektoriavaruuksien rakenne on sama.

\begin{definition}
  Olkoon \(V\) vektoriavaruus ja \(M\) ja \(N\) sen
  aliavaruuksia. Määritellään aliavaruuksien \(M\) ja \(N\)
  \index{summa}
  \index{vektorialiavaruuksien summa}
  \defterm{summa}
  \begin{displaymath}
    M + N := \{ m + n \setsep, m \in M, n \in N \} .
  \end{displaymath}
  Jos \(M \intersection N = \zeroset\), niin sanomme, että \(M + N\)
  \index{suora summa}
  on \defterm{suora summa}, ja sitä merkitään \(M \directsum N\).
\end{definition}

\begin{theorem}
  %% L14052015-1
  Olkoon \(V\) vektoriavaruus ja \(M\) ja \(N\) sen
  aliavaruuksia. Summa \(M + N\) on suora, jos ja vain jos jokaisen
  vektorin \(x \in M + N\) esitys \(x = m + n\), \(m \in M\), \(n \in
  N\), on yksikäsitteinen.
\end{theorem}

\begin{proof}
  \mbox{ }
  \begin{enumerate}
    \item
      Oletetaan ensin, että \(M \intersection N = \zeroset\).
      Olkoon \(x \in M + N\), \(x = m + n\), missä \(m \in M\) ja \(n
      \in N\). Oletetaan, että olisi \(x = m' + n'\), \(m' \in M\),
      \(n' \in N\) ja \(m' \not= m\) tai \(n' \not= n\). Nyt \(x = m +
      n = m' + n'\), joten \(m - m' = n' - n \not= 0\). Toisaalta
      \(m - m' = n' - n \in M \intersection N\). Siis \(M
      \intersection N \not= \zeroset\). Tämä on ristiriita, joten
      esitys \(x = m + n\) on yksikäsitteinen.
    \item
      Oletetaan, että jokaisen
      vektorin \(x \in M + N\) esitys \(x = m + n\), \(m \in M\), \(n \in
      N\), on yksikäsitteinen.
      Oletetaan, että olisi \(y \in M + N\) ja \(y \not= 0\). Olkoon
      \(x \in M + N\) eli \(x = m + n\), \(m \in M\) ja \(n \in
      N\). Nyt \(x = m + n = (m + y) + (n - y)\), missä \(m + y \in
      M\) ja \(n - y \in N\). Siten esitys \(x = m + n\) ei ole
      yksikäsitteinen. Tämä on ristiriita, joten on oltava \(M
      \intersection N = \zeroset\).
  \end{enumerate}
\end{proof}

\begin{definition}
  Olkoot \(V\) ja \(W\) vektoriavaruuksia kerroinkunnalla
  \(K\). Olkoon \(f : V \to W\) lineaarinen funktio. Määritellään
  \index{ydin}
  \index{$\ker$}
  funktion \(f\) \defterm{ydin}
  \begin{displaymath}
    \ker f := \{ x \in V \setsep f(x) = 0 \}
  \end{displaymath}
  \index{kuva-avaruus}
  \index{$\image$}
  ja \defterm{kuva-avaruus}
  \begin{displaymath}
    \image f := \setimage{f}{V} .
  \end{displaymath}
\end{definition}

\index{ominaisarvoyhtälö yleisessä vektoriavaruudessa}

\begin{definition}
  \label{def:oay}
  Olkoon \(V\) vektoriavaruus kerroinkunnalla \(K\) ja \(f: V \to V\)
  lineaarinen funktio. Jos
  \begin{displaymath}
    f(v) = \lambda v
  \end{displaymath}
  missä \(v \in V\), \(v \not= 0\) ja \(\lambda \in K\), niin sanomme, että
  \index{ominaisarvo}
  \(\lambda\) on funktion \(f\) \defterm{ominaisarvo} ja \(v\)
  \index{ominaisvektori}
  sitä vastaava \defterm{ominaisvektori}.
\end{definition}

Vertaa määritelmään \ref{def:matriisien-oay}. Kun matriisit tulkitaan
lineaarisiksi kuvauksiksi, niin molemmat määritelmät ovat ekvivalentteja.

\begin{definition}
  Olkoon \(V\) vektoriavaruus kerroinkunnalla \(K = \realnumbers\) tai
  \(K = \complexnumbers\). Olkoon \(A \subset V\).
  \index{konveksi joukko}
  Sanomme, että \(A\) on \defterm{konveksi}, jos
  \begin{math}
    tA + (1-t)A \subset A
  \end{math}
  jokaiselle \(t \in \intervalcc{0}{1}\).
  \index{balansoitu joukko}
  Sanomme, että \(A\) on \defterm{balansoitu}, jos
  \(\alpha A \subset A\)
  jokaiselle \(\alpha \in K\).
\end{definition}

\begin{definition}
  Olkoon \(V\) vektoriavaruus ei-diskreetillä kerroinkunnalla \(K\) ja
  \index{absorboida}
  \(A, B \subset V\). Sanomme, että \(A\) \defterm{absorboi} joukon
  \(B\), jos on olemassa \(\lambda_0 \in K\) siten, että \(B \subset
  \lambda A\) kaikille \(\abs{\lambda} \geq \abs{\lambda_0}\).
  \index{radiaalinen}
  \index{absorboiva}
  Olkoon \(U \subset V\). Sanomme, että \(U\) on \defterm{radiaalinen
    (absorboiva)}, jos \(U\) absorboi jokaisen avaruuden \(V\)
  \index{rengastettu}
  äärellisen osajoukon. Sanomme, että \(U\) on \defterm{rengastettu},
  jos \(\lambda U \subset U\) kaikille \(\abs{\lambda} \leq 1\).
\end{definition}

\begin{exercises}
  \begin{exercise}{4.1}
    Laske matriisitulo
    \begin{displaymath}
      \begin{pmatrix}
        9 & 6 & 7 \\
        8 & -5 & 4 \\
        0 & -1 & 2
      \end{pmatrix}
      \begin{pmatrix}
        1 \\
        -4 \\
        -5
      \end{pmatrix} .
    \end{displaymath}
  \end{exercise}
  \begin{exercise}{4.2}
    Laske matriisitulo
    \begin{displaymath}
      \begin{pmatrix}
        -10 & 7 & 5 & 8 \\
        5 & 7 & 6 & 9 \\
        0 & 8 & 7 & 4
      \end{pmatrix}
      \begin{pmatrix}
        1 & 10 \\
        -1 & 0 \\
        1 & 8 \\
        -1 & 9
      \end{pmatrix} .
    \end{displaymath}
  \end{exercise}
  \begin{exercise}{4.3}
    %% L13052015-5
    Olkoot \(V\) ja \(W\) vektoriavaruuksia kerroinkunnalla \(K\) ja
    \(f : V \to W\) lineaarinen funktio. Osoita, että \(\ker f\) on
    avaruuden \(V\) aliavaruus.
  \end{exercise}
  \begin{exercise}{4.4}
    %% L13052015-6
    Olkoot \(V\) ja \(W\) vektoriavaruuksia kerroinkunnalla \(K\) ja
    \(f : V \to W\) lineaarinen funktio. Osoita, että \(\image f\) on
    avaruuden \(W\) aliavaruus.
  \end{exercise}
  \begin{exercise}{4.5}
    Olkoot \(n, m, p \in \positiveintegers\). Olkoon \(A\)
    nxp-reaalimatriisi ja \(B\) pxm-reaalimatriisi.
    Olkoon \(f : \realnumbers^p \to \realnumbers^n\) matriisia \(A\)
    vastaava lineaarinen funktio ja \(g : \realnumbers^m \to
    \realnumbers^p\) matriisia \(B\) vastaava lineaarinen
    funktio. Osoita, että
    \begin{displaymath}
      f(g(\rnvec{v})) = AB\rnvec{v}
    \end{displaymath}
    kaikille reaalisille mx1-pystyvektoreille \(\rnvec{v}\).
  \end{exercise}
  \begin{exercise}{4.6}
    Todista lause \ref{th:vekt-av}.
  \end{exercise}
\end{exercises}

\chapter{Topologiset avaruudet}
\label{ch:topologiset-avaruudet}

Tässä luvussa on käytetty lähteitä
\cite{schaefer1971,wpnet,wptopspace, wptopprop}.

\section{Jonot, verkot ja filtterit}
\label{sec:j-v-f}

\begin{definition}
  Olkoon \(A\) epätyhjä joukko. Sanomme, että funktio \(f \in
  \index{jono}
  A^\naturalnumbers\) on \defterm{jono}.
\end{definition}

Jonoa merkitään yleensä \((a_k)_{k=0}^\infty\) tai vain \((a_k)\),
missä \(a_k = f(k)\),
\(k \in \naturalnumbers\). Jos \((a_k)_{k=0}^\infty\) on jono
joukon \(A\) alkioita, niin merkitään \((a_k)_{k=0}^\infty \subset
A\).
Kun \(A\) on kunta, niin jonot
lasketaan yhteen ja vähennetään alkioittain, ja myös vastaluku
lasketaan alkioittain. Jono kerrotaan ja jaetaan kunnan \(A\) alkiolla
alkioittain.
\begin{eqnarray*}
  \seqalt{a}{k} + \seqalt{b}{k} & := &
  \left(a_k+b_k\right)_{k=0}^\infty \\
  \seqalt{a}{k} - \seqalt{b}{k} & := &
  \left(a_k-b_k\right)_{k=0}^\infty \\
  -\seqalt{a}{k} & := &
  \left(-a_k\right)_{k=0}^\infty \\
  c \cdot \seqalt{a}{k} & := &
  \left(ca_k\right)_{k=0}^\infty \\
  \seqalt{a}{k} / c & := &
  \left(\frac{a_k}{c}\right)_{k=0}^\infty, \;\;\;\; c \not= 0
\end{eqnarray*}
missä \(a_k, b_k \in A\), \(k \in \naturalnumbers\),
ja \(c \in A\).

\begin{definition}
  Olkoon \(\seqalt{a}{k} \in \realnumbers^\naturalnumbers\). Määritellään jonon
  \index{jonon raja-arvo}
  \index{raja-arvo}
  \(\seqalt{a}{k}\) \defterm{raja-arvo} \(b\) asettamalla
  \begin{displaymath}
    \lim_{k\to\infty} a_k = b ,
  \end{displaymath}
  jos ja vain jos
  \begin{displaymath}
    \forall \varepsilon \in \realnumbers_+ : \exists n \in
    \naturalnumbers : \forall k \in \naturalnumbers :
    ( k > n \implies \abs{a_k-b} < \varepsilon ) .
  \end{displaymath}
  Tällöin sanomme, että raja-arvo \(\lim_{k\to\infty} a_k\) on
  olemassa.
  Asetetaan
  \begin{displaymath}
    \lim_{k\to\infty} a_k = \infty ,
  \end{displaymath}
  jos ja vain jos
  \begin{displaymath}
    \forall m \in \realnumbers_+ : \exists n \in \naturalnumbers :
    \forall k \in \naturalnumbers : ( k > n \implies \abs{a_k} > m )
  \end{displaymath}
  ja
  \begin{displaymath}
    \lim_{k\to\infty} a_k = - \infty ,
  \end{displaymath}
  jos ja vain jos
  \begin{displaymath}
    \lim_{k\to\infty} (-a_k) = \infty .
  \end{displaymath}
\end{definition}

Olkoon
\begin{displaymath}
  s_k := \sum_{j=0}^k a_k ,
\end{displaymath}
missä \(a_k \in \complexnumbers\), \(k \in \naturalnumbers\).
\index{osasummien jono}
Jonoa \(\seqalt{s}{k}\) sanotaan \defterm{sarjan osasummien jonoksi}.
Osasummien jonon raja-arvoa merkitään
\begin{displaymath}
  \sum_{k=0}^\infty a_k
\end{displaymath}
\index{sarja}
ja tätä merkintää sanotaan \defterm{sarjaksi}.
 
\begin{definition}
  Olkoon \(X\) epätyhjä joukko ja \(A\) suunnattu
  \index{verkko}
  joukko. \defterm{Verkko} on funktio \(f : A \to X\), ja sitä
  merkitään yleensä \((x_\alpha)_{\alpha \in A}\) tai vain
  \((x_\alpha)\), missä \(x_\alpha = f(\alpha)\).
  Merkintä \((x_\alpha) \subset X\) tarkoittaa, että \((x_\alpha)\) on
  verkko joltakin suunnatulta joukolta joukkoon \(X\).
\end{definition}

\begin{examples}
  \begin{example}
    Jokainen totaalisesti järjestetty joukko on suunnattu. Erityisesti
    luonnollisten lukujen joukko on suunnattu, joten jokainen jono on
    myös verkko.
  \end{example}
  \begin{example}
    Olkoon \(X\) topologinen avaruus ja \(x \in X\). Olkoon \(N\)
    kaikkien pisteen \(x\) ympäristöjen joukko ja määritellään \(A
    \geq B \iff B \subset A\) kaikille \(A, B \in N\). Nyt \(N\) on
    verkko.
  \end{example}
\end{examples}

\begin{definition}
  Olkoon \((x_\alpha)\) verkko suunnatulta joukolta \(A\)
  epätyhjään joukkoon \(X\) ja \(P(x)\), \(x \in X\), predikaatti.

  Sanomme, että
  \index{lopulta tosi}
  \(P(x_\alpha)\) on \defterm{lopulta tosi}, jos on olemassa
  \(\alpha \in A\) siten, että jokaiselle \(\beta \in A\), \(\beta
  \geq \alpha\), predikaatti \(P(x_\beta)\) on tosi.

  \index{usein tosi}
  Sanomme, että \((x_\alpha)\) on \defterm{usein tosi}, jos
  jokaiselle \(\alpha \in A\) on olemassa \(\beta \in A\), \(\beta
  \geq \alpha\) siten, että \(P(x_\beta)\) on tosi.

  \index{universaali}
  Sanomme, että \((x_\alpha)\) on \defterm{universaali}, jos
  jokaiselle \(A \subset X\) on \(x_\alpha \in A\) lopulta tosi
  tai \(x_\alpha \in X \setminus A\) lopulta tosi.
\end{definition}

\begin{definition}
  Olkoon \(X\) joukko ja \(I \subset \powerset{X}\). Sanomme, että
  \index{filtteri}
  \(I\) on \defterm{filtteri} avaruudessa \(X\), jos
  seuraavat ehdot ovat voimassa:
  \begin{itemize}
    \item[(i)] \(I \not= \emptyset \land \emptyset \not\in I\).
    \item[(ii)] \(F \in I \land F \subset G \subset X \implies G \in
      I\).
    \item[(iii)] \(F \in I \land G \in I \implies F \intersection G
      \in I\).
  \end{itemize}
\end{definition}

\hyphenation{filtteri-kanta}

\begin{definition}
  Olkoon \(X\) joukko ja \(B \subset \powerset{X}\). Sanomme, että
  \index{filtterikanta}
  \(B\) on \defterm{filtterikanta}, jos
  \begin{itemize}
    \item[(i)] \(B \not= \emptyset \land \emptyset \not\in B\).
    \item[(ii)] Jos \(B_1 \in B\) ja \(B_2 \in B\), niin on olemassa
      \(B_3 \in B\) siten, että \(B_3 \subset B_1 \intersection B_2\).
  \end{itemize}
\end{definition}

Jokainen filtterikanta \(B\) generoi yksikäsitteisen filtterin \(I\)
siten, että \(F \in I\), jos ja vain jos \(A \subset F\) ainakin
\index{kanta}
yhdelle \(A \in B\). Sanomme, että \(B\) on filtterin \(I\) kanta.
Olkoot \(I_1\) ja \(I_2\) filttereitä joukossa \(X\). Jos \(I_1
\index{karkeampi filtteri}
\subset I_2\) sanomme, että \(I_1\) on \defterm{karkeampi} kuin
\index{hienompi filtteri}
\(I_2\) ja \(I_2\) on \defterm{hienompi} kuin \(I_1\).

\section{Yleiset topologiset avaruudet}

\begin{definition}
  Olkoon \(\mytopspace\) epätyhjä joukko ja \(\mytop \subset
  \powerset{\mytopspace}\). Sanomme,
  \index{topologia}
  että \(\mytop\) on avaruuden \(\mytopspace\) \defterm{topologia},
  jos seuraavat aksioomat ovat voimassa:
  \begin{itemize}
    \item[(T1)] \(\emptyset \in \mytop\) ja \(\mytopspace \in \mytop\).
    \item[(T2)] Olkoon \(I \subset \mytop\). Nyt
      \begin{displaymath}
        \bigcup_{C \in I} C \in \mytop .
      \end{displaymath}
    \item[(T3)] Olkoon \(I\) joukon \(\mytop\) \textit{äärellinen}
      osajoukko. Nyt
      \begin{displaymath}
        \bigcap_{C \in I} C \in \mytop .
      \end{displaymath}
  \end{itemize}
  \index{topologinen avaruus}
  Paria \((\mytopspace, \mytop)\) sanotaan \defterm{topologiseksi avaruudeksi}.
\end{definition}

\hyphenation{osa-joukko}

\index{avoin joukko}
Topologian \(\mytop\) alkioita sanotaan \defterm{avoimiksi
  joukoiksi}.
\index{suljettu joukko}
Sanomme, että osajoukko \(B \subset A\) on suljettu, jos
ja vain jos \(A \setminus B\) on avoin, ts. \(A \setminus B \in
\mytop\).

\begin{definition}
  Olkoon \((\mytopspace, \mytop)\) topologinen avaruus. Sanomme, että
  \index{topologian kanta}
  \index{kanta}
  joukko \(B \subset \mytop\) on topologian \(\mytop\) \defterm{kanta}, jos
  \begin{displaymath}
    \forall C \in \mytop : \exists D \in \powerset{B} : C = \bigcup_{J \in
      D} J .
  \end{displaymath}
\end{definition}

\begin{definition}
  Olkoon \((\firsttopspace, \mytop)\) topologinen avaruus ja
  \(\secondtopspace \subset \firsttopspace\), \(\secondtopspace \not=
  \emptyset\). Asetetaan
  \begin{displaymath}
    \secondtop := \{ V \intersection \secondtopspace \setsep V \in
    \mytop\} .
  \end{displaymath}
  Nyt \((\secondtopspace,\secondtop)\) on topologinen avaruus, jotka kutsutaan topologisen
  \index{topologinen aliavaruus}
  avaruuden \defterm{topologiseksi aliavaruudeksi}. Joukkoa \(\secondtop\)
  \index{relatiivinen topologia}
  sanotaan avaruuden \(\secondtopspace\) \defterm{relatiiviseksi topologiaksi}.
\end{definition}

\begin{theorem}
  %% L19052015-1
  Edellisessä määritelmässä \(\secondtop\) on topologia joukossa
  \(\secondtopspace\).
\end{theorem}

\begin{proof}
  Nyt \(\emptyset = \emptyset \intersection \secondtopspace \in
  \secondtop\) ja \(\secondtopspace = \firsttopspace
  \intersection \secondtopspace \in \secondtop\). Olkoon
  \((K_\lambda)_{\lambda \in I} \subset \secondtop\), missä \(I\) on
  mielivaltainen joukko. Nyt \(K_\lambda = L_\lambda \intersection
  \secondtopspace\), \(L_\lambda \in \firsttop\), kaikille \(\lambda
  \in I\). Edelleen
  \begin{displaymath}
    \bigcup_{\lambda \in I} K_\lambda = \bigcup_{\lambda \in I}
    (L_\lambda \intersection \secondtopspace) = \left( \bigcup_{\lambda \in I}
    L_\lambda \right) \intersection \secondtopspace \in \secondtop .
  \end{displaymath}
  Oletetaan sitten, että \(M, N \in \secondtop\). Nyt \(M = P
  \intersection \secondtopspace\)
  ja \(N = Q \intersection \secondtopspace\), missä \(P, Q \in
  \firsttop\). Siten
  \begin{displaymath}
    M \intersection N = (P \intersection \secondtopspace)
    \intersection (Q \intersection \secondtopspace) = (P \intersection
    Q) \intersection \secondtopspace \in \secondtop .
  \end{displaymath}
\end{proof}

\begin{definition}
  Olkoon \(\mytopspace\) joukko ja \(\firsttop\) ja \(\secondtop\)
  joukon \(\mytopspace\)
  \index{hienompi topologia}
  topologioita. Sanomme, että \(\firsttop\) on \defterm{hienompi} kuin
  \(\secondtop\), jos \(\secondtop \subset \firsttop\).
  \index{karkeampi topologia}
  Sanomme, että \(\firsttop\) on \defterm{karkeampi} kuin
  \(\secondtop\), jos \(\firsttop \subset \secondtop\).
\end{definition}

\begin{remark}
  On myös mahdollista, että \(\firsttop\) ei ole hienompi eikä karkeampi
  kuin \(\secondtop\).
\end{remark}

\begin{definition}
  Olkoon \((\mytopspace, \mytop)\) topologinen avaruus ja \(x
  \in \mytopspace\). Sanomme, että \index{ympäristö} joukko \(Y
  \subset \mytopspace\) on pisteen \(x\) \defterm{ympäristö}, jos on
  olemassa avoin joukko \(A \subset \mytopspace\) siten, että \(x \in
  A\) ja \(A \subset Y\). Jos \(Y\) on avoin, niin sanomme sitä
  pisteen \(x\)
  \index{avoin ympäristö}
  \defterm{avoimeksi ympäristöksi}.
\end{definition}

\begin{figure}[htb]
\label{fig:ymparisto}
\caption{Joukko $Y$ on pisteen $x$ ympäristö, jos ja vain jos on
  olemassa avoin joukko $A$ siten, että $x \in A \subset Y$.}
\includegraphics{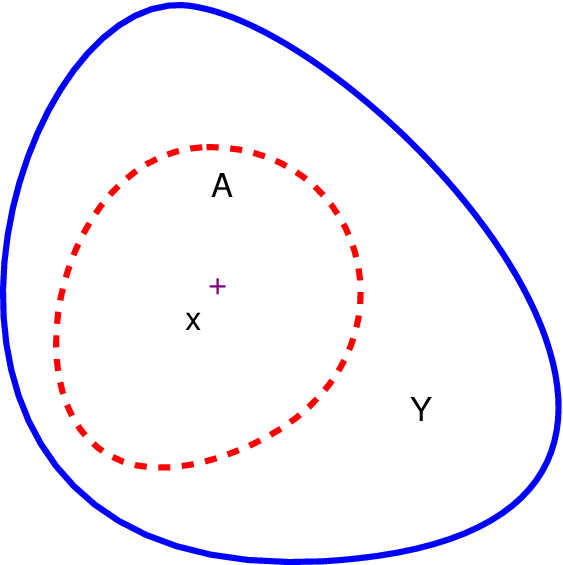}
\end{figure}

Ympäristön käsitettä on havainnollistettu kuvassa \ref{fig:ymparisto}.

\hyphenation{ympäristö-kanta}

\begin{definition}
  Olkoon \((\mytopspace, \mytop)\) topologinen avaruus ja \(x
  \in \mytopspace\).  Pisteen \(x\) kaikki ympäristöt muodostavat
  pisteen \(x\)
  \index{ympäristöfiltteri}
  \defterm{ympäristöfiltterin}, jota merkitään \(U_x\).  Olkoon \(B\)
  jokin filtterin \(U_x\) kanta. Sanomme, että \(B\) on
  \index{ympäristökanta} pisteen \(x\) \defterm{ympäristökanta}.
\end{definition}

\begin{definition}
  Olkoon \((\mytopspace, \mytop)\) topologinen avaruus ja \(x
  \in \mytopspace\).  Olkoon \(I\) filtteri joukossa
  \(\mytopspace\). Sanomme, että \(I\) suppenee pisteeseen \(x\), jos
  \(I\) on hienompi kuin \(U_x\).
\end{definition}

\begin{theorem}
  \label{th:avoin-joukko}
  %% L15052015-1
  Olkoon \(\mytopspace\) topologinen avaruus ja \(A \subset \mytopspace\). Joukko
  \(A\) on avoin, jos ja vain jos jokaisella pisteellä \(x \in A\)
  on avoin ympäristö \(Y\) siten, että \(Y \subset A\).
\end{theorem}

\begin{proof}
  Jos \(A\) on avoin, niin \(A\) on jokaisen pisteen \(x \in A\) avoin
  ympäristö.  Oletetaan, että jokaisella \(x \in A\) on avoin
  ympäristö \(Y\) siten, että \(Y \subset A\).  Olkoon \(Y(x)\) jokin
  pisteen \(x\) avoin ympäristö jokaiselle \(x \in B\). Nyt
  \begin{displaymath}
    \bigcup_{x \in A} Y(x) = A .
  \end{displaymath}
  Avoimien joukkojen mielivaltainen unioni on avoin, joten \(A\) on avoin.
\end{proof}

\begin{definition}
  Olkoon \((\mytopspace, \mytop)\) topologinen avaruus. Sanomme, että
  \((\mytopspace, \mytop)\)
  \index{Hausdorffin erotusaksiooma}
  toteuttaa \defterm{Hausdorffin erotusaksiooman}, jos jokaiselle
  kahdelle eri pisteelle \(x, y \in \mytopspace\) on olemassa pisteen
  \(x\) ympäristö \(Y\) ja pisteen \(y\) ympäristö \(Z\) siten, että
  \(Y \intersection Z = \emptyset\).  Jos \((\mytopspace, \mytop)\)
  toteuttaa Hausdorffin erotusaksiooman, niin sanomme,
  \index{Hausdorffin avaruus}
  että \((\mytopspace, \mytop)\) on \defterm{Hausdorffin avaruus}.
\end{definition}

\begin{figure}[htb]
\label{fig:hausdorff}
\caption{Tilanne, jossa Hausdorffin erotusaksiooma ei ole
  voimassa. Jokainen avoin joukko $A$, joka sisältää pisteen $x$, sisältää
  myös pisteen $y$.}
\includegraphics{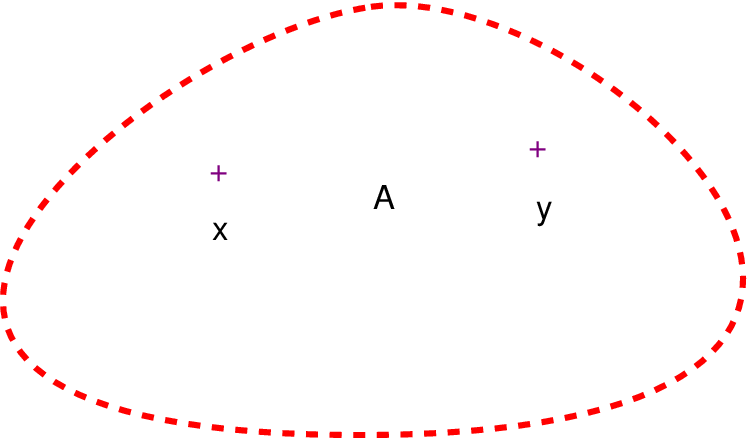}
\end{figure}

\begin{remark}
  Yleisesti raja-arvo ei välttämättä ole yksikäsitteinen, eli
  jonolla tai funktiolla voi olla useampia raja-arvoja samassa pisteessä.
  Hausdorffin avaruuksissa raja-arvot ovat yksikäsitteisiä.
  Kuvassa \ref{fig:hausdorff} on havainnollistettu Hausdorffin
  erotusaksioomaa.
\end{remark}

\begin{definition}
  Olkoon \((x_\alpha)\) verkko topologisessa avaruudessa
  \index{suppeneminen}
  \index{verkon suppeneminen}
  \(\mytopspace\). Sanomme, että verkko \defterm{suppenee kohti
    pistettä \(y \in \mytopspace\)},
  jos jokaiselle pisteen \(y\) ympäristölle \(Y\) verkko
  \((x_\alpha)\) kuuluu lopulta joukkoon \(Y\).
  Tällöin merkitsemme \(x_\alpha \to y\).
  Jos \(\mytopspace\) on Hausdorffin avaruus, niin raja-arvo on
  yksikäsitteinen, ja merkitsemme
  \begin{displaymath}
    \lim x_\alpha = y .
  \end{displaymath}
\end{definition}

\begin{definition}
  \index{kasautumispiste}
  Olkoon \((x_\alpha)\) verkko topologisessa
  avaruudessa \(\mytopspace\).  Sanomme, että \(y \in \mytopspace\) on
  verkon \((x_\alpha)\) \defterm{kasautumispiste}, jos jokaiselle
  pisteen \(y\) ympäristölle \(Y\) verkko \((x_\alpha)\) on usein
  joukossa \(Y\).
\end{definition}

\begin{theorem}
  %% L25052015-2, L26052015-1
  \label{th:suljettu-joukko}
  Olkoon \(\mytopspace\) topologinen avaruus ja \(A \subset \mytopspace\). Tällöin \(A\)
  on suljettu, jos ja vain jos jokaiselle verkolle \((x_\alpha)_{\alpha \in I}
  \subset A\) on
  \begin{displaymath}
    x_\alpha \to a \in \mytopspace \implies a \in A .
  \end{displaymath}
\end{theorem}

\begin{proof}
  Jos \(A = \emptyset\) or \(A = \mytopspace\), niin väite on tosi.
  Oletetaan, että \(A \not= \emptyset\) ja \(A \not= \mytopspace\).
  Olkoot \(A\) suljettu, \((x_\alpha)_{\alpha \in I} \subset A\)
  verkko, ja \(x_\alpha \to a \in \mytopspace\).  Oletetaan, että
  olisi \(a \not\in A\) eli \(a \in \mytopspace \setminus A\)
  (vastaoletus). Joukko \(\mytopspace \setminus A\) on avoin, joten
  lauseen \ref{th:avoin-joukko} nojalla on olemassa pisteen \(a\)
  ympäristö \(Y\) siten, että \(Y \subset \mytopspace \setminus
  A\).  Kaikille \(x_\alpha\), \(\alpha \in I\), on \(x_\alpha \not\in
  Y\), joten \(x_\alpha\) ei voi olla lopulta joukossa \(Y\). Siten
  \(x_\alpha \not\to a\), mikä on ristiriita, joten on oltava \(a \in
  A\).

  Oletetaan sitten, että
  \begin{displaymath}
    x_\alpha \to a \in \mytopspace \implies a \in A .
  \end{displaymath}
  jokaiselle verkolle \((x_\alpha) \subset A\).
  Olkoon \(x \in \mytopspace \setminus A\). Ei ole olemassa verkkoa \((x_\alpha)
  \subset A\) siten, että \(x_\alpha \to x\).
  Oletetaan, että \(Y \intersection A \not= \emptyset\) jokaiselle
  pisteen \(x\) ympäristölle \(Y\) (avaruudessa \(\mytopspace\)). Olkoon
  \begin{displaymath}
    V := \{ Y \setsep Y \;\text{on pisteen}\; x \;\text{ympäristö
      avaruudessa}\; \mytopspace \} .
  \end{displaymath}
  Kun \(Y_1, Y_2 \in V\), niin määritellään
  \begin{displaymath}
    Y_1 \geq Y_2 \iff Y_1 \subset Y_2 .
  \end{displaymath}
  Nyt \((V, \geq)\) on suunnattu joukko.  Kun \(Y \in V\) niin
  valitaan \(x_Y \in Y \intersection A\). Nyt \(x_Y \to x\), mikä on
  ristiriita, joten pisteellä \(x\) on ympäristö \(Y_0\) siten, että
  \(Y_0 \intersection A = \emptyset\). Piste \(x\) oli mielivaltainen,
  joten lauseen \ref{th:avoin-joukko} nojalla \(\mytopspace \setminus
  A\) on avoin, ja siten \(A\) on suljettu.
\end{proof}

\begin{definition}
  \label{def:jonon-raja-arvo}
  Olkoon \(\mytopspace\) topologinen avaruus ja \((x_k)_{k=0}^\infty
  \subset \mytopspace\) jono. Sanomme, että jono \((x_k)_{k=0}^\infty
  \subset \mytopspace\) suppenee kohti pistettä \(y \in X\), jos
  jokaiselle pisteen \(y\) ympäristölle \(Y\) on olemassa \(m \in
  \naturalnumbers\) siten, että \(x_k \in Y\) kaikille \(k \in
  \naturalnumbers\), \(k \geq m\).  Jos \(\mytopspace\) on Hausdorffin
  avaruus, niin raja-arvo on yksikäsitteinen, ja merkitsemme
  \begin{displaymath}
    \lim_{k \to \infty} x_k = y .
  \end{displaymath}
\end{definition}

\begin{definition}
  \label{def:raja-arvo}
  Olkoot \((\firsttopspace, \firsttop)\) ja \((\secondtopspace,
  \secondtop)\) topologisia avaruuksia, \(f
  \in \secondtopspace^\firsttopspace\) ja \(x_0
  \in \firsttopspace\). Sanomme, että piste \(y \in \secondtopspace\)
  on funktion \(f\)
  \index{raja-arvo}
  \index{funktion raja-arvo}
  \defterm{raja-arvo} pisteessä \(x_0\), jos jokaiselle pisteen \(y\)
  ympäristölle \(Z \subset \secondtopspace\) on olemassa pisten \(x_0\)
  ympäristö \(Y \subset \firsttopspace\) siten, että \(\setimage{f}{Y}
  \subset Z\).  Jos \(\secondtopspace\) on
  Hausdorffin avaruus, niin raja-arvo on yksikäsitteinen, ja merkitään
  \begin{displaymath}
    \lim_{x \to x_0} f(x) = y .
  \end{displaymath}
\end{definition}

\begin{definition}
  \label{def:jatkuvuus}
  Olkoot \((\firsttopspace, \firsttop)\) ja \((\secondtopspace,
  \secondtop)\) topologisia avaruuksia, \(f \in
  \secondtopspace^\firsttopspace\) ja \(x \in \firsttopspace\).
  \index{jatkuvuus}
  Sanomme, että f on \defterm{jatkuva} pisteessä
  \(x\), jos jokaiselle pisteen \(f(x)\) ympäristölle \(Z
  \subset \secondtopspace\) on olemassa pisteen \(x\) ympäristö \(Y
  \subset \firsttopspace\) siten, että \(\setimage{f}{Y} \subset
  Z\). Jos \(f\) jatkuva kaikissa pisteissä \(x \in \firsttopspace\),
  niin
  \index{jatkuva funktio}
  sanomme, että \(f\) on
  \defterm{jatkuva}.
\end{definition}

\begin{remark}
  Jatkuvuuden määritelmässä ympäristö voidaan korvata avoimella
  ympäristöllä, ts. jokaiselle pisteen \(f(x)\) avoimelle ympäristölle
  \(Z \subset \secondtopspace\) on olemassa pisteen \(x\) avoin
  ympäristö \(Y \subset \firsttopspace\) siten, että \(\setimage{f}{Y}
  \subset Z\).
\end{remark}

\begin{figure}[htb]
\label{fig:jatkuvuus}
\caption{Funktion jatkuvuus pisteessä $x$. Merkinnät ovat
  määritelmästä \ref{def:jatkuvuus}.}
\includegraphics{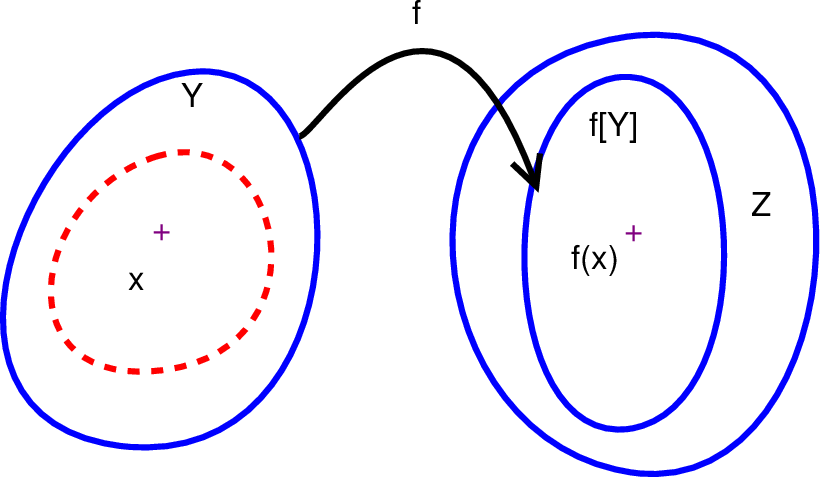}
\end{figure}

\begin{theorem}
  \label{th:top-av-suppeneminen-a}
  Olkoot \(\firsttopspace\) ja \(\secondtopspace\) topologisia
  avaruuksia ja \(f : \firsttopspace \to \secondtopspace\)
  funktio. Funktio \(f\) on jatkuva pisteessä \(x_0
  \in \firsttopspace\), jos ja vain jos jokaiselle verkolle
  \((x_\alpha) \subset \firsttopspace\) on voimassa
  \begin{displaymath}
    x_\alpha \to x_0 \implies f(x_\alpha) \to f(x_0) .
  \end{displaymath}
\end{theorem}

Vertaa lauseeseen \ref{th:top-av-suppeneminen-b}.

\begin{definition}
  Olkoot \(\firsttopspace\) ja \(\secondtopspace\) topologisia
  avaruuksia ja \(f\) kuvaus avaruudelta \(\firsttopspace\) avaruuteen
  \(\secondtopspace\). Sanomme, että \(f\) on
  \index{avoin funktio}
  \defterm{avoin}, jos \(\setimage{f}{A}\) on avoin avaruudessa
  \(\secondtopspace\) aina, kun \(A\) on avoin avaruudessa
  \(\firsttopspace\).
\end{definition}

\begin{definition}
  Olkoon \((\mytopspace, \mytop)\) topologinen avaruus ja \(A \subset \mytopspace\). Olkoon
  \begin{displaymath}
    I := \{ Y \subset \mytopspace \setsep \mytopspace \setminus Y \in \mytop \land A \subset Y \} .
  \end{displaymath}
  \index{sulkeuma}
  Määritellään joukon \(A\) \defterm{sulkeuma} asettamalla
  \begin{displaymath}
    \clos A := \bigcap_{G \in I} G .
  \end{displaymath} 
  \index{reuna}
  Joukon \(A\) \defterm{reuna} määritellään
  \begin{displaymath}
    \partial A := (\clos A) \intersection (\clos (\firsttopspace \setminus A)) .
  \end{displaymath}
  \index{sisäpuoli}
  Joukon \(A\) \defterm{sisäpuoli} määritellään
  \begin{displaymath}
    \interior A := \firsttopspace \setminus \partial A .
  \end{displaymath}
  \index{ulkopuoli}
  Joukon \(A\) \defterm{ulkopuoli} määritellään
  \begin{displaymath}
    \exterior A := \firsttopspace \setminus \clos A .
  \end{displaymath}
\end{definition}

%% \begin{definition}
%%   Olkoon \((\mytopspace, \mytop)\) topologinen avaruus, \(\mytopspace \subset A\) ja \(Y \subset
%%   A\). Sanomme, että joukot \(\mytopspace\) ja \(Y\) \defterm{eivät ole
%%     kytketty}, jos on olemassa avoimet joukot \(M, N \subset A\)
%%   siten, että \(\mytopspace \subset M\), \(Y \subset N\) ja \(M \intersection N
%%   = \emptyset\). Muussa tapauksessa sanomme, että \(\mytopspace\) ja \(Y\) ovat
%%   \index{kytketyt joukot}
%%   \defterm{kytketty}.
%% \end{definition}

\begin{definition}
  Olkoon \(\mytopspace\) topologinen avaruus. Sanomme, että \(\mytopspace\) on
  \index{kytketty topologinen avaruus}
  \defterm{kytketty}, jos se ei ole kahden epätyhjän avoimen joukon
  unioni.
\end{definition}

\begin{figure}[htb]
\label{fig:ei-kytketty-joukko}
\caption{Topologinen avaruus, joka ei ole kytketty.}
\includegraphics{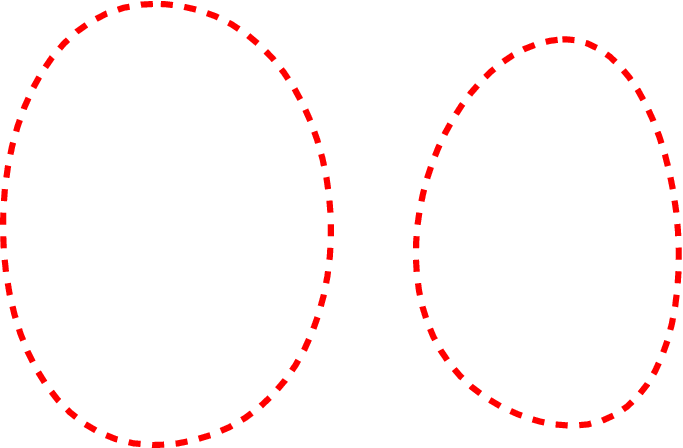}
\end{figure}

Kuva \ref{fig:ei-kytketty-joukko} esittää joukkoa, joka ei ole
kytketty.

\hyphenation{polku-kytketty}

\begin{definition}
  Olkoon \(\mytopspace\) topologinen avaruus. Sanomme, että
  \index{polkukytketty topologinen avaruus} \(\mytopspace\) on
  \defterm{polkukytketty}, jos
  kaikille kaksille pisteille \(x, y \in \mytopspace\) on olemassa jatkuva
  funktio \(p : \intervalcc{0}{1} \to \mytopspace\), jolle \(p(0)=x\) ja
  \(p(1)=y\).
\end{definition}

\begin{remark}
  Polkukytketty avaruus on aina kytketty.
\end{remark}

\begin{definition}
  Olkoon \(\mytopspace\) topologinen avaruus ja \(A, B \subset \mytopspace\). Sanomme,
  \index{tiheä joukko}
  että \(B\) on \defterm{tiheä} avaruudessa \(A\), jos \(B \subset A\)
  ja \(A \subset \clos B\).
\end{definition}

\begin{definition}
  Olkoon \(\mytopspace\) topologinen avaruus ja \(E \subset \mytopspace\). Sanomme, että
  \index{ei missään tiheä}
  \(E\) \defterm{ei ole missään tiheä}, jos \(\interior \clos E =
  \emptyset\).
\end{definition}

\begin{definition}
  Olkoon \(\mytopspace\) topologinen avaruus ja \(A \subset \mytopspace\). Sanomme, että
  \index{ensimmäinen kategoria}
  \(A\) on \defterm{ensimmäistä kategoriaa avaruudessa \(\mytopspace\)}, jos
  \(A\) on korkeintaan numeroituva unioni joukoista, jotka eivät ole
  \index{toinen kategoria}
  missään tiheitä. Sanomme, että \(A\) on \defterm{toista kategoriaa
    avaruudessa \(\mytopspace\)}, jos se ei ole ensimmäistä kategoriaa
  avaruudessa \(\mytopspace\).
\end{definition}

\begin{definition}
  Olkoon \(\mytopspace\) topologinen avaruus. Sanomme, että \(\mytopspace\) on
  \index{separoituva avaruus}
  \defterm{separoituva}, jos on olemassa korkeintaan numeroituva ja
  tiheä (avaruudessa \(\mytopspace\)) avaruuden \(\mytopspace\) osajoukko.
\end{definition}

\begin{definition}
  Olkoon \((\mytopspace, \mytop)\) topologinen avaruus ja \(V
  \subset \mytopspace\). Olkoon \(\Gamma \subset
  \powerset{\mytopspace}\). Sanomme, että \(\Gamma\) on joukon
  \index{peite}
  \(V\) \defterm{peite}, jos
  \begin{displaymath}
    V \subset \bigcup_{G \in \Gamma} G .
  \end{displaymath}
  Jos kaikki joukot \(G \in \Gamma\) ovat avoimia, niin sanomme, että
  \index{avoin peite}
  \(\Gamma\) on \defterm{avoin peite}.
  Jos \(\Xi \subset \Gamma\) ja \(\Xi\) on joukon \(V\) peite, niin
  \index{alipeite}
  sanomme, että \(\Xi\) on peitteen \(\Gamma\) \defterm{alipeite}
  joukolle \(V\).
\end{definition}

\begin{definition}
  Olkoon \((\mytopspace, \mytop)\) topologinen avaruus ja \(A
  \subset \mytopspace\). Sanomme,
  \index{kompakti joukko}
  että \(A\)
  on \defterm{kompakti}, jos jokaisella avoimella peitteellä joukolle
  \(A\) on äärellinen alipeite joukolle \(A\).  Jos \(\mytopspace\) on
  kompakti,
  \index{kompakti topologinen avaruus}
  niin sanomme, että
  \((\mytopspace, \mytop)\) on \defterm{kompakti topologinen avaruus}.
\end{definition}

\begin{definition}  
  Olkoon \((\mytopspace, \mytop)\) Hausdorffin avaruus ja \(A
  \subset \mytopspace\). Sanomme, että \(A\) on
  \index{jonokompakti joukko}
  \defterm{jonokompakti}, jos jokaisella jonolla joukossa \(A\) on
  suppeneva osajono, jonka raja-arvo on joukossa \(A\).
\end{definition}

\begin{definition}
  Olkoon \((\mytopspace, \mytop)\) topologinen avaruus. Sanomme, että
  \((\mytopspace, \mytop)\) on
  \index{lokaalikompakti avaruus}
  \defterm{lokaalikompakti avaruus}, jos jokaisella \(x \in \mytopspace\) on
  kompakti ympäristö.
\end{definition}

\begin{definition}
  Olkoot \((\firsttopspace, \firsttop)\) ja \((\secondtopspace,
  \secondtop)\) topologisia avaruuksia. Sanomme,
  \index{homeomorfismi}
  että funktio \(f : \firsttopspace \to \secondtopspace\) on
  \defterm{homeomorfismi}, jos \(f\) on bijektio ja funktiot \(f\) ja
  \(f^{-1}\) ovat jatkuvia.  Jos on olemassa homeomorfismi avaruudelta
  \((\firsttopspace, \firsttop)\) avaruudelle \((\secondtopspace,
  \secondtop)\), niin sanomme, että \((\firsttopspace, \firsttop)\) ja
  \((\secondtopspace, \secondtop)\) ovat
  \index{homeomorfiset avaruudet}
  \defterm{homeomorfiset}.
\end{definition}

\begin{remark}
  Homeomorfisten avaruuksien topologioita voidaan pitää samoina.
\end{remark}

\begin{definition}
  Olkoot \((\firsttopspace, \firsttop)\) ja \((\secondtopspace,
  \secondtop)\) topologisia avaruuksia. Sanomme,
  \index{upotus}
  että funktio \(f : \firsttopspace \to \secondtopspace\) on
  \defterm{upotus}, jos \(f\) on homeomorfismi avaruudelta
  \(\firsttopspace\) avaruudelle \(\setimage{f}{\firsttopspace}\).
\end{definition}

\hyphenation{kompakti-kantajainen}

\begin{definition}
  \label{def:topologinen-kantaja}
  Olkoon \(\mytopspace\) topologinen avaruus ja \(f : \mytopspace \to
  \complexnumbers\)
  \index{topologinen kantaja}
  \index{kantaja}
  funktio. Määrittelemme funktion \(f\) \defterm{(topologisen)
    kantajan} asettamalla
  \begin{displaymath}
    \supp f := \clos \stsupp f .
  \end{displaymath}
  Jos funktion \(f\) kantaja on kompakti, niin sanomme, että \(f\) on
  \index{kompaktikantajainen funktio}
  \defterm{kompaktikantajainen}.
\end{definition}

\begin{definition}
  Olkoon \(\mytopspace\) topologinen avaruus. Sanomme, että \(\mytopspace\) on
  \defterm{1. lajia laskettava}, jos jokaisella pisteellä \(x \in \mytopspace\)
  on korkeintaan numeroituva paikallinen ympäristökanta. Sanomme, että
  \(\mytopspace\) on \defterm{2. lajia laskettava}, jos sen topologialla on
  korkeintaan numeroituva kanta.
\end{definition}

\begin{remark}
  2. lajia laskettava topologinen avaruus on aina 1. lajia laskettava
  ja separoituva.
\end{remark}

\begin{theorem}
  \label{th:top-av-suppeneminen-b}
  Olkoon \(\firsttopspace\) 1. lajia laskettava topologinen avaruus ja
  \(\secondtopspace\) topologinen avaruus. Olkoon \(f : \firsttopspace
  \to \secondtopspace\) funktio. Funktio \(f\) on jatkuva, jos ja vain
  jos jokaiselle jonolle \((x_k)_{k=0}^\infty \subset \firsttopspace\)
  on voimassa
  \begin{displaymath}
    x_k \to y \in \firsttopspace \implies f(x_k) \to f(y) .
  \end{displaymath}
\end{theorem}

Vertaa lauseeseen \ref{th:top-av-suppeneminen-a}.

\begin{theorem}
  \label{th:suljettu-joukko-b}
  Olkoon \(\firsttopspace\) 1. lajia laskettava topologinen avaruus
  ja \(A \subset \firsttopspace\). Tällöin \(A\) on suljettu, jos ja
  vain jos jokaiselle jonolle
  \((x_k)_{k=0}^\infty \subset A\) on
  \begin{displaymath}
    x_k \to a \in \firsttopspace \implies a \in A .
  \end{displaymath}
\end{theorem}

Vertaa lauseeseen \ref{th:suljettu-joukko}.

\begin{remark}
  Lauseet \ref{th:top-av-suppeneminen-b} ja \ref{th:suljettu-joukko-b}
  ovat erityisesti voimassa, jos \(\firsttopspace\) on metrinen
  avaruus, ks. luku \ref{sec:metr-av}.
\end{remark}

\begin{definition}
  \label{def:tasainen-struktuuri}
  Olkoon \(\mytopspace\) joukko. Kun \(V, W \in \mytopspace \times \mytopspace\), niin merkitään
  \(W^{-1} := \{ (y,x) \setsep (x,y) \in W\}\) ja \(V \circ W := \{
  (x, z) \setsep \exists y \in \mytopspace : (x,y) \in W \land (y,z) \in V \}\).
  Joukkoa \(\Delta := \{ (x,x) \setsep x \in \mytopspace\}\) sanotaan joukon \(\mytopspace
  \times \mytopspace\) diagonaaliksi. Olkoon \(Z\) filtteri avaruudessa \(\mytopspace
  \times \mytopspace\) siten, että
  \begin{itemize}
    \item[(i)] \(\forall W \in Z : \Delta \subset W\).
    \item[(ii)] \(W \in Z \implies W^{-1} \in Z\).
    \item[(iii)] Jokaiselle \(W \in Z\) on olemassa \(V \in Z\) siten,
      että \(V \circ V \in W\).
  \end{itemize}
  Sanomme, että filtteri \(Z\) (tai jokin sen kanta) määrittelee
  \index{tasainen struktuuri}
  \defterm{tasaisen struktuurin} avaruudessa \(\mytopspace\). Jokaista \(W \in
  \index{lähistö}
  Z\) kutsutaan tämän struktuurin \defterm{lähistöksi}.
\end{definition}

Olkoon \(Q := \{ G \in \powerset{\mytopspace} \setsep x \in G \implies \exists W
\in Z : \{ y \setsep (x,y) \in W \} \subset G \}\).
Olkoon \(Y(W, x) := \{ y \setsep (x,y) \in W \}\), missä \(W \in Z\) ja
\(x \in \mytopspace\).
Nyt \(Q\)
määrittelee topologian \(\mytop\) joukossa \(\mytopspace\) siten, että
jokaiselle \(x \in \mytopspace\) joukko \(\{ Y(W, x)\setsep W \in Z \}\)
on pisteen \(x\)
\index{ympäristökanta}
\index{tasainen avaruus}
ympäristökanta. Nyt paria \((\mytopspace, \mytop)\) kutsutaan
\defterm{tasaiseksi avaruudeksi}.
Topologinen avaruus \(\mytopspace\) on
\index{tasoittuva}
\defterm{tasoittuva}, jos sen topologia voidaan generoida yllä
kuvatulla tavalla.

\begin{definition}
  \label{def:taydellisyys}
  Olkoon \(\mytopspace\) tasainen avaruus ja \(I\) filtteri avaruudessa
  \index{Cauchyn filtteri}
  \(\mytopspace\). Sanomme, että \(I\) on \defterm{Cauchyn filtteri}, jos
  jokaiselle lähistölle \(V\) on olemassa \(F \in I\) siten, että \(F
  \times F \subset V\).
  Jos jokainen Cauchyn filtteri suppenee kohti jotakin pistettä \(x
  \index{täydellinen}
  \in \mytopspace\), niin sanomme, että \(\mytopspace\) on \defterm{täydellinen}.
\end{definition}

Jokaiselle tasaiselle avaruudelle \(\mytopspace\) voidaan konstruoida
täydellinen tasainen avaruus \(\tilde{\mytopspace}\) siten, että \(\mytopspace\) on
isomorfinen joukon \(\tilde{\mytopspace}\) jonkin tiheän osajoukon kanssa, ja
siten, että \(\tilde{\mytopspace}\) on separoituva, jos \(\mytopspace\) on.
Jos \(\mytopspace\) on separoituva, niin \(\tilde{\mytopspace}\) on määritelty näiden
ominaisuuksien nojalla isomorfisuuteen asti, ja sitä kutsutaan
\index{täydellistymä}
avaruuden \(\mytopspace\) \defterm{täydellistymäksi}.

\begin{definition}
  Olkoon \(\mytopspace\) epätyhjä joukko ja \((Y_\alpha)_{\alpha \in I}\) perhe
  topologisia avaruuksia. Olkoon \(f_\alpha\) funktio avaruudelta
  \(\mytopspace\) avaruuteen \(Y_\alpha\) kullekin \(\alpha \in
  \index{projektiivinen topologia}
  I\). \defterm{Projektiivinen
    topologia} avaruudella \(\mytopspace\) perheen \((Y_\alpha,
  f_\alpha)_{\alpha \in I}\) suhteen määritellään karkeimmaksi
  topologiaksi (avaruudella \(\mytopspace\)), jossa jokainen \(f_\alpha\) on
  jatkuva.
  \index{induktiivinen topologia}
  \defterm{Induktiivinen topologia} avaruudella \(\mytopspace\) perheen
  \((Y_\alpha, f_\alpha)_{\alpha \in I}\) suhteen määritellään
  hienoimmaksi topologiaksi (avaruudella \(\mytopspace\)), jossa jokainen
  \(f_\alpha\) on jatkuva.
\end{definition}

\begin{definition}
  Olkoon \((Y_\alpha)_{\alpha \in I}\) perhe topologisia avaruuksia,
  \(\mytopspace\) niiden karteesinen tulo ja \(f_\alpha\) projektio avaruudelta
  \(\mytopspace\) avaruudelle \(Y_\alpha\) kullekin \(\alpha \in
  I\). Projektiivista topologiaa avaruudella \(\mytopspace\) perheen
  \((Y_\alpha,f_\alpha)_{\alpha \in I}\) suhteen kutsutaan
  \index{tulotopologia}
  \defterm{tulotopologiaksi} ja avaruutta \(\mytopspace\) kutsutaan perheen
  \index{topologinen tulo}
  \index{tulo}
  \((\mytopspace_\alpha)_{\alpha \in I}\) \defterm{(topologiseksi) tuloksi}.
\end{definition}

\begin{definition}
  Olkoon \((G, \circ)\) ryhmä ja \(T\) topologia joukossa \(G\).
  \index{topologinen ryhmä}
  Sanomme, että \((G, \circ, T)\) on \defterm{topologinen ryhmä},
  jos funktio \(\circ : G \times G \to G\) on jatkuva. Avaruudessa \(G
  \times G\) käytetään tulotopologiaa.
\end{definition}

\begin{definition}
  Olkoon \((K, +, \cdot)\) kunta ja \(T\) topologia joukossa \(K\).
  \index{topologinen kunta}
  Sanomme, että \((K, +, \cdot, T)\) on \defterm{topologinen kunta},
  jos funktiot \(+ : K \times K \to K\) ja \(\cdot : K \times K \to
  K\) ovat jatkuvia. Avaruudessa \(K \times K\) käytetään tulotopologiaa.
\end{definition}

\index{diskreetti topologia}

\begin{definition}
  Olkoon \((\mytopspace, \mytop)\) topologinen avaruus. Sanomme, että
  \(\mytop\) on \defterm{diskreetti}, jos jokainen \(A
  \subset \mytopspace\) on avoin.
\end{definition}

\section{Metriset avaruudet}
\label{sec:metr-av}

Metriset avaruudet ovat joukkoja, joissa on määritelty tietyt
aksioomat toteuttava etäisyys pisteiden välillä.

\begin{figure}[htb]
\label{fig:kolmioepayhtalo}
\caption{Kolmioepäyhtälö: $d(x,y) \leq d(x,z) + d(z,y)$.}
\includegraphics{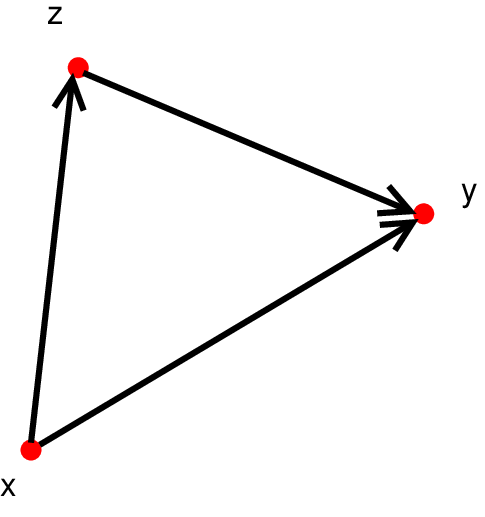}
\end{figure}

\begin{definition}
  \label{def:metr-av}
  Olkoon \(E\) epätyhjä joukko ja \(d: E \times E \to \realnumbers_0\)
  \index{metriikka}
  funktio. Sanomme, että \(d\) on \defterm{metriikka}, jos seuraavat
  aksioomat ovat voimassa:
  \begin{itemize}
    \item[(M1)] \(d(x,y) >= 0\) kaikille \(x,y \in E\).
    \item[(M2)] \(d(x,x) = 0 \iff x = 0\) kaikille \(x \in E\).
    \item[(M3)] \(d(x,y) = d(y,x)\) kaikille \(x,y \in E\).
    \item[(M4)] \(d(x,y) \leq d(x,z) + d(z,y)\) kaikille \(x,y,z \in E\)
  \end{itemize}
  \index{metrinen avaruus}
  Paria \((E, d)\) kutsutaan \defterm{metriseksi avaruudeksi}.
  \index{kolmioepäyhtälö}
  Aksioomaa (M4) kutsutaan \defterm{kolmioepäyhtälöksi}, ja sitä on
  havainnollistettu kuvassa \ref{fig:kolmioepayhtalo}.
\end{definition}

\begin{examples}
  \begin{example}
    Reaalilukujen joukko \(\realnumbers\) on metrinen avaruus, kun
    määritellään \(d(x,y) := \abs{x-y}\) kaikille \(x,y \in
    \realnumbers\).
  \end{example}
  \begin{example}
    Vastaavasti kompleksilukujen joukko \(\complexnumbers\) on
    metrinen avaruus.
  \end{example}
  \begin{example}
    \(\realnumbers^n\) ja \(\complexnumbers^n\), \(n \in
    \positiveintegers\) ovat metrisiä avaruuksia, kun asetetaan
    \begin{displaymath}
      d\left(
      \left(
      \begin{array}{c}
        x_1 \\
        \vdots \\
        x_n
      \end{array}
      \right) ,
      \left(
      \begin{array}{c}
        y_1 \\
        \vdots \\
        y_n
      \end{array}
      \right)
      \right)
      :=
      \sqrt{\sum_{k=1}^n \abs{y_k - x_k}^2} ,
      \;\;\;\;
      \left(
      \begin{array}{c}
        x_1 \\
        \vdots \\
        x_n
      \end{array}
      \right) ,
      \left(
      \begin{array}{c}
        y_1 \\
        \vdots \\
        y_n
      \end{array}
      \right)
      \in \realnumbers \;\textrm{tai}\; \complexnumbers .
    \end{displaymath}
    %% Tod. HT
  \end{example}
  \begin{example}
    %% VIITTAUS FUNKTIOAVARUUKSIIN
    Olkoon \(F\) kaikkien jatkuvien ja neliöllisesti
    integroituvien joukossa \(\realnumbers^n\) määriteltyjen
    kompleksiarvoisten funktioiden joukko. \(F\) on metrinen avaruus,
    kun asetetaan
    \begin{displaymath}
      d(f,g) := \sqrt{\int_{x \in \realnumbers^n} \abs{f(x) - g(x)}^2 d\mu}
    \end{displaymath}
    missä \(f,g \in F\).
  \end{example}
\end{examples}

\begin{definition}
  \index{avoin pallo}
  Olkoon \((E,d)\) metrinen avaruus. Määritellään \defterm{avoin pallo}
  pisteen \(x \in E\) ympärillä asettamalla
  \begin{displaymath}
    \openball{E}{x}{r} := \{ y \in E \setsep d(x,y) < r \}
  \end{displaymath}
  ja \defterm{suljettu pallo}
  \begin{displaymath}
    \closedball{E}{x}{r} := \{ y \in E \setsep d(x,y) \leq r \}.
  \end{displaymath}
\end{definition}

\index{metrisen avaruuden topologia}

Metrisen avaruuden \((E,d)\) topologia \(\mytop\) määritellään topologiana, jonka
kanta on
\begin{equation}
  \label{eq:metr-av-kanta}
  B := \{ \openball{E}{x}{r} \setsep x \in E, r \in \realnumbers_0 \} .
\end{equation}
\index{yhteensopivuus}
Sanomme, että metriikka \(d\) ja topologia \(\mytop\) ovat
\defterm{yhteensopivia} toistensa kanssa.

\begin{theorem}
  %% L19052015-2
  Olkoon \((E,d)\) metrinen avaruus. Kannan \eqref{eq:metr-av-kanta}
  määrittelemä topologia \(\mytop\) toteuttaa topologisen avaruuden
  aksioomat.
\end{theorem}

\begin{proof}
  \mbox{ }
  \begin{itemize}
    \item[(T1):]
      \(\emptyset = \openball{E}{x}{0}\) jollekin \(x \in E\), joten
      \(\emptyset \in \mytop\).
      \(E = \union_{y \in E} \openball{E}{y}{1}\), joten \(E \in \mytop\).
    \item[(T2):]
       Olkoon \((M_\alpha)_{\alpha \in I} \subset \mytop\), missä \(I\) on
       mielivaltainen joukko. Nyt
       \begin{displaymath}
         M_\alpha = \bigcup_{C \in B_\alpha} C,
       \end{displaymath}
       missä \(B_\alpha \subset B\) kaikilla \(\alpha \in I\).
       Edelleen
       \begin{displaymath}
         \bigcup_{\alpha \in I} M_\alpha = \bigcup_{\alpha \in I}
         \bigcup_{C \in B_\alpha} C = \bigcup_{C \in V} C ,
       \end{displaymath}
       missä
       \begin{displaymath}
         V := \bigcup_{\alpha \in I} B_\alpha .
       \end{displaymath}
       Täten \(\union_{\alpha \in I} M_\alpha \in \mytop\).
    \item[(T3):]
      Olkoot \(M, N \in \mytop\). Nyt
      \begin{displaymath}
        M = \bigcup_{C \in B_1} C
      \end{displaymath}
      ja
      \begin{displaymath}
        N = \bigcup_{C \in B_2} C ,
      \end{displaymath}
      missä \(B_1, B_2 \in B\).
      Jos \(M = \emptyset\) tai \(N = \emptyset\), niin (T3) on
      tosi. Voidaan olettaa, että \(M \not= \emptyset\) ja \(N \not=
      \emptyset\). Olkoon \(x_0 \in M \intersection N\). Nyt \(x_0 \in
      P\) jollekin \(P \in B_1\), ja \(x_0 \in Q\) jollekin \(Q \in
      B_2\).
      Meillä on \(P = \openball{E}{x_1}{r_1}\) ja \(Q =
      \openball{E}{x_2}{r_2}\).
      Olkoon \(r_3:= \frac{1}{2} \min \{ r_1, r_2 \}\) ja \(y \in
      \openball{E}{x_0}{r_3}\).
      Olkoot
      \begin{eqnarray*}
        d_1 & := & d(x_0, x_1) < r_1 \\
        d_2 & := & d(x_0, x_2) < r_2 .
      \end{eqnarray*}
      Nyt
      \begin{eqnarray*}
        d(y, x_1) & \leq & d(y, x_0) + d(x_0, x_1) \\
        & < & \frac{1}{2}(r_1-d_1) + d_1 = \frac{1}{2} r_1 +
        \frac{1}{2} d_1 \\
        & < & \frac{1}{2} r_1 + \frac{1}{2} r_1 = r_1
      \end{eqnarray*}
      ja vastaavasti
      \begin{eqnarray*}
        d(y, x_2) & \leq & d(y, x_0) + d(x_0, x_2) \\
        & < & \frac{1}{2}(r_2-d_2) + d_2 = \frac{1}{2} r_2 +
        \frac{1}{2} d_2 \\
        & < & \frac{1}{2} r_2 + \frac{1}{2} r_2 = r_2 .
      \end{eqnarray*}
      Alkio \(y \in \openball{E}{x_0}{r_3}\) oli mielivaltainen, joten
      \[\openball{E}{x_0}{r_3} \subset P\] ja \[\openball{E}{x_0}{r_3}
      \subset Q , \] mistä seuraa, että
      \(\openball{E}{x_0}{r_3} \subset P \intersection Q \subset M
      \intersection N\).
      Siis jokaisella \(x_0 \in M \intersection N\) on olemassa \(Y
      \in B\) siten, että \(x_0 \in Y\) ja \(Y \subset M \intersection
      N\). Täten \(M \intersection N \in \mytop\).
  \end{itemize}
\end{proof}

Kun \((E,d)\) on metrinen avaruus, niin sen topologian generoi
tasainen struktuuri, joka koostuu joukoista \(W_\varepsilon := \{ (x,
y) \setsep x, y \in E \land d(x, y) < \varepsilon \}\), \(\varepsilon \in
\realnumbers_+\). Kohdan (iii) määritelmässä
\ref{def:tasainen-struktuuri} todistamiseksi asetetaan \(V :=
W_{\varepsilon / 2}\). Tällöin \((x,y) \in V\) ja \((y,z) \in
V\). Olkoon \((x,z) \in V \circ V\). Nyt \(d(x,z) \leq d(x,y) + d(y,z)
< \frac{\varepsilon}{2} + \frac{\varepsilon}{2} = \varepsilon\), joten
\((x, z) \in W_\varepsilon\).

\begin{definition}
  Olkoon \((\mytopspace, \mytop)\) topologinen avaruus. Sanomme, että
  topologia \(\mytop\)
  \index{metrisoituva topologia}
  on \defterm{metrisoituva}, jos on olemassa metriikka \(d : A \times
  A \to \realnumbers_0\) siten, että \(d\) generoi topologian \(\mytop\)
  yllä olevien määritelmien mukaisesti.
\end{definition}

\begin{definition}
  Olkoon \((E,d)\) metrinen avaruus.
  Olkoon \(\seq{a} := (a_k)_{k=0}^\infty \subset E\). Sanomme,
  \index{Cauchyn jono}
  että \(\seq{a}\) on \defterm{Cauchyn jono}, jos
  \begin{displaymath}
    \forall \varepsilon \in \realnumbers_+ : \exists N \in
    \naturalnumbers : \forall n, m \in \naturalnumbers :
    n > N \land m > N \implies d(a_n, a_m) < \varepsilon .
  \end{displaymath}
\end{definition}

\begin{theorem}
  \label{th:metr-av-taydellinen}
  %% L06062015-1
  Olkoon \((E,d)\) metrinen avaruus. Avaruus \((E,d)\) on
  \index{täydellinen}
  täydellinen, jos ja vain jos jokainen joukon \(E\) alkioiden Cauchyn
  jono suppenee kohti jotakin avaruuden \(E\) alkiota.
\end{theorem}

Todistus harjoitustehtävänä.

Metrisen avaruuden \((E, d)\) täydellistymä voidaan konstruoida sen
Cauchyn jonojen muodostamista ekvivalenssiluokista seuraavasti
\cite{wpcomplmetrspace}: Kun \(\seq{x}:=(x_k), \seq{y}:=(y_k) \subset
E\), niin määritellään
\begin{displaymath}
  d( \seq{x}, \seq{y} ) = \lim_{k \to \infty} d(x_k, y_k)
\end{displaymath}
Määritellään kaksi Cauchyn jonoa \(\seq{x}\) ja \(\seq{y}\)
ekvivalenteiksi, jos ja vain jos \(d( \seq{x}, \seq{y} ) = 0\).

\begin{definition}
  Olkoot \((E_1, d_1)\) ja \((E_2, d_2)\) metrisiä avaruuksia ja \(f :
  \index{isometria}
  E_1 \to E_2\) funktio. Sanomme, että \(f\) on \defterm{isometria},
  jos
  \begin{displaymath}
    d_2(f(x),f(y)) = d_1(x,y)
  \end{displaymath}
  kaikille \(x,y \in E\).
\end{definition}

\begin{definition}
  Olkoon \((E, d)\) metrinen avaruus ja \(A \subset E\).
  \index{rajoitettu}
  Sanomme, että \(A\) on \defterm{rajoitettu}, jos se kuuluu johonkin
  \index{halkaisija}
  (äärellissäteiseen) palloon. Joukon \(A\) \defterm{halkaisija} on pienimmän
  joukon \(A\) sisältävän suljetun pallon halkaisija.
  \index{etäisyys}
  \index{$\dist$}
  Määritellään pisteen \(f \in E\) \defterm{etäisyys} joukosta \(A\)
  kaavalla
  \begin{displaymath}
    \dist(f,A) := \inf_{g \in A} d(f,g) .
  \end{displaymath}
\end{definition}

Metrisessä avaruudessa funktion raja-arvon määritelmä saa muodon
\begin{align}
  \nonumber
  & \lim_{x \to x_0} f(x) = a \iff \\
  \label{eq:metr-av-raja-arvo}
  & \forall \varepsilon \in \realnumbers_+ : \exists \delta \in
  \realnumbers_+ : \forall x, y \in E_1 :
  d_1(x, y) < \delta \implies d_2(f(x), f(y)) < \varepsilon ,
\end{align}
missä \((E_1, d_1)\) ja \((E_2, d_2)\) ovat metrisiä avaruuksia ja \(f
: E_1 \to E_2\) funktio, ja jonon raja-arvo saa muodon
\begin{align}
  \nonumber
  & \lim_{k \to \infty} x_k = a \iff \\
  \label{eq:metr-av-jonon-raja-arvo}
  & \forall \varepsilon \in \realnumbers_+ : \exists m \in
  \naturalnumbers: \forall k \in \naturalnumbers:
  k \geq m \implies d(x_k, a) < \varepsilon ,
\end{align}
missä \((E, d)\) on metrinen avaruus ja \((x_k)_{k=0}^\infty \subset
E\).

\begin{definition}
  Olkoon \(E\) metrinen avaruus, \(A \subset E\) ja \(\Gamma\) joukon
  \(A\) peite. Jos kaikkien joukkojen \(G \in \Gamma\) halkaisija on
  pienempi kuin \(\varepsilon \in \realnumbers_+\), niin sanomme, että
  \index{\(\varepsilon\)-peite}
  \(\Gamma\) on joukon \(A\) \defterm{\(\varepsilon\)-peite}.
\end{definition}

\begin{definition}
  Olkoon \(E\) metrinen avaruus ja \(A \subset E\). Sanomme, että
  \index{totaalisesti rajoitettu}
  \(A\) on \defterm{totaalisesti rajoitettu}, jos jokaiselle
  \(\varepsilon \in \realnumbers_+\) on olemassa äärellinen joukko
  avoimia palloja, jotka muodostavat \(\varepsilon\)-peitteen joukolle
  \(A\).
\end{definition}

\begin{theorem}
  Olkoon \(E\) metrinen avaruus. Avaruus \(E\) on kompakti, jos ja
  vain jos se on jonokompakti \cite{wpseqcompspace}.
\end{theorem}

\begin{theorem}
  Tämä lause ja todistus pohjautuvat lauseeseen 5.2.8 kirjassa
  \cite{hpc2005}. Ks. myös \cite{wpcompactspace,wpcomplmetrspace}.

  Olkoon \(E\) metrinen avaruus. Avaruus \(E\) on
  kompakti, jos ja vain jos se on täydellinen ja totaalisesti
  rajoitettu. 
\end{theorem}

\hyphenation{jono-kompakti}

\begin{proof}
  Oletetaan ensin, että \(E\) on kompakti. Olkoon \(\varepsilon \in
  \realnumbers_+\) ja
  \begin{displaymath}
    P := \{ \openball{E}{x}{\varepsilon} \setsep x \in E \} .
  \end{displaymath}
  Joukko \(P\) on peite avaruudelle \(E\), joten joukolla \(P\) on
  äärellinen alipeite avaruudelle \(E\). Siten \(E\) on totaalisesti
  rajoitettu.

  Olkoon \((x_i)\) Cauchyn jono avaruudessa \(E\).
  Koska \(E\) on kompakti metrinen avaruus, niin se on jonokompakti,
  ja täten jonolla \((x_i)\) on suppeneva osajono
  \((x_{n(k)})_{k=0}^\infty \subset E\). Olkoon
  \begin{displaymath}
    a = \lim_{k \to \infty} x_{n(k)} .
  \end{displaymath}
  Olkoon \(\varepsilon \in \realnumbers_+\). On olemassa \(k_0 \in
  \naturalnumbers\) siten, että \(d(x_{n(k)}, a) < \varepsilon / 2\)
  kaikilla \(k \in \naturalnumbers\), \(k \geq k_0\).
  Koska \((x_i)\) on Cauchyn jono, voidaan valita \(i_0 \in
  \naturalnumbers\) siten, että \(i_0 > n_{n(k_0)}\) ja \(d(x_i, x_j)
  < \varepsilon / 2\) kaikilla \(i, j \in \naturalnumbers\), \(i, j
  \geq i_0\). Edelleen voidaan valita \(k_1 \in \naturalnumbers\),
  jolle \(n(k_1) > i_0\). Kun \(i > i_0\), niin
  \begin{displaymath}
    d(x_i, a) \leq d(x_i, x_{n(k_1)}) + d(x_{n(k_1)}, a)
    < \frac{\varepsilon}{2} + \frac{\varepsilon}{2}
    = \varepsilon .
  \end{displaymath}
  Täten \((x_i)\) suppenee avaruudessa \(E\) (raja-arvo \(a\)) ja
  \(E\) on täydellinen.

  Oletetaan sitten, että \(E\) on täydellinen ja totaalisesti
  rajoitettu. Olkoon \((x_n)\) ääretön jono avaruudessa \(E\). Koska
  \(E\) on
  totaalisesti rajoitettu, on olemassa avaruuden \(E\) äärellinen
  1-peite. Merkitään tätä \(\Gamma :=
  A(g_1,\frac{1}{2}),\ldots,A(g_N,\frac{1}{2})\). Ainakin yksi näistä
  palloista sisältää jonon \((x_n)\) äärettömän osajonon, jota
  merkitään \((x_{n,1})\). Otetaan seuraavaksi äärellinen
  \(\frac{1}{2}\)-peite avaruudelle \(E\) ja kuten aikaisemminkin
  muodostetaan jonon \((x_{n,1})\) ääretön osajono \((x_{n,2})\), joka
  sisältyy yhteen peitteen palloista. Jatkamme näin löytääksemme
  jokaiselle \(m\) jonon \((x_{n,m-1})\) äärettömän osajonon
  \((x_{n,m})\) siten, että \((x_{n,m})\) sisältyy palloon, jonka
  halkaisija on \(m^{-1}\).
  %% Nyt \((x_{j,j})_{j=n}^\infty\) on jonon \((x_{j,n})_{j=n}^\infty\)
  %% osajono ja sisältyy siten palloon, jonka halkaisija on
  %% \(n^{-1}\).
%%  Nyt kukin jono \((x_{j,j})_{j=n}^\infty\), \(n \in \naturalnumbers\), sisältyy
%%  palloon, jonka halkaisija on \(n^{-1}\).
  Kun \(n, m \in \naturalnumbers\), \(n < m\), niin \(x_{n,n}\) ja
  \(x_{m,m}\) kuuluvat samaan palloon, jonka halkaisija on \(n^{-1}\).
  Täten
  \begin{displaymath}
    d(x_{i,i},x_{j,j}) \leq \frac{1}{\min \{ i, j \}}
  \end{displaymath}
  kaikille \(i, j \in \naturalnumbers\).
  Siis \((x_{n,n})\) on Cauchyn jono ja koska \(E\) on täydellinen,
  jono \((x_{n,n})\) suppenee avaruudessa \(E\), ja siten jonolla
  \((x_n)\) on suppeneva osajono. Täten \(E\) on jonokompakti, mistä
  seuraa, että se on kompakti.
\end{proof}

\begin{corollary}
  Olkoon \(E\) täydellinen metrinen avaruus ja \(A \subset E\). Joukko
  \(A\) on kompakti, jos ja vain jos se on suljettu ja totaalisesti
  rajoitettu.
\end{corollary}

\index{itseisarvo}

\begin{definition}
  Olkoon \(K\) kunta, jossa on määritelty itseisarvo.
  Funktio \((\lambda, \mu) \mapsto \abs{\lambda-\mu}\) on metriikka
  kunnassa \(K\). Varustettuna tällä metriikalla kuntaa \(K\) kutsutaan
  \index{arvotettu kunta}
  \defterm{arvotetuksi kunnaksi}.
  Sanomme, että \(K\) on \defterm{ei-diskreetti}, jos sen topologia ei
  ole diskreetti.
\end{definition}

\begin{exercises}
  \begin{exercise}{5.1}
    %% L15052015-2
    Olkoot \(\firsttopspace\) ja \(\secondtopspace\) topologisia
    avaruuksia ja \(f : \firsttopspace \to \secondtopspace\)
    funktio. Osoita, että \(f\) on jatkuva, jos ja vain jos jokaisen
    avoimen joukon \(R \subset \secondtopspace\) alkukuva
    \(\inverseimage{f}{\secondtopspace} \subset \firsttopspace\) on
    avoin.
  \end{exercise}
  \begin{exercise}{5.2}
    %% L19052015-3
    Olkoot \((E_1, d_1)\) ja \((E_2, d_2)\) metrisiä avaruuksia
    ja \(f : E_1 \to E_2\) funktio. Osoita, että raja-arvon
    määritelmät \ref{def:raja-arvo} ja \eqref{eq:metr-av-raja-arvo}
    ovat tässä tapauksessa ekvivalentteja.
  \end{exercise}
  \begin{exercise}{5.3}
    %% L23052015-1
    Olkoon \((E, d)\) metrinen avaruus ja \((x_k)_{k=0}^\infty \subset
    E\).
    Osoita, että raja-arvon määritelmät \ref{def:jonon-raja-arvo} ja
    \eqref{eq:metr-av-jonon-raja-arvo} ovat tässä tapauksessa
    ekvivalentteja.
  \end{exercise}
  \begin{exercise}{5.4}
    %% L25052015-2, L26052015-2
    Olkoon \(E\) metrinen avaruus ja \(B \subset E\). Osoita, että \(B\)
    on suljettu, jos ja vain jos jokaiselle jonolle \((x_k)_{k=0}^\infty
    \subset B\) on
    \begin{displaymath}
      x_\alpha \to a \in E \implies a \in B .
    \end{displaymath}
  \end{exercise}
  \begin{exercise}{5.5}
    %% L06062015-1
    Todista lause \ref{th:metr-av-taydellinen}.
  \end{exercise}
\end{exercises}

\chapter{Topologiset vektoriavaruudet}
\label{ch:topologiset-vektoriavaruudet}

Tässä luvussa on käytetty kirjaa \cite{schaefer1971}.

\section{Yleiset topologiset vektoriavaruudet}

\begin{definition}
  Olkoon \(V\) vektoriavaruus ja \(\mytop\) topologia avaruudessa
  \(V\). Sanomme,
  \index{translaatioinvariantti topologia}
  että topologia \(\mytop\) on \defterm{translaatioinvariantti}, jos
  funktio \(x \in V \mapsto x + x_0\) on homeomorfismi jokaiselle \(x_0 \in
  V\).
\end{definition}

\begin{definition}
  Olkoon \(V\) vektoriavaruus ja \(d\) metriikka avaruudessa
  \index{translaatioinvariantti metriikka}
  \(V\). Sanomme, että \(d\) on \defterm{translaatioinvariantti}, jos
  \(d(x, y) = d(x+a,y+a)\) kaikille \(x,y,a \in V\).
\end{definition}

\begin{definition}
  Olkoon \(V\) vektoriavaruus ei-diskreetillä topologisella
  kerroinkunnalla \((K, U)\) (esim. \(\realnumbers\) tai
  \(\complexnumbers\)) ja olkoon \(\mytop\) topologia joukolle
  \(V\). Sanomme, että \index{topologinen vektoriavaruus}
  \((V,+,\cdot,\mytop)\) on \defterm{topologinen vektoriavaruus}, jos
  funktiot \(+ : V \times V \to V\) ja \(\cdot : K \times V \to V\)
  ovat jatkuvia. Avaruuksissa \(V \times V\) ja \(K \times V\)
  käytetään tulotopologioita.
\end{definition}

\begin{definition}
  Olkoon \(V\) topologinen vektoriavaruus.
  Sanomme pisteen \(0 \in V\) ympäristökantaa avaruuden
  \index{0-ympäristökanta}
  \(V\) \defterm{0-ympäristökannaksi}.
\end{definition}

\begin{theorem}
  \cite{schaefer1971}
  Olkoon \((V,+,\cdot)\) vektoriavaruus ja \(\mytop\) siinä määritelty
  topologia.
  Avaruus \((V,+,\cdot,\mytop)\) on topologinen vektoriavaruus,
  ts. operaatiot \(+\) ja \(\cdot\) ovat jatkuvia, jos ja vain jos
  \(\mytop\) on translaatioinvariantti ja
  sillä on 0-ympäristökanta \(B\), joka toteuttaa
  seuraavat ehdot:
  \begin{itemize}
    \item[(a)]
      Jokaiselle \(W \in B\) on olemassa \(U \in B\) siten,
      että \(U + U \subset W\).
    \item[(b)]
      Jokainen \(W \in B\) on radiaalinen ja rengastettu.
    \item[(c)]
      On olemassa \(\lambda \in K\), \(0 < \abs{\lambda} < 1\) siten,
      että ehdosta \(W \in B\) seuraa \(\lambda W \in B\).
  \end{itemize}
  Jos \(K\) on Arkhimedeen kunta, kuten esim. \(\realnumbers\) tai
  \(\complexnumbers\), niin ehto (c) voidaan jättää pois.
\end{theorem}

\begin{proof}
  Tämä todistus on otettu kirjasta \cite{schaefer1971}.  Osoitetaan
  ensin, että vaaditunlainen 0-ympäristökanta on olemassa.  Olkoon \(W
  \subset V\) 0-ympäristö. Operaation \(\cdot\) jatkuvuuden nojalla on
  olemassa 0-ympäristö \(U\) ja reaaliluku \(\varepsilon \in
  \realnumbers_+\) siten, että \(\lambda U \subset W\) aina, kun
  \(\abs{\lambda} \leq \varepsilon\). Koska \(K\) on ei-diskreetti,
  \(V := \union \{ \lambda U \setsep \abs{\lambda} \leq \varepsilon
  \}\) on rengastettu 0-ympäristö ja sisältyy joukkoon \(W\).  Olkoon
  \(B\) kaikkien rengastettujen 0-ympäristöjen joukko avaruudessa
  \(V\). Nyt \(B\) on 0-ympäristökanta. Koska funktio \((\lambda, x_0)
  \mapsto \lambda x_0\) on jatkuva pisteessä \(\lambda = 0\)
  jokaiselle \(x_0 \in V\), niin jokainen \(Z \in B\) on
  radiaalinen. Operaation \(+\) jatkuvuuden nojalla \(B\) toteuttaa
  ehdon (a). Kohdan (c) todistamiseksi riittää huomata, että koska
  \(K\) on ei-diskreetti, niin on olemassa \(\lambda \in K\) jolle \(0
  < \lambda < 1\) ja että 0-ympäristö \(\lambda Z\), \(Z \in
  B\), on rengastettu. Lauseen \cite[1.1]{schaefer1971} nojalla \(\mytop\)
  on translaatioinvariantti.

  Oletetaan sitten, että \(\mytop\) on translaatioinvariantti topologia
  avaruudessa \(V\) siten, että topologialla \(\mytop\) on 0-ympäristökanta
  \(B\), joka toteuttaa ehdot (a), (b) ja (c). On osoitettava, että
  operaatiot \(+\) ja \(\cdot\) ovat jatkuvia.
  On selvää, että \(\{x_0 + Z \setsep Z \in B\}\) on ympäristökanta
  pisteessä \(x_0 \in V\). Siten jos \(Z \in B\) on annettu ja \(U \in
  B\) valitaan niin, että \(U + U \subset Z\), niin
  \begin{displaymath}
    x - x_0 \in U \land y - y_0 \in U \implies x + y \in x_0 + y_0 + Z .
  \end{displaymath}
  Täten operaatio \(+\) on jatkuva.

  Olkoot \(\lambda_0 \in K\) ja \(x_0 \in V\) mielivaltaisia. Jos \(Z
  \in B\) on annettu, niin kohdan (a) nojalla on olemassa \(U \in B\)
  siten, että \(U + U \subset Z\). Kohdan (b) nojalla \(U\) on
  radiaalinen, on olemassa reaaliluku \(\varepsilon \in
  \realnumbers_+\) siten, että \((\lambda - \lambda_0) x_0 \in U\)
  aina, kun \(\abs{\lambda-\lambda_0} \leq \varepsilon\).
  Oletetaan, että \(\mu \in K\) toteuttaa ehdon (c). Tällöin on
  olemassa \(n \in \integernumbers\) siten, että \(\abs{\mu^{-n}} =
  \abs{\mu}^{-n} \geq \abs{\lambda_0} + \varepsilon\).
  Määritellään \(W \in B\) asettamalla \(W := \mu^n U\). Koska \(U\)
  on rengastettu, niin relaatioista \(x-x_0 \in W\) ja
  \(\abs{\lambda-\lambda_0} \leq \varepsilon\) seuraa, että \(\lambda
  (x-x_0) \in U\). Täten
  \begin{displaymath}
    \lambda x = \lambda_0 x_0 + (\lambda - \lambda_0) x_0 + \lambda
    (x-x_0) ,
  \end{displaymath}
  mistä seuraa, että
  \(\lambda x \in \lambda_0 x_0 + U + U \subset \lambda_0 x_0 + Z\),
  mikä osoittaa, että operaatio \(\cdot\) on jatkuva.

  Oletetaan lopuksi, että \(K\) on Arkhimedeen kunta. Tällöin
  \(\abs{2} > 1\) alkiolle \(2 \in K\), mistä seuraa, että \(\abs{2^n}
  = \abs{2}^n > \lambda_0 + \varepsilon\) jollekin \(n \in
  \naturalnumbers\) (merkinnät edellisestä kappaleesta).
  Kohdan (b) toistuvalla soveltamisella voidaan valita \(W \in B\)
  siten, että \(2^n W_1 \subset W_1 + \ldots + W_1 \subset U\), missä
  summalla on \(2^n\) yhteenlaskettavaa. Koska \(W_1\) ja siten \(2^n
  W_1\) ovat rengastettuja, \(W_1\) voidaan sijoittaa joukon \(W\)
  paikalle edellisen kappaleen todistuksessa (operaation \(\cdot\)
  jatkuvuus), ja siten ehto (c) voidaan jättää pois tässä
  tapauksessa.
\end{proof}

%% Seuraava määrittely ei välttämättä ole aina ekvivalentti metrisessä
%% avaruudessa määritellyn rajoitettuuden kanssa.

%% \begin{definition}
%%   Olkoon \(V\) topologinen vektoriavaruus kerroinkunnalla
%%   \(K\). Olkoon \(A \subset V\). Sanomme, että \(A\) on
%%   \index{rajoitettu joukko}
%%   \defterm{rajoitettu}, jos jokaiselle \(0\)-ympäristölle \(U\)
%%   avaruudessa \(V\) on olemassa \(\lambda \in K\) siten, että \(A
%%   \subset \lambda U\).
%% \end{definition}

Kun 0-ympäristökanta \(B\) on annettu, niin määritellään topologia
\(\mytop\) avaruudessa \(V\) seuraavasti:
\begin{equation}
  \label{eq:y-kanta}
  A \in \mytop \iff \forall x \in A : \exists y \in V, Y \in B : x \in y
  + Y,
\end{equation}
missä \(y + Y := \{ y + x \setsep x \in Y \}\).

\begin{definition}
  Olkoot \((V, \firsttop)\) ja \((W, \secondtop)\) topologisia
  vektoriavaruuksia. Olkoon
  \index{topologinen isomorfismi}
  \(f : V \to W\) funktio. Sanomme, että \(f\) on \defterm{topologinen
    isomorfismi}, jos \(f\) on algebrallinen isomorfismi
  vektoriavaruudelta \(V\) vektoriavaruudelle \(W\) ja homeomorfismi
  topologiselta avaruudelta \((V, \firsttop)\) topologiselle
  avaruudelle \((W, \secondtop)\).
  Jos on olemassa topologinen isomorfismi avaruudelta \((V, \firsttop)\)
  avaruudelle \((W, \secondtop)\), niin sanomme, että \((V,
  \firsttop)\) ja \((W, \secondtop)\)
  \index{topologisesti isomorfiset avaruudet}
  ovat \defterm{topologisesti isomorfiset}.
\end{definition}

\begin{definition}
  Kun \(V\) ja \(W\) ovat topologisia vektoriavaruuksia, niin
  \index{$\equaltvs$}
  merkintä \(V \equaltvs W\) tarkoittaa, että \(V\) ja \(W\) ovat
  samat topologisina vektoriavaruuksina, ts. \(V\) ja \(W\) ovat sama
  vektoriavaruus ja niiden topologiat ovat samat.
\end{definition}

\begin{definition}
  Olkoon \((V, +, \cdot, \mytop)\) topologinen vektoriavaruus. Määritellään
  \index{duaali}
  \index{topologinen duaali}
  vektoriavaruuden \(V\) \defterm{(topologinen) duaali} asettamalla
  \begin{displaymath}
    V^* := \{ \tilde{f} \in V^\sharp \setsep \tilde{f} \;\textrm{on
      jatkuva} \} .
  \end{displaymath}
\end{definition}

\begin{definition}
  \cite{wpweaktop}
  Olkoon \(V\) topologinen vektoriavaruus ja \(F \subset V^\sharp\).
  \index{heikko topologia}
  Määritellään avaruuden \(V\) \defterm{heikko topologia} joukon \(F\)
  suhteen, merkitään \(\sigma(V, F)\), siten, että se on karkein
  topologia, jossa kaikki joukon \(F\) funktionaalit ovat jatkuvia.
  Topologiaa \(\sigma(V, V^*)\) sanotaan avaruuden \(V\) heikoksi
  topologiaksi.
\end{definition}

Olkoon \(V\) topologinen vektoriavaruus, \(F \subset V^\sharp\) ja
\(x_\alpha \in V\) verkko. Nyt \(x_\alpha \to x\) topologiassa
\(\sigma(V,F)\), jos ja vain jos \(a(x_\alpha) \to a(x)\) kaikilla \(a
\in F\).

Jos \(V\) on vektoriavaruus ja \(x \in V\), niin vektori \(x\) voidaan
samaistaa avaruuden \(V^{**}\) alkioon \(T_x\) asettamalla
\begin{math}
  T_x(\tilde{y}) = \tilde{y}(x)
\end{math}
kaikille \(\tilde{y} \in V^*\).
Täten \(V \subset V^{**}\).

\begin{definition}
  Olkoon \(V\) topologinen vektoriavaruus. Avaruuden \(V^*\) topologiaa
  \(\sigma(V^*,V)\) sanotaan avaruuden \(V^*\)
  \index{heikko-*-topologia}
  \defterm{heikko-*-topologiaksi}.
\end{definition}

Olkoon \(V\) topologinen vektoriavaruus ja \(\tilde{x}_\alpha \in
V^*\) verkko. Nyt \(\tilde{x}_\alpha \to \tilde{x}\) topologiassa
\(\sigma(V^*,V)\), jos ja vain jos \(\tilde{x}_\alpha(y) \to
\tilde{x}(y)\) kaikilla \(y \in V\).

\hyphenation{translaatio-invariantti}

\begin{definition}
  Olkoon \(A\) vektoriavaruus ja \(\mytop\) topologia avaruudessa
  \(A\). Sanomme, että \((A,\mytop)\) on \textit{F}-avaruus, jos se on
  täydellinen ja \(\mytop\) on translaatioinvariantti ja metrisoituva.
\end{definition}

\section{Lokaalikonveksit avaruudet}

Tässä luvussa oletetaan, että kerroinkunta \(K = \realnumbers\) tai
\(K = \complexnumbers\).

\begin{definition}
  \label{def:minkowski}
  Olkoon \(V\) vektoriavaruus ja \(A \subset V\)
  \index{Minkowskin funktionaali}
  absorboiva. Määritellään joukon \(A\) \defterm{Minkowskin
    funktionaali}
  \begin{displaymath}
    \mu_A(x) := \inf \{ t > 0 \setsep t^{-1}x \in A \} .
  \end{displaymath}
\end{definition}

\begin{remark}
  Koska edellisessä määritelmässä \(A\) on absorboiva, niin \(\mu_A(x)
  < \infty\) kaikille \(x \in V\).
\end{remark}

\begin{definition}
  Olkoon \(V\) vektoriavaruus kerroinkunnalla \(K\). Sanomme, että
  \index{seminormi}
  funktio \(\rho : V \to \realnumbers_0\) on \defterm{seminormi}, jos
  seuraavat aksioomat ovat voimassa:
  \begin{itemize}
    \item[(S1)] \(\rho(x) \geq 0\) kaikille \(x \in V\).
    \item[(S2)] \(\rho(ax) = \abs{a} \rho(x)\) kaikille \(a \in K\),
      \(x \in V\).
    \item[(S3)] \(\rho(x-y) \leq \rho(x-z) + \rho(z-y)\) kaikille
      \(x,y,z \in V\).
  \end{itemize}
\end{definition}

\begin{definition}
  Olkoon \(P\) joukko seminormeja vektoriavaruudessa \(V\). Sanomme,
  että \(P\) on \defterm{separoiva}, jos jokaiselle \(x \in V
  \setminus \{0\}\) on olemassa \(p \in P\) siten, että \(p(x)
  \not=0\).
\end{definition}

\begin{theorem}
  \label{th:seminorm-a}
  %% Rudin Th. 1.34(e)
  Olkoon \(V\) vektoriavaruus ja \(p\) seminormi avaruudessa
  \(V\). Nyt \(\abs{p(x)-p(y)} \leq p(x-y)\) kaikille \(x, y \in V\).
\end{theorem}

\begin{proof}
  Tämä todistus on otettu kirjasta \cite{rudin1991}.
  Olkoon \(x,y \in V\). Seminormien subadditiivisuudesta seuraa, että
  \(p(x) = p(x-y+y) \leq p(x-y)+p(y)\), joten \(p(x)-p(y) \leq
  p(x-y)\). Vastaavasti saadaan \(p(y)-p(x) \leq p(y-x) =
  p(x-y)\). Täten \(\abs{p(x)-p(y)} \leq p(x-y)\).
\end{proof}

\begin{theorem}
  %% Rudin Th. 1.35(c)
  \label{th:seminorm-b}
  Olkoon \(V\) vektoriavaruus ja \(A \subset V\) konveksi, absorboiva
  ja balansoitu. Nyt \(\mu_A\) on seminormi.
\end{theorem}

\begin{proof}
  Tämä todistus on otettu kirjasta \cite{rudin1991}.
  Olkoon \(x, y \in V\) ja \(\varepsilon \in \realnumbers_+\). Olkoon
  \(t := \mu_A(x) + \varepsilon\) ja \(s := \mu_A(y) +
  \varepsilon\). Nyt \(x/t \in A\) ja \(y/s \in A\). Siten myös niiden
  konveksi kombinaatio
  \begin{displaymath}
    \frac{x+y}{s+t} = \frac{t}{s+t} \cdot \frac{x}{t} + \frac{s}{s+t}
    \cdot \frac{y}{s}
  \end{displaymath}
  kuuluu joukkoon \(A\). Täten \(\mu_A(x+y) \leq s + t = \mu_A(x) +
  \mu_A(y) + 2 \varepsilon\), mistä seuraa, että \(\mu_A(x+y) \leq
  \mu_A(x) + \mu_A(y) \).
  Minkowskin funktionaalin määritelmästä \ref{def:minkowski} seuraa,
  että \(\mu_A(tx) = t\mu_A(x)\) kaikille \(x \in V\) ja \(t \geq 0\).
  Olkoon \(x \in V\) ja \(s \in \realnumbers\). Meillä on
  \begin{displaymath}
    \mu_A(sx) = \inf \{ t > 0 \setsep t^{-1}sx \in A \} .
  \end{displaymath}
  Oletetaan, että \(s < 0\)-
  Koska \(A\) on balansoitu, niin \(t^{-1}sx =-t^{-1}\abs{s}x \in A
  \iff t^{-1}\abs{s}x \in A\), mistä seuraa, että \(\mu_A(sx) =
  \abs{s} \mu_A(x)\).
\end{proof}

\begin{definition}
  Olkoon \(V\) topologinen vektoriavaruus kerroinkunnalla \(K\) ja
  topologialla \(\mytop\). Sanomme, että
  \index{lokaalikonveksi topologia}
  \(\mytop\) on \defterm{lokaalikonveksi topologia} jos topologialla \(\mytop\)
  on konveksi balansoitu 0-ympäristökanta.
  Jos lisäksi \((V, \mytop)\) on Hausdorffin avaruus, niin sanomme, että
  \index{lokaalikonveksi avaruus}
  \(V\) on \defterm{lokaalikonveksi avaruus}.
\end{definition}

\begin{theorem}
  %% L11062015-1 (Rudin th. 1.36, 1.37)
  Olkoon \(V\) topologinen vektoriavaruus topologialla \(\mytop\).
  Olkoon \(P\) separoiva joukko seminormeja avaruudessa \(V\). Olkoon
  \begin{displaymath}
    V(p,n) := \left\{ x \in V \bigsetsep p(x) < \frac{1}{n} \right\} ,
  \end{displaymath}
  missä \(p \in P\) ja \(n \in \positiveintegers\).
  Olkoon \(B\) joukko kaikista joukkojen \(V(p,n)\) äärellisistä
  leikkauksista. Nyt \(V\) on lokaalikonveksi avaruus, jos ja vain jos
  \(B\) on topologian \(\mytop\) 0-ympäristökanta.
\end{theorem}

\hyphenation{translaatio-invariantin}
\hyphenation{lokaali-konveksi}
\hyphenation{lokaali-konvekseissa}

\begin{proof}
  Tässä todistuksessa on käytetty lauseiden \cite[1.36 ja
    1.37]{rudin1991} todistuksia.
  Oletetaan ensin, että \(V\) on lokaalikonveksi avaruus. Olkoon \(B\)
  avaruuden \(V\) konveksi ja balansoitu 0-ympäristökanta.
  Olkoon \(A \in B\). Jos \(x \in A\), niin \(x/t \in A\) jollakin \(t
  < 1\), koska \(A\) on avoin. Täten \(\mu_A(x) < 1\). Jos \(x \not\in
  A\), niin väitteestä \(x/t \in A\) seuraa \(t \geq 1\), koska \(A\)
  on balansoitu. Täten \(\mu_A(x) \geq 1\). Näin ollen
  \begin{equation}
    \label{eq:c}
    A = \{ x \in V \setsep \mu_A(x) < 1 \} .
  \end{equation}
  Lauseesta \ref{th:seminorm-b} seuraa, että \(\mu_A\) on
  seminormi. Jos \(r > 0\), niin yhtälöstä \eqref{eq:c} ja lauseesta
  \ref{th:seminorm-a} seuraa, että
  \begin{displaymath}
    \abs{\mu_A(x)-\mu_A(y)} \leq \mu_A(x-y) < r ,
  \end{displaymath}
  jos \(x-y \in rA\). Täten \(\mu_A\) on jatkuva. Jos \(x \in V\) ja
  \(x \not= 0\), niin \(x \not\in A_1\) jollekin \(A_1 \in B\) ja
  \(\mu_{A_1}(x) \geq 1\). Siis \(\{\mu_A\}\) on separoiva.

  Oletetaan sitten, että lauseen ehdot seminormijoukolle \(P\) ovat
  voimassa.
  Määritellään joukko \(C \subset A\) avoimeksi, jos ja vain jos \(C\)
  on (mahdollisesti tyhjä) unioni kannan \(B\) alkioiden
  translaatioista. Tämä määrittelee translaatioinvariantin topologian
  \(\tau\) avaruudessa \(A\). Jokainen kannan \(B\) jäsen on konveksi
  ja balansoitu ja \(B\) on topologian \(\tau\) 0-ympäristökanta.

  Olkoon \(x \in A\), \(x \not= 0\). Nyt \(p(x) > 0\) jollakin \(p \in
  P\). Koska \(x \not\in V(p,n)\), jos \(np(x) > 1\), niin \(0\) ei
  ole pisteen \(x\) ympäristössä \(x-V(p,n)\), joten \(x\) ei kuulu
  joukon \(\{0\}\) sulkeumaan. Siten \(\zeroset\) on suljettu joukko,
  ja koska \(\tau\) on translaatioinvariantti jokainen joukon \(A\)
  piste on suljettu.

  Seuraavaksi näytetään, että yhteenlasku ja skalaarilla kertominen
  ovat jatkuvia. Olkoon \(U\) pisteen \(0\) ympäristö avaruudessa
  \(A\). Nyt
  \begin{displaymath}
    V(p_1, n_1) \cap \ldots \cap V(p_m,n_m) \subset U
  \end{displaymath}
  joillekin \(p_1,\ldots,p_m \in P\) ja \(n_1,\ldots,n_m \in
  \positiveintegers\). Asetetaan
  \begin{displaymath}
    V := V(p_1,2n_1) \cap \ldots \cap V(p_m,2n_m) .
  \end{displaymath}
  Koska jokainen \(p \in P\) on subadditiivinen, niin \(V + V \subset
  U\). Tämä osoittaa, että yhteenlasku on jatkuva.
  Oletetaan, että \(x \in A\), \(\alpha \in K\), ja \(U\) ja \(V\)
  ovat kuten määritelty yllä. Nyt \(x \in sV\) jollekin \(s >
  0\). Olkoon
  \begin{displaymath}
    t := \frac{s}{1+\abs{\alpha}s} .
  \end{displaymath}
  Jos \(y \in x + tV\) ja \(\abs{\beta-\alpha} < 1/s\), niin
  \begin{displaymath}
    \beta y - \alpha x = \beta (y-x) + (\beta - \alpha) x
  \end{displaymath}
  joka kuuluu joukkoon
  \begin{displaymath}
    \abs{\beta} t V + \abs{\beta-\alpha} s V \subset V + V \subset U ,
  \end{displaymath}
  koska \(\abs{\beta} t \leq 1\) ja \(V\) on balansoitu. Tämä
  osoittaa, että skalaarilla kertominen on jatkuva.

  Siten \(A\) on lokaalikonveksi avaruus. Joukon \(V(p,n)\)
  määritelmä osoittaa, että jokainen \(p \in P\) on jatkuva pisteessä
  \(0\), mistä seuraa, että \(p\) on jatkuva avaruudessa \(A\).
\end{proof}

\begin{remark}
  Ehto \(P\) on separoiva on yhtäpitävä sen kanssa, että \(\mytop\) on
  Hausdorffin topologia.
\end{remark}

Olkoon \(V\) lokaalikonveksi avaruus seminormiperheellä
\((\rho_\alpha)_{\alpha \in I}\), \((x_\gamma)_{\gamma \in G}\) verkko
avaruudessa \(V\) ja \(x \in V\). Olkoon \(y_{\alpha,\gamma} :=
\rho_\alpha(x_\gamma - x)\) kaikille \(\gamma \in G\) ja \(\alpha \in
I\). Tällöin \(x_\gamma \to x\), jos ja vain jos verkko
\((y_{\alpha,\gamma})_{\gamma \in G}\) suppenee kohti lukua 0 kaikille
\(\alpha \in I\).

Olkoon \(V\) lokaalikonveksi avaruus seminormiperheellä
\((\rho_\alpha)_{\alpha \in I}\) ja
\(W\) lokaalikonveksi avaruus seminormiperheellä
\((\eta_\beta)_{\beta \in J}\). Olkoon \(f : V \to W\) funktio.
Lokaalikonvekseissa avaruuksissa funktion raja-arvo saadaan
seuraavasti: \[\lim_{x \to x_0} f(x) = a,\] jos ja vain jos
jokaiselle \(\varepsilon \in \realnumbers_+\) ja jokaiselle
äärelliselle joukolle seminormeja \(\eta_{\beta(1)}, \ldots,
\eta_{\beta(m)}\) on olemassa \(\delta \in \realnumbers_+\) ja
äärellinen joukko seminormeja \(\rho_{\alpha(1)}, \ldots,
\rho_{\alpha(n)}\) siten, että
\begin{align*}
  & \rho_{\alpha(1)}(x-x_0) < \delta \;\land\; \ldots \;\land\;
  \rho_{\alpha(n)}(x-x_0) < \delta \\
  & \implies
  \eta_{\beta(1)}(f(x)-a) < \varepsilon \;\land\; \ldots \;\land\;
  \eta_{\beta(m)}(f(x)-a) < \varepsilon .
\end{align*}

Jos \(V\) on topologinen vektoriavaruus kerroinkunnalla \(K =
\realnumbers\) tai \(K = \complexnumbers\) ja \(F \subset V^\sharp\),
niin \(\sigma(V,F)\) on lokaalikonveksi topologia, jonka generoivat
seminormit
\begin{displaymath}
  \norm{x}_{\tilde{f}} := \abs{\tilde{f}(x)} , \spaceafter x \in V,
\end{displaymath}
kaikille \(\tilde{f} \in F\).
Erityisesti avaruuden \(V^*\) heikko-*-topologian generoivat
seminormit
\begin{displaymath}
  \norm{\tilde{g}}_x := \abs{\tilde{g}(x)}, \spaceafter \tilde{g} \in V^* ,
\end{displaymath}
kaikille \(x \in V\).

\begin{definition}
  \label{def:vahva-topologia}
  Olkoon \(V\) lokaalikonveksi avaruus. Duaaliavaruuden \(V^*\)
  \index{vahva topologia}
  \defterm{vahva topologia} määritellään seminormien
  \begin{displaymath}
    \norm{f}_B := \sup_{x \in B} \abs{f(x)}, \spaceafter f \in V^*, B
    \in S,
  \end{displaymath}
  generoimana topologiana,
  missä
  \begin{displaymath}
    S := \{ \textrm{kaikki avaruuden}\; V \;\textrm{rajoitetut
      joukot}\} .
  \end{displaymath}
  Avaruutta \(V^*\) varustettuna tällä topologialla kutsutaan
  \index{vahva duaaliavaruus}
  avaruuden \(V\) \defterm{vahvaksi duaaliavaruudeksi}
  \cite{wpstrongtop}.
\end{definition}

\section{Normiavaruudet}

\hyphenation{normi-avaruus}

\begin{definition} 
  Olkoon \(V\) vektoriavaruus kerroinkunnalla \(K = \realnumbers\) tai
  \(K = \complexnumbers\). Olkoon \(\norm{\cdot} : V \to
  \realnumbers_0\) funktio. Sanomme, että \(\norm{\cdot}\) on
  \index{normiavaruus}
  \index{normi}
  \defterm{normi} ja \(V\) \defterm{normiavaruus}, jos seuraavat
  aksioomat ovat voimassa:
  \begin{itemize}
    \item[(N1)] \(\norm{x} = 0 \iff x = 0\) kaikille \(x \in V\).
    \item[(N2)] \(\norm{ax} = \abs{a} \norm{x}\) kaikille \(a \in K\),
      \(x \in V\).
    \item[(N3)] \(\norm{x-y} \leq \norm{x-z} + \norm{z-y}\) kaikille
      \(x,y,z \in V\).
  \end{itemize}
\end{definition}

\index{kolmioepäyhtälö}
Aksiooma (N3) on yksi muoto kolmioepäyhtälöstä. Jo normiavaruus on
\index{Banachin avaruus}
täydellinen, niin sitä kutsutaan \defterm{Banachin avaruudeksi}.

\begin{examples}
  \begin{example}
    Jatkuvien ja rajoitettujen funktioiden \(f : \realnumbers^n \to
    \complexnumbers\) muodostama avaruus \(\Cb{\realnumbers^n}\)
    varustettuna normilla
    \begin{displaymath}
      \norm{f} := \sup_{\rnvec{x} \in \realnumbers^n}
      \abs{f(\rnvec{x})}, \spaceafter
      f \in \Cb{\realnumbers^n},
    \end{displaymath}
    on Banachin avaruus. Ks. luku \ref{sec:Cb}.
  \end{example}
  \begin{example}
    Vektoriavaruus \((\realnumbers,\realnumbers,+,\cdot)\)
    varustettuna normille \(\norm{x} := \abs{x}\), \(x \in
    \realnumbers\), on Banachin avaruus.
  \end{example}
  \begin{example}
    Vektoriavaruus \((\complexnumbers,\complexnumbers,+,\cdot)\)
    varustettuna normille \(\norm{x} := \abs{x}\), \(x \in
    \complexnumbers\), on Banachin avaruus.
  \end{example}
  \begin{example}
    Vektoriavaruus \(\realnumbers^n\), \(n \in \positiveintegers\),
    varustettuna normilla
    \begin{displaymath}
      \norm{\rnvec{x}} := \sqrt{\sum_{k=1}^n x_k^2}, \spaceafter
      \rnvec{x} \in \realnumbers^n,
    \end{displaymath}
    on Banachin avaruus.
  \end{example}
  \begin{example}
    Vektoriavaruus \(\complexnumbers^n\), \(n \in \positiveintegers\),
    varustettuna normilla
    \begin{displaymath}
      \norm{\cnvec{x}} := \sqrt{\sum_{k=1}^n \abs{x_k}^2}, \spaceafter
      \cnvec{x} \in \complexnumbers^n,
    \end{displaymath}
    on Banachin avaruus.
  \end{example}
\end{examples}

\index{metrinen avaruus}
\index{metriikka}
Jokainen normiavaruus \(V\) on metrinen avaruus, kun metriikka
määritellään
\begin{displaymath}
  d(x,y) := \norm{x-y}, \spaceafter x, y \in V .
\end{displaymath}

\begin{definition}
  Olkoon \(V\) vektoriavaruus ja \(\norm{\cdot}_1\) ja
  \(\norm{\cdot}_2\) kaksi avaruudessa \(V\) määriteltyä
  normia. Sanomme, että normit \(\norm{\cdot}_1\) ja
  \index{ekvivalentit normit}
  \(\norm{\cdot}_2\) ovat \defterm{ekvivalentit}, jos on olemassa
  luvut \(a,b \in \realnumbers_+\) siten, että
  \begin{displaymath}
    a \norm{x}_1 \leq \norm{x}_2 \leq b \norm{x}_1
  \end{displaymath}
  kaikille \(x \in V\).
\end{definition}

\begin{theorem}
  %% L24052015-1
  Olkoon \(V\) vektoriavaruus ja \(\norm{\cdot}_1\) ja
  \(\norm{\cdot}_2\) kaksi avaruudessa \(V\) määriteltyä
  normia. Normit \(\norm{\cdot}_1\) ja \(\norm{\cdot}_2\) generoivat
  saman topologian, jos ja vain jos ne ovat ekvivalentit.
\end{theorem}

\begin{proof}
  Olkoot \(\tau_1\) ja \(\tau_2\) normien \(\norm{\cdot}_1\) ja
  \(\norm{\cdot}_2\) generoimat topologiat ja avaruudet \(V_1\) ja
  \(V_2\) normiavaruuksia \(V\) varustettuna kullakin topologialla.

  Oletetaan ensin, että \(\tau_1 = \tau_2\).
  Oletetaan, että olisi \( \lnot \exists a \in \realnumbers_+ : \forall
  x \in V : a \norm{x}_1 \leq \norm{x}_2 \) (vastaoletus 1).
  Jokaisella \(k \in \naturalnumbers\) on olemassa \(x_k \in V\)
  siten, että
  \begin{math}
    \norm{x_k}_1 = 1/(k+1)
  \end{math}
  ja
  \(\norm{x_k}_2 \geq (k+1) \norm{x_k}_1 = 1\).
  Nyt \(x_k \to 0\) topologiassa \(T_1\), mutta \(x_k \not\to 0\)
  topologiassa \(T_2\). Tämä on ristiriita, joten vastaoletus 1 on
  väärä.

  Oletetaan, että olisi \( \lnot \exists b \in \realnumbers_+ : \forall
  x \in V : \norm{x}_2 \leq b \norm{x}_1 \) (vastaoletus 2).
  Jokaisella \(k \in \naturalnumbers\) on olemassa \(x_k \in V\)
  siten, että
  \begin{math}
    \norm{x_k}_2 = 1/(k+1)
  \end{math}
  ja
  \(\norm{x_k}_1 \geq (k+1) \norm{x_k}_2 = 1\).
  Nyt \(x_k \to 0\) topologiassa \(\tau_2\), mutta \(x_k \not\to 0\)
  topologiassa \(\tau_1\). Tämä on ristiriita, joten vastaoletus 2 on
  väärä.

  Oletetaan sitten, että kaikille \(x \in V\) on \(a \norm{x}_1 \leq
  \norm{x}_2 \leq b \norm{x}_1\), missä \(a, b \in \realnumbers_+\)
  ovat vakioita.
  Olkoon \(A \in \tau_1\) ja \(x \in A\). Nyt \(\openball{V_1}{x}{r_1}
  \subset A\) jollakin \(r_1 \in \realnumbers_+\). Olkoon \(y \in
  \openball{V_2}{x}{ar_1}\). Nyt \(\norm{y-x}_2 < ar_1\) ja
  \(\norm{y-x}_1 \leq \frac{1}{a} \norm{y-x}_2 < r_1\). Siis
  \(\openball{V_2}{x}{ar_1} \subset \openball{V_1}{x}{r_1} \subset
  A\). Koska \(x \in A\) oli mielivaltainen, niin lauseen
  \ref{th:avoin-joukko} nojalla \(A \in \tau_2\).

  Oletetaan edelleen, että normit ovat ekvivalentteja. Olkoon \(A \in
  \tau_2\) ja \(x \in A\). Nyt \(\openball{V_2}{x}{r_2}
  \subset A\) jollakin \(r_2 \in \realnumbers_+\). Olkoon \(y \in
  \openball{V_1}{x}{r_2/b}\). Nyt \(\norm{y-x}_1 < r_2 / b\) ja
  \(\norm{y-x}_2 \leq b \norm{y-x}_1 < r_2\). Siis
  \(\openball{V_1}{x}{r_2/b} \subset \openball{V_2}{x}{r_2} \subset
  A\). Koska \(x \in A\) oli mielivaltainen, niin lauseen
  \ref{th:avoin-joukko} nojalla \(A \in \tau_1\).
\end{proof}

Normiavaruuksille rajoitetun joukon määritelmä voidaan kirjoittaa
seuraavaan muotoon:

\begin{definition}
  Kun \(V\) on normiavaruus ja \(A \subset V\), \(A \not= \emptyset\),
  \index{rajoitettu joukko}
  niin sanomme, että \(A\) on \defterm{rajoitettu}, jos on olemassa
  luku \(M \in \realnumbers_+\) siten, että \(\norm{x} \leq M\)
  kaikille \(x \in A\).
\end{definition}

\begin{definition}
  Olkoon \(X\) epätyhjä joukko, \(V\) normiavaruus ja \(f : X \to V\)
  funktio. Sanomme, että \(f\) on
  \index{rajoitettu funktio}
  \defterm{rajoitettu}, jos on olemassa luku \(m \in \realnumbers\)
  siten, että \(\norm{f(x)} \leq m\) kaikilla \(x \in X\).
\end{definition}

\begin{theorem}
  %% L23052015-2
  Olkoot \((V, \norm{\cdot}_V)\) ja \((W, \norm{\cdot}_W)\)
  normiavaruuksia ja \(f : V \to W\) lineaarinen funktio. Tällöin
  \(f\) on jatkuva, jos ja vain jos se on rajoitettu.
\end{theorem}

\begin{proof}
  Oletetaan ensin, että \(f\) on rajoitettu. Nyt \(\norm{f(x)}_W \leq
  m \norm{x}_V\) kaikille \(x \in V\), missä \(m \in \realnumbers_+\)
  on vakio. Oletetaan, että \((x_k) \subset V\) ja \(x_k \to y \in
  V\), kun \(k \to \infty\). Nyt \(\norm{f(x_k) - f(y)}_W =
  \norm{f(x_k-y)}_W \leq m \norm{x_k-y}_V \to 0\), kun \(k \to
  \infty\). Siten lauseen \ref{th:top-av-suppeneminen-b} nojalla \(f\)
  on jatkuva.

  Oletetaan sitten, että \(f\) on jatkuva. Oletetaan, että \(f\) ei
  olisi rajoitettu (vastaoletus). Koska \(f\) ei ole rajoitettu, niin
  jokaisella \(r \in \realnumbers_+\) on olemassa \(x \in V\) siten,
  että \(\norm{x}_V = 1\) ja \(\norm{f(x)}_W > r\). Kun \(k \in
  \naturalnumbers\), niin valitaan \(x_k\) siten, että \(\norm{x_k}_V
  = 1\) ja \(\norm{f(x_k)}_W > k + 1\). Nyt
  \begin{displaymath}
    \sznorm{\frac{x_k}{k+1}}_V = \frac{1}{k+1} \to 0 ,
  \end{displaymath}
  kun \(k \to 0\), mutta
  \begin{displaymath}
    \sznorm{f\left(\frac{x_k}{k+1}\right)}_W
    = \frac{1}{k+1} \norm{f(x_k)}_W
    > \frac{1}{k+1} \cdot ( k + 1)
    = 1 .
  \end{displaymath}
  Siis
  \begin{displaymath}
    f\left(\frac{x_k}{k+1}\right) \not\to 0 = f(0) ,
  \end{displaymath}
  mikä on ristiriidassa lauseen \ref{th:top-av-suppeneminen-b} kanssa,
  joten vastaoletus on väärä ja f on rajoitettu.
\end{proof}

\begin{theorem}
  \cite[lause 5.2.4]{hpc2005}
  Olkoon \(V\) Banachin avaruus ja \(A \subset V\).
  \begin{itemize}
    \item[(i)] Jos \(A\) on kompakti, niin se on suljettu ja rajoitettu.
    \item[(ii)] Jos \(A\) on suljettu ja rajoitettu ja \(V\) on
      äärellisulotteinen, niin \(A\) on kompakti.
  \end{itemize}
\end{theorem}

\begin{definition}
  Olkoot \(V\) ja \(W\) normiavaruuksia. Sanomme, että funktio \(\iota
  \index{isometrinen isomorfismi}
  : V \onto W\) on \defterm{isometrinen isomorfismi}, jos se on
  algebrallinen isomorfismi vektoriavaruudelta \(V\)
  vektoriavaruudelle \(W\) ja lisäksi isometria metriseltä avaruudelta
  \(V\) metriselle avaruudelle \(W\). Jos on olemassa isometrinen
  isomorfismi normiavaruudelta \(V\) normiavaruudelle \(W\), niin
  \index{isometrisesti isomorfiset avaruudet}
  sanomme, että \(V\) ja \(W\) ovat \defterm{isometrisesti
    isomorfiset}.
\end{definition}

\begin{remark}
  Isometrisesti isomorfiset normiavaruudet voidaan samaistaa
  keskenään.
\end{remark}

%% \begin{definition}
%%   Olkoot \(A\) ja \(B\) normiavaruuksia ja \(A \subset B\).
%%   Jos \(\norm{a}_A = \norm{a}_B\) kaikille \(a \in A\), niin sanomme,
%%   \index{isometrinen upotus}
%%   että \(A\) on \defterm{isometrisesti upotettu} avaruuteen \(B\) ja
%%   merkitsemme tätä \(A \isomemb B\).
%% \end{definition}

\begin{definition}
  Olkoot \(V\) ja \(W\) Banachin avaruuksia.
  Jos \(V\) on isometrisesti isomorfinen jonkin avaruuden \(W\)
  osajoukon kanssa, niin sanomme,
  \index{isometrinen upotus}
  että \(V\) on \defterm{isometrisesti upotettu} avaruuteen \(W\) ja
  merkitsemme tätä \(V \isomemb W\).
\end{definition}

%% TARK. määritelmät
\begin{definition}
  Olkoon \(V\) Banachin avaruus kerroinkunnalla \(K\) ja
  \((e_k)_{k=0}^\infty \subset V\).  Sanomme, että \((e_k)\) on
  \index{Schauderin kanta}
  avaruuden \(V\) \defterm{Schauderin kanta}, jos \((e_k)\) on
  lineaarisesti riippumaton ja kaikille \(f \in V\) on olemassa
  skalaarit \(a_k \in K\), \(k \in \naturalnumbers\), siten, että
  \begin{equation}
    \label{eq:schauder}
    f = \sum_{k=0}^\infty a_k e_k .
  \end{equation}
  Jos sarja \eqref{eq:schauder} suppenee riippumatta termien
  järjestyksestä, ts.
  \begin{displaymath}
    \sum_{k=0}^\infty a_k e_k = \sum_{k=0}^\infty a_{\sigma(k)} e_{\sigma(k)}
  \end{displaymath}
  kaikille bijektioille \(\sigma : \naturalnumbers \onto
  \index{ehdoton kanta}
  \naturalnumbers\) ja kaikille \(f \in V\), niin sanomme, että \((e_k)\) on
  \defterm{ehdoton kanta}. Muussa tapauksessa sanomme, että \((e_k)\) on
  \index{ehdollinen kanta}
  \defterm{ehdollinen kanta}.
  %% ITSEINEN KANTA
\end{definition}

\begin{definition}
  Olkoon \(p \in [1, \infty]\).
  Kun \(p < \infty\), määritellään Banachin avaruus
  \index{$l^p$}
  \begin{displaymath}
    l^p := \left\{ \seq{a} \in \complexnumbers^\naturalnumbers \bigsetsep
    \norm{\seq{a}} := \left(\sum_{k=0}^\infty
    \abs{\seqelem{\seq{a}}{k}}^p\right)^{\frac{1}{p}} < \infty
    \right\} .
  \end{displaymath}
  Kun \(p = \infty\), määritellään Banachin avaruus
  \begin{displaymath}
    l^\infty := \left\{ \seq{a} \in \complexnumbers^\naturalnumbers \bigsetsep
    \norm{\seq{a}} := \sup_{k \in \naturalnumbers}
    \abs{\seqelem{\seq{a}}{k}}  < \infty
    \right\} .
  \end{displaymath}
\end{definition}

\begin{definition}
  Määritellään Banachin avaruus \(c_0\) seuraavasti:
  \index{$c_0$}
  \begin{displaymath}
    c_0 := \left\{ \seq{a} \in \complexnumbers^\naturalnumbers \bigsetsep
    \lim_{k \to \infty} \seqelem{\seq{a}}{k} = 0 \right\}
  \end{displaymath}
  ja
  \begin{displaymath}
    \norm{\seq{a} \;\vert\; c_0} := \sup_{k \in \naturalnumbers}
    \abs{\seqelem{\seq{a}}{k}} .
  \end{displaymath}
\end{definition}

\section{Sisätuloavaruudet}

\begin{definition} 
  Olkoon \(V\) vektoriavaruus kerroinkunnalla \(K = \realnumbers\) tai
  \(K = \complexnumbers\). Olkoon \(\ip{\cdot}{\cdot} : V \times V \to
  K\) funktio. Sanomme, että \(\ip{\cdot}{\cdot}\) on
  \index{sisätulo}
  \index{sisätuloavaruus}
  \defterm{sisätulo} ja \(V\) \defterm{sisätuloavaruus}, jos seuraavat
  aksioomat ovat voimassa:
  \begin{itemize}
    \item[(ST1)] \(\ip{x}{x} \in \realnumbers_0\) kaikille \(x \in
      V\).
    \item[(ST2)] \(\ip{x}{x} = 0 \iff x = 0\) kaikille \(x \in
      V\).
    \item[(ST3)] \(\ip{x}{y} = \ip{y}{x}^*\) kaikille \(x, y \in V\).
    \item[(ST4)] \(\ip{x}{ay} = a \ip{x}{y}\) kaikille \(x, y \in V, a
      \in K\).
    \item[(ST5)] \(\ip{x}{y+z} = \ip{x}{y} + \ip{x}{z}\) kaikille
      \(x,y,z \in V\).
  \end{itemize}
\end{definition}

\index{Hilbertin avaruus}
Täydellistä sisätuloavaruutta kutsutaan \defterm{Hilbertin
  avaruudeksi}.

\begin{examples}
  \begin{example}
    Reaalilukujen joukko \(\realnumbers\) varustettuna sisätulolla
    \(\ip{x}{y} := xy\), \(x, y \in \realnumbers\), on Hilbertin avaruus.
  \end{example}
  \begin{example}
    Kompleksilukujen joukko \(\complexnumbers\) varustettuna sisätulolla
    \(\ip{x}{y} := x^*y\), \(x, y \in \complexnumbers\), on Hilbertin avaruus.
  \end{example}
  \begin{example}
    \(\realnumbers^n\) varustettuna sisätulolla
    \begin{displaymath}
      \szip{\columnvector{x_1}{x_n}}{\columnvector{y_1}{y_n}}
      :=
      \sum_{k=1}^n x_k y_k , \spaceafter
      \columnvector{x_1}{x_n}, \columnvector{y_1}{y_n} \in
          \realnumbers^n
    \end{displaymath}
    on Hilbertin avaruus.
  \end{example}
  \begin{example}
    \(\complexnumbers^n\) varustettuna sisätulolla
    \begin{displaymath}
      \szip{\columnvector{x_1}{x_n}}{\columnvector{y_1}{y_n}}
      :=
      \sum_{k=1}^n x_k^* y_k , \spaceafter
      \columnvector{x_1}{x_n}, \columnvector{y_1}{y_n} \in
          \complexnumbers^n
    \end{displaymath}
    on Hilbertin avaruus.
  \end{example}
  \begin{example}
    \(\Lp{p}{\realnumbers^n}\) koostuu Borel-mitallisten funktioiden
    \(f: \realnumbers^n \to \complexnumbers\) määräämistä
    ekvivalenssiluokista:
    funktiot, jotka eroavat vain 0-mitallisessa joukossa,
    samaistetaan.
    Ks. luvut \ref{sec:Lp} ja \ref{ch:mittateoriaa}.
    Avaruudessa \(\Lp{2}{\realnumbers^n}\) sisätulo määritellään
    \begin{displaymath}
      \ip{f}{g} := \int_{\rnvec{x} \in \realnumbers^n} f(\rnvec{x})^*
      g(\rnvec{x}) d\mu
    \end{displaymath}
    missä \(f,g \in \Lp{2}{\realnumbers^n}\).
    Avaruus \(\Lp{2}{\realnumbers^n}\) on Hilbertin avaruus.
  \end{example}
\end{examples}

\index{sisätuloavaruuden normi}
Kun \(S\) on sisätuloavaruus, määritellään
\begin{equation}
  \label{eq:st-normi}
  \norm{x} := \sqrt{\ip{x}{x}}, \spaceafter x \in S .
\end{equation}
Myöhemmin osoitetaan, että tämä on normi avaruudessa \(S\).

\begin{theorem}
  \label{th:schwarzs-inequality}
  \textbf{Schwartzin epäyhtälö.}
  \cite[lause 1.5.5]{hpc2005}

  Olkoon \(V\) sisätuloavaruus ja \(f, g \in V\). Tällöin
  \(\abs{\ip{f}{g}} \leq \norm{f}\norm{g}\).
\end{theorem}

\begin{proof}
  \cite[lause 1.5.5]{hpc2005}
  Jos \(g = 0\), niin väite on tosi. Oletetaan jatkossa, että \(g
  \not= 0\). Riittää osoittaa, että
  \begin{displaymath}
    \abs{\ip{f}{\frac{g}{\norm{g}}}} \leq \norm{f}
  \end{displaymath}
  kaikille \(g \not= 0\), mikä on yhtäpitävää sen kanssa, että
  \(\abs{\ip{f}{h}} \leq \norm{f}\) kaikille yksikkövektoreille
  \(h\). Nyt
  \begin{align*}
    0 & \leq \norm{f - \ip{f}{h} h}^2 = \ip{f - \ip{f}{h} h}{f - \ip{f}{h} h}
    \\
    & = \ip{f}{f} - \ip{f}{h} \ip{f}{h} - \ip{f}{h} \ip{h}{f} +
    \ip{f}{h}^2 \ip{h}{h} \\
    & = \ip{f}{f} - \ip{f}{h} \ip{h}{f} \\
    & = \norm{f}^2 - \abs{\ip{f}{h}}^2 .
  \end{align*}
\end{proof}

\begin{theorem}
  Kaava \eqref{eq:st-normi} määrittelee normin avaruudessa
  \(V\).
\end{theorem}

\hyphenation{ilmi-selviä}

\begin{proof}
  \cite[lause 1.5.6]{hpc2005}
  Kolmioepäyhtälöä lukuunottamatta normin aksioomat ovat
  ilmiselviä. Kolmioepäyhtälö todistetaan seuraavasti:
  \begin{align*}
    \norm{f + g}^2 & = \ip{f+g}{f+g} \\
    & = \norm{f}^2 + \ip{f}{g} + \ip{g}{f} + \norm{g}^2 \\
    & \leq \norm{f}^2 + 2 \norm{f}\norm{g} + \norm{g}^2 \\
    & = (\norm{f} + \norm{g})^2 .
  \end{align*}
\end{proof}

On olemassa myös normiavaruuksia, jotka eivät ole sisätuloavaruuksia.

\begin{definition}
  \cite[määritelmä 1.5.7]{hpc2005}
  Olkoon \(V\) sisätuloavaruus ja \(A \subset V\).
  Määritellään
  \begin{displaymath}
    x \bot y \iff \ip{x}{y} = 0, \spaceafter x, y \in V .
  \end{displaymath}
  \index{vektorien kohtisuoruus}
  \index{kohtisuoruus}
  Jos \(x \bot y\), niin sanomme, että \(x\) ja \(y\) ovat
  kohtisuorassa toisiaan vastaan.  Sanomme, että \(A\) on
  \index{ortogonaalinen vektorijoukko}
  \defterm{ortogonaalinen}, jos \(x \bot y\) kaikille \(x, y \in A\),
  \(x \not= y\).  Jos lisäksi \(\ip{x}{x} = 1\) kaikille \(x \in A\),
  \index{ortonormaali vektorijoukko}
  niin sanomme, että \(A\) on \defterm{ortonormaali}.
  \index{ortogonaalinen komplementti}
  Määritellään joukon \(A\) \defterm{ortogonaalinen komplementti}
  asettamalla
  \begin{displaymath}
    A^\bot := \{ x \in V \setsep \forall y \in A : x \bot y \} .
  \end{displaymath}
  Siis joukon ortogonaalinen komplementti koostuu niistä vektoreista, jotka
  ovat kohtisuorassa kaikkia joukon vektoreita vastaan.
\end{definition}

\begin{theorem}
  \label{th:lahin}
  \cite[lause 1.5.11]{hpc2005}
  Olkoon \(V\) Hilbertin avaruus, \(M\) sen suljettu
  aliavaruus ja \(f \in V\). Tällöin on olemassa yksikäsitteinen
  avaruuden \(M\) piste lähimpänä pistettä \(f\).
\end{theorem}

\begin{proof}
  \cite[lause 1.5.11]{hpc2005}
  Olkoon \(d := \dist(f, M)\) ja muodostetaan jono \((g_n) \subset M\)
  siten, että \(\lim \norm{f-g_n} = d\). Suunnikassäännöstä
  \begin{displaymath}
    2\norm{a}^2 + 2\norm{b}^2 = \norm{a-b}^2 + \norm{a+b}^2
  \end{displaymath}
  seuraa, että
  \begin{align*}
    \norm{g_n - g_m}^2 & = \norm{(g_n - f) - (g_m - f)}^2 \\
    & = 2 \norm{g_n-f}^2 + 2 \norm{g_m-f}^2 - 4
    \sznorm{\frac{1}{2}(g_n+g_m)-f}^2 .
  \end{align*}
  Koska \((g_n+g_m)/2 \in M\), niin viimeisen termin itseisarvo ei ole
  vähemmän kuin \(4d^2\). Siten
  \begin{displaymath}
    \lim_{m,n \to \infty} \norm{g_n-g_m} \leq 2d^2 + 2d^2 - 4d^2 = 0 ,
  \end{displaymath}
  joten \((g_n)\) on Cauchyn jono. Koska \(M\) on suljettu, niin \(g_n
  \to g\) jollekin \(g \in M\). Koska \(\lim \norm{f-g_n} = d\), niin
  \(\norm{f-g} = d\).

  Oletetaan, että olisi toinen alkio \(g' \in M\) siten, että
  \(\norm{f-g'} = d\). Suunnikassäännöstä seuraa, että
  \begin{displaymath}
    \sznorm{f - \frac{1}{2}(g+g')}^2 = d^2 - \norm{g-g'}^2 .
  \end{displaymath}
  Koska \((g+g')/2 \in M\), niin tämän yhtälön vasen puoli ei ole
  vähemmän kuin \(d^2\). Siten \(g=g'\).
\end{proof}

\begin{theorem}
  \cite[lause 1.5.12]{hpc2005}
  Olkoon \(V\) Hilbertin avaruus ja \(M\) sen suljettu
  aliavaruus. Tällöin \(M^\bot\) on myös avaruuden \(H\) suljettu
  aliavaruus,
  \(M \intersection M^\bot = \zeroset\) ja
  \begin{displaymath}
    V = M \oplus M^\bot .
  \end{displaymath}
  Edelleen kehitelmässä \(f = g + h\), \(g \in M\), \(h \in M^\bot\),
  on \(g\) joukon \(M\) alkio, joka on lähimpänä alkiota \(f\).
\end{theorem}

\begin{proof}
  \cite[lause 1.5.12]{hpc2005}
  Olkoon \((x_k) \subset M^\bot\) ja \(x_k \to y \in H\). Olkoon \(m
  \in M\) mielivaltainen. Nyt \(\ip{m}{x_k} = 0\) kaikille \(k \in
  \naturalnumbers\). Toisaalta \(\ip{m}{x_k} \to \ip{m}{y}\), joten
  \(\ip{m}{y} = 0\). Koska \(m \in M\) oli mielivaltainen, niin \(y
  \in M^\bot\). Olkoon \(z \in M \intersection M^\bot\). Koska \(z \in
  M^\bot\), niin \(z\) on kohtisuorassa kaikkia avaruuden \(M\)
  vektoreita vastaan, erityisesti vektoria \(z \in M\). Siis
  \(\ip{z}{z} = 0\), joten \(z = 0\).
  Täten \(M^\bot\) on suljettu ja \(M \intersection M^\bot = \zeroset\).

  Riittää osoittaa, että \(V = M +
  M^\bot\). Jos \(f \in M\), niin tulos on ilmiselvä. Oletetaan, että
  \(f \not\in M\). Olkoon \(g \in M\) avaruuden \(M\) lähinnä alkiota
  \(f\) oleva piste. Tällainen piste on olemassa lauseen
  \ref{th:lahin} nojalla. Nyt on osoitettava, että \(f - g \in
  M^\bot\).  Mille tahansa \(h \in M\) ja \(\alpha > 0\), on \(g +
  \alpha h \in M\), ja siten
  \begin{displaymath}
    \norm{f-g}^2 \leq \norm{f-g-\alpha h}^2 = \norm{f-g}^2 - 2 \ReAlt
    \alpha \ip{h}{f-g} + \alpha^2 \norm{h}^2 .
  \end{displaymath}
  Tästä seuraa, että \(2 \ReAlt \alpha \ip{h}{f-g} \leq \alpha^2
  \norm{h}^2\) ja edelleen \(2 \ReAlt \ip{h}{f-g} \leq \alpha
  \norm{h}^2\). Koska \(\alpha > 0\) oli mielivaltainen, niin \(\ReAlt
  \ip{h}{f-g} \leq 0\). Toistamalla sama päättely kun \(\alpha < 0\)
  saadaan \(\ReAlt  \ip{h}{f-g} \geq 0\). Siten \(\ReAlt
  \ip{h}{f-g} = 0\) ja samoin \(\ImAlt \ip{h}{f-g} = 0\). Siten
  \(\ip{h}{f-g} = 0\) kaikille \(h \in M\), joten \(f-g \in M^\bot\).
\end{proof}

\begin{theorem}
  \cite[lemma 1.5.15]{hpc2005}
  Olkoon \(V\) Hilbertin avaruus.
  Olkoon \((\varphi_n) \subset H\) ortonormaali vektorijoukko ja \((\alpha_n)
  \subset \complexnumbers\). Tällöin \(\sum \alpha_n \varphi_n\)
  suppenee, jos ja vain jos \(\sum \abs{\alpha_n}^2\) suppenee. Jos
  sarja suppenee, niin sen summa on riippumaton termien järjestyksestä
  ja
  \begin{equation}
    \label{eq:orton-a}
    \sznorm{\sum_{n=0}^\infty \alpha_n \varphi_n} =
    \sqrt{\sum_{n=0}^\infty \abs{\alpha_n}^2} .
  \end{equation}
\end{theorem}

\begin{proof}
  Koska \((\varphi_n)\) on ortonormaali, niin
  \begin{equation}
    \label{eq:orton-b}
    \sznorm{\sum_{n=i}^j \alpha_n \varphi_n} =
    \sqrt{\sum_{n=i}^j \abs{\alpha_n}^2} .
  \end{equation}
  Siis jos \(\sum \abs{\alpha_n}^2\) suppenee, niin jono
  \((\sum_{n=0}^m \alpha_n \varphi_n)_{m=0}^\infty\) on Cauchyn jono,
  ja koska \(V\) on täydellinen, jono suppenee. Asettamalla \(i = 0\)
  ja antamalla \(j \to \infty\) kaavassa \eqref{eq:orton-b} saadaan
  kaava \eqref{eq:orton-a}.

  Olkoon \(f := \sum \alpha_{m_n} \varphi_{m_n}\) sarjan \(g := \sum
  \alpha_n \varphi_n\) uudelleenjärjestely. Nyt
  \begin{align}
    \label{eq:uj-a}
    \norm{f}^2 & = \norm{g}^2 = \sum \abs{\alpha_n}^2 \\
    \label{eq:uj-b}
    \norm{f-g}^2 & = \norm{f}^2 - 2 \ReAlt \ip{f}{g} + \norm{g}^2
  \end{align}
  Toisaalta
  \begin{displaymath}
    \ip{f}{g} = \lim_{j \to \infty} \szip{\sum \alpha_{m_n} \varphi_{m_n}}{\sum
      \alpha_n \varphi_n} = \sum \abs{\alpha_n}^2.
  \end{displaymath}
  Siten kaavojen \eqref{eq:uj-a} ja \eqref{eq:uj-b} nojalla
  \(\norm{f-g}=0\).
\end{proof}

\begin{definition}
  Olkoon \(V\) Hilbertin avaruus ja \(L = (\varphi_n)\) ortonormaali
  joukko. Sanomme, että \(L\) on Hilbertin avaruuden \(V\)
  \index{ortonormaali kanta}
  \defterm{ortonormaali kanta}, jos jokaiselle \(f \in V\) on
  \begin{equation}
    \label{eq:orton-c}
    f = \sum_{n=0}^\infty \ip{\varphi_n}{f} \varphi_n .
  \end{equation}
\end{definition}

\begin{remark}
  Ortonormaali kanta ei yleisesti ole Hamelin kanta.
\end{remark}

\begin{theorem}
  \cite[lause 1.5.18]{hpc2005}
  Olkoon \(V\) separoituva Hilbertin avaruus ja olkoon \(L =
  (\varphi_n)\) ortonormaali joukko avaruudessa \(V\). Tällöin
  seuraavat ehdot ovat yhtäpitäviä:
  \begin{itemize}
    \item[(i)] \(L^\bot = \zeroset\) ,
    \item[(ii)] \(\clos \spanop L = V\) ,
    \item[(iii)] \(V\) on ortonormaali kanta ja
    \item[(iv)] (Parsevalin kaava) Jokaiselle \(f \in V\) on
      \begin{displaymath}
        \norm{f}^2 = \sum \abs{\ip{\varphi_n}{f}}^2 .
      \end{displaymath}
  \end{itemize}
\end{theorem}

\begin{proof}
  (i) \(\implies\) (ii). Jos \(\clos \spanop L \not= V\), niin on
  olemassa piste \(f \in V\), \(f \not\in \clos \spanop L\). Siten \(f
  = g + h\) joillekin \(g \in \clos \spanop L\) ja \(h \in
  (\clos \spanop L)^\bot\), \(h \not=0\). Tästä seuraa, että \(L^\bot
  = (\clos \spanop L)^\bot \not= \zeroset\), mikä on ristiriidassa
  ehdon (i) kanssa.

  (ii) \(\implies\) (iii). Sarja \(\sum \ip{\varphi_n}{f}\) suppenee
  kohti funktion \(f\) ortogonaalista projektiota avaruudelle \(\clos
  \spanop L = V\).

  (iii) \(\implies\) (iv). Tämä seuraa yhtälöstä
  \begin{displaymath}
    \sznorm{f - \sum_{n=0}^\infty \ip{\varphi_n}{f} \varphi_n}^2
    =
    \norm{f}^2 - \sum_{n=0}^\infty \abs{\ip{\varphi_n}{f}}^2 .
  \end{displaymath}

  (iv) \(\implies\) (i). Jos \(f \in L^\bot\), niin jokainen termi
  summassa \eqref{eq:orton-c} on 0, mistä seuraa, että \(f = 0\). 
\end{proof}

\hyphenation{orto-normaali}

\begin{remark}
  Jokaisella separoituvalla Hilbertin avaruudella on ortonormaali
  kanta.
%% TOD.?
\end{remark}

\begin{definition}
  \cite{wprieszsequence}
  Olkoon \(V\) Hilbertin avaruus ja \((x_k)_{k=0}^\infty \subset
  \index{Rieszin jono}
  V\). Sanomme, että \((x_k)_{k=0}^\infty\) on \defterm{Rieszin jono},
  jos on olemassa vakiot \(c, d \in \realnumbers_+\) siten, että
  \begin{displaymath}
    c \sum_{k=0}^\infty \abs{a_k}^2
    \leq
    \sznorm{\sum_{k=0}^\infty a_k x_k}^2
    \leq
    d \sum_{k=0}^\infty \abs{a_k}^2
  \end{displaymath}
  kaikille jonoille \((a_k)_{k=0}^\infty \in l^2\).
  Jos \(\clos \spanop \{x_k\} = V\), niin sanomme, että
  \index{Rieszin kanta}
  \((x_k)_{k=0}^\infty\) on avaruuden \(V\) \defterm{Rieszin kanta}.
\end{definition}

%% Kummin päin konjugoinnit?
\begin{definition}
  \label{def:transpoosi-ja-h-konjugaatti}
  Olkoon \(V\) sisätuloavaruus ja \(A : V \to V\) lineaarinen
  \index{transpoosi}
  funktio. Määritellään funktion \(A\) \defterm{transpoosi}
  yhtälöllä
  \begin{displaymath}
    \ip{x}{Ay} = \ip{A^Tx}{y}^*
  \end{displaymath}
  kaikille \(x, y \in V\)
  \index{Hermiten konjugaatti}
  ja funktion \(A\) \defterm{Hermiten konjugaatti}
  yhtälöllä
  \begin{displaymath}
    \ip{x}{Ay} = {\ip{A^\herm x}{y}}
  \end{displaymath}
  kaikille \(x, y \in V\).
\end{definition}

Vertaa vastaavaan määritelmään matriiseille, määritelmä
\ref{def:matr-tr-hk}. Kuin matriisit tulkitaan lineaarisiksi
kuvauksiksi ja pystyvektorien välinen sisätulo määritellään
\begin{displaymath}
  \ip{x}{y} := x^\herm y , \;\;\;\; x, y \in \complexnumbers^{Z(n) \times
    Z(1)} ,
\end{displaymath}
niin molemmat määritelmät ovat ekvivalentteja (harjoitustehtävä).

\begin{definition}
  \label{def:lin-kuv-lajeja}.
  Olkoon \(V\) vektoriavaruus ja \(f : V \to V\)
  lineaarinen funktio. Sanomme, että
  \begin{itemize}
    \index{symmetrinen lineaarifunktio}
  \item \(f\) on \defterm{symmetrinen}, jos \(f = f^T\).
    \index{hermiittinen lineaarifunktio}
  \item \(f\) on \defterm{hermiittinen}, jos \(f = f^\herm\).
    \index{ortogonaalinen lineaarifunktio}
  \item \(f\) on \defterm{ortogonaalinen}, jos \(f^{-1} = f^T\).
    \index{unitaarinen lineaarifunktio}
  \item \(f\) on \defterm{unitaarinen}, jos \(f^{-1} = f^\herm\).
  \end{itemize}
\end{definition}

Vertaa vastaavaan määritelmään matriiseille, määritelmä
\ref{def:matriisilajeja}.

\section{Topologisiin vektoriavaruuksiin liittyviä lauseita}

\begin{definition}
  Olkoon \(V\) vektoriavaruus kerroinkunnalla \(\complexnumbers\).
  Olkoon \(f\) funktio avaruudelta \(V\) joukkoon \(\complexnumbers\)
  siten, että \(f(z+w) = f(z) + f(w)\) kaikille \(z, w \in V\).
  \index{reaalilineaarinen}
  Sanomme, että \(f\) on \defterm{reaalilineaarinen}, jos \(f(\alpha
    z) = \alpha f(z)\) kaikille \(\alpha \in \realnumbers\) ja \(z \in V\).
  \index{kompleksilineaarinen}
  Sanomme, että \(f\) on \defterm{kompleksilineaarinen}, jos \(f(\alpha
    z) = \alpha f(z)\) kaikille \(\alpha \in \complexnumbers\) ja \(z \in
    V\).
\end{definition}

\hyphenation{mukaan-lukien}

\index{Hahn-Banachin lause}

Nimitystä ``Hahn-Banachin lause'' käytetään useammalle eri lauseelle,
mukaanlukien tässä esitetyt seuraavat kolme lausetta.

\begin{theorem}
  \label{th:hb-a}
  Tämä lause ja todistus pohjautuvat lauseeseen 3.2 viitteessä \cite{rudin1991}.
  Oletetaan, että
  \begin{itemize}
    \item[(a)]
      \(V\) on vektoriavaruus kerroinkunnalla \(\realnumbers\) ja
      \(M\) on avaruuden \(V\) aliavaruus.
    \item[(b)]
      Funktio \(p : V \to \realnumbers\) täyttää ehdot
      \begin{enumerate}
        \item \(p(x + y) \leq p(x) + p(y)\) ja
        \item \(p(tx) = tp(x)\)
      \end{enumerate}
      kaikille \(x, y \in V\) ja \(t \in \realnumbers_0\).
    \item[(c)]
      Funktio \(f : M \to \realnumbers\) on lineaarinen ja \(f(x) \leq
      p(x)\) kaikille \(x \in M\).
  \end{itemize}
  Tällöin on olemassa lineaarinen funktio \(g : V \to \realnumbers\)
  siten, että \(g(x) = f(x)\) kaikille \(x \in M\) ja \(-p(-x) \leq
  g(x) \leq p(x)\) kaikille \(x \in V\).
\end{theorem}

\begin{proof}
  Jos \(M \not= V\) valitaan \(x_1 \in V\), \(x_1 \not\in V\), ja
  määritellään
  \begin{displaymath}
    M_1 := \{ x + tx_1 \setsep x \in M, t \in \realnumbers \} .
  \end{displaymath}
  Nyt \(M_1\) on vektoriavaruus. Koska
  \begin{displaymath}
    f(x) + f(y) = f(x+y) \leq p(x+y) \leq p(x-x_1) + p(x_1 + y) ,
  \end{displaymath}
  niin
  \begin{displaymath}
    f(x) - p(x-x_1) \leq p(y+x_1) - f(y)
  \end{displaymath}
  kaikille \(x, y \in M\). Olkoon
  \begin{displaymath}
    \alpha := \sup_{x \in M} (f(x) - p(x-x_1)) .
  \end{displaymath}
  Nyt
  \begin{equation}
    \label{eq:hb-a}
    f(x) - \alpha \leq p(x-x_1)
  \end{equation}
  kaikille \(x \in M\) ja
  \begin{equation}
    \label{eq:hb-b}
    f(y) + \alpha \leq p(y+x_1)
  \end{equation}
  kaikille \(x \in M\).
  Olkoon \(t \in \realnumbers_0\).
  Määritellään \(f_1 : M_1 \to \realnumbers\) yhtälöllä
  \begin{displaymath}
    f_1(x+tx_1) = f(x) + t \alpha ,
  \end{displaymath}
  missä \(x \in M\).
  Epäyhtälöstä \eqref{eq:hb-a} seuraa
  \begin{displaymath}
    f(t^{-1} x) - \alpha \leq p(t^{-1}x-x_1) ,
  \end{displaymath}
  joten
  \begin{displaymath}
    f(x) - \alpha t \leq p(x-tx_1) .
  \end{displaymath}
  Vastaavasti epäyhtälöstä \eqref{eq:hb-b} seuraa
  \begin{displaymath}
    f(y) + \alpha t \leq p(y+tx_1) 
  \end{displaymath}
  kaikille \(y \in M\).
  Edelleen
  \begin{displaymath}
    f_1(x-tx_1) = f(x) - t \alpha \leq p(x-tx_1)
  \end{displaymath}
  ja
  \begin{displaymath}
    f_1(y+tx_1) = f(y) + t \alpha \leq p(y+tx_1) .
  \end{displaymath}
  Jos \(t = 0\), niin meillä on \(f(x) \leq p(x)\) kaikille \(x \in M\).
  Täten
  \begin{displaymath}
    f_1(z+ux_1) \leq p(z+ux_1)
  \end{displaymath}
  kaikille \(u \in \realnumbers\) ja \(z \in M\).
  Siis \(f_1 \leq p\) avaruudessa \(M_1\).

  \hyphenation{maksimaalisuus-teoreeman}

  Olkoon \(P\) joukko kaikista järjestetyistä pareista \((M', f')\),
  missä \(M'\) on avaruuden \(V\) aliavaruus, joka sisältää avaruuden
  \(M\), ja \(f' : M \to \realnumbers\) lineaarinen funktio, joka
  on funktion \(f\) laajennus ja joka toteuttaa epäyhtälön \(f'(x)
  \leq p(x)\) kaikille \(x \in M'\). Järjestetään \(P\) osittain
  asettamalla \((M',f') \leq (M'', f'')\), jos ja vain jos \(M'
  \subset M''\) ja \(f''(x) = f'(x)\) kaikilla \(x \in
  M'\). Hausdorffin maksimaalisuusteoreeman nojalla on olemassa
  joukon \(P\) maksimaalinen totaalisesti järjestetty osajoukko
  \(\Omega\).

  Olkoon \(\Phi\) joukko kaikista \(M'\) siten, että \((M',f') \in
  \Omega\). Nyt \(\Phi\) on totaalisesti järjestetty
  osajoukkorelaatiolla, ja siten kaikkien joukon \(\Phi\) jäsenten
  unioni \(\tilde{M}\) on avaruuden \(V\) aliavaruus.
  Jos \(x \in \tilde{M}\), niin \(x \in M'\) jollekin \(M' \in \Phi\)
  ja määritellään \(g(x) := f'(x)\), missä \(f'\) on parin \((M',f')
  \in \Omega\) jälkimmäinen jäsen.

  Nyt \(g\) on hyvin määritelty joukossa \(\tilde{M}\), \(g\) on
  lineaarinen ja \(g(x) \leq f'(x)\) kaikille \(x \in \tilde{M}\). Jos
  \(\tilde{M}\) olisi avaruuden \(V\) aito aliavaruus, niin tämän
  todistuksen ensimmäinen osa antaisi uuden laajennuksen funktiolle
  \(g\), mikä olisi ristiriidassa joukon \(\Omega\)
  maksimaalisuuden kanssa. Siten \(\tilde{M} = V\).

  Edelleen epäyhtälöstä \(g \leq p\) seuraa, että
  \begin{displaymath}
    -p(-x) \leq -g(-x) = g(x)
  \end{displaymath}
  kaikille \(x \in V\).
\end{proof}

\begin{theorem}
  Tämä lause ja todistus on otettu viitteestä \cite[lause
    3.3]{rudin1991}.  Olkoon \(V\) on vektoriavaruus kerroinkunnalla
  \(K = \realnumbers\) tai \(K = \complexnumbers\), \(M\) avaruuden
  \(V\) aliavaruus, \(p\) seminormi avaruudessa \(V\), \(f : M \to K\)
  lineaarinen funktio siten, että \(\abs{f(x)} \leq p(x)\) kaikille
  \(x \in M\).
  Tällöin \(f\) voidaan laajentaa lineaariseksi funktionaaliksi \(g :
  V \to K\), jolle \(\abs{g(x)} \leq p(x)\) kaikille \(x \in V\).
\end{theorem}

\begin{proof}
  Jos \(K = \realnumbers\), niin tämä lause seuraa lauseesta
  \ref{th:hb-a}, koska nyt \(p(-x) = p(x)\). Oletetaan, että \(K =
  \complexnumbers\). Määritellään \(u = \ReAlt f\). Lauseen
  \ref{th:hb-a} nojalla on olemassa reaalilineaarinen \(U : V \to
  \realnumbers\) siten, että \(U(x) = u(x)\) kaikilla \(x \in M\) ja
  \(U(y) \leq p(y)\) kaikilla \(y \in V\).
  Olkoon \(g\) kompleksilineaarinen funktio avaruudessa \(V\) siten,
  että \(\ReAlt g = U\). Nyt \(g(x) = f(x)\) kaikilla \(x \in M\).
  Jokaiselle \(x \in V\) on olemassa \(\alpha \in \complexnumbers\),
  \(\abs{\alpha} = 1\) siten, että \(\alpha g(x) = \abs{g(x)}\). Täten
  \begin{displaymath}
    \abs{g(x)} = g(\alpha x) = U(\alpha x) \leq p(\alpha x) = p(x).
  \end{displaymath}
\end{proof}

\begin{corollary}
  Ks. \cite[luku II.3.2]{schaefer1971}.
  Olkoon \(V\) normiavaruus, \(M\) avaruuden \(V\) normialiavaruus ja
  \(f \in M^*\), jolle \(\abs{f(x)} \leq \norm{x}\) kaikille \(x
  \in M\).
  Tällöin funktionaalilla \(f\) on lineaarinen laajennus \(f_1\)
  avaruudelle \(V\) siten, että \(\abs{f_1(x)} \leq \norm{x}\)
  kaikille \(x \in V\).
\end{corollary}

\begin{theorem}
  \textbf{Avoimen kuvauksen lause.}
  Tämä lause ja todistus pohjautuvat lauseeseen I.2.11 viitteessä
  \cite{rudin1991}.

  Oletetaan, että
  \begin{itemize}
    \item[(a)] \(V\) on \textit{F}-avaruus,
    \item[(b)] \(W\) on topologinen vektoriavaruus,
    \item[(c)] \(f : V \to W\) on jatkuva ja lineaarinen ja
    \item[(d)] \(\setimage{f}{V}\) on toista kategoriaa avaruudessa
      \(W\).
  \end{itemize}
  Tällöin
  \begin{itemize}
    \item[(i)] \(\setimage{f}{V} = W\),
    \item[(ii)] \(f\) on avoin ja
    \item[(iii)] \(W\) on \textit{F}-avaruus.
  \end{itemize}
\end{theorem}

\begin{proof}
  Kohdasta (ii) seuraa (i), koska \(W\) on itsensä ainoa epätyhjä
  avoin aliavaruus.
  Kohdan (ii) todistamiseksi olkoon \(V\) pisteen \(0\) ympäristö
  avaruudessa \(V\). On osoitettava, että \(\setimage{f}{V}\) sisältää
  pisteen \(0\) ympäristön avaruudessa \(W\).
  Olkoon \(d\) translaatioinvariantti metriikka avaruudessa \(V\)
  siten, että \(d\) on yhteensopiva avaruuden \(V\) topologian kanssa.
  Olkoon
  \begin{displaymath}
    V_n := \{ x \in V \setsep d(x, 0) < 2^{-n}r \} , \spaceafter n \in
    \naturalnumbers ,
  \end{displaymath}
  missä \( r > 0 \) on niin pieni, että \(V_0 \subset V\).
  Osoitamme, että pisteen \(0\) eräs ympäristö \(W\) avaruudessa \(W\)
  toteuttaa ehdon
  \begin{equation}
    \label{eq:av-k-a}
    \setimage{f}{V} \subset \clos \setimage{f}{V_1} \subset W .
  \end{equation}
  Koska \(V_2 - V_2 \subset V_1\), niin
  \begin{displaymath}
    (\clos \setimage{f}{V_2})
    -
    (\clos \setimage{f}{V_2})
    \subset
    \clos (\setimage{f}{V_2} - \setimage{f}{V_2})
    \subset
    \clos \setimage{f}{V_1} .
  \end{displaymath}
  Väitteen \eqref{eq:av-k-a} ensimmäisen osan todistamiseksi riittää,
  että joukolla \(\clos \setimage{f}{V_2}\) on epätyhjä
  sisäpuoli. Mutta
  \begin{displaymath}
    \setimage{f}{V} = \bigcup_{k=1} k \setimage{f}{V_2}
  \end{displaymath}
  koska \(V_2\) on pisteen \(0\) ympäristö.
  Ainakin yksi \(k \setimage{f}{V_2}\) on siten toista kategoriaa
  avaruudessa \(W\). Koska \(y \mapsto ky\) on homeomorfismi
  avaruudelta \(W\) avaruudelle \(W\), niin \(\setimage{f}{V_2}\) on
  toista kategoriaa avaruudessa \(W\), joten joukolla
  \(\setimage{f}{V_2}\) on epätyhjä sisäpuoli.

  Todistaaksemme toisen osan väitteestä \eqref{eq:av-k-a} kiinnitetään
  \(y_1 \in \clos \setimage{f}{V_1}\). Oletetaan, että \(n \geq 1\) ja
  \(y_n \in \clos \setimage{f}{V_n}\) on valittu. Edellä oleva
  todistus joukolle \(V_1\) pätee yhtä hyvin joukolle
  \(V_{n+1}\), joten \(\clos \setimage{f}{V_{n+1}}\) sisältää pisteen
  \(0\) ympäristön. Täten
  \begin{displaymath}
    \left(
    y_n - \clos \setimage{f}{V_{n+1}}
    \right)
    \intersection \setimage{f}{V_n}
    \not= \emptyset .
  \end{displaymath}
  Siten olemassa \(x_n \in V_n\) siten, että
  \begin{displaymath}
    f(x_n) \in y_n - \clos \setimage{f}{V_{n+1}} .
  \end{displaymath}
  Asetetaan \(y_{n+1} := y_n - f(x_n)\). Nyt
  \begin{math}
    y_{n+1} \in \clos \setimage{f}{V_{n+1}}
  \end{math}.
  Koska \(d(x_n,0) < 2^{-n}r\) kaikille \(n \in \positiveintegers\),
  summat \(x_1 + \ldots + x_n\) muodostavat Cauchyn jonon joka
  avaruuden \(V\) täydellisyyden nojalla suppenee kohti jotakin \(x
  \in V\), missä \(d(x,0) < r\). Täten \(x \in V\). Koska
  \begin{displaymath}
    \sum_{n=1}^m f(x_n) = \sum_{n=1}^m (y_n - y_{n+1}) = y_1 - y_{m+1} ,
  \end{displaymath}
  ja koska funktion \(f\) jatkuvuuden nojalla \(y_{m+1} \to 0\) kun
  \(m \to \infty\), niin \(y_1 = f(x) \in \setimage{f}{V}\). Tämä
  osoittaa toisen osan väitteestä \eqref{eq:av-k-a}, joten (ii) on
  todistettu.

  Olkoon \(N := \ker f\).
  Voidaan osoittaa, että \(V / N\) on \textit{F}-avaruus. Siten
  kohdan (iii) todistamiseksi riittää, jos löydämme isomorfismin \(g :
  V / N \onto W\), joka on myös homeomorfismi. Määritellään
  \begin{displaymath}
    g(x + N) := f(x), \spaceafter x \in V .
  \end{displaymath}
  Nyt \(g\) on isomorfismi ja \(f(x)=g(\pi(x))\), missä \(\pi : V \to
  V / N\) on jäännösluokkafunktio \(\pi(x) = x + N\), \(x \in V\).

  Jos \(V\) on avoin avaruudessa \(W\), niin
  \begin{displaymath}
    \setimage{g^{-1}}{V} = \setimage{\pi}{\setimage{f^{-1}}{V}}
  \end{displaymath}
  on avoin, koska \(f\) on jatkuva ja \(\pi\) on avoin. Täten \(g\) on
  jatkuva. Jos \(E\) on avoin avaruudessa \(V / N\), niin
  \begin{displaymath}
    \setimage{g}{E} = \setimage{f}{\setimage{\pi^{-1}}{E}}
  \end{displaymath}
  on avoin, koska \(\pi\) on jatkuva ja \(f\) on avoin. Täten \(g\) on
  homeomorfismi.
\end{proof}

\begin{definition}
  Jos \(V\) ja \(W\) ovat epätyhjiä joukkoja ja \(f\) on funktio
  joukolta \(V\) joukkoon \(W\), niin määrittelemme funktion \(f\)
  \index{kuvaaja}
  \index{funktion kuvaaja}
  \defterm{kuvaajan} joukoksi
  \begin{displaymath}
    \{ (x, f(x)) \setsep x \in V\ \} \subset V \times W .
  \end{displaymath}
\end{definition}

\begin{theorem}
  Tämä lause ja todistus on otettu lauseesta I.2.14 viitteessä
  \cite{rudin1991}.

  Jos \(V\) on topologinen avaruus, \(W\) on Hausdorffin avaruus ja
  \(f : V \to W\) on jatkuva, niin funktion \(f\) kuvaaja on suljettu.
\end{theorem}

\begin{proof}
  Olkoon \(G\) funktion \(f\) kuvaaja ja \(\Omega :=  \cp{V}{W}
  \setminus G\). Olkoon \((x_0, y_0) \in \Omega\). Nyt \(y_0 \not=
  f(x_0)\).  Täten pisteillä \(y_0\) ja \(f(x_0)\) on erilliset
  ympäristöt \(V\) ja \(W\) avaruudessa \(W\). Koska \(f\) on jatkuva,
  pisteellä \(x_0)\) on ympäristö \(U\) siten, että \(\setimage{f}{U}
  \subset W\).  Pisteen \((x_0, y_0)\) ympäristö \(U \times V\)
  sisältyy siten joukkoon \(\Omega\). Tämä osoittaa, että \(\Omega\)
  on avoin.
\end{proof}

\begin{theorem}
  \textbf{Suljetun kuvaajan lause.}
  Tämä lause ja todistus on otettu lauseesta I.2.15 viitteessä
  \cite{rudin1991}.

  Oletetaan, että
  \begin{itemize}
    \item[(a)] \(V\) ja \(W\) ovat \textit{F}-avaruuksia,
    \item[(b)] \(f : V \to W\) on lineaarinen ja
    \item[(c)] \(G := \{ (x, f(x)) \setsep x \in V \}\) on suljettu
      avaruudessa \(V \times W\).
  \end{itemize}
  Tällöin \(f\) on jatkuva.
\end{theorem}

\begin{proof}
  Joukko \(V \times W\) on vektoriavaruus, kun yhteenlasku ja
  skalaarilla kertominen määritellään alkioittain:
  \begin{displaymath}
    \alpha (x_1, y_1) + \beta (x_2, y_2) = (\alpha x_1 + \beta x_2,
    \alpha y_1 + \beta y_2 ) .
  \end{displaymath}
  On olemassa täydelliset translaatioinvariantit metriikat \(d_V\) ja
  \(d_W\) avaruuksissa \(V\) ja \(W\), siten että nämä metriikat
  generoivat joukkojen \(V\) ja \(W\) topologiat.
  Asetetaan
  \begin{displaymath}
    d((x_1, y_1), (x_2, y_2)) := d_V(x_1, x_2) + d_W(y_1, y_2) .
  \end{displaymath}
  Nyt \(d\) on translaatioinvariantti metriikka avaruudessa \(V \times
  W\), \(d\) ja avaruuden \(V \times W\) tulotopologia ovat
  yhteensopivia, ja \(d\) tekee avaruudesta \(V \times W\)
  \textit{F}-avaruuden.

  Koska \(f\) on lineaarinen, niin \(G\) on avaruuden \(V \times W\)
  aliavaruus. Täydellisen metrisen avaruuden suljetut osajoukot ovat
  täydellisiä. Siten \(G\) on \textit{F}-avaruus.
  Määritellään \(\pi_1 : G \to V\) ja \(\pi_2 : V \times W \to W\)
  asettamalla \(\pi_1(x, f(x)) := x\) ja \(\pi_2(x, y) := y\).
  Nyt \(\pi_1\) on jatkuva lineaarinen bijektio \textit{F}-avaruudelta
  \(G\) \textit{F}-avaruuteen \(V\). Avoimen kuvauksen lauseen
  perusteella \(\pi_1^{-1} : V \to G\) on jatkuva. Nyt \(f = \pi_2
  \circ \pi_1^{-1}\) ja \(\pi_2\) on jatkuva, joten \(f\) on jatkuva.
\end{proof}

\section{Differentiaalilaskentaa}

\begin{definition}
  Olkoon \(A \subset \realnumbers\), \(A \not= \emptyset\) ja \(f: A
  \index{derivaatta}
  \to \realnumbers\) funktio. Funktion \(f\) \defterm{derivaatta}
  pisteessä \(x \in A\) määritellään
  \begin{displaymath}
    f'(x) := \lim_{h \to 0} \frac{f(x+h)-f(x)}{h} .
  \end{displaymath}
  Jos yllä oleva raja-arvo on olemassa, niin sanomme, että \(f\) on
  \defterm{derivoituva} pisteessä \(x\). Jos \(f\) on derivoituva
  kaikilla \(x \in A\), niin sanomme, että \(f\) on
  \index{derivoituvuus}
  \index{derivaattafunktio}
  \defterm{derivoituva}. Tällöin \defterm{derivaattafunktio} \(f'\) on
  funktio joukolta \(A\) joukkoon \(\realnumbers\).
\end{definition}

Avaruuden \(\naturalnumbers^n\), \(n \in \positiveintegers\), alkioita
\index{multi-indeksi}
kutsutaan \defterm{multi-indekseiksi}. Kun \(\alpha \in
\naturalnumbers^n\), niin määritellään
\begin{displaymath}
  \abs{\alpha} = \abs{(\alpha_1,\ldots,\alpha_n)} := \sum_{k=1}^n
  \alpha_k .
\end{displaymath}

\begin{definition}
  Olkoon \(n \in \positiveintegers\), \(A \subset \realnumbers^n\),
  \(A \not= \emptyset\) ja \(f : A \to \complexnumbers\)
  funktio. Kun \(k \in Z(n)\), niin määritellään funktion \(f\) 
  \index{osittaisderivaatta}
  \defterm{osittaisderivaatta} pisteessä \(\rnvec{x} \in A\) asettamalla
  \begin{displaymath}
    (\partial^k f)(\rnx) = \frac{\partial f}{\partial x_k}\left(\rnx\right)
    := \lim_{h \to 0} \frac{f(\rnx+h\unitvec{n}{k}) - f(\rnx)}{h} .
  \end{displaymath}
  Määritellään lisäksi
  \index{$\mideriv{\alpha}$}
  \begin{displaymath}
    \mideriv{\alpha} f := \frac{\partial^{\abs{\alpha}}
      f}{\partial_{x_1}^{\alpha_1} \cdots \partial_{x_n}^{\alpha_n}}
    = \frac{\partial^{\alpha_1}}{\partial x_1^{\alpha_1}} \cdots
    \frac{\partial^{\alpha_n}}{\partial x_n^{\alpha_n}} f
  \end{displaymath}
  kaikille \(\alpha \in \naturalnumbers^n\) ja \(\rnx \in A\).
  Tässä asetetaan
  \begin{displaymath}
    \frac{\partial^0}{\partial x_k^0} f = f .
  \end{displaymath}
  Määritellään lisäksi funktion \(f\) \defterm{gradientti} pisteessä
  \(\rnx \in A\) kaavalla
  \begin{displaymath}
    (\nabla f)(\rnx) :=
    \begin{pmatrix}
      \left(\frac{\partial f}{\partial x_1}\right) \left( \rnx \right)
      \\
      \vdots \\
      \left(\frac{\partial f}{\partial x_n}\right) \left( \rnx \right)
    \end{pmatrix} .
  \end{displaymath}
\end{definition}

\begin{definition}
  Olkoon \(n \in \positiveintegers\) ja \(U \subset \rn\), \(U \not=
  \emptyset\). Sanomme, että funktion \(f : U \to \complexnumbers\) on
  \index{sileä funktio}
  \defterm{sileä}, jos \(f\) on äärettömän monta kertaa derivoituva.
\end{definition}

Oletetaan, että funktio \(f : \intervalcc{x_0}{x} \to \realnumbers\),
missä \(x_0, x \in \realnumbers\),
ja sen \(m+1\) ensimmäistä derivaattaa ovat jatkuvia välillä
\(\intervalcc{x_0}{x}\). Tällöin on
\index{Taylorin kehitelmä}
voimassa \defterm{Taylorin kehitelmä} \cite{grossman}
\begin{equation}
  \label{eq:taylor-one-dim}
  f(x) = f(x_0) + \sum_{k=1}^m \frac{f^{(k)}(x_0)}{k!}(x-x_0)^k + R_m(x)
\end{equation}
missä jäännöstermi on
\begin{displaymath}
  R_m(x) = \frac{f^{(m+1)}(c)}{(m+1)!}(x-x_0)^{m+1}
\end{displaymath}
jollekin \(c \in \intervalcc{x_0}{x}\).

\begin{definition}
  Kun \(\rnx, \rny \in \rn\) niin määritellään pisteitä \(\rnx\) ja
  \(\rny\)
  \index{viivasegmentti}
  yhdistävä \defterm{viivasegmentti} kaavalla
  \begin{displaymath}
    V(\rnx, \rny) := \{ \rnz \in \rn \setsep \rnz = \rnx + t(\rny -
    \rnx), t \in \intervalcc{0}{1} \} .
  \end{displaymath}
\end{definition}

Olkoon \(\Omega \subset \rn\) on epätyhjä ja avoin, \(m \in
\positiveintegers\), \(\rnx_0, \rnx \in \rn\) ja \(V(\rnx_0, \rnx)
\subset \Omega\).  Olkoon \(f : \rn \to \realnumbers\) jatkuva
funktio, jonka kaikki korkeintaan \(m+1\) kertaluvun
osittaisderivaatat ovat jatkuvia.  \index{Taylorin kehitelmä} Tällöin
on voimassa \defterm{Taylorin kehitelmä} \cite{grossman}
\begin{equation}
  \label{eq:taylor-mdim}
  f(\rnx) = f(\rnx_0) + \sum_{k=1}^m \sum_{i_1,\ldots,i_k = 1}^n
  \frac{\partial}{\partial x_{i_1}} \cdots \frac{\partial}{\partial
    x_{i_k}}
  f(\rnx_0) \cdot \prod_{\alpha = 1}^k
  \left(x_{i_\alpha}-x^{(0)}_{i_\alpha}\right)
  + R_m(\rnx)
\end{equation}
missä jäännöstermi on
\begin{align}
  \label{eq:taylor-residual-mdim}
  R_m(\rnx) = \frac{1}{(m+1)!} \sum_{i_1,\ldots,i_{m+1}=1}^n &
  \frac{\partial}{\partial x_{i_1}} \cdots \frac{\partial}{\partial
    x_{i_{m+1}}} f \left( \rnx_0 + c(\rnx - \rnx_0) \right) \\
  & \cdot \prod_{\alpha = 1}^{m+1}
  \left(x_{i_\alpha}-x^{(0)}_{i_\alpha}\right)
\end{align}
jollekin \(c \in \intervaloo{0}{1}\).

Taylorin kehitelmien voimassaolo on todistettu kirjassa
\cite{grossman}.

\hyphenation{tensori-tulot}

\section{Hilbertin ja Banachin avaruuksien tensoritulot}

\hyphenation{kerroin-kunnalla}

Tässä luvussa on käytetty kirjaa \cite{ryan2002} ja viitettä
\cite{wptoptp}.
\index{tensoritulo}
Kun \(X\) ja \(Y\) ovat vektoriavaruuksia
kerroinkunnalla \(K\), niin määritellään
\begin{displaymath}
  \bil{X}{Y} := \{ \text{kaikki bilineaariset funktiot}\;X \times Y
  \to K \}
\end{displaymath}
ja alkioiden välinen tensoritulo
\(x \otimes y \in \bil{X}{Y}^\sharp\) asettamalla
\begin{displaymath}
  (x \otimes y)(A) := \langle A , x \otimes y \rangle := A(x, y)
\end{displaymath}
kaikille \(A \in \bil{X}{Y}\).
Määritellään vektoriavaruuksien \(X\) ja \(Y\) välinen
\index{algebrallinen tensoritulo}
\defterm{algebrallinen tensoritulo} kaavalla
\begin{displaymath}
  X \otimes Y := \spanop \{ x \otimes y \setsep x \in X, y \in Y \}
\end{displaymath}
Kun \(u \in X \otimes Y\), \(u = \sum_{k=1}^n a_k x_k \otimes
y_k\), \(a_k \in K\),
niin \(u = 0\), jos ja vain jos
\begin{displaymath}
  \sum_{k=1}^n a_k \varphi(x_k) \psi(y_k) = 0
\end{displaymath}
kaikille \(\varphi \in X^\sharp\) ja \(\psi \in Y^\sharp\).

Olkoot \(H_1\) ja \(H_2\) Hilbertin avaruuksia sisätuloilla
\(\ip{\cdot}{\cdot}_1\) ja \(\ip{\cdot}{\cdot}_2\)
\cite{wptphilbert}. Voimme saada edellä määritellystä algebrallisesta
tensoritulosta sisätuloavaruuden määrittelemällä
\begin{displaymath}
  \ip{\varphi_1 \otimes \varphi_2}{\psi_1 \otimes \psi_2}
  :=
  \ip{\varphi_1}{\psi_1}_1 \ip{\varphi_2}{\psi_2}_2
\end{displaymath}
kaikille \(\varphi_1, \psi_1 \in H_1\) ja \(\varphi_2, \psi_2 \in
H_2\) ja laajentamalla määrittely lineaarisuuden nojalla koko
avaruuteen \(H_1 \otimes H_2\). Muodostamalla täydellistymä näin
saadusta sisätuloavaruudesta saadaan Hilbertin avaruuksien \(H_1\) ja
\(H_2\) välinen tensoritulo, jota merkitään myös \(H_1 \otimes H_2\).

\hyphenation{tensori-tuloksi}
\hyphenation{normi-avaruudessa}

\index{järkevä ristinormi}
Olkoot \(A\) ja \(B\) Banachin avaruuksia
\cite{wptoptp}. \defterm{Järkevä ristinormi} \(\alpha\) on normi
algebrallisessa tensoritulossa \(A \otimes B\) siten, että
\begin{itemize}
  \item \(\alpha(a \otimes b) = \norm{a}\norm{b}\).
  \item Kaikille \(\varphi \in X^*\) ja \(\psi \in Y^*\) lineaarinen
    funktionaali \(\varphi \otimes \psi\) normiavaruudessa \((A
    \otimes B, \alpha)\) on rajoitettu, ja \(\norm{\varphi \otimes
      \psi} = \norm{\varphi} \norm{\psi}\).
\end{itemize}
Normiavaruutta, joka saadaan varustamalla vektoriavaruus \(A \otimes
B\) normilla \(\alpha\) merkitään \(A \otimes_\alpha B\).
On olemassa suurin järkevä ristinormi \(\pi\), jota kutsutaan
\index{projektiivinen ristinormi}
\defterm{projektiiviseksi ristinormiksi} ja joka määritellään
\begin{displaymath}
  \pi(x) := \inf \left\{ \sum_{i=1}^n \norm{a_i} \norm{b_i} \bigsetsep x
  = \sum_{i=1}^n a_i \otimes b_i \right\}
\end{displaymath}
missä \(x \in A \otimes B\).
On olemassa pienin järkevä ristinormi \(\varepsilon\), jota kutsutaan
\index{injektiivinen ristinormi}
\defterm{injektiiviseksi ristinormiksi}, ja joka määritellään
\begin{displaymath}
  \varepsilon(x) := \sup \{ (\varphi \otimes \psi)(x) \setsep \varphi
  \in A^*, \psi \in B^*, \norm{\varphi} = \norm{\psi} = 1 \}
\end{displaymath}
missä \(x \in A \otimes B\).
Albegrallisten tensoritulojen täydellistymiä näiden kahden normin
\index{projektiivinen tensoritulo}
\index{injektiivinen tensoritulo}
suhteen kutsutaan \defterm{projektiiviseksi} ja \defterm{injektiiviseksi
  tensorituloksi} ja merkitään \(A \hat{\otimes}_\pi B\) ja \(A
\hat{\otimes}_\varepsilon B\).
Jos \(A\) ja \(B\) ovat Hilbertin avaruuksia, niin Hilbertin
avaruuksien tensoritulon normi ei yleisesti ole kumpikaan normeista
\(\pi\) tai \(\varepsilon\). Hilbertin avaruuksien tensoritulon normia
merkitään joskus \(\sigma\).

\hyphenation{risti-normiksi}
\hyphenation{jäännös-luokka-operaattoreita}

\index{tasainen ristinormi}
\defterm{Tasainen ristinormi} \(\alpha\) liittää jokaiseen pariin
\((X, Y)\) Banachin avaruuksia järkevän ristinormin \(\alpha_{X,Y}\)
avaruudessa \(X \otimes Y\) siten, että seuraava ehto täyttyy:
Olkoot \(X\), \(W\), \(Y\) ja \(Z\) Banachin avaruuksia. Tällöin
kaikille jatkuville lineaarisille funktioille \(S : X \to W\) ja \(T :
Y \to Z\) funktio \(S \otimes T : X \otimes_\alpha Y \to W
\otimes_\alpha Z\) on jatkuva ja \(\norm{S \otimes T} \leq \norm{S}
\norm{T}\). Vektoriavaruutta \(X \otimes Y\) varustettuna normilla
\(\alpha_{X,Y}\) merkitään \(X \otimes_\alpha Y\) ja sen
täydellistymää \(X \hat{\otimes}_\alpha Y\).
\index{äärellisesti generoitu}
Sanomme, että tasainen ristinormi \(\alpha\) on \defterm{äärellisesti
  generoitu}, jos jokaiselle parille \((X, Y)\) Banachin avaruuksia ja
jokaiselle \(u \in X \otimes Y\) on
\begin{displaymath}
  \alpha_{X,Y}(u) = \inf \{ \alpha_{M,N}(u) \setsep u \in  M \otimes
  N, \dim M < \infty, \dim N < \infty \} .
\end{displaymath}
\index{kofiniittisesti generoitu}
Sanomme, että tasainen ristinormi \(\alpha\) on \defterm{kofiniittisesti
  generoitu}, jos jokaiselle parille \((X, Y)\) Banachin avaruuksia ja
jokaiselle \(u \in X \otimes Y\) on
\begin{displaymath}
  \alpha_{X,Y}(u) = \sup \{ \alpha_{(X/E) \otimes
    (Y/F)}\left(\left(Q_E \otimes Q_F\right)u\right)
  \setsep
  \dim X / E < \infty, \dim Y / F < \infty \} ,
\end{displaymath}
missä \(Q_E : X \to X / E\) ja \(Q_F : Y \to Y / F\) ovat
jäännösluokkaoperaattoreita.
\index{tensorinormi}
\defterm{Tensorinormi} määritellään äärellisesti generoiduksi
tasaiseksi ristinormiksi. Projektiivinen ja injektiivinen
ristinormi ovat molemmat tensorinormeja, ja niitä kutsutaan
\index{projektiivinen tensorinormi}
\index{injektiivinen tensorinormi}
projektiiviseksi ja injektiiviseksi tensorinormiksi.

Jos \(E\) ja \(F\) ovat äärellisulotteisia Banachin avaruuksia ja
\(\alpha\) tensorinormi, niin
\begin{displaymath}
  E \otimes_{\alpha'} F = ( E^* \otimes_\alpha F^*)^* ,
\end{displaymath}
missä
\begin{displaymath}
  \alpha'(u) = \sup \{ \abs{\langle u, v \rangle} \setsep v \in  E^*
  \otimes F^*, \alpha(v) \leq 1 \}
\end{displaymath}
Ääretönulotteisessa tapauksessa avaruus \(X \otimes Y\) ei ole enää
sama kuin \((X^* \otimes_\alpha Y^*)^*\), mutta on olemassa kanoninen
algebrallinen upotus siihen avaruudelta \(X \otimes Y\).
Siten avaruuden \((X^* \otimes_\alpha Y^*)^*\) duaalinormi määrittelee
normin avaruudessa \(X \otimes Y\), jota kutsumme \defterm{Schattenin
  duaalinormiksi} ja merkitsemma \(\alpha^s\). Siis meillä on upotus
\begin{displaymath}
  X \otimes_{\alpha^s} Y \isomemb (X^* \otimes_\alpha Y^*)^*
\end{displaymath}
ja
\begin{displaymath}
  \alpha^s(u) = \sup \{ \abs{\langle u, v \rangle} \setsep v \in  X^*
  \otimes Y^*, \alpha(v) \leq 1 \} .
\end{displaymath}
Meillä on \(\pi^s = \varepsilon\), mutta \(\varepsilon^s <
\pi\). Ongelma Schattenin duaalissa on, että se ei ole äärellisesti
generoitu. Tämän ratkaisemiseksi määritellään \defterm{duaalinormi}
\begin{displaymath}
  \alpha'(u) := \inf \{ \alpha'_{E,F}(u) \setsep u \in E \otimes F,
  \dim E < \infty, \dim F < \infty \}
\end{displaymath}
missä \(E\) ja \(F\) käyvät läpi kaikki avaruuksien \(X\) ja \(Y\)
äärellisulotteiset aliavaruudet joille \(u \in E \otimes F\).
Nyt meillä on \(\pi' = \varepsilon\), \(\varepsilon' = \pi\), ja
yleisesti \((\alpha')' = \alpha\) jokaiselle tensorinormille
\(\alpha\).

\index{neliöjärjestys}
Neliöjärjestys \(\sqordfunction\) määritellään kuten kirjoissa
\cite{ryan2002} ja \cite{singer1970} seuraavassa:

\begin{definition}
  Määritellään funktio \(\sqordfunction : \naturalnumbers \to
  \naturalnumbers^2\)
  asettamalla
  \begin{displaymath}
    \sqord{0} := (0, 0)
  \end{displaymath}
  ja
  \begin{displaymath}
    \sqord{k} :=
    \left\{
    \begin{array}{ll}
      (i, n) ; \;\; & k = n^2 + 1 \;\land\; i \in \{0, \ldots, n\}
      \;\land\; n \in \naturalnumbers \\
      (n, n - i) ; \;\; & k = n^2 + n + i \;\land\; i \in \{1,\ldots,n\}
      \;\land\; n \in \naturalnumbers .
    \end{array}
    \right.
  \end{displaymath}
\end{definition}

Funktio \(\sqordfunction\) on bijektio joukolta \(\naturalnumbers\)
joukolle \(\naturalnumbers^2\). Tensoritulokanta
määritellään kuten kirjassa \cite{ryan2002}:

\begin{definition}
  Olkoot \(E\) ja \(F\) Banachin avaruuksia, joilla on Schauderin
  kannat \((a_k)_{k=0}^\infty\) ja \((b_k)_{k=0}^\infty\),
  vastaavasti. Jonoa
  \begin{math}
    (a_{\sqordfirst{k}} \otimes b_{\sqordsecond{k}})_{k=0}^\infty
    \subset E \otimes F
  \end{math}
  sanotaan kantojen \((a_k)_{k=0}^\infty\) ja \((b_k)_{k=0}^\infty\)
  \index{tensoritulokanta}
  generoimaksi \defterm{tensoritulokannaksi}. 
\end{definition}

Tensoritulokanta on Schauderin kanta sekä avaruudelle \(E
\otimes_\varepsilon F\) että \(E \otimes_\pi F\).

\begin{exercises}
  \begin{exercise}{6.1}
    Osoita, että kun matriisit tulkitaan lineaarisiksi funktioksi,
    Hermiten konjugaatin määritelmät
    \ref{def:transpoosi-ja-h-konjugaatti} ja \ref{def:matr-tr-hk} ovat
    ekvivalentteja.
  \end{exercise}
  \begin{exercise}{6.2}
    Määritellään funktio \(f : l^1 \to l^\infty\), \(f(\seq{x}) :=
    \seq{x}\) kaikille \(\seq{x} \in l^1\). Osoita, että \(f\) on
    jatkuva lineaarinen injektio, mutta käänteisfunktio \(f^{-1} :
    \setimage{f}{l^1} \to l^1\) ei ole jatkuva.
  \end{exercise}
  \begin{exercise}{6.3}
    Kirjoita auki yhden muuttujan Taylorin kehitelmä
    \eqref{eq:taylor-one-dim}, kun \(m = 4\).
  \end{exercise}
  \begin{exercise}{6.4}
    Kirjoita auki monen muuttujan Taylorin kehitelmä
    \eqref{eq:taylor-mdim}, kun muuttujien lukumäärä \(n = 3\) ja
    \(m = 4\).
  \end{exercise}
\end{exercises}

\chapter{Distribuutiot}
\label{ch:distribuutiot}

\section{Fourier-muunnos ja Fourier'n sarjat}
\label{sec:fourier-analyysi}

Olkoot \(a_k, b_k \in \realnumbers\), \(k \in \positiveintegers\)
ja \(a_0 \in \realnumbers\).
Sarjaa
\begin{displaymath}
  S(x) := a_0 + \sum_{k=1} \left(a_k \cos(x) + b_k \sin(x)\right)
\end{displaymath}
\index{Fourier'n sarja}
kutsutaan \defterm{Fourier'n sarjaksi}.

\index{skalaaritulo}
Vektorien \(\rnx, \rny \in \rn\) \defterm{skalaaritulo} määritellään
kaavalla
\begin{displaymath}
  \rnx \cdot \rny := \sum_{k=1}^n x_k y_k .
\end{displaymath}

%% Määrittelyalue?
Olkoon \(n \in \positiveintegers\) ja \(f \in
\index{Fourier-muunnos}
\Lp{1}{\rn}\). Määritellään funktion \(f\) \defterm{Fourier-muunnos}
asettamalla \cite[luku 1.2.1]{triebel1983}
\index{$\ft$}
\begin{displaymath}
  (\ft f)(\rny) := \left(2\pi\right)^{-\frac{n}{2}} \int_{\rnx \in
    \rn} e^{-\imagunit \rny \cdot \rnx}
  f(\rnx) d\mu
\end{displaymath}
\index{Fourier-käänteismuunnos}
ja \defterm{Fourier-käänteismuunnos}
\index{$\ft^{-1}$}
\begin{displaymath}
  (\ft^{-1} f)(\rnx) := \left(2\pi\right)^{-\frac{n}{2}} \int_{\rny
    \in \rn} e^{\imagunit \rnx \cdot \rny}
  f(\rny) d\mu .
\end{displaymath}

\section{Tavalliset distribuutiot}

Tämä luku pohjautuu viitteeseen \cite{wpdistr}.
Tässä luvussa oletetaan, että \(n \in \positiveintegers\) ja \(U
\subset \rn\), \(U \not= \emptyset\), \(U\) avoin.

\begin{definition}
  Sanomme, että \(\varphi : U \to \realnumbers\) on
  \index{testifunktio}
  \defterm{testifunktio}, jos se on äärettömän monta kertaa derivoituva
  ja kompaktikantajainen.
  \index{testifunktioavaruus}
  Määritellään \defterm{testifunktioavaruus}
  \begin{displaymath}
    \testfunctionspace{U} := \{ f : U \to \realnumbers \setsep
    f \; \text{on testifunktio}\} .
  \end{displaymath}
  Testifunktioavaruuden kerroinkunta on \(\realnumbers\).
\end{definition}

\begin{definition}
  Olkoon \(X\) epätyhjä joukko.
  Kun \(f \in \realnumbers^X\), määritellään
  \begin{displaymath}
    \norminfty{f} := \sup_{x \in X} \abs{f(x)} .
  \end{displaymath}
  Kun \((f_k)_{k=0}^\infty \subset \realnumbers^X\), sanomme,
  \index{tasainen suppeneminen}
  että \(f_k \to g \in \realnumbers^X\) tasaisesti, jos
  \begin{displaymath}
    \lim_{k \to \infty} \norminfty{g - f_k} = 0 .
  \end{displaymath}
\end{definition}

Määritellään topologia avaruudessa \(\testfunctionspace{U}\)
seuraavasti:
Kun \((\varphi_k)_{k=0}^\infty \subset \testfunctionspace{U}\) niin
\(\varphi_k \to \varphi \in \testfunctionspace{U}\), jos ja vain jos
\begin{enumerate}
  \item
    On olemassa \(K \subset U\), \(K\) kompakti siten, että
    \begin{displaymath}
      \bigcup_{k \in \naturalnumbers} \supp \varphi_k \subset K .
    \end{displaymath}
  \item
    Jokaiselle \(\alpha \in \naturalnumbers^n\)
    \begin{displaymath}
      \mideriv{\alpha} \varphi_k \to \mideriv{\alpha} \varphi \;\;
      \text{tasaisesti} .
    \end{displaymath}
\end{enumerate}

\(\testfunctionspace{U}\) on täydellinen lokaalikonveksi avaruus.

Yhtäpitävästi avaruuden \(\testfunctionspace{U}\) topologia voidaan
määritellä seuraavasti:
Olkoon
\begin{displaymath}
  U = \bigcup_{i \in \naturalnumbers} U_i ,
\end{displaymath}
missä \(U_i \in \rn\) avoin, \(K_i := \clos U_i\) kompakti, \(i \in
\naturalnumbers\).
Nyt meillä on
\begin{displaymath}
  \testfunctionspace{U} = \bigcup_{i \in \naturalnumbers} D_{K_i} ,
\end{displaymath}
missä
\begin{displaymath}
  D_{k_i} := \{ f : U \to \realnumbers \setsep f \;\text{on sileä
    ja}\; \supp f \subset K_i \}.
\end{displaymath}
Jokaisessa \(D_{K_i}\) määritellään topologia seminormeilla
\begin{displaymath}
  \norm{\varphi}_\alpha := \max_{x \in K_i} \abs{\mideriv{\alpha}
    \phi}, \spaceafter \alpha \in \naturalnumbers^n ,
\end{displaymath}
ts. mielivaltaisen asteisten derivaattojen tasaisen suppenemisen
topologia.
Jokainen \(D_{K_i}\) on Fr{\'e}chet'n avaruus.
Olkoon \(\iota_i : D_{K_i} \to \testfunctionspace{U}\) funktio \(x \in
D_{K_i} \mapsto x \in \testfunctionspace{U}\).
\index{induktiivinen rajatopologia}
Induktiivinen rajatopologia \(\tau\) avaruudessa
\(\testfunctionspace{U}\) on hienoin lokaalikonveksin vektoriavaruuden
topologia, jossa kaikki funktiot \(\iota_i\), \(i \in
\naturalnumbers\), ovat jatkuvia.

\index{distribuutioavaruus}
Määritellään \defterm{distribuutioavaruus} \(\distrspace{U} :=
(\testfunctionspace{U})^*\). Avaruuden \(\distrspace{U}\) alkioita
\index{distribuutio}
kutsutaan \defterm{distribuutioiksi}.
Siis kun \(T \in \distrspace{U}\), niin \(T\) on jatkuva, jos ja vain jos
\begin{displaymath}
  \lim_{k \to \infty} T(\varphi_k) = T\left( \lim_{k \to \infty}
  \varphi_k \right)
\end{displaymath}
Riippumatta distribuutioavaruuteen valitusta topologiasta jono
distribuutioita suppenee, jos ja vain jos se suppenee pisteittäin.
Tämän takia distribuutioavaruuden topologiaksi valitaan joskus
\index{heikko-*-topologia}
heikko-*-topologia, mutta usein käytetään myös rajoitetun suppenemisen
topologiaa, joka tässä tapauksessa on sama kuin tasaisen suppenemisen
topologia kompakteissa joukoissa.
Merkitään
\begin{displaymath}
  \distrappl{T}{\varphi} := T(\varphi)
\end{displaymath}
kaikille \(T \in \distrspace{U}\) ja \(\varphi
\in \testfunctionspace{U}\).

\begin{definition}
  Olkoon \(f : U \to \realnumbers\) funktio. Sanomme, että \(f\) on
  \index{lokaalisti integroituva funktio}
  \defterm{lokaalisti integroituva}, jos se on Lebesgue-integroituva
  kaikilla \(K \subset U\), \(K\) kompakti.
\end{definition}

Kaikki jatkuvat funktiot, kaikki testifunktiot ja kaikki
\(L^p\)-funktiot (ks. luku \ref{sec:Lp}) ovat lokaalisti
integroituvia.  Lokaalisti integroituva funktio \(f\) määrittelee
distribuution \(T_f\), jolle
\begin{displaymath}
  \distrappl{T_f}{\varphi} := \int_{\rnx \in U} f(\rnx) \varphi(\rnx)
  d\mu
\end{displaymath}
kaikille \(\varphi \in \testfunctionspace{U}\).
Usein merkitään
\begin{displaymath}
  \distrappl{f}{\varphi} := \distrappl{T_f}{\varphi} .
\end{displaymath}
Jos \(f\) ja \(g\) ovat lokaalisti integroituvia funktioita, niin
\(T_f = T_g\), jos ja vain jos \(f = g\) melkein kaikkialla.
Jos \(R \in \distrspace{U}\) ja \(R = T_h\) jollekin lokaalisti
integroituvalle funktiolle \(h\), niin sanomme, että distribuutio
\index{säännöllinen distribuutio}
\(T\) on \defterm{säännöllinen}.
Avaruus \(\testfunctionspace{U}\) on tiheä avaruudessa
\(\distrspace{U}\). Jokaiselle \(T \in \distrspace{U}\) on olemassa
\((\varphi_k) \subset \testfunctionspace{U}\) siten, että
\begin{displaymath}
  \distrappl{\varphi_k}{\psi} \to \distrappl{T}{\psi}
\end{displaymath}
kaikille \(\psi \in \testfunctionspace{U}\).
Tämä seuraa Hahn-Banachin lauseesta, koska avaruuden
\(\distrspace{U}\) duaali sen heikko-*-topologialla on avaruus
\(\testfunctionspace{U}\).

\begin{examples}
  \begin{example}
    Jos \(0 \in U\), niin
    \index{Diracin delta}
    \defterm{Diracin delta} määritellään
    \begin{displaymath}
      \delta(\varphi) := \varphi(0), \spaceafter \varphi
      \in \tfs{U} .
    \end{displaymath}
    Diracin delta on distribuutio.
  \end{example}
  \begin{example}
    \index{Cauchyn pääarvo}
    \defterm{Cauchyn pääarvo} määritellään kaavalla
    \begin{displaymath}
      \left( \text{p.v.} \frac{1}{x} \right)(\varphi)
      :=
      \lim_{\varepsilon \to 0+} \int_{\abs{x} \geq \varepsilon}
        \frac{\varphi(x)}{x} d\mu, \spaceafter \varphi
      \in \tfs{U} .
    \end{displaymath}
    Cauchyn pääarvo on distribuutio.
  \end{example}
\end{examples}

Jos \(A : \testfunctionspace{U} \to \testfunctionspace{U}\) on
lineaarinen funktio, joka on jatkuva heikko-*-topologiassa, niin \(A\)
on mahdollista laajentaa funktioksi \(A : \distrspace{U}
\to \distrspace{U}\).
Jos \(A : \testfunctionspace{U} \to \testfunctionspace{U}\) on
\index{transpoosi}
\index{distribuution transpoosi}
lineaarinen ja jatkuva, niin sen \defterm{transpoosi} määritellään funktioksi
\(A^t : \testfunctionspace{U} \to \testfunctionspace{U}\), jolle
\begin{displaymath}
  \int_{\rnx \in U} A\varphi(\rnx) \cdot \psi(\rnx) d\mu = \int_{\rnx \in U}
  \varphi(\rnx) \cdot A^t \psi(\rnx) d\mu
\end{displaymath}
kaikille \(\varphi, \psi \in \testfunctionspace{U}\).
Jos \(A^t\) on olemassa ja jatkuva avaruudessa \(\tfs{U}\), niin \(A\)
voidaan laajentaa avaruudelle \(\distrspace{U}\) määrittelemällä
\begin{displaymath}
  \distrappl{AT}{\varphi} := \distrappl{T}{A^t \varphi}
\end{displaymath}
kaikille \(\varphi \in \tfs{U}\).

Olkoon \(A : \testfunctionspace{U} \to \testfunctionspace{U}\),
\begin{displaymath}
  A \varphi := \frac{\partial \varphi}{\partial x_k} .
\end{displaymath}
Jos \(\varphi, \psi \in \tfs{U}\), niin osittaisintegrointi antaa
\begin{displaymath}
  \int_U \frac{\partial \varphi}{\partial x_k} \psi d\mu = - \int_U
  \varphi \frac{\partial \psi}{\partial x_k} d\mu ,
\end{displaymath}
joten \(A^t = -A\).
\index{osittaisderivaatta}
\index{distribuution osittaisderivaatta}
Siis jos \(T \in \distrspace{U}\), niin distribuution \(T\)
osittaisderivaatta koordinaatin \(x_k\) suhteen määritellään kaavalla
\begin{displaymath}
  \szdistrappl{\frac{\partial T}{\partial x_k}}{\varphi}
  :=
  - \szdistrappl{T}{\frac{\partial \varphi}{\partial x_k}}
\end{displaymath}
kaikille \(\varphi \in \tfs{U}\).
Yleisesti, jos \(\alpha \in \naturalnumbers^n\), niin distribuution
\(T \in \distrspace{U}\) osittaisderivaatta \(\mideriv{\alpha} T\)
määritellään
\begin{displaymath}
  \distrappl{\mideriv{\alpha} T}{\varphi} := (-1)^{\abs{\alpha}}
  \distrappl{T}{\mideriv{\alpha} \varphi}
\end{displaymath}
kaikille \(\varphi \in \tfs{U}\).
Jokainen distribuutio on äärettömän monta kertaa derivoituva, ja
osittaisderivointi \(\partial^\alpha\) on lineaarinen ja jatkuva
operaatio avaruudessa \(\distrspace{U}\). Useimmat muut derivaatan
määritelmät eivät ole jatkuvia.

Jos \(m : U \to \realnumbers\) on äärettömän monta kertaa derivoituva
funktio ja \(T\) distribuutio, niin määritellään
\begin{displaymath}
  \distrappl{mT}{\varphi} := \distrappl{T}{m\varphi}
\end{displaymath}
kaikille \(\varphi \in \tfs{U}\).
Sileillä funktioilla kertomisen suhteen \(\distrspace{U}\) on moduuli
yli renkaan \(C^\infty(U)\).

\begin{definition}
  Olkoot \(U, V \subset \rn\), \(U\) ja \(V\) epätyhjiä ja
  avoimia. Olkoon \(F : V \to U\) funktio. Sanomme, että \(F\) on
  \index{submersio}
  \defterm{submersio}, jos Jacobin derivaatta \(dF(x)\) on lineaarinen
  surjektio jokaiselle \(x \in V\).
\end{definition}

Olkoot \(T \in \distrspace{U}\), \(U, V \in \rn\), \(U\) ja \(V\)
epätyhjiä ja avoimia ja \(F : V \to U\) submersio.
Nyt voidaan määritellä \(T \circ F \in \distrspace{V}\). Joskus
merkitään
\begin{displaymath}
  F^\sharp : T \mapsto F^\sharp T = T \circ F .
\end{displaymath}

Olkoot \(U, V \subset \rn\), \(U\) ja \(V\) epätyhjiä ja avoimia ja \(V \subset U\).
Määritellään funktio \(E_{VU} : \tfs{V} \to \tfs{U}\) asettamalla
\begin{displaymath}
  \left( E_{VU}(f) \right) \left( \rnx \right)
  :=
  \left\{
  \begin{array}{ll}
    f(\rnx) ; & \rnx \in V \\
    0 ; & \rnx \in U \setminus V
  \end{array}
  \right.
\end{displaymath}
kaikille \(\rnx \in U\). Määritellään \(\rho_{VU} := (E_{VU})^t\).
Jokaiselle \(T \in \distrspace{U}\) rajoittuma \(\rho_{VU}(T)\) on
distribuutio avaruudessa \(\distrspace{V}\):
\begin{displaymath}
  \distrappl{\rho_{VU}(T)}{\varphi} = \distrappl{T}{E_{VU}(\varphi)}
\end{displaymath}
kaikille \(\varphi \in \tfs{V}\).

\begin{definition}
  Olkoon \(U \subset \rn\), \(U\) epätyhjä ja avoin. Olkoon \(T
  \index{häviävä distribuutio}
  \in \distrspace{U}\). Sanomme, että \(T\) \defterm{häviää} avoimessa
  joukossa \(V \subset U\), jos
  \begin{displaymath}
    T \in \ker \rho_{VU} .
  \end{displaymath}
  Toisin sanoen, \(T\) häviää joukossa \(V\), jos ja vain jos
  \begin{displaymath}
    \distrappl{T}{\varphi} = 0
  \end{displaymath}
  kaikilla \(\varphi \in \tfs{U}\) joille \(\supp \varphi \subset V\).
\end{definition}

\begin{definition}
  Olkoon \(U \subset \rn\), \(U\) epätyhjä ja avoin. Olkoon \(T
  \in \distrspace{U}\).
  \index{kantaja}
  \index{distribuution kantaja}
  Määritellään distribuution \(T\) \defterm{kantaja} kaavalla
  \begin{displaymath}
    \supp T := U \setminus \{ V \setsep \rho_{VU}(T) = 0 \}. 
  \end{displaymath}
  Jos \(\supp T\) on kompakti, niin sanomme, että \(T\) on
  \index{kompaktikantajainen distribuutio}
  \defterm{kompaktikantajainen}.
\end{definition}

\begin{definition}
  Olkoot \(a \in \realnumbers_+\), \(b \in \rn\)
  ja \(T \in \distrspace{\rn}\).
  Distribuution \(T\) \defterm{\(a\)-dilataatio ja
    \(b\)-translaatio}, jota merkitään \(T(a \cdot - b)\),
  määritellään \cite{cl1996}
  \begin{displaymath}
    \distrappl{T(a \cdot - b)}{f} := 
    \frac{1}{a^n}
    \szdistrappl{T}{f \left( \frac{\cdot + b}{a} 
    \right)}
  \end{displaymath}
  kaikille
  \begin{math}
    f \in \tfs{\rn}
  \end{math}.
\end{definition}

\section{Temperoidut distribuutiot}

Tämä jakso pohjautuu kirjaan \cite{triebel1983} ja viitteeseen
\cite{wpdistr}.
Olkoon \(n \in \positiveintegers\).
Olkoon
\begin{displaymath}
  p_{\alpha,\beta}(\varphi) = \sup_{\rnx \in \rn} \abs{x^\alpha
    \mideriv{\beta} \varphi(\rnx)} , \spaceafter \alpha, \beta \in
  \naturalnumbers^n ,
\end{displaymath}
missä \(\varphi : \rn \to \realnumbers\) on äärettömän monta kertaa
\index{Schwartzin avaruus}
\index{$S(\rn)$}
derivoituva. Määritellään \defterm{Schwartzin avaruus}
\begin{align*}
  S(\rn) = \{ \varphi : \rn \to \realnumbers \setsep & \varphi \;\text{on äärettömän
    monta kertaa derivoituva ja} \\
  & p_{\alpha,\beta}(\varphi) < \infty
  \;\text{kaikille}\; \alpha, \beta \in \naturalnumbers^n \} .
\end{align*}
\index{Schwartzin funktio}
Schwartzin avaruuden alkioita sanotaan \defterm{Schwartzin funktioiksi}.
Seminormiperhe \(p_{\alpha,\beta}\) määrittelee lokaalikonveksin
topologian avaruudessa \(S(\rn)\). Schwartzin avaruus on metrisoituva
ja täydellinen. Fourier-muunnos on bijektio avaruudelta \(S(\rn)\)
avaruudelle \(S(\rn)\).

\index{$\tempdistrspace{\rn}$}
Olkoon \(\tempdistrspace{\rn}\) avaruuden \(S(\rn)\) topologinen duaali
varustettuna vahvalla topologialla, ks. määritelmä
\ref{def:vahva-topologia}. Avaruuden \(\tempdistrspace{\rn}\) alkioita sanotaan
\index{temperoitu distribuutio}
\defterm{temperoiduiksi distribuutioiksi}.

\begin{remark}
  Avaruuden \(\tempdistrspace{\rn}\) alkiot eivät välttämättä ole funktioita.
\end{remark}

Nyt \(\tempdistrspace{\rn} \subset \distrspace{\rn}\). Kaikilla
temperoiduilla distribuutioilla on Fourier-muunnos, mutta kaikilla
avaruuden \(\distrspace{\rn}\) ei ole.

Distribuutio \(T\) on temperoitu distribuutio, jos ja vain jos
\begin{displaymath}
  \lim_{k \to \infty} T(\varphi_k) = 0
\end{displaymath}
on tosi aina, kun
\begin{displaymath}
  \lim_{k \to \infty} p_{\alpha,\beta}(\varphi_k) = 0
\end{displaymath}
kaikille
\((\varphi_k) \subset S(\rn)\) ja
\(\alpha, \beta \in \naturalnumbers^n\).

Temperoidun distribuution derivaatta on temperoitu
distribuutio. Kaikki kompaktikantajaiset distribuutiot ja kaikki
neliöllisesti integroituvat funktiot ovat temperoituja
distribuutioita. Kaikki funktiot \(P(\rnx)f(\rnx)\), missä \(P(\rnx)\)
on polynomi ja \(f \in \Lp{p}{\rn}\), \(p \geq 1\), ovat temperoituja
distribuutioita.

%% Temperoituvat distribuutiot ovat \defterm{hitaasti kasvavia}, mikä on
%% duaali sille, että Schwartzin funktiot ovat \defterm{nopeasti
%%   väheneviä}, eli
%% \begin{displaymath}
%%   \varphi \sim \norm{x}^n \exp \left(-x^2\right)
%% \end{displaymath}
%% kaikille \(\varphi \in S(\rn)\).

Laskettaessa Fourier-muunnoksia on parasta tarkastella
kompleksiarvoisia testifunktioita ja kompleksisia lineaarisia
distribuutioita.
\index{Fourier-muunnos}
\index{temperoidun distribuution Fourier-muunnos}
Kun \(T \in \tempdistrspace{\rn}\), niin määritellään temperoidun
distribuution \(T\) \defterm{Fourier}-muunnos asettamalla
\begin{displaymath}
  (\ft T)(\psi) := T(\ft \psi)
\end{displaymath}
kaikille \(\psi \in S(\rn)\).
Meillä on
%% useammat ulottuvuudet?
\begin{displaymath}
  \ft \frac{dT}{dx} = \imagunit x \ft T .
\end{displaymath}
Fourier-muunnos on bijektio avaruudelta \(\tempdistrspace{\rn}\)
avaruudelle \(\tempdistrspace{\rn}\).

\section{Konvoluutio}

\begin{definition}
  Olkoon \(\psi \in \S(\rn)\). Sanomme, että \(\psi\) on
  \index{hitaasti kasvava funktio}
  \defterm{hitaasti kasvava}, jos kaikki funktion \(\psi\) derivaatat
  kasvavat korkeintaan yhtä nopeasti kuin polynomit.
\end{definition}

Olkoon \(T \in \tempdistrspace{\rn}\) ja olkoon \(\psi : \rn \to
\complexnumbers\) hitaasti kasvava äärettömän monta kertaa derivoituva
funktio joukossa \(\rn\). Nyt \(\psi T \in \tempdistrspace{\rn}\) ja
\begin{displaymath}
  \ft(\psi T) = (\ft \psi)*(\ft T) ,
\end{displaymath}
mikä on distribuutioiden \(\ft T\) ja \(\ft \psi\) (mahdollisesti
\index{konvoluutio}
funktio) \defterm{konvoluutio}.

Meillä on \(\ft 1 = \delta\), missä \(1\) on vakiofunktio \(f(\rnx)=1\)
ja \(\delta\) on Diracin \(\delta\)-distribuutio.
Kun \(f \in \complexnumbers^{\rn}\), niin määritellään
\begin{displaymath}
  \check{f}(\rnx) := f(-\rnx), \spaceafter \rnx \in \rn ,
\end{displaymath}
ja konvoluutio
\begin{displaymath}
  \distrappl{f * T}{\varphi} := \distrappl{T}{\check{f} * \varphi},
  \spaceafter f, \varphi \in \tfs{\rn} .
\end{displaymath}

On myös mahdollista määritellä kahden distribuution (avaruudella
\(\rn\)) konvoluutio, jos toinen niistä on kompaktikantajainen.
Olkoot \(S, T \in \distrspace{\rn}\) ja olkoon \(T\)
kompaktikantajainen. Määritellää konvoluutio
\begin{displaymath}
  (S*T)*\varphi := S*(T*\varphi), \spaceafter \varphi \in \tfs{\rn} .
\end{displaymath}

Määritellään
\begin{displaymath}
  (\tau_{\rnx} \varphi)(\rny) := \varphi(\rny - \rnx), \spaceafter
  \rnx, \rny \in \rn
  \;\text{ja}\; \varphi \in \tfs{\rn}
\end{displaymath}
ja
\begin{displaymath}
  \psi(\rnx) := \distrappl{T}{\tau_{-\rnx}\varphi}, \spaceafter \varphi \in
  \tfs{\rn} .
\end{displaymath}
Vaihtoehtoinen määritelmä konvoluutiolle on
\begin{displaymath}
  \distrappl{S*T}{\varphi} := \distrappl{S}{\psi} .
\end{displaymath}

Olkoot \(f \in \distrspace{\rn}\) ja \(K \subset \rn\) kompakti. On
olemassa jatkuva ja kompaktikantajainen funktio \(F : \rn \to
\realnumbers\) ja multi-indeksi \(\alpha \in \naturalnumbers^n\)
%% avaruuden määritelmä?
siten, että \(f = \mideriv{\alpha} F\) on sileä ja kompaktikantajainen.
Jos \(\supp f\) sisältää vain yhden pisteen \(x\), niin
\begin{displaymath}
  f = \sum_{\abs{\alpha} \leq m} a_\alpha \mideriv{\alpha}(\tau_x
  \delta) .
\end{displaymath}

\begin{exercises}
  \begin{exercise}{7.1}
    %% L05062009-6
    Olkoot \(a, c \in \realnumbers_+\) ja \(\rnb, \rnd \in \rn\). Kun \(T
    \in \distrspace{\rn}\), määritellään \(\tau_{a,\rnb}T := T(a \cdot -
    \rnb)\). Osoita, että \(\tau_{a,\rnb}(\tau_{c,\rnd} T) = \tau_{ac,c
      \rnb + \rnd} T\).
  \end{exercise}
  \begin{exercise}{7.2}
    %% L18052015-5
    Olkoon \(n \in \positiveintegers\). Määritellään \(\delta :
    \tfs{\rn} \to \realnumbers\) (Diracin delta) asettamalla
    \(\delta(\varphi) := \varphi(0)\) kaikille \(\varphi \in
    \tfs{\rn}\). Osoita, että \(\delta \in \distrspace{\rn}\).
  \end{exercise}
\end{exercises}

\chapter{Mittateoriaa}
\label{ch:mittateoriaa}

Tässä luvussa on käytetty lähteitä \cite{knapp,wpborelset,wpmeasfunc}.

\section{$\sigma$-algebra ja $\sigma$-rengas}

\begin{definition}
  Olkoon \(X\) epätyhjä joukko. Sanomme, että joukko \(\mca
  \index{joukkoalgebra}
  \subset \powerset{X}\) on \defterm{joukkoalgebra}, jos seuraavat
  aksioomat ovat voimassa:
  \begin{itemize}
    \item[(1)] \(\emptyset \in \mca\) ja \(X \in \mca\).
    \item[(2)] Jos \(E, F \in \mca\), niin \(E \union F \in
      \mca\).
    \item[(3)] Jos \(E \in \mca\), niin \(X \setminus E \in
      \mca\).
  \end{itemize}
\end{definition}

\begin{definition}
  Olkoon \(X\) epätyhjä joukko ja \(\mca \subset \powerset{X}\)
  joukkoalgebra. Sanomme, että \(\mca\) on
  \index{$\sigma$-algebra}
  \defterm{\(\sigma\)-algebra}, jos
  \begin{displaymath}
    \bigcup_{k=0}^\infty J_k \in \mca
  \end{displaymath}
  kaikille \(\{ J_k \setsep k \in \naturalnumbers \} \subset \mca\).
\end{definition}

\begin{definition}
  Olkoon \(X\) epätyhjä joukko ja \(\mcr \subset
  \index{rengas}
  \powerset{X}\). Sanomme, että \(\mcr\) on \defterm{rengas}, jos
  \begin{itemize}
    \item[(1)] \(\union_{k=1}^n A_k \in \mcr\) kaikille
      \(A_1,\ldots,A_n \in \mcr\), \(n \in \naturalnumbers\).
    \item[(2)] Jos \(A, B \in \mcr\), niin \(A \setminus B \in \mcr\).
  \end{itemize}
  Jos lisäksi
  \begin{itemize}
    \item[(3)] \(\union_{k=0}^\infty A_k \in \mcr\) kaikille
      numeroituville kokoelmille \(A_k \in \mcr\), \(k \in
      \naturalnumbers\), niin sanomme, että \(\mcr\) on
      \index{$\sigma$-rengas}
      \defterm{\(\sigma\)-rengas}.
  \end{itemize}
\end{definition}

Seuraavissa esimerkeissä \(X\) oletetaan epätyhjäksi joukoksi.

\hyphenation{osa-joukot}

\begin{examples}
  \begin{example}
    \(\mca := \{ \emptyset, X \}\) on \(\sigma\)-algebra.
  \end{example}
  \begin{example}
    \(\powerset{X}\) on \sigmaalgebra.
  \end{example}
  \begin{example}
    Kaikki joukon \(X\) äärelliset
    osajoukot muodostavat \(\sigma\)-algebran.
  \end{example}
  \begin{example}
    Kun \(X\) on ääretön joukko, kaikki joukon \(X\) korkeintaan
    numeroituvat osajoukot muodostavat \(\sigma\)-renkaan. Jos näiden
    joukkojen komplementit otetaan mukaan, niin saadaan \(\sigma\)-algebra.
  \end{example}
  \begin{example}
    Jos \(\mathcal{B} \subset \powerset{X}\), niin on olemassa
    yksikäsitteinen pienin joukkoalgebra \(\mca\) joka sisältää joukon
    \(\mathcal{B}\).
  \end{example}
\end{examples}

\section{Mitta}

\begin{definition}
  Olkoon \(X\) epätyhjä joukko ja \(\mcr \subset \powerset{X}\)
  rengas. Sanomme, että funktio \(\rho : \mcr \to
  \index{joukkofunktio}
  \extrealnumbers\) on \defterm{joukkofunktio}. Määritellään, että
  \begin{itemize}
    \index{ei-negatiivinen joukkofunktio}
    \item \(\rho\) on \defterm{ei-negatiivinen}, jos \(\rho(E) \geq
      0\) kaikilla \(E \in \mcr\).
      \index{additiivinen joukkofunktio}
    \item \(\rho\) on \defterm{additiivinen}, jos \(\rho(\emptyset) =
      0\) ja \(\rho(E \union F) = \rho(E) + \rho(F)\) kaikille \(E, F
      \in \mcr\), \(E \intersection F = \emptyset\).
      \index{täysin additiivinen joukkofunktio}
    \item \(\rho\) on \defterm{täysin additiivinen}, jos
      \(\rho(\emptyset) = 0\) ja
      \begin{displaymath}
        \rho\left( \bigcup_{k=0}^\infty E_k \right) =
        \sum_{k=0}^\infty \rho(E_k)
      \end{displaymath}
      aina, kun
      \begin{itemize}
        \item[(1)] \(E_k \in \mcr, \; k \in \naturalnumbers\).
        \item[(2)] \(E_i \intersection E_j = \emptyset\) kaikille \(i,
          j \in \naturalnumbers\), \(i \not= j\).
        \item[(3)] \(\union_{k=0}^\infty E_k \in \mcr\).
      \end{itemize}
  \end{itemize}
\end{definition}

\begin{remark}
  Täysin additiivisuudesta seuraa additiivisuus, koska
  \(\rho(\emptyset)\) = 0.
  %% Tod. HT
\end{remark}

\begin{definition}
  Olkoon \(E\) joukko ja \((E_k)_{k=0}^\infty \subset E\). Sanomme,
  \index{kasvava jono joukkoja}
  että \((E_k)_{k=0}^\infty\) on \defterm{kasvava}, jos \(E_k \subset
  E_l\) kaikille \(k, l \in \naturalnumbers, \; k < l\).  Sanomme,
  \index{vähenevä jono joukkoja}
  että \((E_k)_{k=0}^\infty\) on \defterm{vähenevä}, jos \(E_k \subset
  E_l\) kaikille \(k, l \in \naturalnumbers, \; k > l\).
\end{definition}

\begin{theorem}
  \label{th:additiivisuus}
  \cite[lause 5.2]{knapp}
  Olkoon \(\rho\) additiivinen joukkofunktio joukkorenkaassa
  \(\mcr\). Jos \(\rho\) on täysin additiivinen, niin
  \begin{displaymath}
    \rho(E) = \lim_{k \to \infty} \rho(E_k)
  \end{displaymath}
  aina kun \((E_k)_{k=0}^\infty\) on
  kasvava jono renkaan \(\mcr\) alkioita siten, että
  \begin{displaymath}
    E = \bigcup_{k=0}^\infty E_k \in \mcr .
  \end{displaymath}
  Kääntäen, jos
  \begin{displaymath}
    \rho(E) = \lim_{k \to \infty} \rho(E_k)
  \end{displaymath}
  kaikille tällaisille jonoille \((E_k)_{k=0}^\infty\), niin \(\rho\)
  on täysin additiivinen.
\end{theorem}

\begin{proof}
  \cite[lause 5.2]{knapp}
  Oletetaan ensin, että \(\rho\) on täysin additiivinen. Määritellään
  \(F_0 := E_0\) ja \(F_k = E_k \setminus E_{k-1}\) kun \(k \in
  \positiveintegers\). Nyt
  \begin{displaymath}
    E_n := \sum_{k=0}^n F_k
  \end{displaymath}
  ja joukot \(F_k\) ovat erillisiä. Additiivisuuden nojalla
  \begin{displaymath}
    \rho(E_n) = \bigcup_{k=0}^n \rho(F_k) .
  \end{displaymath}
  Edelleen
  \begin{displaymath}
    E = \bigcup_{k=0}^\infty F_k
  \end{displaymath}
  ja täysin additiivisuudesta seuraa
  \begin{displaymath}
    \rho(E) = \sum_{k=0}^\infty F_k = \lim_{n \to \infty} \sum_{k=0}^n
    \rho(F_k) = \lim_{n \to \infty} E_n .
  \end{displaymath}

  Oletetaan sitten, että \((F_n)\) on erillinen jono renkaan \(\mcr\)
  alkioita ja \(\union_{n=0}^\infty F_n \in \mcr\). Määritellään
  \begin{displaymath}
    E_n := \bigcup_{k=0}^n F_k, \spaceafter n \in \naturalnumbers .
  \end{displaymath}
  Nyt \((E_n)\) on kasvava jono renkaan \(\mcr\) alkioita ja
  \begin{displaymath}
    F := \bigcup_{n=0}^\infty E_n \in \mcr .
  \end{displaymath}
  Oletuksen nojalla \(\rho(F) = \lim_{n \to \infty} \rho(E_n)\) ja
  \begin{displaymath}
    \rho(E_n) = \sum_{k=0}^n \rho(F_k)
  \end{displaymath}
  kaikille \(n \in \naturalnumbers\).
  Täten
  \begin{displaymath}
    \rho(F) = \sum_{k=0}^\infty \rho(F_k) .
  \end{displaymath}
\end{proof}

\begin{corollary}
  \cite[seurauslause 5.3]{knapp}
  Olkoon \(X\) epätyhjä joukko ja \(\mca \subset \powerset{X}\)
  joukkoalgebra. Olkoon \(\rho : \mca \to \extrealnumbers\)
  additiivinen joukkofunktio, jolle \(\rho(X) < +\infty\). Jos
  \(\rho\) on täysin additiivinen, niin
  \begin{displaymath}
    \rho(E) = \lim_{k \to \infty} \rho(E_k)
  \end{displaymath}
  aina, kun \((E_k)_{k=0}^\infty\) on vähenevä ja
  \begin{displaymath}
    \bigcap_{k=0}^\infty E_k \in \mca .
  \end{displaymath}
  Kääntäen, jos
  \begin{displaymath}
    \lim_{k \to \infty} \rho(E_k) = \rho(E)
  \end{displaymath}
  kaikille tällaisille jonoille \((E_k)_{k=0}^\infty\), niin \(\rho\)
  on täysin additiivinen.
\end{corollary}

\begin{definition}
  Olkoon \(X\) epätyhjä joukko ja \(\mcr \subset \powerset{X}\)
  \(\sigma\)-rengas. Sanomme, että joukkofunktio \(\rho: \mcr \to
  \index{mitta}
  \extrealnumbers\) on \defterm{mitta}, jos \(\rho\) on
  ei-negatiivinen ja täysin additiivinen.
\end{definition}

Kun mittoja käytetään integraalien kanssa työskentelyyn, niin
\(\sigma\)-rengas \(\mcr\) oletetaan \(\sigma\)-algebraksi. Jos
integrointi määriteltäisiin \(\sigma\)-renkaassa, joka ei ole
\(\sigma\)-algebra, niin nollasta eroavat vakiofunktiot eivät olisi
mitallisia. Oletus, että \(\mcr\) on \(\sigma\)-algebra,
ei vähennä tulosten yleisyyttä: on olemassa kanoninen tapa laajentaa
mitta \(\sigma\)-renkaasta pienimmälle \(\sigma\)-algebralle, joka
sisältää ko. \(\sigma\)-renkaan.

\index{todennäköisyyskenttä}
\index{otosavaruus}
\index{tapahtuma-avaruus}
\index{todennäköisyysmitta}
\hyphenation{toden-näköisyys-mitta}
\defterm{Todennäköisyyskenttä} määritellään mitta-avaruutena
\cite{wpprobaxioms} \((\Omega, F, P)\), missä \(\Omega\) on
\defterm{otosavaruus}, \(\sigma\)-rengas \(F \subset \powerset{\Omega}\) on
\defterm{tapahtuma-avaruus} ja \(P : F \to \intervalcc{0}{1}\) on
\defterm{todennäköisyysmitta}. Todennäköisyysmitan on täytettävä ehto
\(P(\Omega) = 1\).

\section{Borelin joukot}

Olkoon \(X\) metrinen avaruus. Kun \(T \subset \powerset{X}\),
määritellään
\begin{itemize}
  \item[(1)] \(T_\sigma\): kaikki joukon \(T\) alkioiden korkeintaan
    numeroituvat unionit
  \item[(2)] \(T_\delta\): kaikki joukon \(T\) alkioiden korkeintaan
    numeroituvat leikkaukset
  \item[(3)] \(T_{\delta\sigma} := (T_\delta)_\sigma\)
\end{itemize}

\hyphenation{ordinaali-luku}

Nyt määritellään transfiniittisella induktiolla jono \(G^m\), missä
\(m\) on ordinaaliluku, seuraavasti:
\begin{itemize}
  \item[(1)] Olkoon \(G^0\) kaikkien avaruuden \(X\) avointen
    joukkojen joukko.
  \item[(2)] Jos \(i\) ei ole rajaordinaali, niin ordinaalilla \(i\)
    on välittömästi edeltävä ordinaali \(i-1\). Asetetaan
    \begin{displaymath}
      G^i := (G^{i-1})_{\delta\sigma} .
    \end{displaymath}
  \item[(3)] Jos \(i\) on rajaordinaali, asetetaan
    \begin{displaymath}
      G^i := \bigcup_{j<i} G^j .
    \end{displaymath}
\end{itemize}

\index{Borelin algebra}
Nyt \defterm{Borelin algebra} määritellään joukoksi \(G^{\omega_1}\), missä
\(\omega_1\) on ensimmäinen ylinumeroituva ordinaaliluku. Siis Borelin
algebra voidaan generoida avointen joukkojen joukolta iteroimalla
operaatiota \(G \mapsto G_{\delta\sigma}\) ensimmäiseen
ylinumeroituvaan ordinaalilukuun asti.
\index{Borelin joukko}
Borelin algebran alkioita sanotaan \defterm{Borelin joukoiksi}.

\begin{examples}
  \begin{example}
    Borelin algebra reaalilukujen joukossa on pienin
    \(\sigma\)-algebra, joka sisältää kaikki reaalilukuvälit.
  \end{example}
\end{examples}

\begin{definition}
  Jos \(X\) on epätyhjä joukko ja \(\mca\) on \(\sigma\)-algebra
  joukossa \(X\), niin sanomme paria \((X, \mca)\)
  \index{mitallinen avaruus}
  \defterm{mitalliseksi avaruudeksi}.
\end{definition}

\begin{definition}
  Olkoon \(E\) metrinen avaruus ja \(B\) sen Borelin joukkojen
  \index{Borelin avaruus}
  joukko. Sanomme paria \((E, B)\) \defterm{Borelin avaruudeksi}.
\end{definition}

\begin{definition}
  Sanomme täydellistä, separoituvaa ja metrisoituvaa topologista
  \index{puolalainen avaruus}
  avaruutta \defterm{puolalaiseksi avaruudeksi}.
\end{definition}

\begin{theorem}
  Olkoon \(X\) puolalainen avaruus. Tällöin \(X\) on Borelin
  avaruutena isomorfinen jollekin seuraavista:
  \begin{itemize}
    \item[(1)] \(\realnumbers\)
    \item[(2)] \(\integernumbers\) tai
    \item[(3)] äärelliseen avaruuteen.
  \end{itemize}
\end{theorem}

%% TOD.

\begin{definition}
  Kutsumme puolalaiseen avaruuteen liittyvää Borelin avaruutta
  \index{standardi Borelin avaruus}
  \defterm{standardiksi Borelin avaruudeksi}.
\end{definition}

\section{Mitalliset funktiot ja Lebesguen integraali}

\begin{definition}
  Olkoot \((X, \Sigma)\) ja \((Y, T)\) mitallisia avaruuksia. Sanomme,
  \index{mitallinen funktio}
  että funktio \(f : X \to Y\) on \defterm{mitallinen}, jos
  \(\inverseimage{f}{E} \in \Sigma\) jokaiselle \(E \in T\).
\end{definition}

Mitallisuus riippuu \(\sigma\)-algebroista \(\Sigma\) ja \(T\). Jos
\((X, \Sigma)\) ja \((Y, T)\) ovat Borelin avaruuksia, mitallista
\index{Borelin funktio}
funktiota \(f: X \to Y\) kutsutaan myös \defterm{Borelin funktioksi}
\index{Borel-mitallinen funktio}
tai \defterm{Borel-mitalliseksi funktioksi}.

\index{Lebesgue-mitallinen funktio}

\defterm{Lebesgue-mitallinen funktio} on funktio \(f : (\realnumbers,
\mathcal{L}) \to (\complexnumbers, \mathcal{B}_\complexnumbers)\),
missä \(\mathcal{L}\) on Lebesgue-mitallisten joukkojen muodostama
\(\sigma\)-algebra ja \(\mathcal{B}_\complexnumbers\) on
kompleksilukujen \(\complexnumbers\) Borelin
algebra. Lebesgue-mitalliset funktiot ovat kiinnostavia analyysissä,
koska niitä voidaan integroida. Tapauksessa \(f : X \to \realnumbers\)
funktio \(f\) on Lebesgue-mitallinen, jos ja vain jos joukko \(\{f >
\alpha\}\) on mitallinen kaikille \(\alpha \in \realnumbers\). Tämä on
yhtäpitävää sen kanssa, että jokin joukoista \(\{f \geq \alpha\}\),
\(\{f < \alpha\}\) tai \(\{f \leq \alpha\}\) on mitallinen kaikilla
\(\alpha \in \realnumbers\). Jatkuvat funktiot, monotoniset funktiot,
porrasfunktiot ja Riemann-integroituvat funktiot ovat kaikki
Lebesgue-mitallisia. Funktio \(f : X \to \complexnumbers\) on
mitallinen, jos ja vain jos \(\ReAlt f\) ja \(\ImAlt f\) ovat
mitallisia.

\begin{definition}
  Olkoon \(X\) epätyhjä joukko ja \(E \subset X\). Määritellään joukon
  \index{karakteristinen funktio}
  \(E\) \defterm{karakteristinen funktio} asettamalla
  \begin{displaymath}
    \chi_E(x) :=
    \left\{
    \begin{array}{rr}
      1 ; & x \in E \\
      0 ; & x \not\in E
    \end{array}
    \right.
  \end{displaymath}
\end{definition}

\begin{definition}
  Olkoon \(X\) epätyhjä joukko. Sanomme funktiota \(s: X \to
  \index{yksinkertainen funktio}
  \extrealnumbers\) \defterm{yksinkertaiseksi funktioksi}, jos
  \(\setimage{s}{X} \subset \realnumbers\) on äärellinen.
\end{definition}

\index{yksinkertaisen funktion kanoninen kehitelmä}

Jokaisella yksinkertaisella funktiolla \(s\) on \defterm{kanoninen
  kehitelmä}
\begin{displaymath}
  s = \sum_{k=1}^n c_k \chi_{E_k} ,
\end{displaymath}
missä reaaliluvut \(c_k, \; k = 1,\ldots,n\) ovat keskenään erisuuria
ja joukot \(E_k, \; k = 1,\ldots,n\) ovat epätyhjiä ja erillisiä. Itse
asiassa \(\setimage{s}{X} = \{c_k \setsep k=1,\ldots,n\}\) ja \(E_k\)
on se joukko, jossa \(s\) saa arvon \(c_k\). Olkoon \(V :=
\left]c,+\infty\right]\). Joukko \(\inverseimage{s}{V}\) on unioni
    joukoista \(E_k\), joille \(c < c_k\), mistä seuraa, että \(s\) on
    mitallinen funktio, jos ja vain jos kaikki joukot \(E_k\)
    kanonisessa kehitelmässä ovat mitallisia.

Seuraavassa \((X, \mca)\) on mitallinen avaruus ja \(\mu : \mca \to
\extrealnumbers\) mitta. Olkoon \(s\) yksinkertainen funktio ja \(s(x)
\geq 0\) kaikilla \(x \in X\). Olkoon \(E\) mitallinen joukko ja
\begin{displaymath}
  s = \sum_{k=1}^n c_k \chi_{A_k}
\end{displaymath}
funktion \(s\) kanoninen kehitelmä. Määritellään
\begin{displaymath}
  \mathcal{I}_E(s) :=\sum_{k=1}^n c_k \mu(A_k \intersection E) .
\end{displaymath}
Kun \(f : X \to \realnumbers\) ja \(g : X \to \realnumbers\) ovat
funktioita, niin merkitään \(f \leq g\), jos ja vain jos \(f(x) \leq
g(x)\) kaikilla \(x \in X\). Määritellään vastaavasti \(f \geq g\).

\begin{definition}
  Kun \(f : E \to \realnumbers\) on mitallinen funktio, \(0 \leq f\)
  ja \(E \subset X\) on mitallinen joukko, määritellään funktion \(f\)
  \index{Lebesguen integraali}
  \defterm{Lebesguen integraali} joukossa \(E\) mitan \(\mu\) suhteen
  asettamalla
  \begin{eqnarray*}
    \int_E f d\mu & = & \int_E f(x) d\mu(x) \\
    & := & \sup \{ \mathcal{I}_E(s)
    \setsep s \;\textrm{on yksinkertainen funktio ja}\; 0 \leq s
    \leq f \} \in \extrealnumbers .
  \end{eqnarray*}
\end{definition}

Yleiselle mitalliselle funktiolle \(f: E \to \realnumbers\), joka ei
välttämättä ole ei-negatiivinen, integraali voidaan tai ei voida
määritellä. Kirjoitetaan \(f = f^+ - f^-\), missä \(f^+ \geq 0\) ja
\(f^- \geq 0\). Funktiot \(f^+\) ja \(f^-\) ovat mitallisia, joten
\(\int_E f^+ d\mu\) ja \(\int_E f^- d\mu\) ovat hyvin määriteltyjä
joukon \(\extrealnumbers\) alkioita. Jos \(\int_E f^+ d\mu\) ja
\(\int_E f^- d\mu\) eivät ole molemmat äärettömiä, määritellään
\begin{displaymath}
  \int_E f d\mu = \int_E f(x) d\mu(x) := \int_E f^+ d\mu - \int_E f^-
  d\mu .
\end{displaymath}
Tämä määrittely on yhteensopiva erikoistapauksen \(f \geq 0\) kanssa,
koska sellaiselle \(f\) on \(f^- = 0\) ja siten \(\int_E f^- d\mu =
0\). Sanomme, että \(f : E \to \realnumbers\) on
\index{integroituva funktio}
\index{Lebesgue-integroituva funktio}
\defterm{integroituva}, jos \(\int_E f^+ d\mu\) ja \(\int_E f^- d\mu\)
ovat molemmat äärellisiä.

%% yleisempi avaruus kuin \realnumbers^n ?
\begin{definition}
  Olkoon \(n \in \positiveintegers\), \(E\) Borelin avaruus ja \(f :
  E \to \realnumbers\) Borel-mitallinen funktio. Määritellään
  \index{$\esssup$}
  \begin{displaymath}
    \esssup_{x \in E} f(x) := \inf \{ r \in \realnumbers
    \setsep \mu(\{f(x) \setsep f(x) > r\}) = 0 \}
  \end{displaymath}
  ja
  \index{$\essinf$}
  \begin{displaymath}
    \essinf_{x \in E} f(x) := \sup \{ r \in \realnumbers
    \setsep \mu(\{f(x) \setsep f(x) < r\}) = 0 \} .
  \end{displaymath}
\end{definition}

%% \begin{definition}
%%   Olkoon \((X,\mca)\) mitallinen avaruus ja \(\mu : \mca \to
%%   \extrealnumbers\) mitta. Olkoot \(f : X \to
%%   \complexnumbers\) ja \(g : X \to \complexnumbers\)
%%   \index{melkein kaikkialla}
%%   funktioita. Sanomme, että \(f = g\) \defterm{melkein kaikkialla},
%%   jos
%%   \begin{displaymath}
%%     \mu(\{x \in X \setsep f(x) \not= g(x)\}) = 0 .
%%   \end{displaymath}
%% \end{definition}

\begin{definition}
  Olkoon \((X,\mca)\) mitallinen avaruus ja \(\mu : \mca \to
  \extrealnumbers\) mitta. Olkoon \(P(x)\), \(x \in X\), predikaatti.
  \index{melkein kaikkialla}
  Sanomme, että \(P(x)\) on tosi \defterm{melkein kaikkialla},
  jos
  \begin{displaymath}
    \mu(\{x \in X \setsep \lnot P(x) \}) = 0 .
  \end{displaymath}
\end{definition}

\begin{exercises}
  \begin{exercise}{8.1}
    %% L18052015-4
    Olkoon \(X\) epätyhjä joukko ja \(\mcr \subset \powerset{X}\)
    joukkorengas. Olkoon \(\rho : \mcr \to \extrealnumbers\) täysin
    additiivinen joukkofunktio. Osoita, että \(\rho\) on additiivinen.
  \end{exercise}
  \begin{exercise}{8.2}
    %% L18052015-3
    Olkoon \(\mca := \powerset{\naturalnumbers}\) ja
    \begin{displaymath}
      \mu(A) = \# A, \spaceafter A \in \mca .
    \end{displaymath}
    Osoita, että \(\mu\) on mitta.
  \end{exercise}
\end{exercises}

\chapter{Yleisiä funktioavaruuksia}
\label{ch:yleisia-funktioavaruuksia}

\section{Tarvittavia määritelmiä}
\label{sec:yf-maaritelmia}

\begin{definition}
  Olkoon \(n \in \positiveintegers\) ja \(f : \rn \to
  \complexnumbers\) funktio.  Kun \(m \in \positiveintegers\), \(m
  \index{differenssi}
  \geq 2\), määritellään \defterm{differenssit}
  \begin{eqnarray*}
    (\diff{\rnh} f)(x) & := & (\higherdiff{1}{\rnh} f)(\rnx) :=
    f(\rnx+\rnh)-f(\rnx) \\
    (\higherdiff{m}{\rnh} f)(\rnx) & := & (\diff{\rnh}
    (\higherdiff{m-1}{\rnh} f))(\rnx)
  \end{eqnarray*}
  kaikille \(\rnh, \rnx \in \rn\).
  \index{jatkuvuusmoduuli}
  Jos \(f\) on Borel-mitallinen, määritellään \defterm{jatkuvuusmoduuli}
  \begin{displaymath}
    \genmodcont{m}{p}{f}{t} := \sup \{ \norm{\higherdiff{m}{\rnh} f}_p
    \setsep \rnh \in \closedball{\rn}{0}{t}\} .
  \end{displaymath}
  Määritellään edelleen
  \begin{displaymath}
    \modcont{f}{t} := \genmodcont{1}{\infty}{f}{t}, \spaceafter t \in
    \realnumbers_0 .
  \end{displaymath}
\end{definition}

Merkitsemme jatkuvaa Fourier-muunnosta \(\ft\) ja
Fourier-käänteismuunnosta \(\ft^{-1}\). Ks. luku
\ref{sec:fourier-analyysi}.

\section{$\Lp{p}{\realnumbers^n}$}
\label{sec:Lp}

Olkoon \(n \in \positiveintegers\) ja
\(\bor(\realnumbers^n,\complexnumbers)\) kaikkien
Borel-mitallisten funktioiden \(f : \realnumbers^n \to
\complexnumbers\) joukko.  Määritellään ekvivalenssirelaatio \(\sim\)
joukossa \(\bor(\realnumbers^n, \complexnumbers)\)
asettamalla
\begin{displaymath}
  f \sim g \iff \mu(\{\rnx \in \rn \setsep f(\rnx) \not= g(\rnx)\}) =
  0 .
\end{displaymath}

Kun \(p \in [1, \infty[\), määritellään
\begin{equation}
  \label{eq:Lpnorm}
  \norm{C}_p := \norm{f}_p := \left( \int_{\rnx \in \rn}
  \abs{f(\rnx)}^p d\mu
  \right)^{\frac{1}{p}} ,
\end{equation}
missä \(C \in \bor(\realnumbers^n,\complexnumbers) / \sim\) ja \(f \in
C\).
Kun \(p = \infty\), määritellään
\begin{equation}
  \label{eq:Linftynorm}
  \norm{C}_\infty := \norm{f}_\infty := \esssup_{\rnx \in \rn} \abs{f(\rnx)}
\end{equation}
missä \(C \in \bor(\realnumbers^n,\complexnumbers) / \sim\) ja \(f \in
C\).

\begin{remark}
  Kaavojen \eqref{eq:Lpnorm} ja \eqref{eq:Linftynorm} arvo ei riipu
  funktion \(f \in C\) valinnasta.
\end{remark}

Kun \(p \in [1, \infty]\), määritellään Banachin avaruus 
\index{$\Lp{p}{\rn}$}
\begin{displaymath}
  \Lp{p}{\realnumbers^n} := \{ C \in
    \bor(\realnumbers^n,\complexnumbers) / \sim \setsep \norm{C}_p <
    \infty \}
\end{displaymath}

Joukko \(\Lp{2}{\rn}\) on myös Hilbertin avaruus, kun sisätulo
määritellään
\begin{displaymath}
  \ip{C}{D} := \int_{\rnx \in \rn} f(\rnx)^* g(\rnx) d\mu
\end{displaymath}
missä \(C, D \in \Lp{2}{\rn}\), \(f \in C\) ja \(g \in D\).

\begin{examples}
  \begin{example}
    Vakiofunktio \(f(\rnx) = c \in \complexnumbers\), \(\rnx \in \rn\),
    kuuluu avaruuteen \(\Lp{\infty}{\rn}\).
  \end{example}
  \begin{example}
    Jokainen jatkuva ja kompaktikantajainen funktio kuuluu kaikkiin
    avaruuksiin \(\Lp{p}{\rn}\), \(p \in \intervalcc{1}{\infty}\).
  \end{example}
  \begin{example}
    Kohdan 1 vakiofunktio ei kuulu mihinkään avaruuteen
    \(\Lp{p}{\rn}\), \(p < \infty\).
  \end{example}
\end{examples}

\section{$\Cb{T}$}
\label{sec:Cb}

\begin{definition}
  Olkoon \(X\) epätyhjä joukko ja \(f : X \to \complexnumbers\)
  \index{rajoitettu funktio}
  funktio. Sanomme, että f on \defterm{rajoitettu}, jos on olemassa
  luku \(M \in \realnumbers_+\) siten, että \(\abs{f(x)} \leq M\)
  kaikilla \(x \in X\).
\end{definition}

\begin{definition}
  Olkoon \(T\) topologinen avaruus. Määritellään Banachin avaruus
  \index{$\Cb{T}$}
  \begin{displaymath}
    \Cb{T} := \{ f \setsep f \;\textrm{on jatkuva ja
      rajoitettu funktio avaruudelta}\; T \;\textrm{joukolle}\;
    \complexnumbers \}
  \end{displaymath}
  normilla
  \begin{displaymath}
    \norminspace{f}{\Cb{T}} := \sup_{x \in T} \abs{f(x)} .
  \end{displaymath}
\end{definition}

\section{$\Co{T}$}

Ks. myös \cite{lacey1974}.

\begin{definition}
  Olkoon \(T\) lokaalisti kompakti Hausdorffin avaruus. Sanomme, että
  \index{äärettömyydessä häviävä funktio}
  funktio \(f : T \to \complexnumbers\) \defterm{häviää
    äärettömyydessä}, jos joukko \(\{t \in T \setsep \abs{f(t)} \geq
  \varepsilon \}\) on kompakti jokaiselle \(\varepsilon \in
  \realnumbers_+\).
  Määritellään Banachin avaruus
  \index{$\Co{T}$}
  \begin{displaymath}
    \Co{T} := \{ f : T \to \complexnumbers \setsep f \;\textrm{häviää
      äärettömyydessä} \}
  \end{displaymath}
  normilla
  \begin{displaymath}
    \norminspace{f}{\Co{T}} := \sup_{x \in T} \abs{f(x)} .
  \end{displaymath}
\end{definition}
 
\section{$\Cu{E}$}

\begin{definition}
  Olkoon \((E, d)\) metrinen avaruus ja \(f : E \to \complexnumbers\)
  \index{tasaisesti jatkuva}
  funktio. Sanomme, että \(f\) on \defterm{tasaisesti jatkuva}, jos
  \begin{displaymath}
    \forall \varepsilon \in \realnumbers_+ : \exists \delta \in
    \realnumbers_+ : \forall x, y \in E : \left( d(x,y) < \delta
    \implies \abs{f(x)-f(y)} < \varepsilon \right) .
  \end{displaymath}
\end{definition}

\begin{definition}
  Olkoon \(E\) metrinen avaruus. Määritellään Banachin avaruus
  \index{$\Cu{E}$}
  \begin{displaymath}
    \Cu{E} := \{ f : E \to \complexnumbers \setsep f \; \textrm{on rajoitettu ja tasaisesti jatkuva} \}
  \end{displaymath}
  normilla
  \begin{displaymath}
    \norminspace{f}{\Cu{E}} := \sup_{x \in E} \abs{f(x)} .
  \end{displaymath}
\end{definition}

\section{$\Ccom{T}$}

Olkoon \(T\) lokaalisti kompakti Hausdorffin avaruus. Määritellään
\index{$\Ccom{T}$}
\begin{displaymath}
  \Ccom{T} := \{ f \in \Cb{T} \setsep \supp f \;\textrm{on kompakti}\}
\end{displaymath}
normilla
\begin{displaymath}
  \norminspace{f}{\Ccom{T}} := \sup_{x \in T} \abs{f(x)} .
\end{displaymath}
\(\Ccom{T}\) ei ole täydellinen ja se on tiheä avaruudessa \(\Co{T}\)
\cite{donoho1992}.

\section{$C^m(\rn)$}

Tässä \(m \in \naturalnumbers\).
Ks. \cite[luku 2.2.2]{schaefer1971}.

\index{$C^m(\rn)$}
\begin{align*}
  & C^m(\rn) := \{ f : \rn \to \complexnumbers \setsep \mideriv{\alpha} f \in
  \Cu{\rn} \;\textrm{kaikille}\; \abs{\alpha} \leq m \} \\
  & \norminspace{f}{C^m(\rn)} := \sum_{\abs{\alpha} \leq
    m} \norminspace{\mideriv{\alpha} f}{\Lp{\infty}{\rn}}, \spaceafter f \in
  C^m(\rn) 
\end{align*}

\begin{remark}
  \(C^0(\rn) \equaltvs C_u(\rn)\).
\end{remark}

\section{Hölderin avaruus $\holderspace{s}{\rn}$}

Ks. \cite[luku 2.2.2]{schaefer1971}.
Tässä \(s > 0\) ja \(s \not\in \integernumbers\).
Asetetaan
\begin{displaymath}
  s = [s] + \{s\} ,
\end{displaymath}
missä \([s] \in \integernumbers\) ja \(0 \leq \{s\} < 1\)
ja määritellään
\index{Hölderin avaruus}
\index{$\holderspace{s}{\rn}$}
\defterm{Hölderin avaruus}
\begin{displaymath}
  \holderspace{s}{\rn} := \left\{ f \in C^{[s]}(\rn)
  \bigsetsep \norminspace{f}{\holderspace{s}{\rn}} < \infty \right\} ,
\end{displaymath}
missä normi määritellään
\begin{displaymath}
  \norminspace{f}{\holderspace{s}{\rn}}
  := \norminspace{f}{C^{[s]}(\rn)} + \sum_{\abs{\alpha} = [s]} \sup
  \left\{ \frac{\abs{(\mideriv{\alpha} f)(\rnx) - (\mideriv{\alpha}
      f)(\rny)}}{\norm{\rnx-\rny}^{\{s\}}} \bigsetsep \rnx, \rny \in
  \rn, x \not= y \right\} .
\end{displaymath}

\section{Zygmundin avaruus $\zygmundspace{s}{\rn}$}

Ks. \cite[luku 2.2.2]{schaefer1971}.
Tässä \(s > 0\).
Asetetaan
\begin{displaymath}
  s = [s]^- + \{s\}^+ ,
\end{displaymath}
missä \([s]^- \in \integernumbers\) ja \(0 < \{s\}^+ \leq 1\)
ja määritellään
\index{Zygmundin avaruus}
\index{$\zygmundspace{s}{\rn}$}
\defterm{Zygmundin avaruus}
\begin{displaymath}
  \zygmundspace{s}{\rn} := \left\{ f \in C^{[s]^-}(\rn)
  \bigsetsep \norminspace{f}{\zygmundspace{s}{\rn}} < \infty \right\} ,
\end{displaymath}
missä normi määritellään
\begin{align*}
  \norminspace{f}{\zygmundspace{s}{\rn}}
  := & \norminspace{f}{C^{[s]^-}(\rn)} \\
  & + \sum_{\abs{\alpha} = [s]^-}
  \sup \left\{
  \norm{\rnh}^{-\{s\}^+} \norminspace{\higherdiff{2}{\rnh}\mideriv{\alpha}
    f}{\Cu{\rn}} \bigsetsep \rnh \in \rn \setminus \{0\} \right\} .
\end{align*}

\begin{remark}
  Kun \(s \not\in \integernumbers\), niin \(\zygmundspace{s}{\rn}
  \equaltvs \holderspace{s}{\rn}\) \cite[huomautus
    2.2.2/3]{triebel1983}.
\end{remark}

\section{Sobolevin avaruus $\sobolevspace{m}{p}{\rn}$}

\hyphenation{distribuutio-derivaatta}

Ks. \cite[luku 2.2.2]{schaefer1971}.
Kun \(1 < p < \infty\) ja \(m \in \positiveintegers\), määritellään
\index{Sobolevin avaruus}
\index{$\sobolevspace{m}{p}{\rn}$}
\defterm{Sobolevin avaruus}
\begin{displaymath}
  \sobolevspace{m}{p}{\rn} :=
  \left\{
  f \in \Lp{p}{\rn}
  \setsep
  \norminspace{f}{\sobolevspace{m}{p}{\rn}} < \infty
  \right\} ,
\end{displaymath}
missä normi määritellään
\begin{displaymath}
  \norminspace{f}{\sobolevspace{m}{p}{\rn}} :=
  \sum_{\abs{\alpha} \leq m} \norminspace{\mideriv{\alpha} f}{\Lp{p}{\rn}} .
\end{displaymath}
Asetetaan lisäksi \(\sobolevspace{0}{p}{\rn} := \Lp{p}{\rn}\).
Tässä derivaatta \(\mideriv{\alpha} f\) on distribuutioderivaatta.

\section{Slobodeckij'n avaruus $\slobspace{s}{p}{\rn}$}

Ks. \cite[luku 2.2.2]{schaefer1971}.
Kun \(1 \leq p < \infty\) ja \(s > 0\), \(s \not\in \integernumbers\),
määritellään
\index{Slobodeckij'n avaruus}
\index{$\slobspace{s}{p}{\rn}$}
\defterm{Slobodeckij'n avaruus}
\begin{displaymath}
  \slobspace{s}{p}{\rn} :=
  \left\{
  f \in \sobolevspace{[s]}{p}{\rn}
  \setsep
  \norminspace{f}{\slobspace{s}{p}{\rn}} < \infty
  \right\} ,
\end{displaymath}
missä normi määritellään
\begin{align*}
  \norminspace{f}{\slobspace{s}{p}{\rn}} := &
  \norminspace{f}{\sobolevspace{[s]}{p}{\rn}} \\
  & +
  \sum_{\abs{\alpha} = [s]}
  \left(
  \int_{(\rnx,\rny) \in \rn \times \rn} \frac{\abs{(\mideriv{\alpha} f)(\rnx) -
      (\mideriv{\alpha} f)(\rny)}^p}{\norm{\rnx-\rny}^{n+\{s\}p}} d\mu
  \right)^{\frac{1}{p}} .
\end{align*}

Tässä derivaatta \(\mideriv{\alpha} f\) on distribuutioderivaatta.

\section{Bessel-potentiaaliavaruus $\besselpotspace{s}{p}{\rn}$}

Ks. \cite[luku 2.2.2]{schaefer1971}.
Kun \(s \in \realnumbers\) ja \(1 < p < \infty\) määritellään
\index{Bessel-potentiaaliavaruus}
\index{$\besselpotspace{s}{p}{\rn}$}
\defterm{Bessel-potentiaaliavaruus}
\begin{displaymath}
  \besselpotspace{s}{p}{\rn} := \{ f \in \tempdistrspace{\rn}
  \setsep \norminspace{f}{\besselpotspace{s}{p}{\rn}} < \infty \}
\end{displaymath}
missä normi määritellään
\begin{displaymath}
  \norminspace{f}{\besselpotspace{s}{p}{\rn}}
  :=
  \norminspace{\ft^{-1}\left(1+\norm{\rnx}^2\right)^{s/2}\ft
    f}{\Lp{p}{\rn}} .
\end{displaymath}

\section{Paikallinen Hardyn avaruus $\localhardyspace{p}{\rn}$}

Ks. \cite[luku 2.2.2]{schaefer1971}.
Olkoon \(0 < p < \infty\) ja \(\varphi\) testifunktio, jolle
\(\varphi(0) = 1\).
Määritellään
\(\varphi_t(x) := \varphi(tx)\) kaikille \(x \in \rn\) ja \(t >
\index{paikallinen Hardyn avaruus}
\index{$\localhardyspace{p}{\rn}$}
0\). Määritellään \defterm{paikallinen Hardyn avaruus} asettamalla
\begin{displaymath}
  \localhardyspace{p}{\rn}
  :=
  \left\{
  f \in \tempdistrspace{\rn}
  \setsep
  \norminspace{f}{\localhardyspace{p}{\rn}}^\varphi < \infty
  \right\}
\end{displaymath}
missä normi määritellään
\begin{displaymath}
  \norminspace{f}{\localhardyspace{p}{\rn}}^\varphi
  :=
  \bignorminspace{\sup_{0 < t < 1} \abs{\ft^{-1}\varphi_t\ft
      f}}{\Lp{p}{\rn}} .
\end{displaymath}

\section{$\bmospace{\rn}$}

Ks. \cite[luku 2.2.2]{schaefer1971}.

\begin{definition}
  Olkoon \(n \in \positiveintegers\). Olkoot \(a_k, b_k \in
  \realnumbers\), \(a_k < b_k\), \(k = 1,\ldots,n\). Sanomme joukkoa
  \begin{displaymath}
    Q := \{ (x_1,\ldots,x_n) \in \rn \setsep a_k \leq x_k \leq b_k, k
    = 1,\ldots,n\}
  \end{displaymath}
  \index{kuutio}
  \defterm{kuutioksi}.
  Määritellään lisäksi
  \begin{displaymath}
    \abs{Q} := \prod_{k=1}^n (b_k-a_k) .
  \end{displaymath}
\end{definition}

Jos \(f : \rn \to \complexnumbers\) on Lebesgue-integroituva funktio
ja \(Q\) on kuutio avaruudessa \(\rn\), niin asetetaan
\begin{displaymath}
  f_Q := \frac{1}{\abs{Q}} \int_{\rnx \in Q} f(\rnx) d\mu .
\end{displaymath}
Määritellään
\index{$\bmospace{\rn}$}
\begin{displaymath}
  \bmospace{\rn}
  :=
  \left\{
  f
  \setsep
  f \;\textrm{on paikallisesti Lebesgue-integroituva joukossa}\; \rn,
  \norminspace{f}{\bmospace{\rn}} < \infty
  \right\} ,
\end{displaymath}
missä normi määritellään
\begin{displaymath}
  \norminspace{f}{\bmospace{\rn}}
  :=
  \sup_{\abs{Q} \leq 1} \left( \frac{1}{\abs{Q}} \int_{\rnx \in Q}
  \abs{f(\rnx)-f_Q} d\mu \right)
  +
  \sup_{\abs{Q} > 1} \left( \frac{1}{\abs{Q}} \int_{\rnx \in Q}
  \abs{f(\rnx)} d\mu \right) .
\end{displaymath}

\section{Besovin ja Triebel-Lizorkinin avaruudet}

%%Tässä luvussa on käytetty kirjaa \cite{schaefer1971}.

\begin{definition}
  \cite[luku 1.2.2]{schaefer1971}
  Olkoon \(n \in \positiveintegers\) ja \(p, q \in
  \intervaloc{0}{\infty}\). Olkoon \((f_k)_{k=0}^\infty\) jono
  Borel-mitallisia funktioita \(f : \rn \to
  \complexnumbers\). Määritellään
  \index{$\norm{\cdot \;\vert\; L^p(l^q)}$}
  \begin{displaymath}
    \norm{f_k \;\vert\; L^p(l^q)} := \bignorm{\norm{f_k(\cdot) \vert l^q}
      \big\vert L^p}
    =
    \left(
    \int_{\rnx \in \rn} \left( \sum_{k=0}^\infty
    \abs{f_k(\rnx)}^q\right)^{p/q} d\mu
    \right)^{1/p}
  \end{displaymath}
  ja
  \index{$\norm{\cdot \;\vert\; l^q(L^p)}$}
  \begin{displaymath}
    \norm{f_k \;\vert\; l^q(L^p)} := \bignorm{\norm{f_k(\cdot) \vert L^p}
      \big\vert l^q}
    =
    \left(
    \sum_{k=0}^\infty \left( \int_{\rnx \in \rn} \abs{f_k(\rnx)}^p d\mu \right)^{q/p}
    \right)^{1/q}
\end{displaymath}
\end{definition}

\begin{remark}
  Yllä oleviin integraalikaavoihin tulee muutoksia, jos \(p=\infty\)
  tai \(q = \infty\).
\end{remark}

\begin{definition}
  \cite[luku 1.2.2]{schaefer1971}
  Olkoon \(A\) (reaalinen tai kompleksinen) vektoriavaruus. Funktiota
  \index{kvasinormi}
  \(\norminspace{\cdot}{A}\) sanotaan \defterm{kvasinormiksi}, jos se
  toteuttaa tavalliset normin aksioomat (N1), (N2) ja (N3)
  lukuunottamatta kolmioepäyhtälöä (N4), joka korvataan aksioomalla
  \begin{displaymath}
    \norminspace{a_1+a_2}{A} \leq c(\norminspace{a_1}{A}
    + \norminspace{a_2}{A})
  \end{displaymath}
  jollekin vakiolle \(c \in \realnumbers_+\), joka ei riipu
  vektoreista \(a_1\) ja \(a_2\).
\end{definition}

\begin{definition}
  \index{$\Phi(\rn)$}
  \cite[luku 2.3.1]{schaefer1971}.
  Olkoon \(\Phi(\rn)\) kaikkien sellaisten systeemien \(\varphi :=
  (\varphi_j)_{j=0}^\infty \subset S(\rn)\) joukko, joille
  \begin{itemize}
    \item[(1)] \(\supp \varphi_0 \subset \{\rnx \in \rn \setsep
      \norm{\rnx} \leq 2 \}\).
    \item[(2)] \(\supp \varphi_j \subset \{\rnx \in \rn \setsep
      2^{j-1} \leq \norm{\rnx} \leq 2^{j+1}\}, \spaceafter j \in
        \positiveintegers\).
  \end{itemize}
\end{definition}

\begin{definition}
  \cite[luku 2.3.1]{schaefer1971}.
  Olkoon \(s \in \realnumbers\) ja \(0 < q \leq \infty\). Olkoon
  \(\varphi := (\varphi_j)_{j=0}^\infty \in \Phi(\rn)\).
  \begin{itemize}
  \item[(1)]
    \index{Besovin avaruus}
    \index{$\besovspace{s}{p}{q}{\rn}$}
    Jos \(0 < p \leq \infty\), niin määritellään \defterm{Besovin avaruus}
    \begin{displaymath}
      \besovspace{s}{p}{q}{\rn} := \{ f \in \tempdistrspace{\rn}
      \setsep \norminspace{f}{\besovspace{s}{p}{q}{\rn}}^\varphi < \infty \}
    \end{displaymath}
    missä
    \begin{displaymath}
      \norminspace{f}{\besovspace{s}{p}{q}{\rn}}^\varphi
      := \norminspace{2^{sj}\ft^{-1}\varphi_j\ft{}f}{l^q(L^p(\rn))} .
    \end{displaymath}
  \item[(2)]
    Jos \(0 < p < \infty\), niin määritellään
    \index{Triebel-Lizorkinin avaruus}
    \index{$\tlspace{s}{p}{q}{\rn}$}
    \defterm{Triebel-Lizorkinin avaruus}
    \begin{displaymath}
      \tlspace{s}{p}{q}{\rn} := \{ f \in \tempdistrspace{\rn}
      \setsep \norminspace{f}{\tlspace{s}{p}{q}{\rn}}^\varphi <
      \infty \}
    \end{displaymath}
    missä
    \begin{displaymath}
      \norminspace{f}{\tlspace{s}{p}{q}{\rn}}^\varphi
      := \norminspace{2^{sj}\ft^{-1}\varphi_j\ft{}f}{L^p(\rn,l^q)} .
    \end{displaymath}
  \end{itemize}      
\end{definition}

\begin{remark}
  \(\ft^{-1} \varphi_j \ft f = \ft^{-1}(\varphi_j \ft f)\) on
  analyyttinen funktio avaruudessa \(\rn\).
  Kvasinormit \(\norminspace{f}{\besovspace{s}{p}{q}{\rn}}^\varphi\)
  ja \(\norminspace{f}{\tlspace{s}{p}{q}{\rn}}^\varphi\) riippuvat
  tietysti valitusta systeemistä \(\varphi \in \Phi(\rn)\), mutta ne
  ovat kaikki ekvivalentteja kvasinormeja kvasinormiavaruuksissa
  \(\besovspace{s}{p}{q}{\rn}\) ja \(\tlspace{s}{p}{q}{\rn}\).
  Jos \(p,q \in [1, \infty]\), niin \(\besovspace{s}{p}{q}{\rn}\)
  Banachin avaruus.
  Jos \(p \in [1, \infty[\) ja \(q \in [1, \infty]\), niin
      \(\tlspace{s}{p}{q}{\rn}\) on Banachin avaruus \cite[theorem
        2.3.3]{triebel1983}.
\end{remark}

%% Määritelmä Triebelin kirjasta.
%% Pitääkö sekä upotuksen että sen käänteisfunktion olla jatkuvia?
\index{$\contemb$}
Jos \(A_1\) ja \(A_2\) ovat kvasinormiavaruuksia, niin \(A_1 \contemb
\index{jatkuva upotus}
A_2\) tarkoittaa, että avaruus \(A_1\) on \defterm{jatkuvasti
  upotettu} avaruuteen \(A_2\), ts. on olemassa vakio \(t \in
\realnumbers_+\) siten, että \(\norminspace{a}{A_1} \leq
t \norminspace{a}{A_2}\) kaikille \(a \in A_1\) \cite[luku
  2.3.2]{schaefer1971}.

\begin{theorem}
  \cite[lause 2.3.2/2]{triebel1983}
  \begin{itemize}
    \item[(1)] Olkoon \(0 < q_0 \leq q_1 \leq \infty\), \(0 < p \leq
      \infty\) ja \(s \in \realnumbers\). Tällöin
      \begin{displaymath}
        \besovspace{s}{p}{q_0}{\rn}
        \contemb \besovspace{s}{p}{q_1}{\rn} .
      \end{displaymath}
    \item[(2)] Olkoon \(0 < q_0 \leq \infty\), \(0 < q_1 \leq
      \infty\), \(0 < p \leq \infty\), \(s \in
      \realnumbers\) ja \(\varepsilon > 0\). Tällöin
      \begin{displaymath}
        \besovspace{s+\varepsilon}{p}{q_0}{\rn}
        \contemb \besovspace{s}{p}{q_1}{\rn} .
      \end{displaymath}
    \item[(3)] Olkoon \(0 < q \leq \infty\), \(0 < p < \infty\) ja \(s
      \in \realnumbers\). Tällöin
      \begin{displaymath}
        \besovspace{s}{p}{\min\{p,q\}}{\rn}
        \contemb
        \tlspace{s}{p}{q}{\rn}
        \contemb
        \besovspace{s}{p}{\max\{p,q\}}{\rn} .
      \end{displaymath}
  \end{itemize}
\end{theorem}

\begin{theorem}
  \cite[lause 2.3.3]{triebel1983}
  Olkoon \(n \in \positiveintegers\). Kun \(p, q \in \intervaloc{0}{\infty},
  s \in \realnumbers\) niin meillä on voimassa jatkuvat upotukset
  \begin{displaymath}
    S(\rn) \contemb \besovspace{s}{p}{q}{\rn} \contemb \tempdistrspace{\rn} .
  \end{displaymath}
  Kun \(p, q \in \intervaloo{0}{\infty},
  s \in \realnumbers\) niin meillä on voimassa jatkuvat upotukset
  \begin{displaymath}
    S(\rn) \contemb \tlspace{s}{p}{q}{\rn} \contemb \tempdistrspace{\rn} .
  \end{displaymath}
  Kun \(s \in \realnumbers\) ja \(p, q \in \intervaloo{0}{\infty}\),
  niin \(S(\rn)\) on tiheä avaruuksissa \(\besovspace{s}{p}{q}{\rn}\)
  ja \(\tlspace{s}{p}{q}{\rn}\).
\end{theorem}

Kun \(s \in \intervaloc{0}{\infty}\) ja \(p,q \in [1, \infty]\), niin Besovin
avaruudelle on voimassa ekvivalentti määritelmä
\begin{displaymath}
  \besovspace{s}{p}{q}{\rn} := \{ f \in \Lp{p}{\rn} \setsep
  \besovaltnorm{s}{p}{q}{\rn}{m}{f} < \infty \}
\end{displaymath}
missä \(m \in \positiveintegers\), \(m > s\), ja ekvivalentti normi on
\begin{displaymath}
  \besovaltnorm{s}{p}{q}{\rn}{m}{f} := \norminspace{f}{\Lp{p}{\rn}}
  + \norminspace{t \in \intervaloc{0}{1} \mapsto 
    t^{-s-1/q} \genmodcont{m}{p}{f}{t}}{\Lp{q}{\intervaloc{0}{1}}}
\end{displaymath}

Monet yleisesti käytetyt funktioavaruudet voidaan esittää Besovin ja
Triebel-Lizorkinin avaruuksien avulla. Olkoon \(n \in
\positiveintegers\). Meillä on \cite[huomautus 2.2.2/3 ja luku
  2.3.5]{triebel1983}
\begin{displaymath}
  \begin{array}{llll}
    \index{$\zygmundspace{s}{\rn}$}
    \zygmundspace{s}{\rn} & \equaltvs & \besovspace{s}{\infty}{\infty}{\rn}
    , & \textrm{jos}\; s > 0, \\
    \index{$\sobolevspace{s}{p}{\rn}$}
    \sobolevspace{s}{p}{\rn} & \equaltvs & \besovspace{s}{p}{p}{\rn},
    & \textrm{jos}\; 1 \leq p < \infty, s > 0 \;\textrm{ja}\;
    s \not\in \integernumbers, \\
    \index{$\besselpotspace{s}{p}{\rn}$}
    \besselpotspace{s}{p}{\rn} & \equaltvs & \tlspace{s}{p}{2}{\rn},
    & \textrm{jos}\; 1 < p < \infty, s \in \realnumbers, \\
    \index{$\localhardyspace{p}{\rn}$}
    \localhardyspace{p}{\rn} & \equaltvs & \tlspace{0}{p}{2}{\rn}, & 0 <
    p < \infty, \\
    \index{$\bmospace{\rn}$}
    \bmospace{\rn} & \equaltvs & \tlspace{0}{\infty}{2}{\rn} . & 
  \end{array}
\end{displaymath}

\begin{exercises}
  \begin{exercise}{9.1}
    %% L21082009-2
    Olkoon \(n, m \in \positiveintegers\), \(f \in
    \complexnumbers^{\rn}\), \(a \in \realnumbers \setminus
    \zeroset\) ja \(\rnb \in \rn\).
    Olkoon \(g(x) = f(a\rnx-\rnb)\), \(\rnx \in \rn\).
    Osoita, että
    \begin{displaymath}
      \left( \higherdiff{m}{\rnh} g\right) \left( \rnx \right) =
      \left( \higherdiff{m}{a\rnh} f\right)\left(a \rnx - \rnb\right)
    \end{displaymath}
    kaikille \(\rnx, \rnh \in \rn\).
  \end{exercise}
  \begin{exercise}{9.2}
    %% L21082009-3
    Olkoon \(n, m \in \positiveintegers\), \(p \in
    \intervalcc{1}{\infty}\), \(f \in \Lp{p}{\rn}\), \(a
    \in \realnumbers \setminus
    \zeroset\), \(\rnb \in \rn\) ja \(t \in \realnumbers_+\).
    Olkoon \(g(\rnx) = f(a \rnx - \rnb)\), \(\rnx \in \rn\).
    Osoita, että
    \begin{displaymath}
      \genmodcont{m}{p}{g}{t} = \abs{a}^{-\frac{n}{p}}
      \genmodcont{m}{p}{f}{\abs{a}t} .
    \end{displaymath}
  \end{exercise}
  \begin{exercise}{9.3}
    %% L10062009-1 ja L18052015-6
    Olkoon \(n \in \positiveintegers\) ja \(f : \rn \to
    \complexnumbers\) Borel-mitallinen. Osoita, että \(f\) on
    tasaisesti jatkuva, jos ja vain jos
    \begin{displaymath}
      \lim_{t \to 0} \modcont{f}{t} = 0 .
    \end{displaymath}
  \end{exercise}
  \begin{exercise}{9.4}
    %% L05062009-3
    Olkoon \(E\) metrinen avaruus. Osoita, että \(\Ccom{E} \subset
    \Cu{E}\).
  \end{exercise}
  \begin{exercise}{9.5}
    %% L05062009-4
    Olkoon \(E\) metrinen avaruus. Osoita, että \(\Co{E} \subset
    \Cu{E}\).
  \end{exercise}
  \begin{exercise}{9.6}
    %% L05062009-11
    Olkoon \(n, m \in \positiveintegers\). Osoita, että
    \(C^m(\rn) \subset \zygmundspace{m}{\rn}\).
  \end{exercise}
  \begin{exercise}{9.7}
    %% L18052015-7
    Olkoon \(p \in \intervalco{1}{\infty}\) ja \(f_k \in \Lp{p}{\rn}
    \intersection \Lp{\infty}{\rn}\) kaikille \(k \in
    \naturalnumbers\). Olkoon \(g \in \Lp{p}{\rn}\) ja \(h \in
    \Lp{\infty}{\rn}\) siten, että \(\norm{f_k-g}_p \to 0\) ja
    \(\norminfty{f_k-h} \to 0\), kun \(k \to \infty\). Osoita, että
    \(g = h\) melkein kaikkialla.
  \end{exercise}
\end{exercises}

\appendix

\chapter{Tehtävien ratkaisut}

\begin{solution}{1.1}
  \begin{itemize}
  \item[(a)]
    \begin{enumerate}
    \item refleksiivisyys:
      \begin{displaymath}
        x \sim x \iff x^2 = x^2
      \end{displaymath}
      tosi kaikilla \(x \in \realnumbers\).
    \item symmetrisyys:
      \begin{displaymath}
        x \sim y \iff x^2 = y^2 \iff y^2 = x^2 \iff y \sim x
      \end{displaymath}
      kaikilla \(x, y \in \realnumbers\).
    \item transitiivisuus:
      \begin{eqnarray*}
        x \sim y \land y \sim z & \iff & x^2 = y^2 \land y^2 = z^2 \\
        & \implies & x^2 = z^2 \\
        & \iff & x \sim z
      \end{eqnarray*}
      kaikilla \(x,y,z \in \realnumbers\).
    \end{enumerate}
  \item[(b)] jäännösluokka-avaruus:
    \(x^2 = y^2 \iff x = \pm y\), joten
    \begin{displaymath}
      \quotientspace{\realnumbers}{\sim} = \{ \{ -a, a \} \setsep a
      \in \realnumbers \} .
    \end{displaymath}
  \item[(c)]
    \(1 \sim 1\) ja \(1 \sim -1\), joten \(\sim\) ei ole funktio.
  \end{itemize}
\end{solution}

\begin{solution}{2.1}
  Olkoon \(G\) ryhmä ja \(a \in G\). Oletetaan, että \(e' \circ a = a\)
  jollekin \(e' \in G\). Nyt \(e \circ a = e' \circ a\), mistä seuraa,
  että \(e = e \circ a \circ a^{-1} = e' \circ a \circ a^{-1} = e'\).
\end{solution}

\begin{solution}{2.2}
  Olkoon \(G\) ryhmä ja \(a \in G\). Oletetaan, että \(a \circ b = a
  \circ c = a\) joillekin \(b, c \in G\). Nyt \(b = a^{-1} \circ a
  \circ b = a^{-1} \circ a \circ c = c\).
\end{solution}

\begin{solution}{2.3}
  Olkoon \(p, q \in P_n\) ja \(k \in \{1,\ldots,n\}\).
  Nyt \((p \circ q)(k) = p(q(k))\). Koska \(k \in \{1,\ldots,n\}\), niin
  \(q(k) \in \{1,\ldots,n\}\), mistä seuraa, että \(p(q(k)) \in
  \{1,\ldots,n\}\).

  Olkoot \(p, q, r \in P_n\). Nyt
  \begin{eqnarray*}
    \left(\left(p \circ q\right) \circ r\right)\left( k \right)
    & = & \left( p \circ q \right)\left( r \left( k \right) \right)
    = p\left(q\left(r\left(k\right)\right)\right)
    = p\left(\left(q \circ r\right)\left(k\right)\right) \\
    & = & \left(p \circ \left(q \circ r)\right)\right)\left(k\right)
  \end{eqnarray*}
  kaikille \(k \in \{1,\ldots,n\}\).

  Identtinen kuvaus \(\id{P_n}\) on neutraalialkio joukossa \(P_n\).
  Koska \( p : \{1,\ldots,n\} \onto \{1,\ldots,n\} \) on bijektio,
  niin sillä on olemassa käänteisfunktio \( p^{-1} : \{1,\ldots,n\}
  \onto \{1,\ldots,n\} \in P_n \).
\end{solution}

\begin{solution}{2.4}
  Jos \(A\) on ryhmän \(G\) aliryhmä, niin \(A\) on itsekin ryhmä,
  joten \(a \circ b \in A\) ja \(a^{-1} \in A\) kaikille \(a, b \in
  A\).

  Oletetaan sitten, että \(A\) on suljettu operaatioiden \(a \circ b\)
  ja \(a^{-1}\) suhteen. Aksiooma (G1) on ilmiselvä. Olkoon \(a\)
  jokin joukon \(A\) mielivaltainen alkio. Nyt oletuksen nojalla
  \(a^{-1} \in A\) ja \(e = a^{-1} \circ a \in A\).
\end{solution}

\begin{solution}{3.1}
  Olkoon \(\varepsilon \in \rationalnumbers_+\). Asetetaan
  \begin{displaymath}
    n := \max \left\{ \szceil{\frac{3}{\varepsilon}}, 1 \right\} .
  \end{displaymath}
  Nyt
  \begin{displaymath}
    n \geq \frac{3}{\varepsilon} \implies \frac{3}{n} \leq \varepsilon .
  \end{displaymath}
  Olkoot \(j, k \in \naturalnumbers\), \(j, k > n\) ja \(j < k\).
  Nyt
  \begin{displaymath}
    \szabs{\frac{1}{j+1} - \frac{1}{k+1}} \leq \frac{1}{j+1} +
    \frac{1}{k+1}
    \leq \frac{1}{j} + \frac{1}{k} \leq \frac{1}{n} + \frac{1}{n} =
    \frac{2}{n} < \varepsilon .
  \end{displaymath}
\end{solution}

\begin{solution}{3.2}
  \((a,b) = (a,0) + (b,0) \cdot (0,1) = a + b \imagunit\).
\end{solution}

\begin{solution}{4.1}
  \begin{displaymath}
    \begin{pmatrix}
      -50 \\
      8 \\
      -6
    \end{pmatrix}
  \end{displaymath}
\end{solution}

\begin{solution}{4.2}
  \begin{displaymath}
    \begin{pmatrix}
      -20 & 12 \\
      -5 & 179 \\
      -5 &  92
    \end{pmatrix}
  \end{displaymath}
\end{solution}

\begin{solution}{4.3}
  Olkoon \(x, y \in \ker f\). Nyt \(f(x) = f(y) = 0\). Koska \(f\) on
  lineaarinen, niin \(f(x+y) = f(x) + f(y) = 0\), mistä seuraa, että
  \(x + y \in \ker f\).
  Olkoon \(a \in K\) ja \(z \in \ker f\). Nyt \(f(z)=0\) ja koska
  \(f\) on lineaarinen, niin \(f(az) = af(z) = 0\), joten \(az \in
  \ker f\). Lauseen \ref{th:vektorialiavaruus} nojalla \(\ker
  f \vectorsubspace V\).
\end{solution}

\begin{solution}{4.4}
  Olkoon \(x, y \in \image f\). Nyt \(x = f(a)\) jollekin \(a \in V\)
  ja \(y = f(b)\) jollekin \(b \in V\). Funktion \(f\)
  lineaarisuudesta seuraa, että \(x + y = f(a) + f(b) = f(a+b) \in
  \image f\).  Olkoon \(a \in K\) ja \(z \in \ker f\). Nyt \(z =
  f(c)\) jollekin \(c \in V\). Funktion \(f\) lineaarisuudesta seuraa,
  että \(az = af(c) = f(ac) \in \image f\). Lauseen
  \ref{th:vektorialiavaruus} nojalla \(\image f \vectorsubspace W\).
\end{solution}

\begin{solution}{4.5}
  Olkoon \(\cnvec{v} \in \complexnumbers^{Z(m)\times{}Z(1)}\).
  \begin{displaymath}
    (g(\cnvec{v}))_k = (B\cnvec{v})_k = \sum_{j=1}^m B_{k,j}
    \vecelem{\cnvec{v}}{j}
  \end{displaymath}
  ja
  \begin{eqnarray*}
    (f(g(\cnvec{v})))_q & = & \sum_{i=1}^p A_{q,i} (B\cnvec{v})_i
    = \sum_{j=1}^m \sum_{i=1}^p A_{q,i} B_{i,j} \vecelem{\cnvec{v}}{j}
    = \sum_{j=1}^m (AB)_{q,j} \vecelem{\cnvec{v}}{j} \\
    & = & (AB\cnvec{v})_q .
  \end{eqnarray*}
\end{solution}

\begin{solution}{4.6}
  \begin{itemize}
    \item[(i)]
      Olkoon \(v \in V\). Nyt \(0_K + 0_K = 0_K\), joten \((0_K +
      0_K)v = 0_Kv\).
      Toisaalta (V6):n nojalla \((0_K + 0_K)v = 0_Kv + 0_Kv\).
      Siten \(0_Kv + 0_Kv = 0_Kv\). Lisäämällä puolittain \(-(0_Kv)\)
      saadaan \(0_Kv + 0_V = 0_V\). (V3):n nojalla \(0_Kv = 0_V\).
    \item[(ii)]
      Olkoon \(v \in V\). Nyt (i)-kohdan nojalla \(v+(-1_K)v =
      (1_K+(-1_K))v = 0_Kv = 0_V\) ja (V4):n nojalla
      \(v+(-v)=0_V\). Siis \(v + (-1_K)v = v + (-v)\). Lisäämällä
      puolittain \(-v\) saadaan \((-1_K)v = -v\).
  \end{itemize}
\end{solution}

\begin{solution}{5.1}
  Oletetaan ensin, että \(f\) on jatkuva. Olkoon \(R
  \subset \secondtopspace\) avoin. Jos \(\inverseimage{f}{R} =
  \emptyset\), niin se on avoin. Oletetaan, että \(\inverseimage{f}{R}
  \not= \emptyset\). Olkoon \(x \in \inverseimage{f}{R}\). Nyt \(f(x)
  \in R\) ja \(R\) on pisteen \(f(x)\) avoin ympäristö. Koska \(f\) on
  jatkuva, niin on olemassa pisteen \(x\) avoin ympäristö \(P\) siten,
  että \(\setimage{f}{P} \subset R\), mistä seuraa, että \(P \subset
  \inverseimage{f}{R}\). Lauseen \ref{th:avoin-joukko} nojalla
  \(\inverseimage{f}{R}\) on avoin.

  Oletetaan sitten, että
  \begin{displaymath}
    \forall R \subset \secondtopspace : ( R \;\text{on avoin}\; \implies
    \inverseimage{f}{R} \;\text{on avoin} ) .
  \end{displaymath}
  Olkoon \(x \in \firsttopspace\), \(y := f(x)\), ja \(R\) pisteen
  \(y\) avoin ympäristö.
  Nyt \(x \in \inverseimage{f}{R}\) ja oletuksen nojalla
  \(\inverseimage{f}{R}\) on avoin. Joukko \(\inverseimage{f}{R}\) on
  pisteen \(x\) avoin ympäristö, ja
  \(\setimage{f}{\inverseimage{f}{R}} = R \subset R\).
  Piste \(x\) ja ympäristö \(R\) olivat mielivaltaisia, joten \(f\) on
  jatkuva.
\end{solution}

\begin{solution}{5.2}
  Oletetaan ensin, että määritelmä \ref{def:raja-arvo} toteutuu.
  Olkoon \(\varepsilon \in \realnumbers_+\) ja \(x \in E_1\) ja olkoon
  \(Y\) pisteen \(f(x)\) ympäristö siten, että
  \(\openball{E_2}{f(x)}{\varepsilon} \subset Y\).
  Oletuksen nojalla on olemassa pisteen \(x\) ympäristö \(Z \subset
  E_1\), jolle \(\setimage{f}{Z} \subset
  \openball{E_2}{f(x)}{\varepsilon}\)
  Nyt \(\openball{E}{x}{\delta} \subset Z\) jollekin \(\delta \in
  \realnumbers_+\) ja
  \begin{displaymath}
    d_1(y,x) < \delta \implies y \in Z \implies d_2(f(y), f(x)) <
    \varepsilon ,
  \end{displaymath}
  missä \(y \in E_1\) on mielivaltainen. Koska myös \(x \in E_1\) oli
  mielivaltainen, niin ehto \eqref{eq:metr-av-raja-arvo} toteutuu.

  Oletetaan sitten, että ehto \eqref{eq:metr-av-raja-arvo} toteutuu.
  Olkoon \(\varepsilon \in \realnumbers_+\) ja \(x \in E_1\). Olkoon
  \(Y\) pisteen \(f(x)\) ympäristö siten, että
  \(\openball{E_2}{f(x)}{\varepsilon} \subset Y\).
  Olkoon \(y \in E_1\) siten, että \(d_1(x, y) < \delta\).
  Oletuksen nojalla \(d_2(f(x), f(y)) < \varepsilon\),
  joten \(f(y) \in \openball{E_2}{f(x)}{\varepsilon}\).
  Siis
  \begin{displaymath}
    \setimage{f}{\openball{E_1}{x}{\delta}} \subset
    \openball{E_2}{f(x)}{\varepsilon} ,
  \end{displaymath}
  ja \(\openball{E_1}{x}{\delta}\) on pisteen \(x\) ympäristö. Piste
  \(x\) oli mielivaltainen, joten määritelmä \ref{def:raja-arvo}
  toteutuu.
\end{solution}

\begin{solution}{5.3}
  Oletetaan ensin, että ehto \eqref{eq:metr-av-jonon-raja-arvo}
  toteutuu. Olkoon \(Y\) pisteen \(a\) jokin ympäristö. Nyt
  \(\openball{E}{a}{\varepsilon} \subset Y\) jollakin \(\varepsilon
  \in \realnumbers_+\). Oletuksen nojalla on olemassa \(m \in
  \naturalnumbers\) siten, että \(k \geq m \implies d(x_k, a) <
  \varepsilon\), ts. \(k \geq m \implies x_k \in
  \openball{E}{a}{\varepsilon} \subset Y\).

  Oletetaan sitten, että määritelmän \ref{def:jonon-raja-arvo} ehdot
  ovat voimassa. Olkoon \(\varepsilon \in \realnumbers_+\). Nyt \( Y
  := \openball{E}{y}{\varepsilon}\) on pisteen \(y\)
  ympäristö. Oletuksen nojalla on olemassa \(m \in \naturalnumbers\)
  siten, että \(x_k \in Y\) kaikilla \(k \geq m\). Ts. \(k \geq m
  \implies d(x_k, y) < \varepsilon\).
\end{solution}

\begin{solution}{5.4}
  Jos \(B = \emptyset\) or \(B = E\), niin väite on tosi.  Oletetaan,
  että \(B \not= \emptyset\) ja \(B \not= E\).
  Jos oletetaan, että \(B\) on suljettu, niin oikeanpuoleinen väite
  seuraa lauseesta \ref{th:suljettu-joukko}, koska jono on verkon
  erikoistapaus.

  Oletetaan sitten, että
  \begin{displaymath}
    x_k \to a \in E \implies a \in B .
  \end{displaymath}
  jokaiselle verkolle \((x_k)_{k=0}^\infty \subset B\).
  Olkoon \(x \in E \setminus B\). Ei ole olemassa jonoa \((x_k)
  \subset B\) siten, että \(x_k \to x\).
  Oletetaan, että \(Y \intersection B \not= \emptyset\) jokaiselle
  pisteen \(x\) ympäristölle \(Y\) (avaruudessa \(E\)). Olkoon
  \begin{displaymath}
    V_k := \szopenball{E}{x}{\frac{1}{k+1}}, \spaceafter k \in
    \naturalnumbers .
  \end{displaymath}
  Kun \(k \in \naturalnumbers\) niin valitaan \(V_k \in Y
  \intersection B\). Nyt \(x_k \to x\), mikä on ristiriita, joten
  pisteellä \(x\) on ympäristö \(Y_0\) siten, että \(Y_0 \intersection
  B = \emptyset\). Piste \(x\) oli mielivaltainen, joten lauseen
  \ref{th:avoin-joukko} nojalla \(A \setminus B\) on avoin, ja siten
  \(B\) on suljettu.
\end{solution}

\hyphenation{viipale-filtteri}
\hyphenation{viipale-filtteri-kannan}

\begin{solution}{5.5}
  Oletetaan ensin, että \(E\) on täydellinen.  Olkoon \(Z\) avaruuden
  \(E\) sellainen tasainen struktuuri, että se generoi avaruuden \(E\)
  normaalin topologian.  Olkoon \(\seq{x} := (x_k)_{k=0}^\infty
  \subset E\) Cauchyn jono. Jonon \(\seq{x}\) viipalefiltterikannan
  muodostavat viipaleet \( \{ x_n \setsep n \geq n_0\}\), missä \(n_0
  \in \naturalnumbers\) vaihtelee. Olkoon \(x \in E\) ja \(W \in
  Z\). Koska \(\{\openball{E}{x}{r} \setsep r \in \realnumbers_+\}\)
  on ympäristökanta pisteessä \(x\) ja \(Y(W, x)\) on pisteen \(x\)
  ympäristö, niin \(\openball{E}{x}{r_0} \subset Y(W, x)\) jollekin \(r_0
  \in \realnumbers_+\).  Nyt
  \begin{equation}
    \label{eq:wx}
    d(x,y) < r_0 \implies y \in Y(W, x) \iff
    (x,y) \in W
  \end{equation}
  Olkoon \(\varepsilon \in \intervaloo{0}{r_0}\).  Koska \(\seq{x}\)
  on Cauchyn jono, on olemassa \(N \in \naturalnumbers\) siten, että
  \(d(x_n,x_m) < \varepsilon\) kaikile \(n, m \geq N\).  Kun \(F = \{
  x_n \setsep n \geq N \}\) on jokin viipalefiltterikannan alkio, niin \(F
  \times F = \{ (x_n, x_m) \setsep n \geq N, m \geq N\}\). Kun \((x,
  y) \in F \times F\), niin \(d(x,y) < \varepsilon\), joten yhtälön
  \eqref{eq:wx} nojalla \((x,y) \in W\). Täten jonon \((x_n)\)
  viipalefiltteri on Cauchyn filtteri, ja oletuksen nojalla se
  suppenee.

  Oletetaan sitten, että jokainen Cauchyn jono suppenee avaruudessa
  \(E\). Olkoon \((F)\) Cauchyn filtteri avaruudessa \(E\).
  Olkoon \(x \in E\) ja \(r \in \realnumbers_+\). Koska \(x\)-keskiset
  pallot ovat ympäristökanta pisteessä \(x\), niin voidaan valita
  \(W\) siten, että \(\openball{E}{x}{r} \subset Y(W, x)\). Merkitään
  tätä ympäristöä \(W_r(x)\). Valitaan \(F_r\) siten, että \(F_r
  \times F_r \subset W_r\). Nyt \(x, y \in F_r \implies d(x,y) <
  r\). Valitaan \(z_m \in F_{1/(m+1)}\) jokaiselle \(m \in
  \naturalnumbers\).
  Jono \((z_m)\) on Cauchyn jono, joten oletuksen nojalla se suppenee
  kohti jotakin avaruuden \(E\) alkiota \(z\).
  On osoitettava että \((F)\) suppenee kohti pistettä \(z\).
  Avoimet \(z\)-keskiset pallot ovat ympäristökanta pisteessä \(z\),
  joten on osoitettava, että jokaiselle \(z\)-keskiselle avoimelle
  pallolle on olemassa filtterin alkio, joka sisältyy palloon.
  Olkoon \(\varepsilon \in \intervaloo{0}{2}\). Voidaan valita \(m_0 \in
  \naturalnumbers\) siten, että \(d(z,z_m) < \frac{\varepsilon}{2}\)
  kaikilla \(m \in \naturalnumbers\) ja \(m \geq m_0\), missä \(m_0 \in
  \naturalnumbers\).
  Olkoon
  \begin{displaymath}
    m := \max \left\{ m_0, \left\lceil \frac{2-\varepsilon}{\varepsilon}
    \right\rceil \right\} + 1 .
  \end{displaymath}
  Nyt \(d(x,z_m) < \frac{1}{m+1} < \frac{\varepsilon}{2}\), mistä
  seuraa, että \(d(x,z) \leq d(x,z_m) + d(z_m,z) <
  \frac{\varepsilon}{2} + \frac{\varepsilon}{2} = \varepsilon\) .
  Siis \(x \in \openball{E}{z}{\varepsilon}\). Täten \(F_{1/(m+1)}
  \subset \openball{E}{z}{\varepsilon}\) ja filtteri \((F)\) suppenee
  kohti pistettä \(z\).
\end{solution}

\begin{solution}{6.1}
  Olkoon \(f : \rn \to \rn\) lineaarinen funktio ja \(M(f)\) sen
  matriisi.  Olkoot \(\cnvec{x},\cnvec{y} \in \complexnumbers^n\). Nyt
  \(\ip{\cnvec{x}}{f(\cnvec{y})} = \cnvec{x}^\herm M(f) \cnvec{y} =
  (\cnvec{y}^\herm (M(f))^\herm \cnvec{x})^* =
  \ip{\cnvec{y}}{(M(f))^\herm \cnvec{x}}^* = \ip{(M(f))^\herm
    \cnvec{x}}{\cnvec{y}}\). Toisaalta \(\ip{\cnvec{x}}{f(\cnvec{y})}
  = \ip{f^\herm(\cnvec{x})}{\cnvec{y}}\). Siten \(M(f^\herm) =
  (M(f))^\herm\).
\end{solution}

\begin{solution}{6.2}
  Funktio \(f\) on määritelmän nojalla lineaarinen injektio. Olkoon
  \(x \in l^1\). Nyt \(\norm{x}_1 \geq \norm{x}_\infty\), joten
  \(\norm{f(x)}_\infty \leq \norm{x}_1\). Siis \(f\) on jatkuva.

  Olkoon
  \begin{displaymath}
    \seqelem{\seq{y}_n}{k}
    :=
    \left\{
    \begin{array}{ll}
      1 ; & k < n \\
      0; & k \geq n
    \end{array}
    \right.
  \end{displaymath}
  missä \(k, n \in \naturalnumbers\).
  Nyt \(\norm{\seq{y}_n}_1 = n\) ja \(\norm{\seq{y}_n}_\infty = 1\) kaikille \(n
  \in \naturalnumbers\).
  Edelleen
  \begin{displaymath}
    \sup \{ \norm{f^{-1}(z)} \setsep \seq{z} \in l^\infty \land
    \norm{\seq{z}} \leq 1 \} = \sup \{ \norm{\seq{z}}_1 \setsep
    \seq{z} \in l^\infty \land \norm{\seq{z}}_\infty \leq 1 \} =
    \infty .
  \end{displaymath}
  Siis \(f^{-1}\) ei ole rajoitettu eikä siten myöskään jatkuva.
\end{solution}

\begin{solution}{6.3}
  \begin{align*}
    f(x) = & f(x_0)
    + \frac{f'(x_0)}{1!}(x-x_0)
    + \frac{f''(x_0)}{2!}(x-x_0)^2 \\
    & + \frac{f^{(3)}(x_0)}{3!}(x-x_0)^3
    + \frac{f^{(4)}(x_0)}{4!}(x-x_0)^4
    + R_4(x)
  \end{align*}
\end{solution}

\begin{solution}{6.4}
  \begin{align*}
  f(\rnx) = & f(\rnxo) \\
  & + \sum_{i_1=1}^3 \pd{x_{i_1}} f(\rnxo) (x_{i_1} - x^{(0)}_{i_1}) \\
  & + \sum_{i_1,i_2=1}^3 \pd{x_{i_1}} \pd{x_{i_2}} f(\rnxo) (x_{i_1} -
  x^{(0)}_{i_1}) (x_{i_2} - x^{(0)}_{i_2}) \\
  & + \sum_{i_1,i_2,i_3=1}^3 \pd{x_{i_1}} \pd{x_{i_2}} \pd{x_{i_3}} f(\rnxo) (x_{i_1} -
  x^{(0)}_{i_1}) (x_{i_2} - x^{(0)}_{i_2}) (x_{i_3} - x^{(0)}_{i_3}) \\
  & + \sum_{i_1,i_2,i_3,i_4=1}^3 \pd{x_{i_1}} \pd{x_{i_2}}
  \pd{x_{i_3}} \pd{x_{i_4}} f(\rnxo) (x_{i_1} -
  x^{(0)}_{i_1}) (x_{i_2} - x^{(0)}_{i_2}) (x_{i_3} - x^{(0)}_{i_3})
  (x_{i_4} - x^{(0)}_{i_4}) \\
  & + R_4(\rnx)
  \end{align*}
\end{solution}

\begin{solution}{7.1}
  Olkoon \(\varphi \in \tfs{\rn}\).
  Olkoon
  \begin{displaymath}
    g := \varphi\left( \frac{\cdot + \rnb}{a} \right) .
  \end{displaymath}
  Nyt
  \begin{eqnarray*}
    \distrappl{\tau_{a,\rnb}\tau_{c,\rnd}T}{\varphi}
    & = & \frac{1}{\abs{a}^n} \distrappl{\tau_{c,\rnd} T}{g} \\
    & = & \frac{1}{\abs{a}^n} \frac{1}{\abs{c}^n}
    \szdistrappl{T}{g\left(\frac{\cdot+\rnd}{c}\right)}
  \end{eqnarray*}
  ja
  \begin{displaymath}
    g\left(\frac{\rnx+\rnd}{c}\right) = \varphi\left(
    \frac{\frac{\rnx+\rnd}{c}+\rnb}{a}\right) =
    \varphi\left(\frac{\rnx+\rnd+c\rnb}{ac}\right) , \spaceafter \rnx \in \rn .
  \end{displaymath}
  Siis
  \begin{displaymath}
    \distrappl{\tau_{a,\rnb}\tau_{c,\rnd}T}{\varphi} = \frac{1}{\abs{a}^n}
    \frac{1}{\abs{c}^n}
    \szdistrappl{T}{\varphi\left(\frac{\cdot+\rnd+c\rnb}{ac}\right)}
    = \distrappl{\tau_{ac,c\rnb+\rnd}T}{\varphi} .
  \end{displaymath}
\end{solution}

\begin{solution}{7.2}
  Kun \(\varphi, \psi \in \tfs{\rn}\), niin
  \begin{displaymath}
    \delta(\varphi+\psi) = \left(\varphi+\psi\right)\left(0\right) =
    \varphi(0) + \psi(0) = \delta(\varphi) + \delta(\psi) .
  \end{displaymath}
  Kun \(a \in \realnumbers\) ja \(\varphi \in \tfs{\rn}\), niin
  \begin{displaymath}
    \delta(a\varphi) = \left(a\varphi\right)(0) = a \varphi(0) = a
    \delta(\varphi) .
  \end{displaymath}
  Siis \(\delta\) on lineaarinen.

  Oletetaan, että \((\varphi_k)_{k=0}^\infty \subset \tfs{\rn}\) ja
  \(\varphi_k \to \psi \in \tfs{\rn}\).
  Nyt \(\varphi_k \to \psi\) tasaisesti, eli
  \begin{equation}
    \label{eq:tas-supp-a}
    \lim_{k \to \infty} \norminfty{\psi-\varphi_k} = 0 .
  \end{equation}
  Edelleen
  \begin{equation}
    \label{eq:tas-supp-b}
    \abs{\psi(0) - \varphi_k(0)} \leq \sup_{\rnx \in \rn}
    \abs{\psi(\rnx) - \varphi_k(\rnx)}
    = \norminfty{\psi-\varphi_k}
  \end{equation}
  kaikille \(k \in \naturalnumbers\).
  Yhtälöistä \eqref{eq:tas-supp-a} ja \eqref{eq:tas-supp-b} seuraa,
  että
  \begin{displaymath}
    \lim_{k \to \infty} \abs{\psi(0)-\varphi_k(0)} = 0 ,
  \end{displaymath}
  eli
  \begin{displaymath}
    \lim_{k \to \infty} \abs{\delta(\psi)-\delta(\varphi_k)} = 0 ,
  \end{displaymath}
  mistä seuraa, että
  \begin{displaymath}
    \delta(\varphi_k) \to \delta(\psi) , \spaceafter k \to \infty .
  \end{displaymath}
  Siis
  \begin{displaymath}
    \lim_{k \to \infty} \delta(\varphi_k) = \delta(\psi) =
    \delta\left( \lim_{k \to \infty} \varphi_k \right) ,
  \end{displaymath}
  joten \(\delta \in \distrspace{\rn}\).
\end{solution}

\begin{solution}{8.1}
  Olkoot \(E, F \in \mcr\) ja \(E \intersection F =
  \emptyset\). Määritellään
  \begin{itemize}
    \item \(E_0 := E\).
    \item \(E_1 := F\).
    \item \(E_k := \emptyset\) kaikille \(k \in \naturalnumbers\), \(k
      \geq 2\).
  \end{itemize}
  Nyt
  \begin{displaymath}
    \rho \left( \bigcup_{k=0}^\infty E_k \right) = \rho(E \union F) =
    \rho(E) + \rho(F) .
  \end{displaymath}
\end{solution}

\begin{solution}{8.2}
  \(\mca\) on \(\sigma\)-algebra, ja \(\mu\) on joukkofunktio. Funktio
  \(\mu\) on määritelmänsä nojalla ei-negatiivinen. Meillä on
  \(\mu(\emptyset) = \# \emptyset = 0\).

  Oletetaan, että
  \begin{enumerate}
    \item \(E_k \in \mca\) kaikille \(k \in \naturalnumbers\).
    \item \(E_i \intersection E_j = \emptyset\) kaikille \(i, j \in
      \naturalnumbers\), \(i \not= j\).
    \item \(\cup_{k=0}^\infty E_k \in \mca \).
  \end{enumerate}
  Olkoon
  \begin{displaymath}
    F := \bigcup_{k=0}^\infty E_k .
  \end{displaymath}
  Jos \(F\) on äärellinen, niin
  \begin{displaymath}
    F := \bigcup_{k=0}^n E_{K(k)} ,
  \end{displaymath}
  missä \(n \in \naturalnumbers\), \(K(k) \in \naturalnumbers\)
  kaikille \(k = 1, \ldots, n\).
  Koska \(F\) on äärellinen ja joukot \(E_k\) ovat erillisiä, niin
  \begin{displaymath}
    \left(\forall k \in \{ 1, \ldots, n \} : i \not= K(k)\right) \implies E_i =
    \emptyset .
  \end{displaymath}

  Oletetaan sitten, että \(F\) on ääretön. Olkoon
  \begin{displaymath}
    G_k := \bigcup_{i=1}^k E_i, \spaceafter k \in \naturalnumbers .
  \end{displaymath}
  Nyt
  \begin{displaymath}
    F = \bigcup_{k=0}^\infty G_k .
  \end{displaymath}
  Olkoon \(m \in \naturalnumbers\). Koska \(F\) on ääretön, niin on
  olemassa \(k_0 \in \naturalnumbers\) siten, että
  \begin{displaymath}
    \mu\left(G_{k_0}\right) = \# G_{k_0} > m .
  \end{displaymath}
  Nyt \(\mu(G_k) > m \) kaikille \(k \geq k_0\).
  Täten lauseen \ref{th:additiivisuus} nojalla
  \begin{displaymath}
    \infty = \mu(F) = \lim_{k \to \infty} \mu\left(G_k\right) .
  \end{displaymath}
  Siis \(\mu\) on ei-negatiivinen ja täysin additiivinen ja siten mitta.
\end{solution}

\begin{solution}{9.1}
  Todistetaan induktiolla.  Nyt
  \begin{eqnarray*}
    (\higherdiff{}{\rnh} g)(\rnx) & = & g(\rnx + \rnh) - g(\rnx)
    = f(a\rnx+a\rnh-\rnb) - f(a\rnx-\rnb) \\
    & = & (\higherdiff{1}{a\rnh})(a\rnx-\rnb),
  \end{eqnarray*}
  kaikille \(\rnx, \rnh \in \rn\).
  Siis väite on tosi, kun \(m = 1\).

  Oletetaan, että lause on tosi jollakin \(m_1 \in
  \positiveintegers\). Nyt
  \begin{eqnarray*}
    (\higherdiff{m_1+1}{\rnh} g)(\rnx) & = & \left( \higherdiff{1}{\rnh}
    \left(\higherdiff{m_1}{\rnh} g\right)\right)(\rnx)
    = (\higherdiff{m_1}{\rnh} g)(\rnx + \rnh) - (\higherdiff{m_1}{\rnh} g)(\rnx) \\
    & = & (\higherdiff{m_1}{a\rnh} f)(a\rnx+a\rnh-\rnb) - (\higherdiff{m_1}{a\rnh}
    f)(a\rnx-\rnb) \\
    & = & \left( \higherdiff{1}{a\rnh} \left( \higherdiff{m_1}{a\rnh} f
    \right)\right)(a\rnx-\rnb) \\
    & = & (\higherdiff{m_1+1}{a\rnh} f)(a\rnx-\rnb) .
  \end{eqnarray*}
  Siis väite on tosi, kun \(m = m_1 + 1\).
\end{solution}

\begin{solution}{9.2}
  Käytetään tehtävää 9.1.
  \begin{eqnarray*}
    \genmodcont{m}{p}{g}{t} & = & \sup \{ \norm{\higherdiff{m}{h} g}_p
      \setsep h \in \closedball{\rn}{0}{t} \} \\
      & = & \sup \{ \norm{(\higherdiff{m}{ah} f)(a\cdot-\rnb)}_p
      \setsep h \in \closedball{\rn}{0}{t} \} \\
      & = & \abs{a}^{-\frac{n}{p}} \sup \{ \norm{\higherdiff{m}{ah}
        f}_p \setsep h \in \closedball{\rn}{0}{t} \} \\
      & = & \abs{a}^{-\frac{n}{p}} \sup \{ \norm{\higherdiff{m}{k}
        f}_p \setsep k \in \closedball{\rn}{0}{\abs{a}t} \} \\
      & = & \abs{a}^{-\frac{n}{p}} \genmodcont{m}{p}{f}{\abs{a}t} .
  \end{eqnarray*}
\end{solution}

\begin{solution}{9.3}
  Oletetaan ensin, että \(f\) on tasaisesti jatkuva.  Olkoon
  \(\varepsilon \in \realnumbers_+\).  Koska \(f\) on tasaisesti
  jatkuva, on olemassa \(t_0 \in \realnumbers_+\) siten, että kaikille
  \(\rnx, \rny \in \rn\), \(\norm{\rny-\rnx}_2 < t_0\) on \(\abs{f(y)
    - f(x)} < \frac{\varepsilon}{2}\).
  Olkoon \(t \in \intervaloo{0}{t_0}\) ja \(\rnh \in
  \closedball{\rn}{0}{t}\).
  Nyt
  \begin{displaymath}
    \szabs{\left(\diff{\rnh} f\right)\left(\rnx\right)}
    =
    \szabs{f(\rnx+\rnh) - f(\rnx)} < \frac{\varepsilon}{2}
  \end{displaymath}
  kaikille \(\rnx \in \rn\).
  Täten
  \begin{displaymath}
    \norminfty{\diff{\rnh} f} \leq \frac{\varepsilon}{2} < \varepsilon .
  \end{displaymath}
  Koska \(\rnh\) oli mielivaltainen, niin
  \begin{displaymath}
    \modcont{f}{t} = \sup \{ \norminfty{\diff{\rnh} f} \setsep \rnh
    \in \closedball{\rn}{0}{t}\} \leq \varepsilon .
  \end{displaymath}

  Oletetaan sitten, että
  \begin{displaymath}
    \lim_{t \to 0} \modcont{f}{t} = 0 .
  \end{displaymath}
  Olkoon \(\varepsilon \in \realnumbers_+\). On olemassa \(\delta \in
  \realnumbers_+\) siten, että \(\modcont{f}{\delta} < \varepsilon\).
  Olkoon \(\rnx, \rny \in \rn\), \(\norm{\rny - \rnx}_2 < \delta\)
  ja \(\rnh := \rny - \rnx\). Nyt
  \begin{displaymath}
    \abs{f(\rny) - f(\rnx)} = \abs{\left(\diff{\rnh}
      f\right)\left(\rnx\right)}
    \leq \norminfty{\diff{\rnh} f}
    \leq \modcont{f}{\delta} < \varepsilon .
  \end{displaymath}
\end{solution}

\begin{solution}{9.4}
  Olkoon \(f \in \Ccom{E}\). Olkoon \(A := \supp f\) ja \(h \in
  \realnumbers_+\). Jokaisella \(x \in A\) on olemassa \(t(x) \in
  \realnumbers_+\) siten, että
  \begin{displaymath}
    \setimage{f}{\openball{E}{x}{t(x)}} \subset
    \szopenball{\complexnumbers}{f(x)}{\frac{h}{2}}
  \end{displaymath}
  Olkoon
  \begin{displaymath}
    K(x) := \szopenball{E}{x}{\frac{1}{2}t(x)} \intersection A,
    \spaceafter x \in A .
  \end{displaymath}
  Nyt \(P := \{ K(x) \setsep x \in A \}\) on topologisen avaruuden
  \(A\) avoin peite. Koska \(A\) on kompakti niin peitteellä \(P\) on
  äärellinen alipeite \(P'\),
  \begin{displaymath}
    P' := \{ K(x_k) \setsep k = 1, \ldots, m \}
  \end{displaymath}
  missä \(m \in \positiveintegers\) ja \(x_k \in A\), \(k = 1, \ldots,
  n\).
  Olkoon
  \begin{displaymath}
    r := \min \{ t(x_k) \setsep k = 1,\ldots,m \}.
  \end{displaymath}
  Olkoot \(x, y \in A \land d(x,y) < \frac{r}{2}\). Nyt \(x \in
  K(x_{k_0})\) jollekin \(k_0 \in \{ 1, \ldots, m \}\) ja \(d(x,
  x_{k_0}) < \frac{1}{2} t(x_{k_0})\). Nyt
  \begin{equation}
    \label{eq:ccom-a}
    \abs{f(x)-f(x_{k_0})} < \frac{h}{2} .
  \end{equation}
  Edelleen
  \begin{eqnarray*}
    d(y, x_{k_0}) & \leq & d(y, x) + d(x, x_{k_0}) \\
    & < & \frac{r}{2} + \frac{1}{2} t(x_{k_0}) \\
    & < & \frac{1}{2} t(x_{k_0}) + \frac{1}{2} t(x_{k_0}) \\
    & = & t(x_{k_0}) ,
  \end{eqnarray*}
  joten \(y \in \openball{E}{x_{k_0}}{t(x_{k_0})}\), mistä seuraa,
  että
  \begin{equation}
    \label{eq:ccom-b}
    \abs{f(y)-f(x_{k_0})} < \frac{h}{2} .
  \end{equation}
  Yhtälöistä \eqref{eq:ccom-a} ja \eqref{eq:ccom-b} seuraa, että
  \begin{equation}
    \label{eq:ccom-c}
    \abs{f(x)-f(y)} \leq \abs{f(x)-f(x_{k_0})} + \abs{f(x_{k_0}) -
      f(y)} < \frac{h}{2} + \frac{h}{2} = h .
  \end{equation}
  Oletetaan sitten, että \(x \in A\), \(y \in E \setminus A\) ja
  \(d(x, y) < \frac{r}{2}\). Nyt \(x \in K(x_{k_1})\) jollakin \(k_1 \in
  \{1, \ldots, m\}\) ja \(d(x, x_{k_1}) < \frac{1}{2} t(x_{k_1})\). Edelleen
  \begin{eqnarray*}
    d(y, x_{k_1}) & \leq & d(y, x) + d(x, x_{k_1}) \\
    & < & \frac{r}{2} + \frac{1}{2} t(x_{k_1}) \\
    & < & \frac{1}{2} t(x_{k_1}) + \frac{1}{2} t(x_{k_1}) \\
    & = & t(x_{k_1}) ,
  \end{eqnarray*}
  joten \(y \in \openball{E}{x_{k_1}}{t(x_{k_1})}\), mistä seuraa,
  että
  \begin{equation}
    \label{eq:ccom-d}
    \abs{f(x)-f(y)} \leq \abs{f(x)-f(x_{k_1})} + \abs{f(x_{k_1}) -
      f(y)} < \frac{h}{2} + \frac{h}{2} = h .
  \end{equation}
  Tasainen jatkuvuus seuraa yhtälöistä \eqref{eq:ccom-c} ja
  \eqref{eq:ccom-d}.
\end{solution}

\begin{solution}{9.5}
  Olkoon \(f \in \Co{E}\). Olkoon \(h \in \realnumbers_+\).
  Määritellään
  \begin{displaymath}
    A := \{ x \in E \setsep \abs{f(x)} \geq \frac{h}{3} \} .
  \end{displaymath}
  Koska \(f \in \Co{E}\), niin \(A\) on kompakti.
  Jokaisella \(x \in A\) on olemassa \(t(x) \in
  \realnumbers_+\) siten, että
  \begin{displaymath}
    \setimage{f}{\openball{E}{x}{t(x)}} \subset
    \szopenball{\complexnumbers}{f(x)}{\frac{h}{2}}
  \end{displaymath}
  Olkoon
  \begin{displaymath}
    K(x) := \szopenball{E}{x}{\frac{1}{2}t(x)} \intersection A,
    \spaceafter x \in A .
  \end{displaymath}
  Nyt \(P := \{ K(x) \setsep x \in A \}\) on topologisen avaruuden
  \(A\) avoin peite. Koska \(A\) on kompakti niin peitteellä \(P\) on
  äärellinen alipeite \(P'\),
  \begin{displaymath}
    P' := \{ K(x_k) \setsep k = 1, \ldots, m \}
  \end{displaymath}
  missä \(m \in \positiveintegers\) ja \(x_k \in A\), \(k = 1, \ldots,
  n\).
  Olkoon
  \begin{displaymath}
    r := \min \{ t(x_k) \setsep k = 1,\ldots,m \}.
  \end{displaymath}
  \begin{enumerate}
    \item
      Olkoot \(x, y \in A \land d(x,y) < \frac{r}{2}\). Nyt \(x \in
      K(x_{k_0})\) jollekin \(k_0 \in \{ 1, \ldots, m \}\) ja
      \begin{eqnarray*}
        d(y, x_{k_0}) & \leq & d(y,x) + d(x,x_{k_0}) \\
        & < & \frac{r}{2} + \frac{1}{2} t(x_{k_0}) \\
        & < & \frac{1}{2} t(x_{k_0}) + \frac{1}{2} t(x_{k_0}) \\
        & = & t(x_{k_0}) .
      \end{eqnarray*}
      Täten \(y \in \openball{E}{x_{k_0}}{t(x_{k_0})}\), joten
      \(\abs{f(y)-f(x_{k_0})} < \frac{h}{2}\). Tästä seuraa, että
      \begin{equation}
        \label{eq:co-a}
        \abs{f(x)-f(y)} \leq \abs{f(x)-f(x_{k_0})} + \abs{f(x_{k_0}) -
          f(y)}
        < \frac{h}{2} + \frac{h}{2} = h .
      \end{equation}
    \item
      Olkoot \(x, y \in E \setminus A\). Nyt
      \begin{equation}
        \label{eq:co-b}
        \abs{f(x) - f(y)} \leq \abs{f(x)} + \abs{f(y)} \leq
        \frac{h}{3} + \frac{h}{3} < h
      \end{equation}
    \item
      Olkoot \(x \in A\), \(y \in E \setminus A\) ja \(d(x, y) <
      \frac{r}{2}\).
      Nyt \(x \in K(x_{k_1})\) jollakin \(k_1 \in \{1, \ldots, m\}\)
      ja \(d(x, x_{k_1}) < \frac{1}{2}t(x_{k_1})\). Edelleen
      \begin{eqnarray*}
        d(y,x_{k_1}) & \leq & d(y, x) + d(x, x_{k_1}) \\
        & < & \frac{r}{2} + \frac{1}{2} t(x_{k_1}) \\
        & \leq & \frac{1}{2} t(x_{k_1}) + \frac{1}{2} t(x_{k_1}) \\
        = t(x_{k_1}) ,
      \end{eqnarray*}
      joten
      \(y \in \openball{E}{x_{k_1}}{t(x_{k_1})}\), mistä seuraa, että
      \begin{equation}
        \label{eq:co-c}
        \abs{f(x)-f(y)} \leq \abs{f(x)-f(x_{k_1})} + \abs{f(x_{k_1}) -
          f(y)}
        < \frac{h}{2} + \frac{h}{2} = h .
      \end{equation}
  \end{enumerate}
  Tasainen jatkuvuus seuraa yhtälöistä \eqref{eq:co-a},
  \eqref{eq:co-b} ja \eqref{eq:co-c}.
\end{solution}

\begin{solution}{9.6}
  Olkoon \(f \in C^m(\rn)\). Nyt \(f \in C^{m-1}(\rn)\).
  Meillä on
  \begin{displaymath}
    \label{eq:zygmund-norm}
    \norminspace{f}{\zygmundspace{m}{\rn}}
    =
    \norminspace{f}{C^{m-1}(\rn)} + \sum_{\abs{\alpha} = m - 1} \sup
    \left\{ \frac{1}{t} \genmodcont{2}{\infty}{\mideriv{\alpha} f}{t}
    \bigsetsep t \in \realnumbers_+ \right\} .
  \end{displaymath}
  Kaikilla \(\alpha \in \naturalnumbers^n\), \(\abs{\alpha} = m-1\) on
  \(\mideriv{\alpha} f \in C^1(\rn)\).
  Olkoot \(\rnx,\rny \in \rn\), \(\rnx \not= \rny\) ja \(\rnh := \rny - \rnx\).
  Taylorin kehitelmästä \eqref{eq:taylor-mdim},
  \eqref{eq:taylor-residual-mdim}
  saadaan
  \begin{displaymath}
    (\mideriv{\alpha} f)(\rnx + \rnh)
    =
    (\mideriv{\alpha} f)(\rnx)
    +
    \sum_{k=1}^n \left(\mideriv{{\unitvec{n}{k}}} \mideriv{\alpha}
    f\right)
    \left(\rnx + c \rnh \right) \rnh[k] ,
  \end{displaymath}
  missä \(c \in \intervaloo{0}{1}\).
  Tästä seuraa, että
  \begin{displaymath}
    (\mideriv{\alpha} f)(\rny) - (\mideriv{\alpha} f)(\rnx)
    =
    \sum_{k=1}^n \left(\mideriv{{\unitvec{n}{k}}} \mideriv{\alpha}
    f\right)
    \left(\rnx + c\rnh\right) \rnh[k] .
  \end{displaymath}
  Nyt
  \begin{eqnarray*}
    \abs{(\mideriv{\alpha} f)(\rny) - (\mideriv{\alpha} f)(\rnx)} & \leq &
    n \norminspace{\mideriv{\alpha} f}{C^1(\rn)} \cdot \norm{\rnh}_\infty
    \\
    & \leq & n \norminspace{f}{C^m(\rn)} \norm{\rnh}_\infty
  \end{eqnarray*}
  mistä seuraa, että
  \begin{displaymath}
    \norm{\diff{\rnh}(\mideriv{\alpha} f)}_\infty \leq
    n \norminspace{f}{C^m(\rn)} \norm{\rnh}_\infty ,
  \end{displaymath}
  joten
  \begin{displaymath}
    \frac{1}{\norm{\rnh}_\infty} \norm{\higherdiff{2}{\rnh}(\mideriv{\alpha}
      f)}_\infty \leq 2n \norminspace{f}{C^m(\rn)} .
  \end{displaymath}
  Täten
  \begin{equation}
    \label{eq:zygmund-sup}
    \sup \left\{ \frac{1}{\norm{\rnh}_\infty}
    \norm{\higherdiff{2}{\rnh}(\mideriv{\alpha} f)}_\infty \bigsetsep \rnh \in
    \rn \setminus \zeroset \right\}
    \leq
    2n \norminspace{f}{C^m(\rn)} .
  \end{equation}
  Yhtälöistä \eqref{eq:zygmund-norm} ja \eqref{eq:zygmund-sup} seuraa,
  että
  \begin{displaymath}
    \norminspace{f}{\zygmundspace{m}{\rn}} \leq
    2n \norminspace{f}{C^m(\rn)} < \infty .
  \end{displaymath}
\end{solution}

\begin{solution}{9.7}
  Olkoon \(B_k := \closedball{\rn}{0}{k+1}\) kaikille \(k \in
  \naturalnumbers\).  Olkoon \(k_0 \in \naturalnumbers\) ja \(c \in
  \realnumbers_+\). On olemassa \(m \in \naturalnumbers\) siten, että
  \begin{displaymath}
    \norm{g-f_j}_p < \frac{c}{2} \;\;\land\;\;
    \norminfty{h-f_j} < \frac{c}{2} \left(\mu(B_{k_0})\right)^{-\frac{1}{p}}
  \end{displaymath}
  kaikille \(j \in \naturalnumbers\), \(j \geq m\).
  Nyt
  \begin{displaymath}
    \norm{(g-f_j) \restriction B_{k_0}} < \frac{c}{2}
  \end{displaymath}
  kaikille \(j \in \naturalnumbers\), \(j \geq m\) ja
  \begin{displaymath}
    \norm{(h-f_j) \restriction B_{k_0}} < \frac{c}{2}
    \left(\mu\left(B_{k_0}\right)\right)^{-\frac{1}{p}}
  \end{displaymath}
  kaikille \(j \in \naturalnumbers\), \(j \geq m\).
  Olkoon \(r \in \naturalnumbers\), \(r \geq m\).
  Nyt
  \begin{displaymath}
    \norm{(h-f_r) \restriction B_{k_0}}_p < \left(\left(\frac{c}{2}
    \left(\mu\left(B_{k_0}\right)\right)^{-\frac{1}{p}}\right)^p
    \mu(B_{k_0}) \right)^\frac{1}{p}
    = \frac{c}{2} .
  \end{displaymath}
  Täten
  \begin{eqnarray*}
    \norm{(g-h) \restriction B_{k_0}}_p & \leq & \norm{(g-x_r) \restriction B_{k_0}}_p +
    \norm{(x_r-h) \restriction B_{k_0}}_p \\
    & < & \frac{c}{2} + \frac{c}{2} = c .
  \end{eqnarray*}
  Koska \(c \in \realnumbers_+\) oli mielivaltainen, nii
  \begin{equation}
    \label{eq:ae}
    g \restriction B_{k_0} = h \restriction B_{k_0} \spaceafter
    \text{melkein kaikkialla.}
  \end{equation}
  Olkoot
  \begin{eqnarray*}
    C_k & := & \left\{ \rnx \in B_k \setsep g(\rnx) \not= h(\rnx) \right\} ,
    \\
    C  & := & \left\{ \rnx \in \rn \setsep g(\rnx) \not= h(\rnx) \right\} .
  \end{eqnarray*}
  Nyt
  \begin{displaymath}
    C = \bigcup_{k=0}^\infty C_k
  \end{displaymath}
  ja \((C_k)_{k=0}^\infty\) on kasvava jono joukkoja. Yhtälön
  \eqref{eq:ae} nojalla \(\mu(C_k) = 0\) kaikilla \(k \in
  \naturalnumbers\).
  Lauseen \ref{th:additiivisuus} nojalla
  \begin{displaymath}
    \mu(C) = \lim_{k \to \infty} \mu(C_k) = \lim_{k \to \infty} 0 = 0 .
  \end{displaymath}
  Siis \(g = h\) melkein kaikkialla.
\end{solution}

\backmatter

\printindex

\bibliographystyle{abbrv}

\bibliography{jfa}

\end{document}